\documentclass[11pt]{article}

\usepackage[pdftex]{graphicx}
\usepackage{amsmath}
\usepackage{amssymb}
\usepackage{caption}
\usepackage[subrefformat=parens]{subcaption}
\usepackage{here}
\newtheorem{theorem}{Theorem}
\newtheorem{corollary}{Corollary}

\usepackage{version}
\excludeversion{en-text}
\usepackage{comment}
\includecomment{jp-text}
\excludecomment{en-text}
\def\be{\begin{equation}}
\def\ee{\end{equation}}
\def\bea{\begin{eqnarray}}
\def\eea{\end{eqnarray}}
\def\beas{\begin{eqnarray*}}
\def\eeas{\end{eqnarray*}}

\def\bN2{\bar{N}_2}
\def\dlnt{\delta}

\def\be{\begin{equation}}
\def\ee{\end{equation}}
\def\bea{\begin{eqnarray}}
\def\eea{\end{eqnarray}}
\def\beas{\begin{eqnarray*}}
\def\eeas{\end{eqnarray*}}

\newif\ifcol
\colfalse

\newcommand{\argsup}{\mathop{\rm arg~sup}\limits}
\newcommand{\arginf}{\mathop{\rm arg~inf}\limits}

\begin{document}
\title{Parametric estimation for a parabolic linear SPDE model \\
based on sampled data
}
\author{
${}^{1}${Yusuke Kaino}
 and
 ${}^{1,2,3}${Masayuki Uchida}
 \\
        ${}^1${\small {Graduate School of Engineering Science, Osaka University}}\\
        ${}^2${\small {Center for Mathematical Modeling and Data Science (MMDS),  
         Osaka University} and }\\
          ${}^3${\small {JST CREST}}, 
        {\small {Toyonaka, Osaka 560-8531, Japan}}
}
\maketitle
\noindent
{\bf Abstract.}
We consider parametric estimation for a parabolic linear second order stochastic partial differential equation (SPDE) 
from high frequency data which are  observed in time and space. 
By using thinned data obtained from the high frequency data,
adaptive estimators of the coefficient parameters including 
the volatility parameter of a parabolic linear SPDE model are proposed. 
Moreover, we give {some} examples and simulation results of the adaptive estimators of the SPDE model 
based on the high frequency data.

\begin{en-text} 
We consider parameter estimation of a parabolic second-order linear stochastic partial differential equation model using discrete observations.
There are four parameters to be estimated, and three of them are used to construct normalized volatility parameter and curvature parameter.
We propose two methods to estimate normalized volatility parameter and curvature parameter.
The first method uses data from time 0 to time 1, and the second method uses data from time 0 to time $T$.
A coordinate process is generated from discrete observations and the estimator of the curvature parameter, and we estimate its drift and diffusion coefficient.
The parameters of the original model are estimated from this estimator and the estimator of the normalized volatility parameter.
We verify the asymptotic behavior of the estimators by simulation.
\end{en-text}

\vspace{0.5cm} 


%




\section{Introduction}

We consider a linear parabolic stochastic partial differential equation (SPDE) with one space dimension.
\begin{eqnarray}
& &dX_{t}(y) = \left( \theta_2 \frac{\partial^2 X_{t}(y)}{\partial y^2} + \theta_1 \frac{\partial X_{t}(y)}{\partial y} + \theta_0 X_{t}(y) \right)dt 
+ \sigma dB_t(y), \quad  (t,y) \in [0, T]\times[0, 1],  \quad
\label{spde0} \\
& & X_t(0) = X_t(1) = 0, \quad t \in [0,T], \qquad X_0(y) = \xi = 0, \quad y \in [0,1],
\nonumber
\end{eqnarray}
where
$T>0$, $B_t$ is defined as a cylindrical Brownian motion in the Sobolev space on $[0,1]$,
the initial condition $\xi = 0$,
an unknown parameter $\theta= (\theta_0, \theta_1, \theta_2, \sigma)$  
and $\theta_0,\theta_1 \in \mathbb{R}, \theta_2, \sigma > 0$,
and the parameter space $\Theta$ is a compact convex subset of $\mathbb{R}^2 \times [0, \infty]^2$.
Moreover, the true value of parameter $\theta^*= (\theta_0^*, \theta_1^*, \theta_2^*, \sigma^*)$
and we assume that $\theta^* \in \mbox{Int}(\Theta)$.
The data are discrete observations ${\bf X}_{N,M}= \left\{ X_{t_{i:N}}(y_{j:M}) \right\}_{i = 1, ..., N, j = 1, ..., M}$, 
$t_{i:N} = i \frac{T}{N}$,
${ y_{j:M} = \frac{j}{M}}$.
{
For the characteristics of the parameters
$\theta_0$, $\theta_1$, $\theta_2$ and $\sigma$ of the SPDE (\ref{spde0}),
see the Appendix
below.
}

Statistical inference for SPDE models based on sampled data has been developed
by many researchers, see for example,
Markussen (2003), Cont (2005), Cialenco and Glatt-Holtz (2011),
Cialenco and Huang (2017), 
Bibinger and Trabs (2017),
Cialenco et.\ al.\ (2018, 2019),
Cialenco (2018),  Chong (2019) and references therein.
Recently, Bibinger and Trabs (2017) studied the parabolic linear second order SPDE model 
based on high frequency data observed on a fixed region and 
proved the asymptotic properties of minimum contrast estimators $\check{\sigma}_0^2$ and $\check{\eta}$ 
for both the normalized volatility parameter $\sigma_0^2 = \frac{\sigma^2}{\sqrt{\theta_2}}$ 
and the curvature parameter $\eta = \frac{\theta_1}{\theta_2}$.

In this paper, we propose adaptive maximum likelihood (ML) type estimator of the coefficient parameter
$\theta= (\theta_0, \theta_1, \theta_2, \sigma)$ 
of the parabolic linear second order SPDE model (\ref{spde0}).
For $k \in \mathbb{N}$, 
the coordinate process  $x_k(t)$ of the SPDE model (\ref{spde0}) is that 
\begin{eqnarray}
x_k(t) 
&=& \int^1_0 X_t(y) \sqrt{2}\sin (\pi k y) \exp \left( \frac{\eta y}{2} \right) dy, \label{coordinate}
\end{eqnarray}
which satisfies that 
$$
dx_k(t) = -\lambda_k x_k(t) dt + \sigma dw_k(t), \quad x_k(0) = 0,
$$
where 
\begin{eqnarray*}
\lambda_k &=& - \theta_0 + \frac{\theta_1^2}{4 \theta_2} + \pi^2 k^2 \theta_2.
\end{eqnarray*}
Note that the coordinate process (\ref{coordinate}) 
is  the Ornstein-Uhlenbeck process.
Using the minimum contrast estimator $\check{\eta}$ proposed by Bibinger and Trabs (2017), 
we obtain the approximate coordinate process 
\begin{eqnarray*}
\check{x}_k(t) = \frac{1}{M}\sum^M_{j = 1} X_{t}(y_{j:M}) \sqrt{2}\sin (\pi ky_{j:M}) 
\exp \left( \frac{\check{\eta} y_{j:M}}{2} \right)
\end{eqnarray*}
and the adaptive estimator is constructed by using  the property that 
the coordinate process (\ref{coordinate}) is a diffusion process. 
{
It is also shown that the adaptive ML type estimators have asymptotic normality under some regularity conditions. 
Furthermore, in order to verify asymptotic performance of the adaptive ML type estimators of the coefficient parameters of the parabolic linear second order SPDE model based on high-frequency data, some examples and simulation results of the adaptive ML type estimators are given. 
}
For details of statistical inference for diffusion type processes and stochastic differential equations,
{ see}
Prakasa Rao (1983,1988), 
Kutoyants (1994, 2004), 
Florens-Zmirou (1989),  
Yoshida (1992, 2011), 
Bibby and S{\o}rensen (1995),  
Kessler (1995, 1997), 
Uchida (2010), 
Uchida and Yoshida (2012, 2014),
De Gregorio and Iacus (2013),
Kamatani and Uchida (2015),
Nakakita and {Uchida (2019)}
 for ergodic diffusion processes,
 {and}
Shimizu and Yoshida (2006), 
Shimizu (2006), 
Ogihara and Yoshida (2011), 
Masuda (2013a, 2013b) for jump diffusion processes and L{\'e}vy type processes, 
 {and}
Dohnal (1987), Genon-Catalot and Jacod (1993, 1994), 
Uchida and Yoshida (2013),
Ogihara and Yoshida (2014),
Ogihara (2018), 
Kaino and Uchida (2018)
for non-ergodic diffusion processes.
{
For  adaptive ML type estimators and thinned data
for diffusion type processes, see for example, Uchida and Yoshida (2012) 
and Kaino and Uchida (2018). 
}

This paper is organized as follows. 
In Section 2, we consider 
{the adaptive estimator} 
of the SPDE model based on 
{the sampled data 
in the fixed region
$[0,1] \times [0,1]$.}
The adaptive estimator is constructed by using the minimum contrast estimators of $\sigma_0^2$ and $\eta$ proposed by 
Bibinger and Trabs (2017). It is shown that the adaptive estimator has asymptotic normality.
In Section 3, we deal with the SPDE model based on sampled data which are observed in the region
$[0,T] \times [0,1]$ when $T$ is large. The quasi log likelihood function is obtained 
by using the approximate coordinate process and we prove that the adaptive ML type estimator
of  $\theta= (\theta_0, \theta_1, \theta_2, \sigma)$ has 
{asymptotic normality.}
In Section 4, concrete examples are given and 
the asymptotic behavior of the estimators proposed in Sections 2 and 3 is verified by simulations.
Section 5 is devoted to the proofs of the results presented in Sections 2 and 3.
{
The Appendix contains the sample paths
with different values of the parameters 
to understand the characteristics of the parameters
$\theta_0$, $\theta_1$, $\theta_2$ and $\sigma$ of the SPDE (\ref{spde0}).}

\section{The case that T is fixed}

In this section, we treat the linear parabolic SPDE (\ref{spde0}) with $T=1$, which is defined as
\begin{eqnarray}
& &dX_{t}(y) = \left( \theta_2 \frac{\partial^2 X_{t}(y)}{\partial y^2} + \theta_1 \frac{\partial X_{t}(y)}{\partial y} + \theta_0 X_{t}(y) \right)dt 
+ \sigma dB_t(y), \quad  (t,y) \in [0, 1]\times[0, 1],  \quad
\label{spde1} \\
& & X_t(0) = X_t(1) = 0, \quad t \in [0,1], \qquad X_0(y) = \xi = 0, \quad y \in [0,1].
\nonumber
\end{eqnarray}
The data are discrete observations ${{\bf \bar{X}}_{N, \bar{M}}}= \left\{ X_{t_{i:N}}({\bar{y}_{j:M}}) \right\}_{i = 1, \ldots, N, j = 1, \ldots, {\bar{M}}}$ with
{$t_{i:N} =  \frac{i}{N}$},
{$\bar{y}_{j:M} = \delta+ \frac{j-1}{M}$ and
$\bar{M}=[(1-2\delta)M]$ for $\delta \in (0,1/2) $.
Note that $\delta \leq \bar{y}_{j:M} \leq 1-\delta$ for $ j = 1, \ldots, \bar{M}$,
and that ${\bf \bar{X}}_{N, \bar{M}}$
is generated by the full data  ${\bf X}_{N,M}= \left\{ X_{t_{i:N}}(y_{j:M}) \right\}_{i = 1, ..., N, j = 1, ..., M}$
with $t_{i:N} = \frac{i}{N}$ and $y_{j:M} = \frac{j}{M}$}.

Setting that {the differential operator}  
\begin{eqnarray*}
A_\theta := \theta_0  + \theta_1\frac{\partial}{\partial y} + \theta_2\frac{\partial^2}{\partial y^2},
\end{eqnarray*}
one has that  for $k \in \mathbb{N}$,
$$
A_\theta  e_k =  -\lambda_k e_k,
$$
where the eigenfunctions $e_k$ of $A_\theta$ and the corresponding eigenvalues $-\lambda_k$ are given by
\begin{eqnarray*}
e_k(y) &=& \sqrt{2} \sin(\pi k y)\exp\left( -\frac{\theta_1}{2\theta_2} y \right), \;\;\; y\in [0,1], \\
\lambda_k &=& - \theta_0 + \frac{\theta_1^2}{4 \theta_2} + \pi^2 k^2 \theta_2.
\end{eqnarray*}
Let 
$H_\theta := \{ f:[0,1] \to \mathbb{R} : \| f \|_\theta < \infty, f(0) = f(1) = 0  \}$  with 
$\langle f,g \rangle_\theta := \int^1_0 e^{y\theta_1/\theta_2}f(y)g(y) dy$ 
and $ \| f \|_\theta := \sqrt{\langle f,f \rangle_\theta}$.
The initial condition $\xi(y) = 0$  and $\xi \in H_\theta$.

\begin{en-text}
The cylindrical Brownian motion $(B_t)_{t \ge 0}$ in (\ref{spde1}) can be defined as 
\begin{eqnarray*}
\langle B_t,f \rangle_\theta = \sum_{k \ge 1} \langle f, e_k \rangle_\theta W^k_t, \quad 
f \in H_\theta, \quad t \ge 0
\end{eqnarray*}
for independent real-valued Brownian motions $(W^k_t)_{t\ge0}$, $k\ge1$. 
$X_t(y)$ is called a mild solution of (\ref{spde1}) on $[0,1]$ if it satisfies that for any $t \in [0,1]$,  
\begin{eqnarray*}
X_t = e^{tA_\theta}\xi + \int^t_0 e^{(t-s)A_\theta}\sigma dB_s \;\;\; a.s.
\end{eqnarray*}
\end{en-text}

\begin{en-text}
\noindent
$[A2-1]$
Suppose we observe a mild solution X of (1) on discrete grid $(t_{i:N}, y_{j:M}) \in [0,1]^2$ with 
\begin{eqnarray*}
t_{i:N} = \frac{i}{N} \;\; for \;\; i = 0,...,n \;\; and \;\; \delta \le y_{1:M} < y_{2:M} < ... < y_{M:M} \le 1 - \delta
\end{eqnarray*}
where  $N,M \in \mathbb{N}$ and $\delta > 0$. Let $M max_{j = 2,...,M}|y_{j:M} - y_{j-1:M}|$ be uniformly bounded.
\end{en-text}

\subsection{Estimation of $\sigma_0^2$ and $\eta$}
Set $\sigma_0^2 = \frac{\sigma^2}{\sqrt{\theta_2}}$, $\eta = \frac{\theta_1}{\theta_2}$.
Let $m \leq {\bar{M}}$
and
$\tilde{y}_{j:m} = {\delta+} \left[ \frac{{\bar{M}}}{m} \right] {\frac{j-1}{M}} $ 
for $j = 1,...,m$.
Assume that 
$m \leq N^\rho$ for $\rho \in (0, 1/2)$.
Set 
\begin{eqnarray*}
Z_{j:m} = \frac{1}{N \sqrt{{t}_{1:N}}} \sum_{i=1}^{N} (X_{\tilde{t}_{i:N}}(\tilde{y}_{j:m} )- X_{\tilde{t}_{i-1:N}}(\tilde{y}_{j:m}))^2.
\end{eqnarray*}
The contrast function is defined as
\begin{eqnarray*}
U_{N,m}(\sigma_0^2,\eta) = \frac{1}{m} \sum^{m}_{j = 1}  \left( Z_{j:m} - \frac{1}{\sqrt{\pi}} \sigma_0^2 \exp(- \eta \tilde{y}_{j:m}) \right)^2.
\end{eqnarray*}
The minimum contrast estimator of $\sigma_0^2$ and $\eta$ are given by 
\begin{eqnarray*}
(\check{\sigma}_0^2,\check{\eta}) = \arginf_{\sigma_0^2,\eta}U_{N,m}(\sigma_0^2,\eta).
\end{eqnarray*}
Let
\begin{eqnarray*}
U(\zeta^*)  &=& 
\begin{pmatrix}
\int^{1-\delta}_\delta e^{-4 \eta^* y}dy & -(\sigma^*_0)^2\int^{1-\delta}_\delta ye^{-4 \eta^* y}dy\\
 -(\sigma_0^*)^2\int^{1-\delta}_\delta ye^{-4 \eta^* y}dy& (\sigma^*_0)^4\int^{1-\delta}_\delta y^2e^{-4 \eta^* y}dy
\end{pmatrix}, \\
V(\zeta^*) &=& 
\begin{pmatrix}
\int^{1-\delta}_\delta e^{-2 \eta^* y}dy & -(\sigma^*_0)^2\int^{1-\delta}_\delta ye^{-2 \eta^2 y}dy\\
 -(\sigma^*_0)^2\int^{1-\delta}_\delta ye^{-2 \eta^* y}dy& (\sigma^*_0)^4\int^{1-\delta}_\delta y^2e^{-2 \eta^* y}dy
\end{pmatrix}, \\
\Gamma &=& \frac{1}{\pi} \sum^{\infty}_{r = 0} I(r)^2 + \frac{2}{\pi}
\ \mbox{with} \ I(r) = 2\sqrt{r+1} - \sqrt{r+2} - \sqrt{r}.
\end{eqnarray*}

\begin{theorem}[Bibinger and Trabs(2017)] \label{thm1}
Let $\zeta^* = ((\sigma_0^*)^2, \eta^*) \in \Xi$ for {a} compact subset $\Xi \subseteq (0,\infty) \times [0,\infty)$.
Assume that $m \to \infty$ and 
$m = \mathcal{O}(N^\rho)$ for some $\rho \in (0,1/2)$.
Then, as $N \to \infty$, 
\begin{eqnarray*}
\sqrt{m N}(\check{\sigma}_0^2 - (\sigma_0^*)^2,  \check{\eta}-\eta^* )
\stackrel{d}{\longrightarrow } N(0,(\sigma_0^*)^4 \Gamma \pi V(\zeta^*)^{-1}U(\zeta^*)V(\zeta^*)^{-1}).
\end{eqnarray*}
\end{theorem}

\subsection{Estimation of $\theta_1,\theta_2$ and $\sigma^2$}
Next we consider the estimation of $\sigma,\theta_1,\theta_2$ using {$\check{\sigma}_0^2,\check{\eta}$}. 
Let $k \in \mathbb{N}$.
The coordinate process is that 
\begin{eqnarray*}
x_k(t) &=& \langle X_t, e_k \rangle_\theta 
= \int^1_0 \exp \left( \frac{\theta_1}{\theta_2} y \right)  X_t(y) e_k(y) dy \\
&=& \int^1_0 X_t(y) \sqrt{2}\sin (\pi k y) \exp \left( \frac{\eta y}{2} \right) dy .
\end{eqnarray*}
Note that
$$
dx_k(t) = -\lambda_k x_k(t) dt + \sigma dw_k(t), \quad x_k(0) = 0, 
$$
which means that $x_k(t)$ is  the Ornstein-Uhlenbeck process.
Furthermore, the random field $X_t(y)$ is that
\begin{eqnarray*}
X_t(y) = \sum_{k=1}^\infty x_k(t) e_k(y).
\end{eqnarray*}

Let $N_2 \leq N$ and ${s}_{i:N_2}=\left[\frac{N}{N_2}\right] t_{i:N}=i \left[\frac{N}{N_2}\right] \frac{1}{N}$ for $i = 1,...,N_2$．
As an approximation of $x_k(t)$, we consider 
\begin{eqnarray*}
\check{x}_k({s}_{i:N_2}) = \frac{1}{M}\sum^M_{j = 1} X_{{s}_{i:N_2}}(y_{j:M}) \sqrt{2}\sin (\pi ky_{j:M}) 
\exp \left( \frac{\check{\eta} y_{j:M}}{2} \right), 
\quad i = 1,...,N_2 .
\end{eqnarray*}
{By using the thinned data based on the approximate coordinate process 
$\left\{ \check{x}_k(s_{i:N_2}) \right\}_{i=1,\ldots, N_2}$, }
the estimator of $\sigma^2$ is defined as 
\begin{eqnarray*}
\check{\sigma}_k^2 = \sum^{N_2}_{i = 1} ( \check{x}_k({s}_{i:N_2}) - \check{x}_k({s}_{i-1:N_2}) )^2.
\end{eqnarray*}
Set  $k=1$ and $\check{\sigma}^2=\check{\sigma}_1^2$．
The estimators of $\theta_1$ and $\theta_2$ are defined as 
\begin{eqnarray*}
\check{\theta}_2 = \left( \frac{ \check{\sigma}^2}{\check{\sigma}_0^2} \right)^2 , \quad  
\check{\theta}_1 = \check{\eta}  \check{\theta_2}.
\end{eqnarray*}

\begin{theorem} \label{thm2}
Assume the same conditions as Theorem \ref{thm1}. 
Moreover, assume that $\frac{N_2^{3/2}}{mN} \rightarrow 0$  
and $\frac{N_2^{3/2}}{M^{1-\rho_1}} \rightarrow 0$ for $\rho_1 \in (0,1)$.
Then, as $N_2 \rightarrow \infty$, \\
$$
\begin{pmatrix}
\sqrt{N_2}(\check{\sigma}^2 - (\sigma^{*})^2) \\
\sqrt{N_2}(\check{\theta}_2 - \theta_{2}^*) \\
\sqrt{N_2}(\check{\theta}_1 - \theta_{1}^*)
\end{pmatrix}
\stackrel{d}{\longrightarrow } N
\left(
\begin{pmatrix}
0 \\
0 \\
0 
\end{pmatrix}
,
\begin{pmatrix}
2(\sigma^*)^4 &4\theta_2^*(\sigma^*)^2 & 4\theta_1^*(\sigma^*)^2\\
4\theta_2^*(\sigma^*)^2 &8(\theta_2^*)^2 & 8\theta_1^*\theta_2^*\\
 4\theta_1^*(\sigma^*)^2 &  8\theta_1^*\theta_2^*& 8(\theta_1^*)^2
\end{pmatrix}
\right)
$$

\end{theorem}

\section{The case that $T$ is large}

In this section, we deal with the linear parabolic SPDE (\ref{spde0}) when $T$ is large, which is defined as
\begin{eqnarray}
& &dX_{t}(y) = \left( \theta_2 \frac{\partial^2 X_{t}(y)}{\partial y^2} + \theta_1 \frac{\partial X_{t}(y)}{\partial y} + \theta_0 X_{t}(y) \right)dt 
+ \sigma dB_t(y), \quad  (t,y) \in [0, T] \times[0, 1],  \quad
\label{spde2} \\
& & X_t(0) = X_t(1) = 0, \quad t \in [0,T], \qquad X_0(y) = \xi = 0, \quad y \in [0,1].
\nonumber
\end{eqnarray}
The data are discrete observations 
${{\bf \check{X}}_{N,\bar{M}}}= \left\{ X_{t_{i:N}}({\bar{y}_{j:M}}) \right\}_{i = 1, ..., N, j = 1, ..., {\bar{M}}}$
with 
$t_{i:N} =  i h_{N:T}$, $h_{N:T}= \frac{T}{N}$,
{$\bar{y}_{j:M} = \delta+ \frac{j-1}{M} $
and  
$\bar{M}=[(1-2\delta)M]$ for $\delta \in (0,1/2) $.
Note that $\delta \leq \bar{y}_{j:M} \leq 1-\delta$ for $ j = 1, \ldots, \bar{M}$,
and that ${\bf \check{X}}_{N, \bar{M}}$
is generated by the full data  ${\bf X}_{N,M}= \left\{ X_{t_{i:N}}(y_{j:M}) \right\}_{i = 1, ..., N, j = 1, ..., M}$
with $t_{i:N} = i h_{N:T}$ and $y_{j:M} = \frac{j}{M}$}.

\begin{en-text}
We consider the second-order linear stochastic partial differential equation model as follows,
\begin{eqnarray*}
dX_{t}(y) = \left( \theta_2 \frac{\partial^2 X_{t}(y)}{\partial y^2} + \theta_1 \frac{\partial X_{t}(y)}{\partial y} + \theta_0 X_{t}(y) \right)dt 
+ \sigma dB_t(y), \\
(t,y) \in [0, T]\times[0, 1],  \quad X_t(0) = X_t(1) = 0, \quad X_0(y) = \xi = 0. 
\end{eqnarray*}
where
$B_t$ is defined as a cylindrical Brownian motion in the Sobolev space on $[0,1]$,
the initial condition $\xi = 0$,
unknown parameter is
$\theta_0,\theta_1 \in \mathbb{R}, \theta_2, \sigma > 0$.
The data are discrete observations ${\bf X}_{N:T, M} 
= \left( X_{t_{i:N:T}}(y_{j:M}) \right)_{i = 1, ..., N, j = 1, ..., M}$,
$t_{i:N:T} = i \Delta_{N:T}$, $\Delta_{N:T} = \frac{T}{N}$, 
$y_{j:M} = (1-2\delta)\frac{j-1}{M+1} + \delta$, $0 < \delta < \frac{1}{2}$.

Differential operator is defined as
\begin{eqnarray*}
A_\theta := \theta_0  + \theta_1\frac{\partial}{\partial y} + \theta_2\frac{\partial^2}{\partial y^2}.
\end{eqnarray*}
The eigenfunctions $e_k$ of $A_\theta$ with corresponding eigenvalues $-\lambda_k$ are given by
\begin{eqnarray*}
e_k(y) &=& \sqrt{2} \sin(\pi k y)\exp\left( -\frac{\theta_1}{2\theta_2} y \right), \;\;\; y\in [0,1]\\
\lambda_k &=& - \theta_0 + \frac{\theta_1^2}{4 \theta_2} + \pi^2 k^2 \theta_2, \;\;\; k \in \mathbb{N}
\end{eqnarray*}

\begin{eqnarray*}
H_\theta := \{ f:[0,1] \to \mathbb{R} : \| f \| < \infty, f(0) = f(1) = 0  \} \;\;\; with \\
\langle f,g \rangle_\theta := \int^1_0 e^{y\theta_1/\theta_2}f(y)g(y) dy \;\;\; and \;\;\; \| f \|^2_\theta := \langle f,f \rangle.
\end{eqnarray*}
Initial condition $\xi(y) = 0$ is $\xi \in H_\theta$.
The cylindrical Brownian motion$(B_t)_t\ge0$ in (1) can be defined via 
\begin{eqnarray*}
\langle B_t,f \rangle_\theta = \sum_{k \ge 1} \langle f, e_k \rangle_\theta W^k_t, \;\;\; 
f \in H_\theta,t \ge 0,
\end{eqnarray*}
for independent real-valued Brownian motions $(W^k_t)_t\ge0,k\ge1$.$X_t(y)$ is called a mild solution of (1) on $[0,1]$ if it satisfies the integral equation for any $t \in [0,1]$
\begin{eqnarray*}
X_t = e^{tA_\theta}\xi + \int^t_0 e^{(t-s)A_\theta}\sigma dB_s \;\;\; a.s.
\end{eqnarray*}

\end{en-text}

\subsection{Estimation of $\sigma_0^2, \eta$}
Let $\bar{m} \leq {\bar{M}}$
and
$\tilde{y}_{j: \bar{m}} = {\delta+} \left[ \frac{{\bar{M}}}{\bar{m}} \right] {\frac{j-1}{M}} $ 
for $j = 1,...,\bar{m}$.
Set 
\begin{eqnarray*}
\bar{{\cal Z}}_{j:\bar{m}} = \frac{1}{N \sqrt{h_{N:T}}} \sum_{i=1}^{N} 
(X_{t_{i:N:T} }(\tilde{y}_{j:\bar{m}}) - X_{ t_{i-1:N:T}}(\tilde{y}_{j:\bar{m}}))^2.
\end{eqnarray*}
The contrast function is defined as
\begin{eqnarray*}
\bar{{\cal U}}_{N,\bar{m}}(\sigma_0^2,\eta) = \frac{1}{\bar{m}} \sum^{\bar{m}}_{j = 1} 
\left( \bar{{\cal Z}}_{j:\bar{m}} - \frac{1}{\sqrt{\pi}} \sigma_0^2 \exp(-\eta \tilde{y}_{j:\bar{m}}) \right)^2. 
\end{eqnarray*}
The minimum contrast estimator of $\sigma_0^2$ and $\eta$ are given by 
\begin{eqnarray*}
(\hat{\sigma}_0^2,\hat{\eta}) = \arginf_{\sigma_0^2,\eta}\bar{{\cal U}}_{N,\bar{m}}(\sigma_0^2,\eta). 
\end{eqnarray*}

\begin{theorem} \label{thm3}
Assume that 
$N h_{N:T}^2 \rightarrow 0$,  
$\bar{m} \rightarrow \infty$ and 
$\bar{m} = O(h_{N:T}^{-\rho})$ 
for
$\rho \in (0, 1/2)$.
Then, as $N \to \infty$, 
\begin{eqnarray*}
\sqrt{\bar{m} N}(\hat{\sigma}_0^2 - (\sigma_0^*)^2,  {\hat{\eta}}-\eta^* )
\stackrel{d}{\longrightarrow } N(0,(\sigma_0^*)^4 \Gamma \pi V(\zeta^*)^{-1}U(\zeta^*)V(\zeta^*)^{-1}).
\end{eqnarray*}
\end{theorem}

\begin{en-text}
where
\begin{eqnarray*}
U(\zeta^*) := 
\begin{pmatrix}
\int^{1-\delta}_\delta e^{-4 \eta^* y}dy & -(\sigma^*_0)^2\int^{1-\delta}_\delta ye^{-4 \eta^* y}dy\\
 -(\sigma_0^*)^2\int^{1-\delta}_\delta ye^{-4 \eta^* y}dy& (\sigma^*_0)^4\int^{1-\delta}_\delta y^2e^{-4 \eta^* y}dy
\end{pmatrix}
,\\
V(\zeta^*) := 
\begin{pmatrix}
\int^{1-\delta}_\delta e^{-2 \eta^* y}dy & -(\sigma^*_0)^2\int^{1-\delta}_\delta ye^{-2 \eta^* y}dy\\
 -(\sigma^*_0)^2\int^{1-\delta}_\delta ye^{-2 \eta^* y}dy& (\sigma^*_0)^4\int^{1-\delta}_\delta y^2e^{-2 \eta^* y}dy
\end{pmatrix}
,\\
\end{eqnarray*}
\begin{eqnarray*}
\Gamma := \frac{1}{\pi} \sum^{\infty}_{r = 0} I(r)^2 + \frac{2}{\pi}\;\;\; with \;\;\; I(r) = 2\sqrt{r+1} - \sqrt{r+2} - \sqrt{r}.
\end{eqnarray*}
\end{en-text} 

Next we consider the case that $\bar{m}$ is fixed.
Let
\begin{eqnarray*}
U_{\bar{m}}(\zeta^*) &=& 
\begin{pmatrix}
\sum_{j=1}^{\bar{m}} e^{-4 \eta^* y_{j:\bar{m}}} & -(\sigma^*_0)^2\sum_{j=1}^{\bar{m}} y_{j:\bar{m}} e^{-4 \eta^* y_{j:\bar{m}}}\\
-(\sigma^*_0)^2\sum_{j=1}^{\bar{m}} y_{j:\bar{m}} e^{-4 \eta^* y_{j:\bar{m}}}& (\sigma^*_0)^4\sum_{j=1}^{\bar{m}} y_{j:\bar{m}} e^{-4 \eta^* y_{j:\bar{m}}^2}
\end{pmatrix},
\\
V_m(\zeta^*) &=& 
\begin{pmatrix}
\sum_{j=1}^{\bar{m}} e^{-2 \eta^* y_{j:\bar{m}}} & -(\sigma^*_0)^2\sum_{j=1}^{\bar{m}} y_{j:\bar{m}} e^{-2 \eta^* y_{j:\bar{m}}}\\
-(\sigma^*_0)^2\sum_{j=1}^{\bar{m}} y_{j:\bar{m}} 
e^{-2 \eta^* y_{j:\bar{m}}}& (\sigma^*_0)^4\sum_{j=1}^{\bar{m}} y_{j:\bar{m}}^2 e^{-2 \eta^* y_{j:\bar{m}}^2}
\end{pmatrix}.
\end{eqnarray*}

\begin{corollary} \label{cor4}
Let $\bar{m} \geq 2$. 
Assume that 
$N h_{N:T}^2 \rightarrow 0$. 
Then, as $N \to \infty$,
\begin{eqnarray*}
\sqrt{N}(\hat{\sigma}_0^2 - (\sigma_0^*)^2,  {\hat{\eta}}-\eta^* )
\stackrel{d}{\longrightarrow } N(0,(\sigma_0^*)^4 \Gamma \pi V_{\bar{m}}(\zeta^*)^{-1}U_{\bar{m}}(\zeta^*)V_{\bar{m}}(\zeta^*)^{-1}).
\end{eqnarray*}
\end{corollary}

\subsection{Estimation of four parameters}
Next we consider the estimation of $\sigma,\theta_1,\theta_2,$ and $\theta_0$ using $\hat{\sigma}_0^2 $ and $\hat{\eta}$.
Let $k \in \mathbb{N}$.
As an approximation of $x_k(t)$, we consider 
\begin{eqnarray*}
\bar{x}_k(t_{i:N:T}) = \frac{1}{M}\sum^{M}_{j = 1} 
X_{t_{i:N:T}}(y_{j:M}) \sqrt{2}\sin (\pi k y_{j:M}) 
\exp \left( \frac{\hat{\eta} y_{j:M}}{2} \right), 
\quad i = 1,...,N. 
\end{eqnarray*}
Let $\bar{N}_2 \leq N$,
$\delta_{\bar{N}_2:T} = \left[  \frac{N}{\bar{N}_2} \right] h_{N:T} = \left[\frac{N}{\bar{N}_2} \right] \frac{T}{N}$
and $s_{i:\bar{N}_2:T} = i \delta_{\bar{N}_2:T}$ for $i = 1,...,\bar{N}_2$.
{The quasi log-likelihood function with 
the thinned data based on the approximate coordinate process }
${\bf \bar{x}}_k =\{\bar{x}_k(s_{i:\bar{N}_2:T})\}_{i=1,\ldots, \bar{N}_2}$ 
is defined as
\begin{eqnarray*}
l_{\bar{N}_2}(\lambda_k,\sigma^2 \ | \ {\bf \bar{x}}_k  ) 
&=&
 -\sum^{\bar{N}_2}_{i=1}\left\{ \frac{1}{2} \log 
 \left( \frac{\sigma^2(1 - \exp (-2\lambda_k \delta_{\bar{N}_2:T} ))}{2\lambda_k} \right) \right. \\
& & \left. 
+ \frac{(\bar{x}_k( s_{i:\bar{N}_2:T} ) - \exp (-\lambda_k \delta_{\bar{N}_2:T} )
\bar{x}_k(s_{i-1:\bar{N}_2:T}))^2}{\frac{2\sigma^2(1-\exp (-2\lambda_k \delta_{\bar{N}_2:T}))}{2\lambda_k}} \right\}.
\end{eqnarray*}
The quasi maximum likelihood estimator of $(\lambda_k, \sigma^2)$ is defined as 
\begin{eqnarray*}
(\hat{\lambda}_k, \hat{\sigma}^2_k) = \argsup_{\lambda,   \sigma^2} l_{\bar{N}_2}(\lambda,\sigma^2  \ | \ {\bf \bar{x}}_k  ). 
\end{eqnarray*}
Set $k = 1$, $\hat{\lambda} = \hat{\lambda}_1,  \hat{\sigma}^2 =  \hat{\sigma}^2_1$.
Note that $\lambda_1^* = - \theta_0^* + \frac{(\theta_1^*)^2}{4 \theta_2^*} + \pi^2 \theta_2^*$.
The estimators of $\theta_1$, $\theta_2$ are defined as 
\begin{eqnarray*}
\hat{\theta}_2 = \left( \frac{\hat{\sigma}^2}{\hat{\sigma}_0^2} \right)^2 , \quad  \hat{\theta}_1 = \hat{\eta} \hat{\theta_2}.
\end{eqnarray*}
The estimator of $\theta_0$ is defined as 
\begin{eqnarray*}
\hat{\theta}_0 = -\hat{\lambda} + \frac{\hat{\theta}_1^2}{4\hat{\theta}_2} + \pi^2\hat{\theta}_2. 
\end{eqnarray*}

\begin{theorem} \label{thm4}
Assume the same conditions as Theorem \ref{thm3}. 
Moreover, assume that 
{
$\frac{\bar{N}_2^{\frac{5}{2}}}{T^{\frac{3}{2}} N \bar{m}} \to 0$
and 
$\frac{\bar{N}_2^{3}}{T^2 M^{1-\rho_1}} \to 0$}
for
$\rho_1 \in (0,1)$.
Then, 
as $\bar{N}_2 \to \infty$, \\
$$
\begin{pmatrix}
\sqrt{N_2}((\hat{\sigma})^2 - (\sigma^{*})^2) \\
\sqrt{N_2}(\hat{\theta}_2 - \theta_{2}^*) \\
\sqrt{N_2}(\hat{\theta}_1 - \theta_{1}^*) \\
\sqrt{T}(\hat{\theta}_0 - \theta_{0}^*)
\end{pmatrix}
\stackrel{d}{\longrightarrow } N
\left(
\begin{pmatrix}
0 \\
0 \\
0 \\
0 
\end{pmatrix}
,
\begin{pmatrix}
2(\sigma^*)^4 &4\theta_2^*(\sigma^*)^2 & 4\theta_1^*(\sigma^*)^2 & 0\\
4\theta_2^*(\sigma^*)^2 &8(\theta_2^*)^2 & 8\theta_1^*\theta_2^* & 0\\
 4\theta_1^*(\sigma^*)^2 &  8\theta_1^*\theta_2^*& 8(\theta_1^*)^2 & 0 \\
0 & 0 & 0 & 2\lambda^*_1
\end{pmatrix}
\right).
$$
\end{theorem}

\begin{corollary}
Let $\bar{m} \geq 2$. 
Assume the same conditions as Corollary \ref{cor4}. 
Moreover, assume that 
{
$\frac{\bar{N}_2^{\frac{5}{2}}}{T^{\frac{3}{2}} N } \to 0$
and 
$\frac{\bar{N}_2^{3}}{T^2 M^{1-\rho_1}} \to 0$}
for
$\rho_1 \in (0,1)$.
Then, 
as $\bar{N}_2 \to \infty$, \\
$$
\begin{pmatrix}
\sqrt{N_2}((\hat{\sigma})^2 - (\sigma^{*})^2) \\
\sqrt{N_2}(\hat{\theta}_2 - \theta_{2}^*) \\
\sqrt{N_2}(\hat{\theta}_1 - \theta_{1}^*) \\
\sqrt{T}(\hat{\theta}_0 - \theta_{0}^*)
\end{pmatrix}
\stackrel{d}{\longrightarrow } N
\left(
\begin{pmatrix}
0 \\
0 \\
0 \\
0 
\end{pmatrix}
,
\begin{pmatrix}
2(\sigma^*)^4 &4\theta_2^*(\sigma^*)^2 & 4\theta_1^*(\sigma^*)^2 & 0\\
4\theta_2^*(\sigma^*)^2 &8(\theta_2^*)^2 & 8\theta_1^*\theta_2^* & 0\\
 4\theta_1^*(\sigma^*)^2 &  8\theta_1^*\theta_2^*& 8(\theta_1^*)^2 & 0 \\
0 & 0 & 0 & 2\lambda^*_1
\end{pmatrix}
\right).
$$
\end{corollary}


\section{Examples and simulation results}

In the same way as Bibinger and Trabs (2017), the numerical solution of the SPDE (\ref{spde0}) is generated by
\begin{eqnarray}
\tilde{X}_{t_{i:N}}(y_{j:M}) = \sum^{K}_{k = 1}x_k(t_{i:N})e_k(y_{j:M}), \quad i = 1, ..., N, j = 1, ..., M,
\label{appro-SPDE}
\end{eqnarray}
where
\begin{eqnarray*}
x_k(t_{i:N}) = \exp \left( -\lambda_k \frac{T}{N} \right) x_k(t_{i-1:N}) + \sqrt{\frac{\sigma^2(1-\exp(-2 \lambda_k \frac{T}{N}))}{2 \lambda_k}}N(0,1) ,\quad i = 1, ..., N.
\end{eqnarray*}
When the tuning parameter $K$ is not large enough, 
the approximation (\ref{appro-SPDE}) 
does not work well and the estimators in Theorem \ref{thm1} and Theorem \ref{thm3} have considerable biases.
Therefore, the tuning parameter $K$ needs $10^5$ or more.

\subsection{The case that $T$ is fixed}

In this subsection, we consider the linear parabolic SPDE (\ref{spde1}) where 
the true value $(\theta_0, \theta_1, \theta_2, \sigma) = (0,0.5,0.1,1)$.
We set that $N = 10^4$, $M = 10^4$, $K = 10^5$, $T=1$.
{
When $N=M=10^4$, the size of data ${\bf X}_{N,M}$ is about {1 GB}.
We used R language to compute {the estimators of Theorems 1 and  2}.
The personal computer with Intel Gold 6128 (3.40GHz) was used for this simulation.
It takes about 4 hours to generate one sample path
{of the SPDE (\ref{spde1}) with $N = 10^4$, $M = 10^4$, $K = 10^5$, $T=1$}.
}

\begin{figure}[h] 
\begin{center}
\includegraphics[width=9cm]{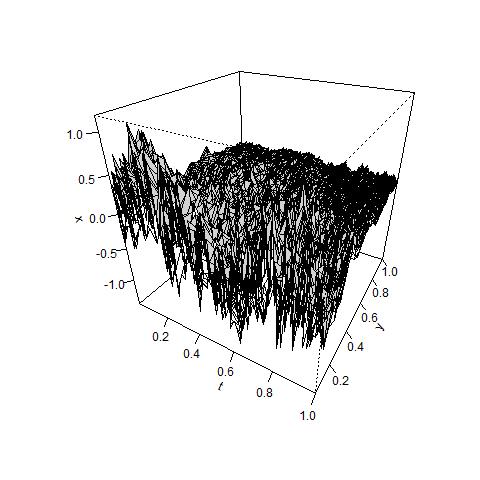} 
\caption{Sample path with $\theta = (0,0.5,0.1,1)$  \label{fig1}}
\end{center}
\end{figure}

Figure \ref{fig1} is {a} sample path of $X_t(y)$ for $(t,y) \in [0,1]\times [0,1]$
when $(\theta_0, \theta_1, \theta_2, \sigma) = (0,0.5,0.1,1)$.




\begin{figure}[h] 
\begin{center}
\includegraphics[width=5cm,pagebox=cropbox,clip]{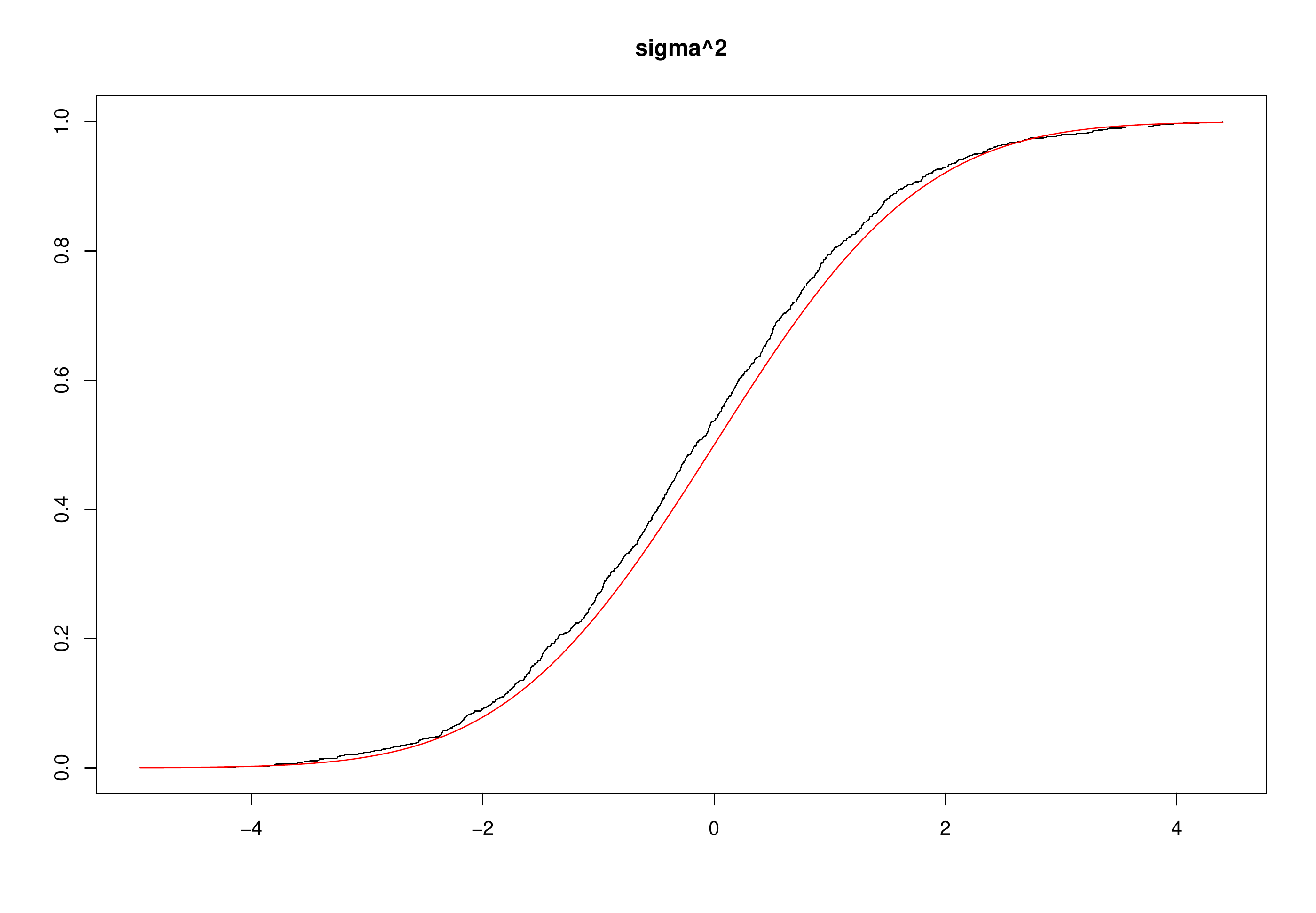}
\includegraphics[width=5cm,pagebox=cropbox,clip]{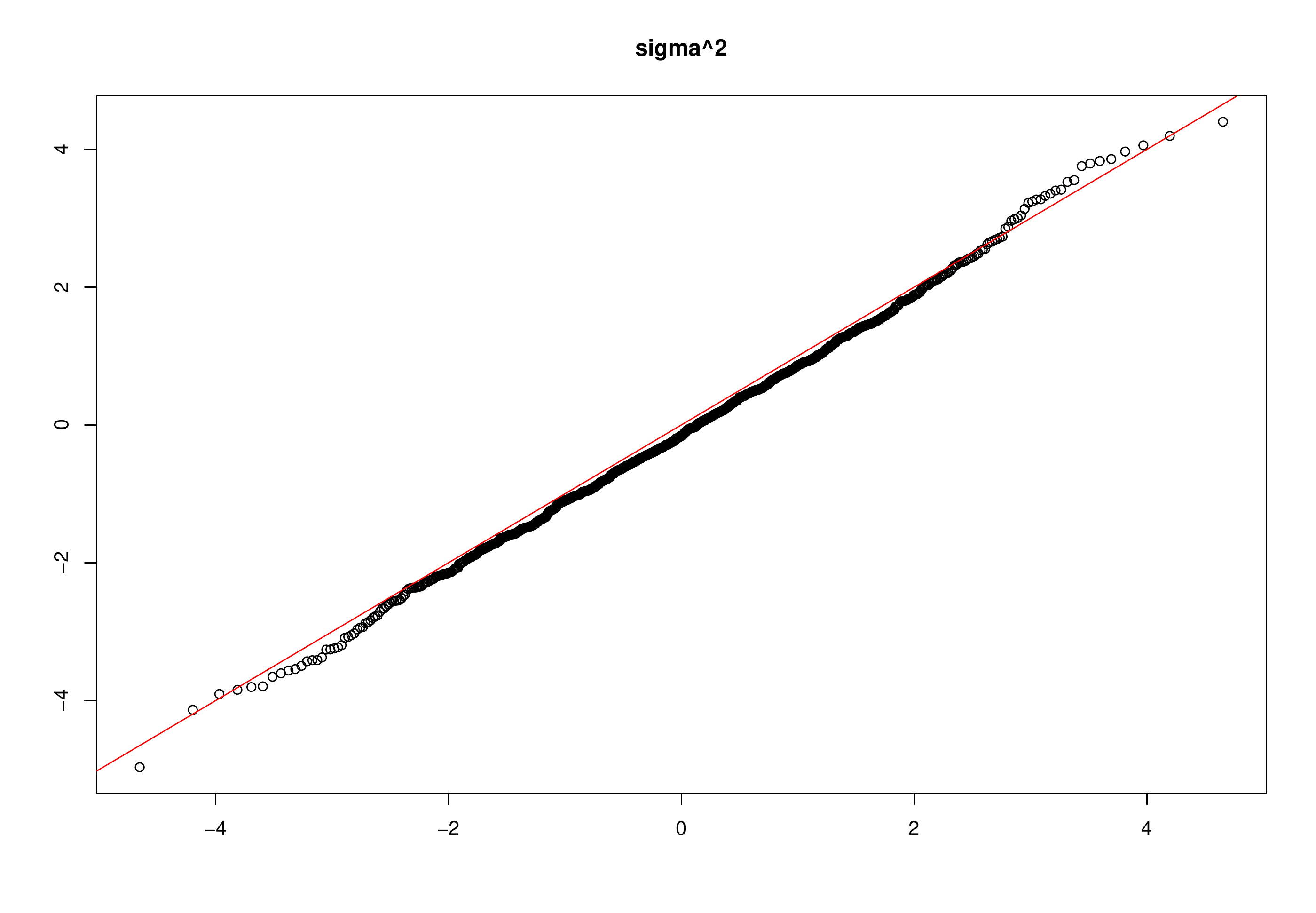}
\includegraphics[width=5cm,pagebox=cropbox,clip]{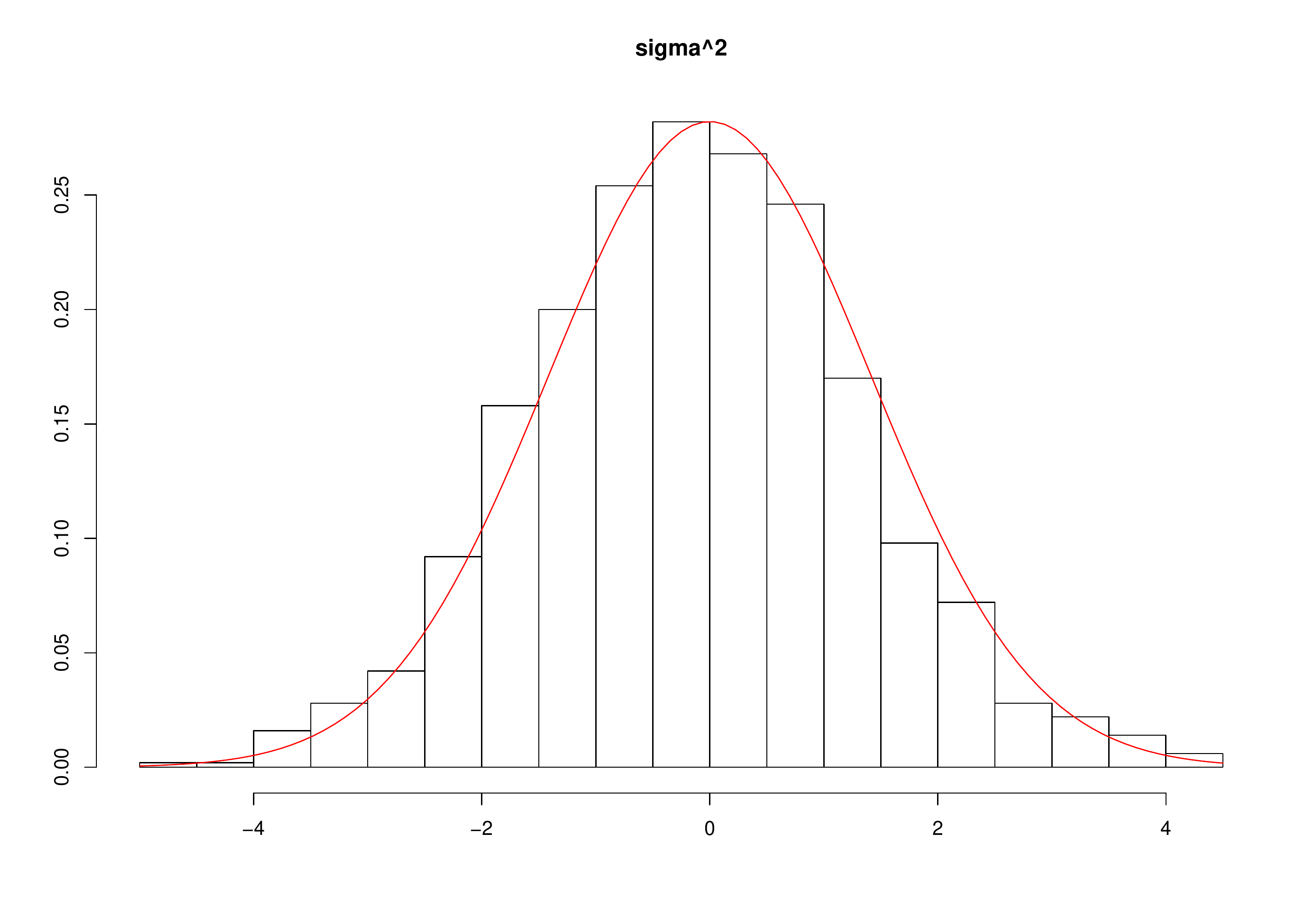}
\caption{Simulation results of $\check{\sigma}^2$ \label{fig2}}
\end{center}
\end{figure}

\begin{figure}[h]
\begin{center}
\includegraphics[width=5cm,pagebox=cropbox,clip]{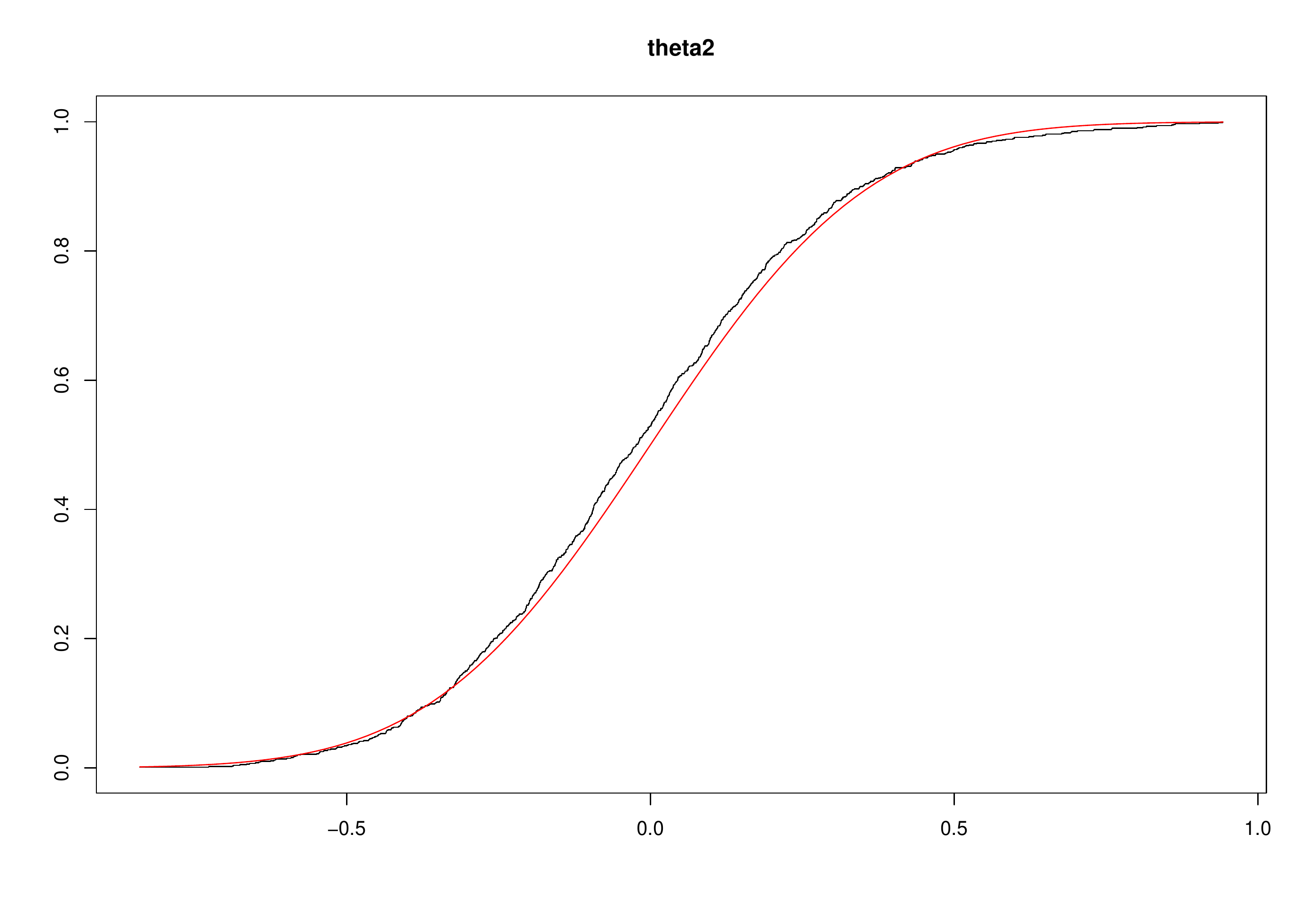}
\includegraphics[width=5cm,pagebox=cropbox,clip]{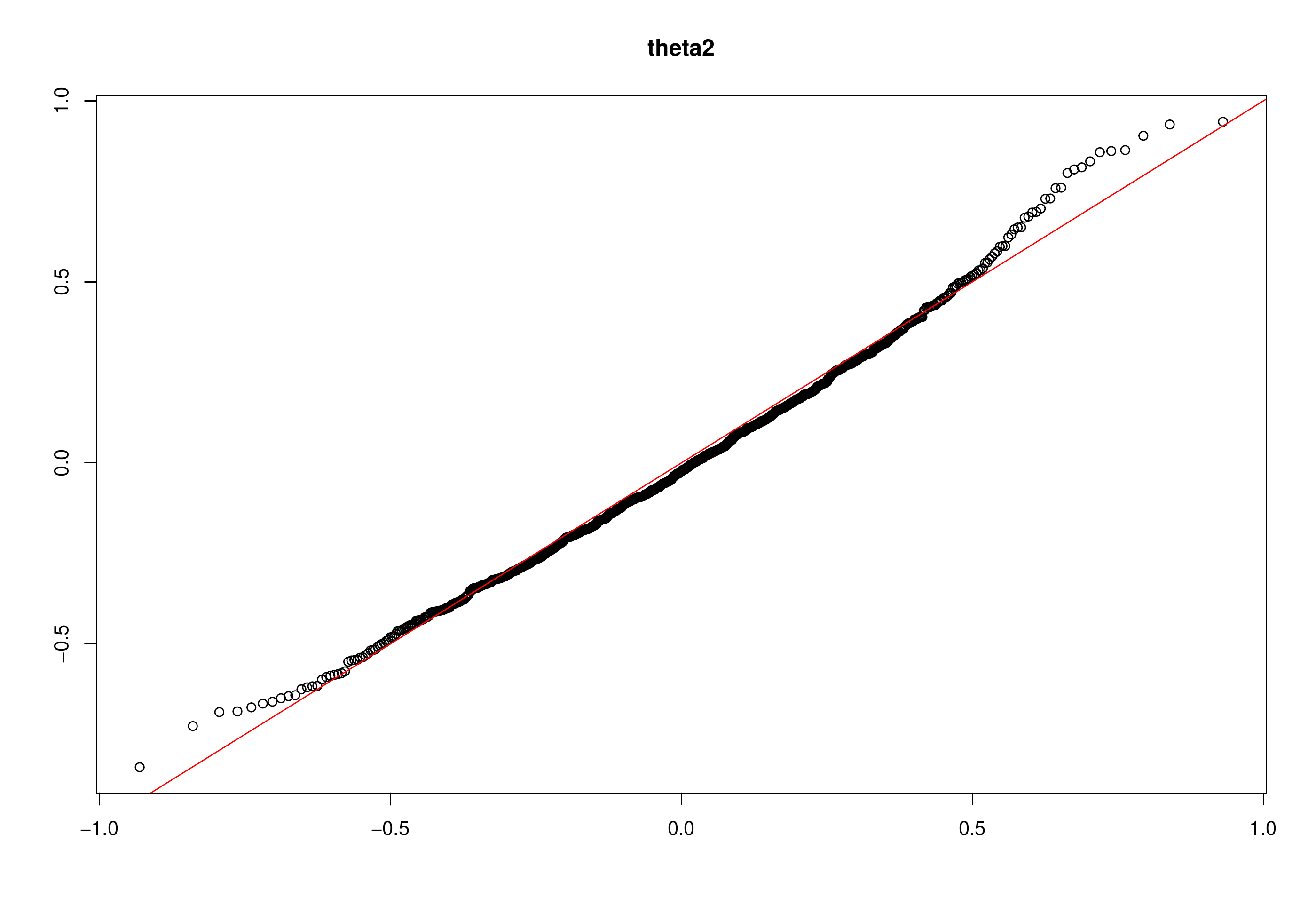}
\includegraphics[width=5cm,pagebox=cropbox,clip]{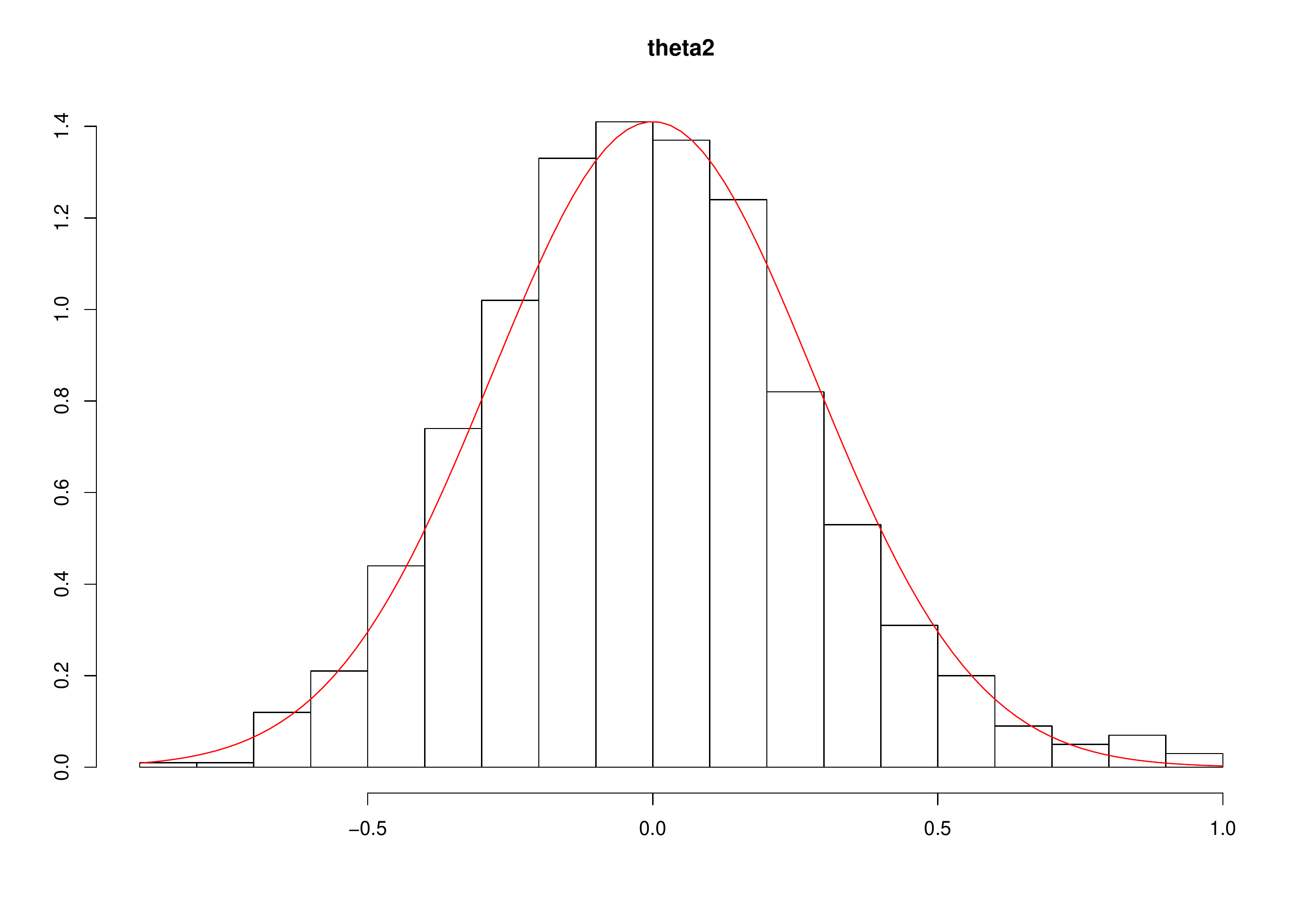}\\
\caption{Simulation results of $\check{\theta}_2$  \label{fig3}}
\end{center}
\end{figure}

\begin{figure}[t] 
\begin{center}
\includegraphics[width=5cm,pagebox=cropbox,clip]{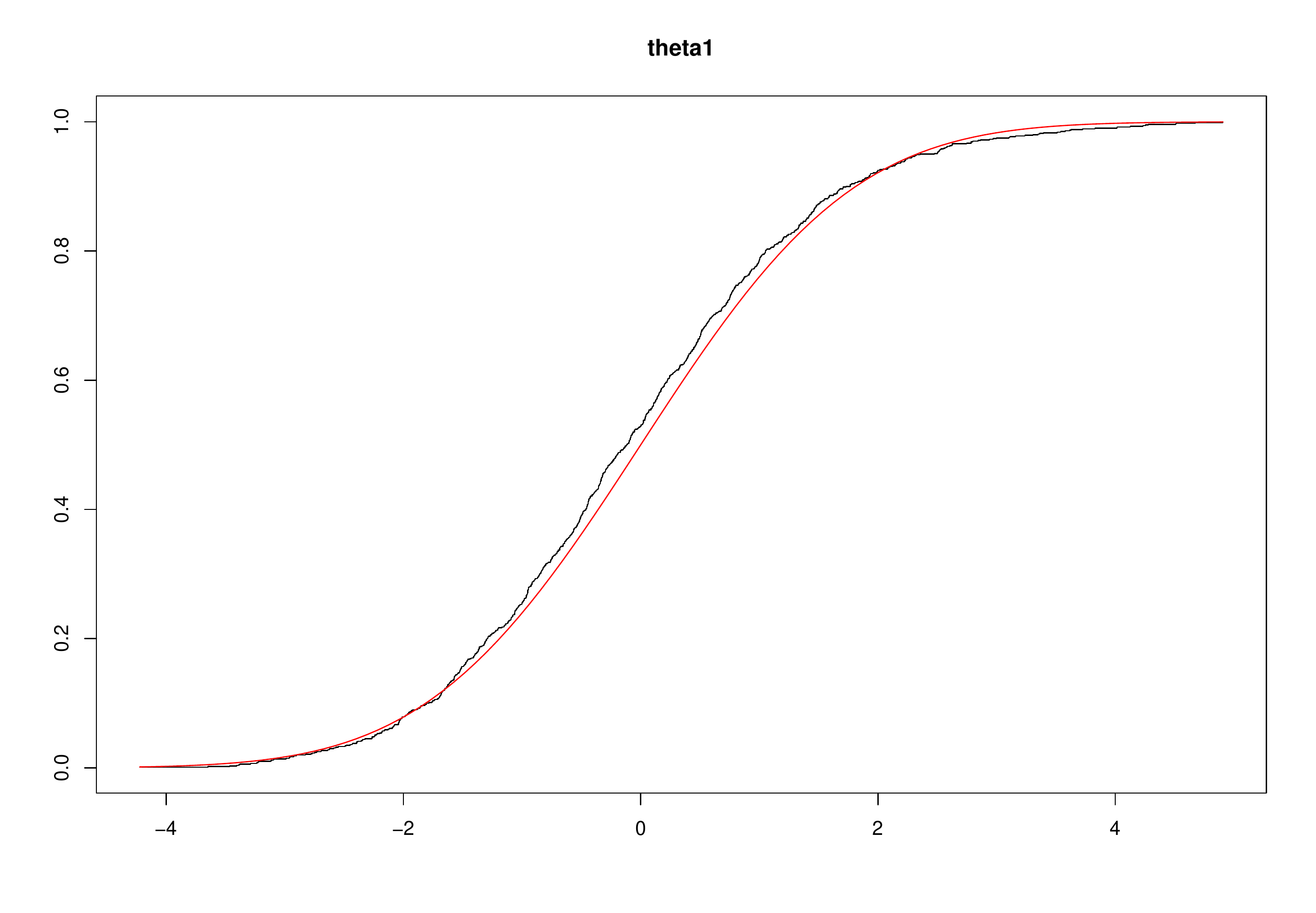}
\includegraphics[width=5cm,pagebox=cropbox,clip]{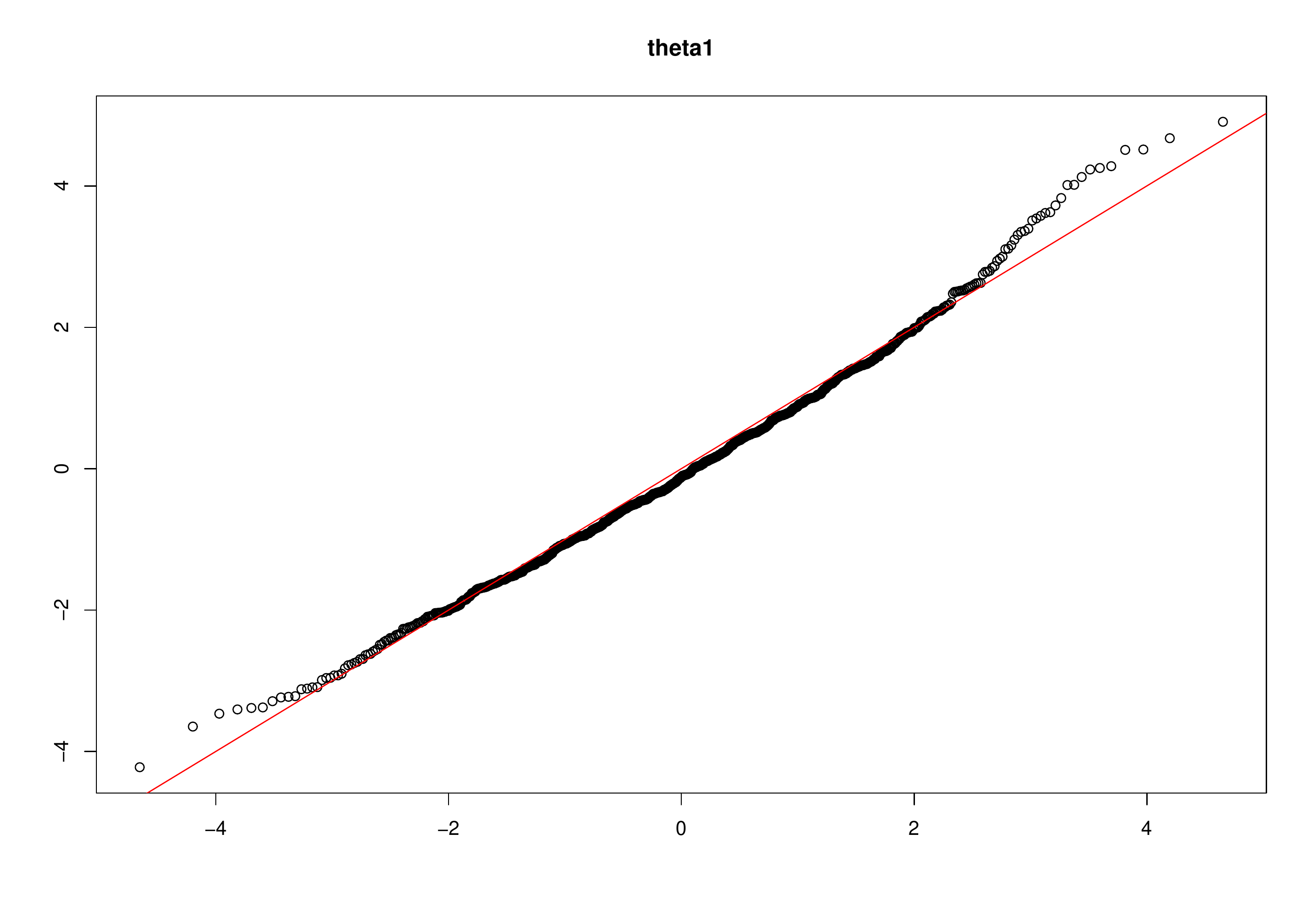}
\includegraphics[width=5cm,pagebox=cropbox,clip]{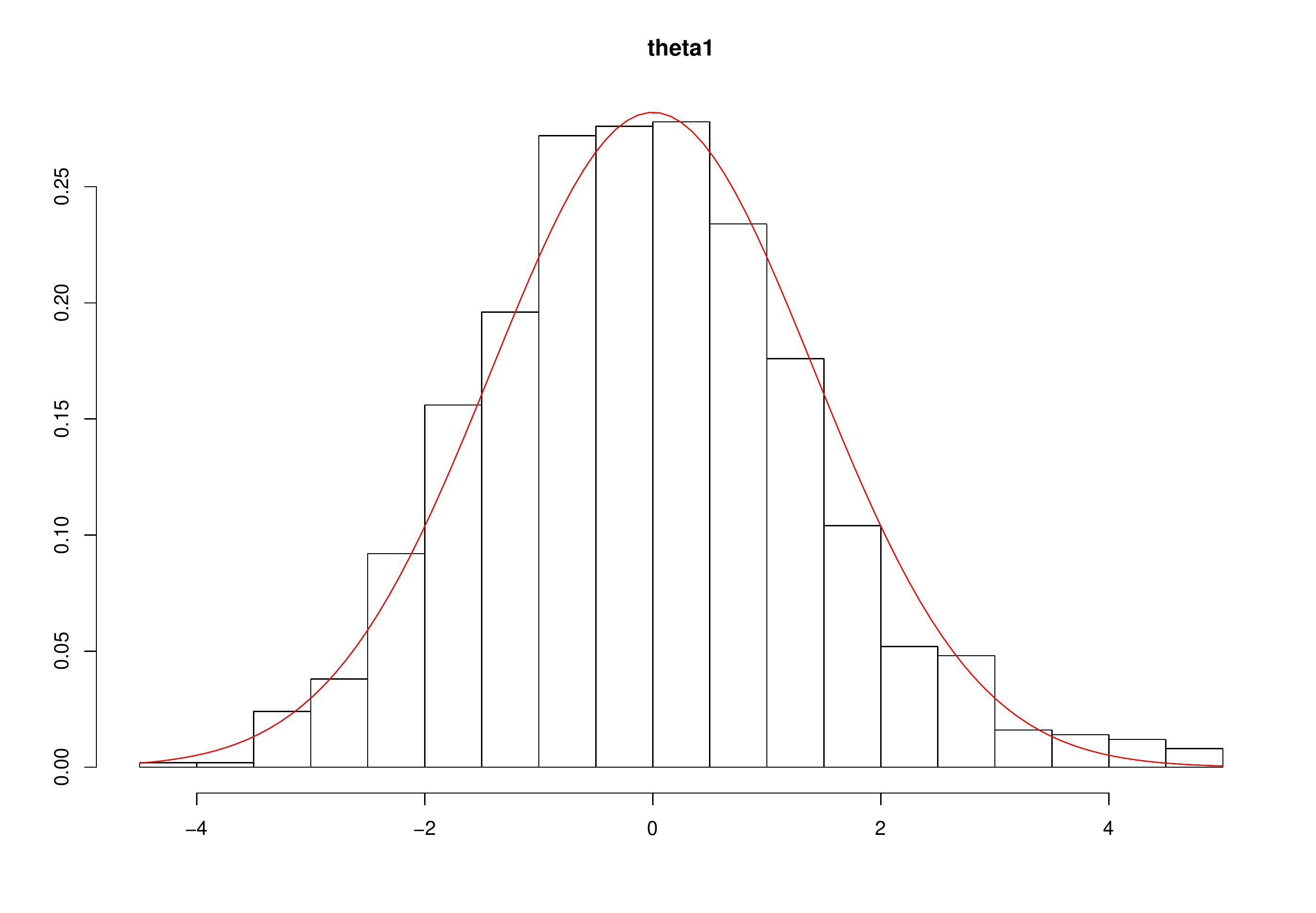}\\
\caption{Simulation results of $\check{\theta}_1$ \label{fig4}}
\end{center}
\end{figure}

Figures \ref{fig2}-\ref{fig4} are the simulation results of $\check{\sigma}^2$, $\check{\theta}_2$  and $\check{\theta}_1$
in Theorem \ref{thm2}
with $(N, m, N_2) = (10^4, 99, 333)$.
Note that $\frac{N_2^{3/2}}{m N} \approx 0.006$, $\frac{N_2^{3/2}}{M^{0.99}} \approx 0.66$
and the number of iteration is $1000$.
The left side of Figure \ref{fig2} is the plot of the empirical distribution function of 
$\sqrt{N_2}(\check{\sigma}^2-(\sigma^*)^2)$ (black line) and the distribution function of $N(0,2(\sigma^*)^4)$ (red line).
The center of Figure \ref{fig2} is the Q-Q plot of $\sqrt{N_2}(\check{\sigma}^2-(\sigma^*)^2)$ and $N(0,2(\sigma^*)^4)$.
The right side of Figure \ref{fig2} is the plot of the histogram of $\sqrt{N_2}(\check{\sigma}^2-(\sigma^*)^2)$ 
and the density function of $N(0,2(\sigma^*)^4)$ (red line).
Figures \ref{fig3} and \ref{fig4} are the plots of the empirical distribution functions, the Q-Q plots 
and the histograms of $\sqrt{N_2}(\check{\theta}_2 - \theta_{2}^*)$ 
and $\sqrt{N_2}(\check{\theta}_1 - \theta_{1}^*)$, respectively.
From Figures \ref{fig2}-\ref{fig4}, 
we can see that the distributions of the estimators in Theorem \ref{thm2} almost 
correspond with the asymptotic distribution and these estimates have good performance.

\vspace{0.5cm}

\clearpage

\subsection{The case that $T$ is large}

In this subsection, we deal with the linear parabolic SPDE (\ref{spde2}) where 
the true value $(\theta_0, \theta_1, \theta_2, \sigma) = (0,0.2,0.2,1)$
and we set $N = 10^5$, $M = 10^5$, $K = 10^5$, $T = 100$.
{
When $N=M=10^5$, 
{${\bf X}_{N,M}$ is about 80 GB,
which is too large for R to handle.}
{We use the Python Programming Language to compute 
the estimators of  Theorems 3 and 4}.
It takes about 30 hours to generate one sample path
{of the SPDE (\ref{spde2}) with $N = 10^5$, $M = 10^5$, $K = 10^5$, $T = 100$}.
}

\begin{figure}[h]
\begin{center}
\includegraphics[width=9cm]
{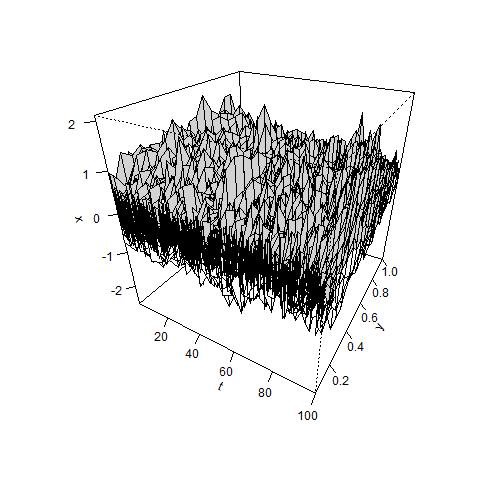}
\caption{Sample path with $\theta = 
(0,0.2,0.2,1)$ 
\label{fig5}}
\end{center}
\end{figure}

Figure \ref{fig5} is {a} sample path of $X_t(y)$ for $(t,y) \in [0,100]\times [0,1]$
when $(\theta_0, \theta_1, \theta_2, \sigma) = 
(0,0.2,0.2,1)$.

\begin{figure}[h]
  \begin{center}
\includegraphics[width=6cm,height = 5cm,angle=0,pagebox=cropbox,clip]
{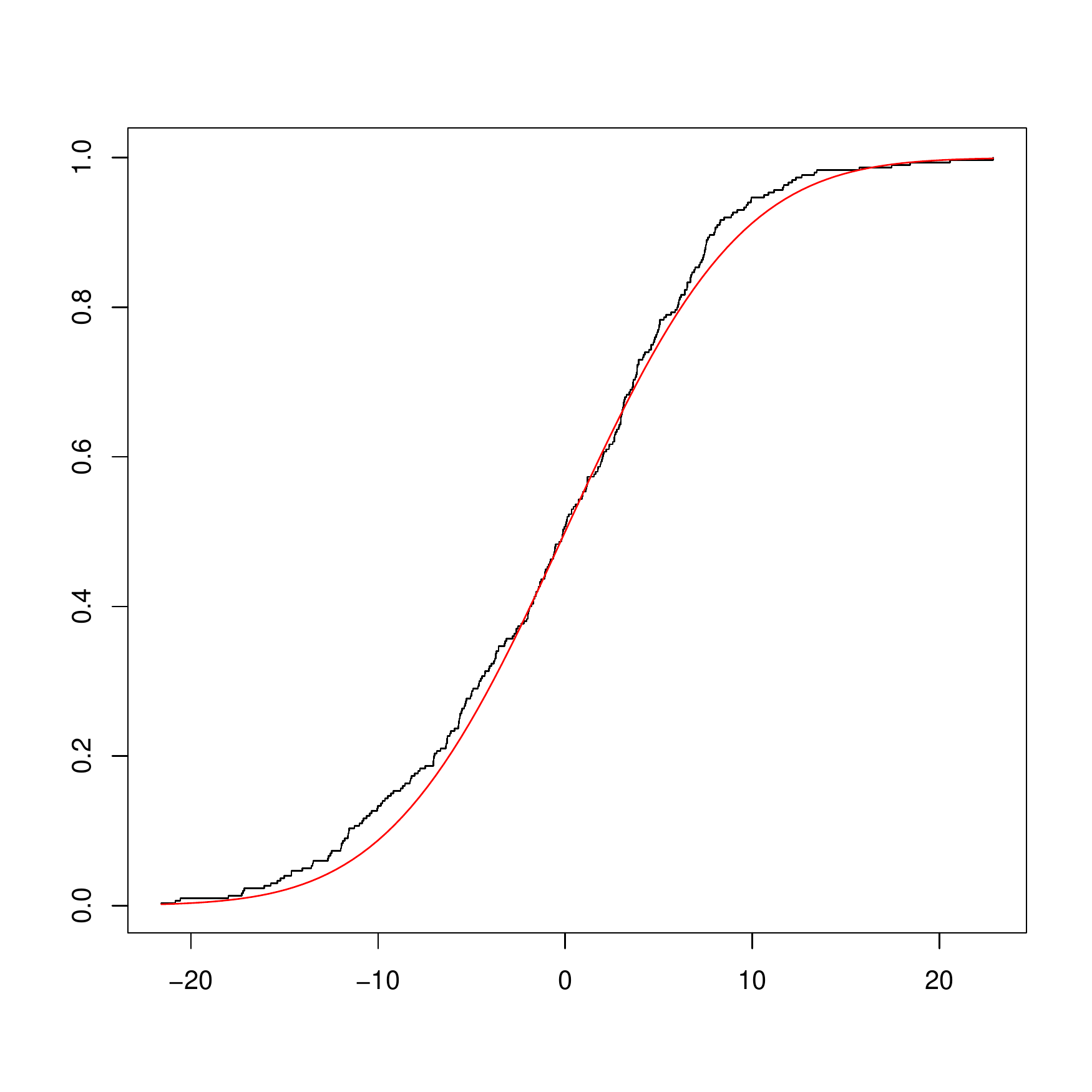}
\includegraphics[width=6cm,height = 5cm,angle=0,pagebox=cropbox,clip]
{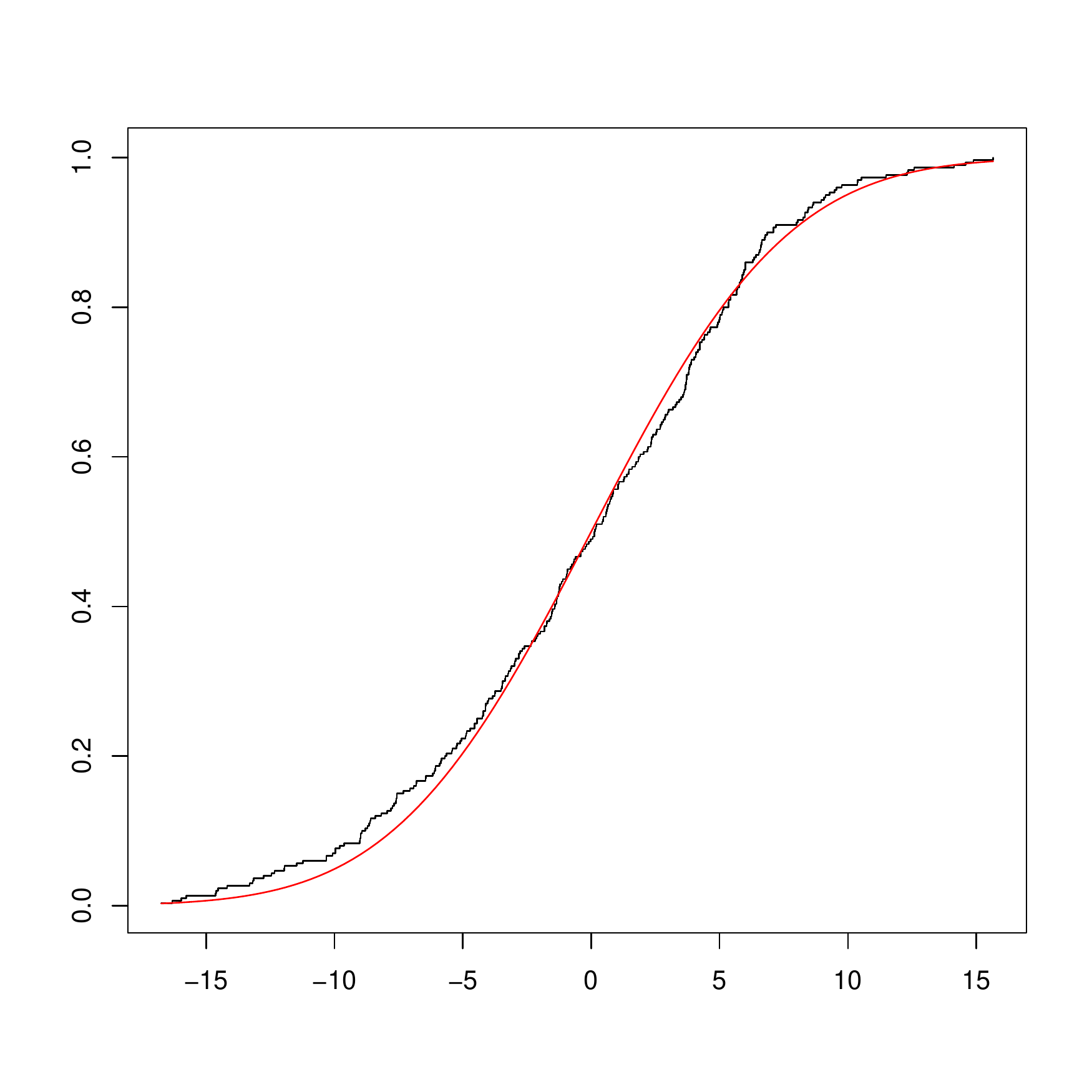}

  \caption{Simulation results of $\hat{\sigma}_0^2$ (left) and $\hat{\eta}$ (right) in Theorem \ref{thm3} \label{fig6}}
    \end{center}
\end{figure}

We set $\bar{m} =  24$ and the number of iteration is $300$.
Figure \ref{fig6} is  the plots of the empirical distribution functions of 
$\sqrt{\bar{m} N}(\hat{\sigma}_0^2 - (\sigma_0^*)^2)$ 
and $\sqrt{\bar{m} N}(\bar{\eta}-\eta^* )$ in Theorem \ref{thm3}, respectively.
The red lines in Figure \ref{fig6} are the plots of the distribution functions of the corresponding asymptotic distributions.
From Figure \ref{fig6}, we can see that the empirical distributions of the estimators in Theorem \ref{thm3} almost fit 
the asymptotic distribution and these estimates have good behavior.
This result indicates that 
the estimators work well when $N = 10^5$ and $\bar{m} = 24$.

\begin{figure}[h]
\begin{center}
\includegraphics[width=5cm,pagebox=cropbox,clip]{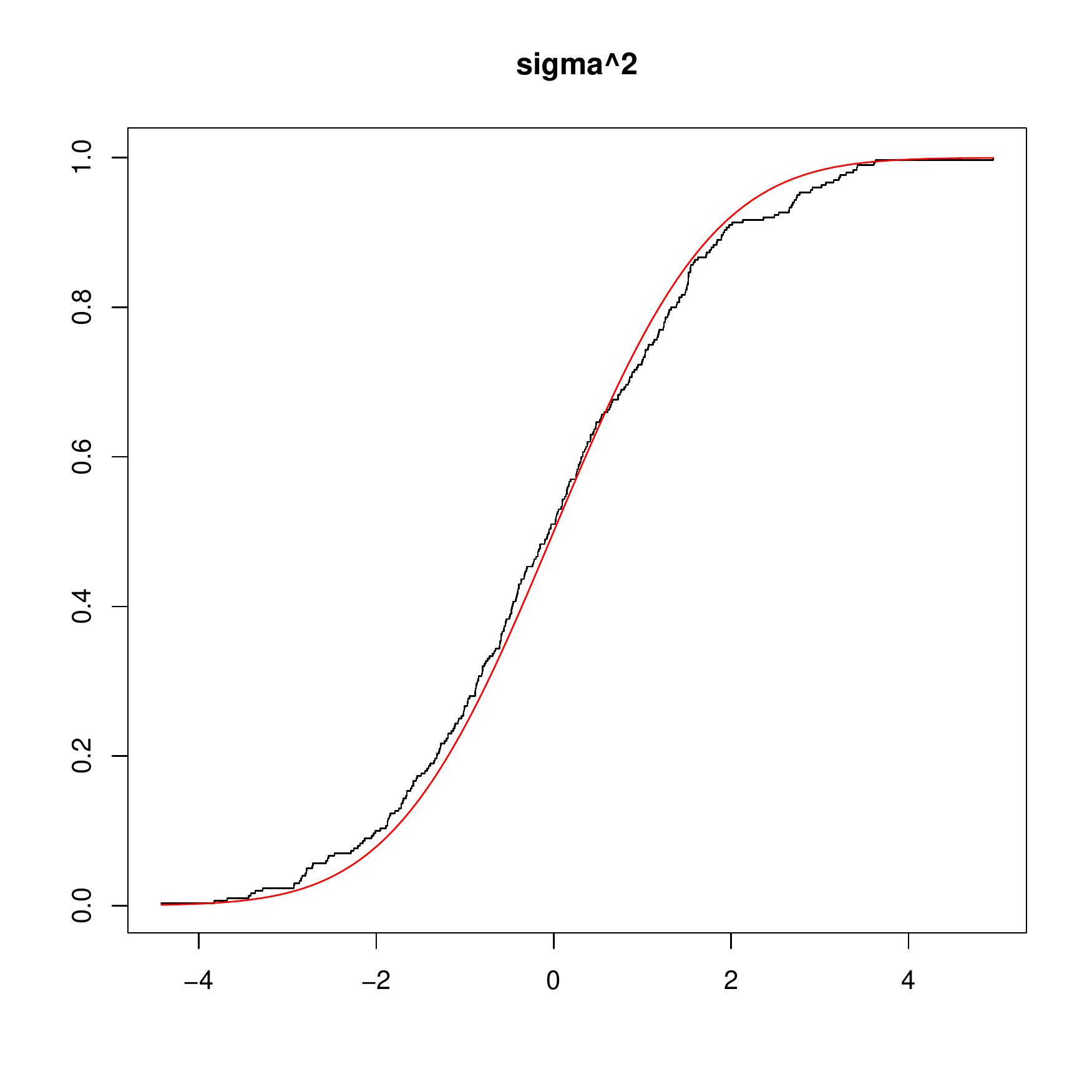}
\includegraphics[width=5cm,pagebox=cropbox,clip]{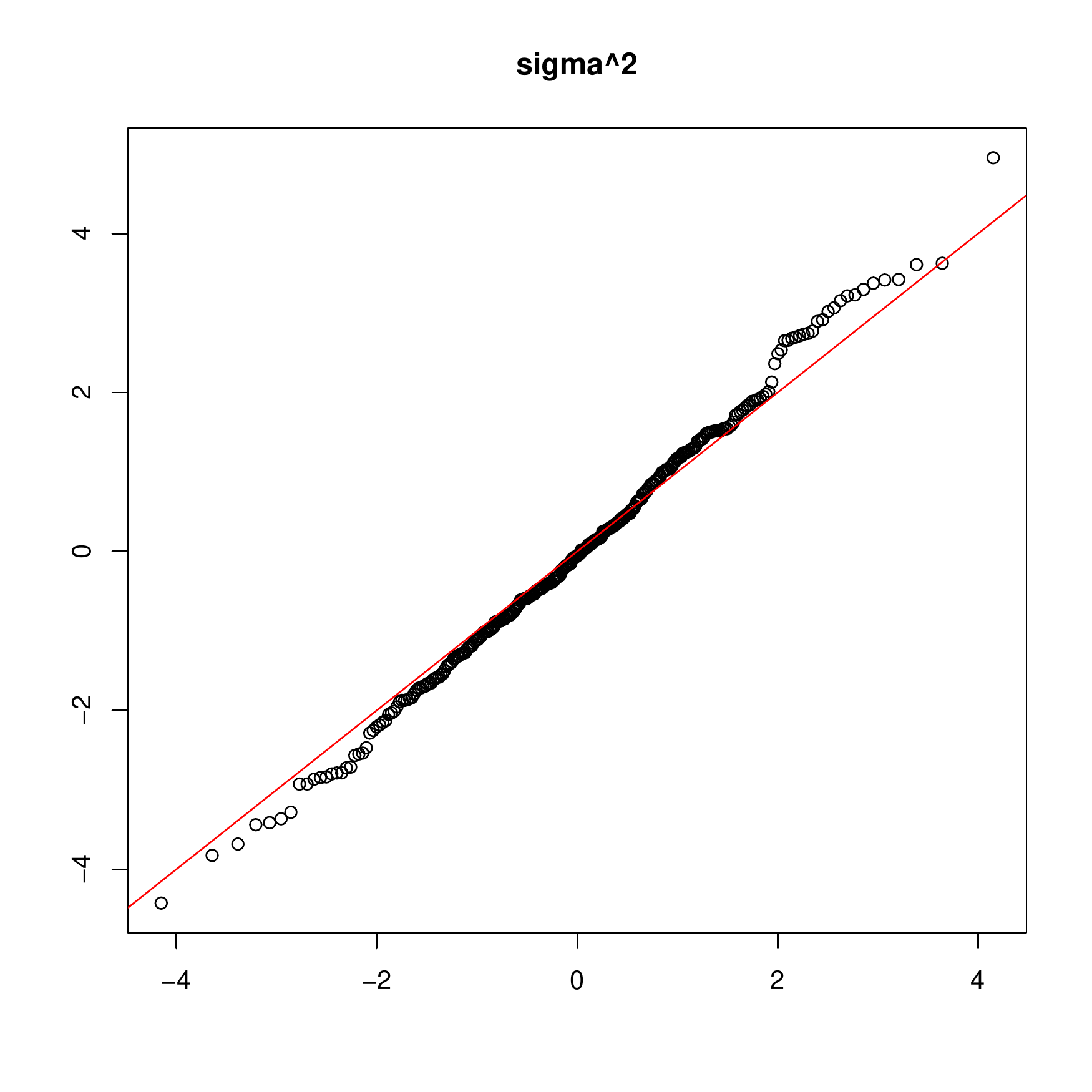}
\includegraphics[width=5cm,pagebox=cropbox,clip]{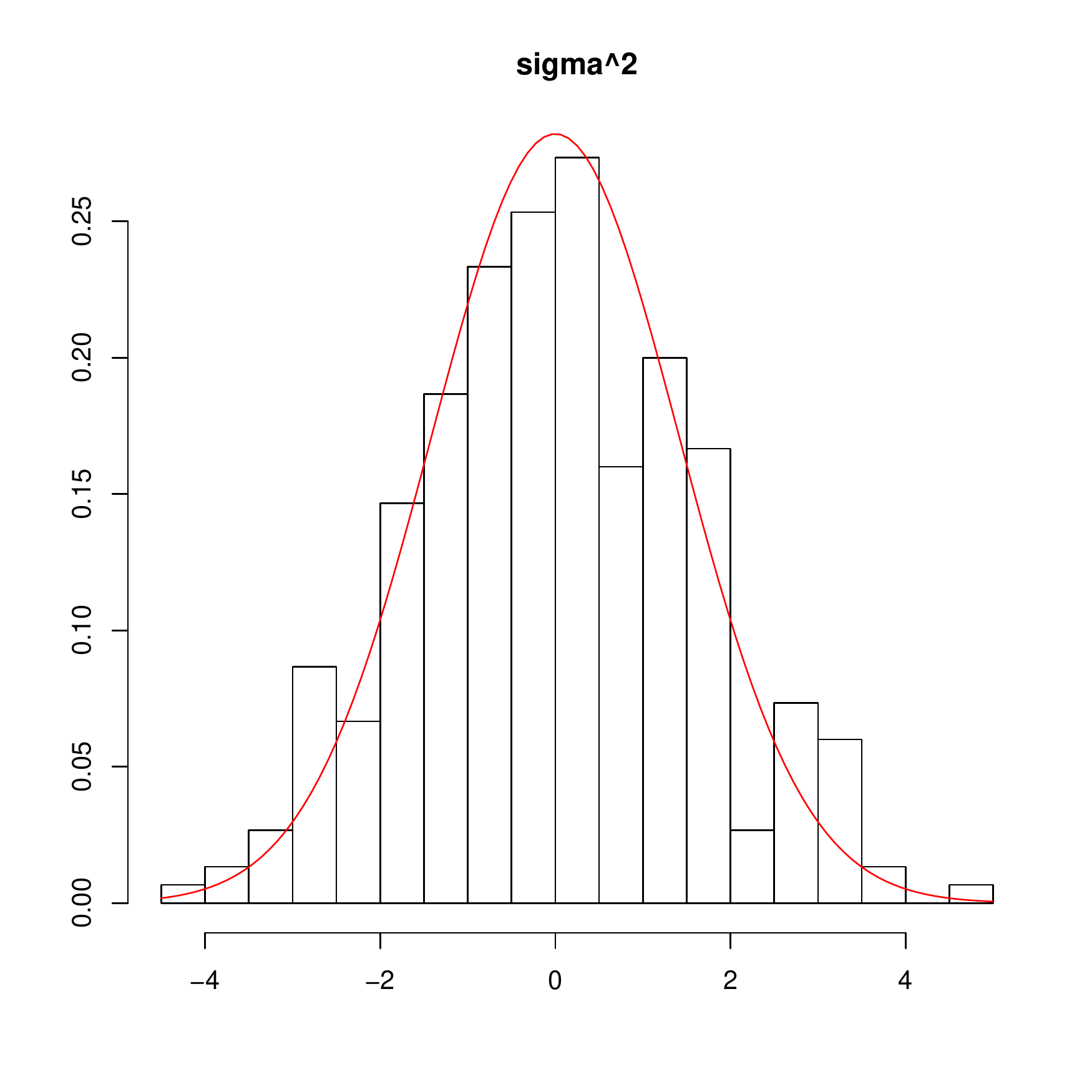}\\
\caption{Simulation results of $\hat{\sigma}^2$ \label{fig7}}
\end{center}
\end{figure}

\begin{figure}[h]
\begin{center}
\includegraphics[width=5cm,pagebox=cropbox,clip]
{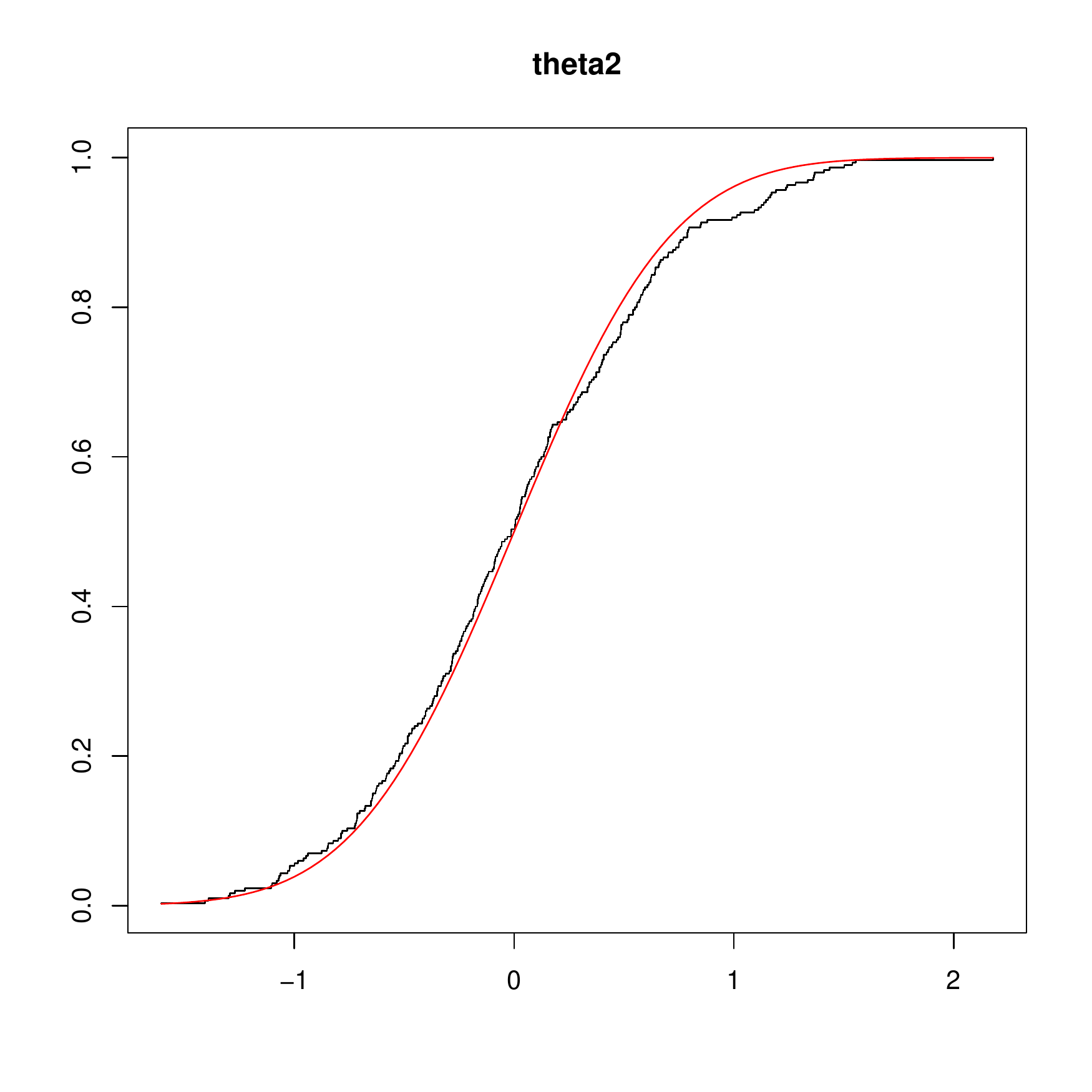}
\includegraphics[width=5cm,pagebox=cropbox,clip]{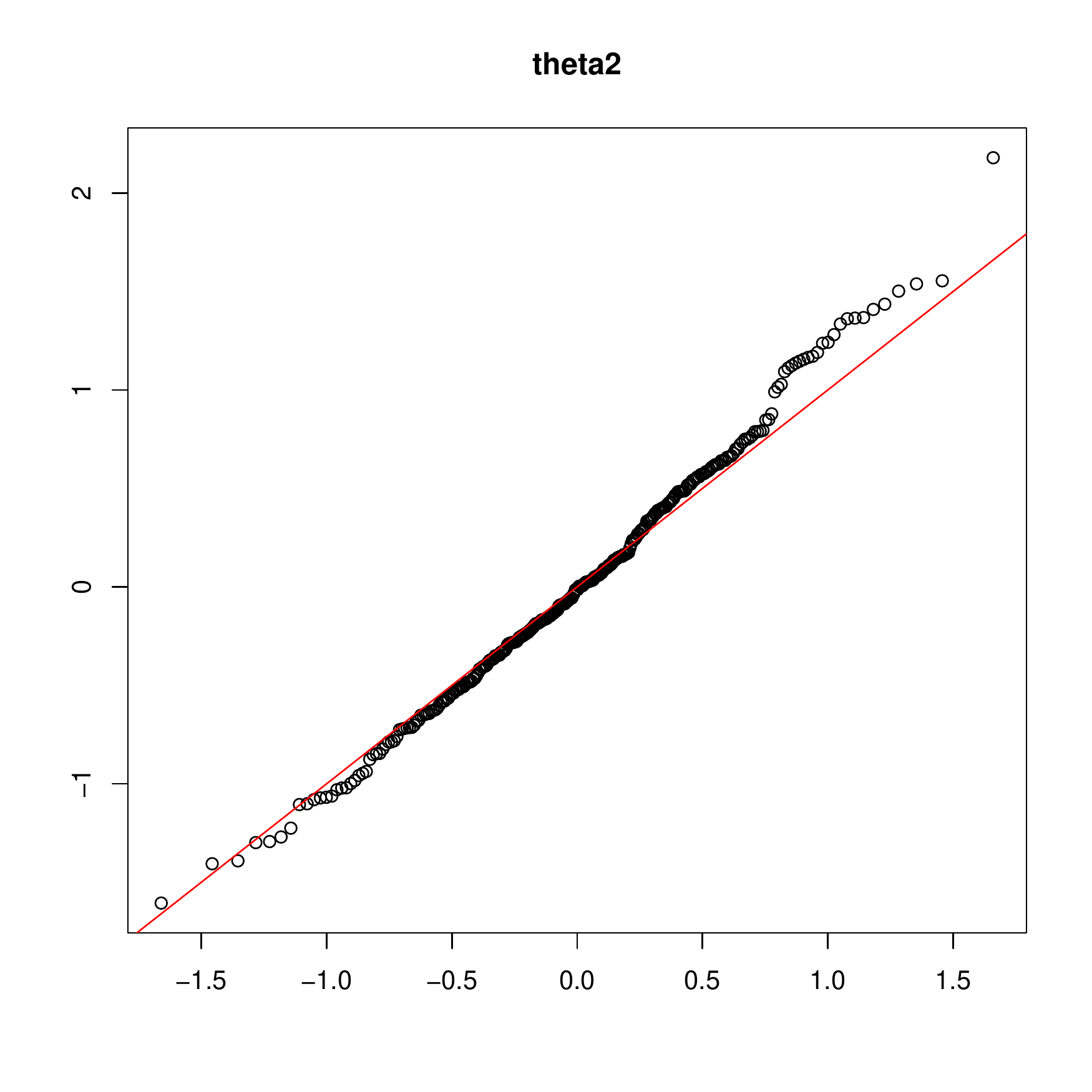}
\includegraphics[width=5cm,pagebox=cropbox,clip]{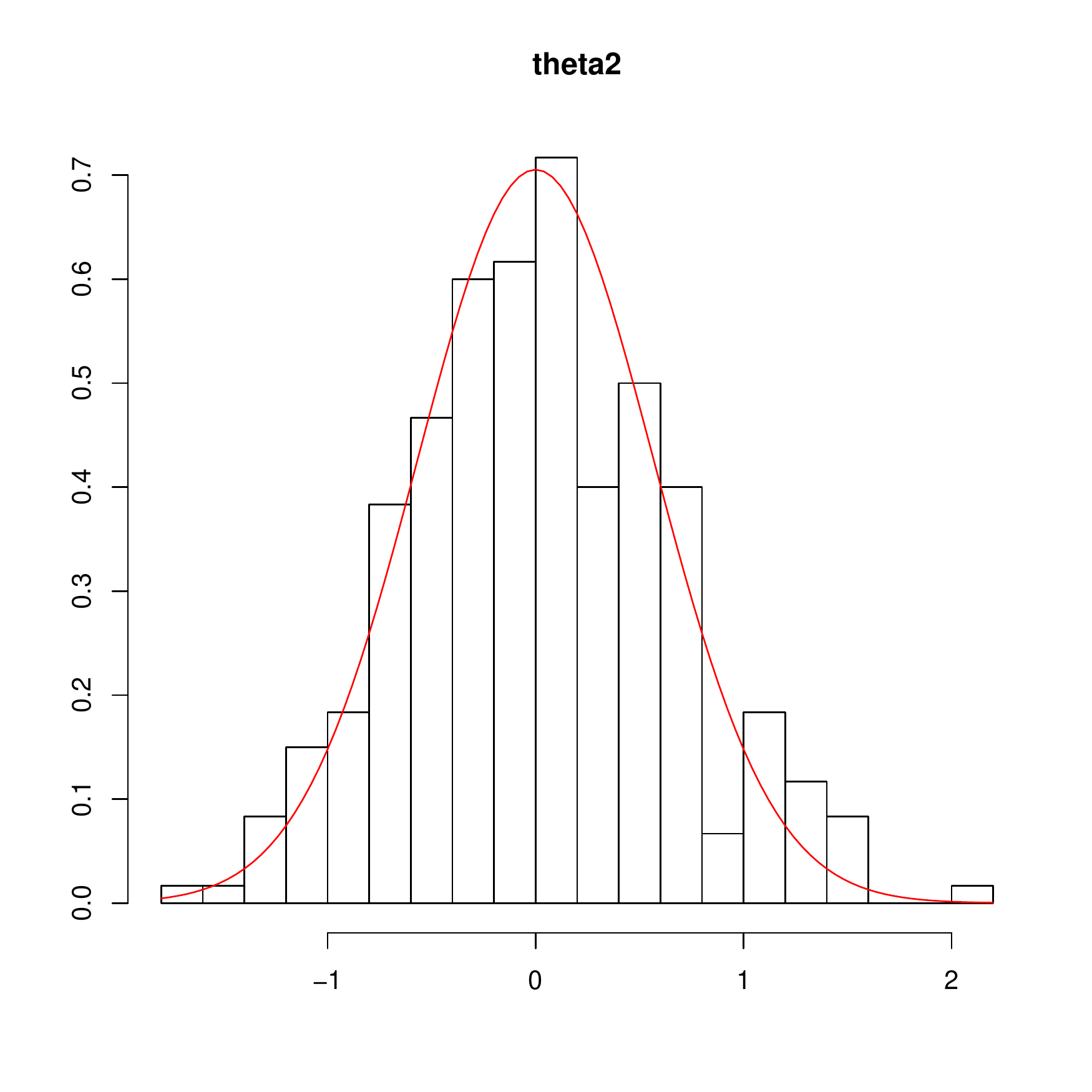}\\
\caption{Simulation results of $\hat{\theta}_2$ \label{fig8}}
\end{center}
\end{figure}

\begin{figure}[h]
\begin{center}
\includegraphics[width=5cm,pagebox=cropbox,clip]
{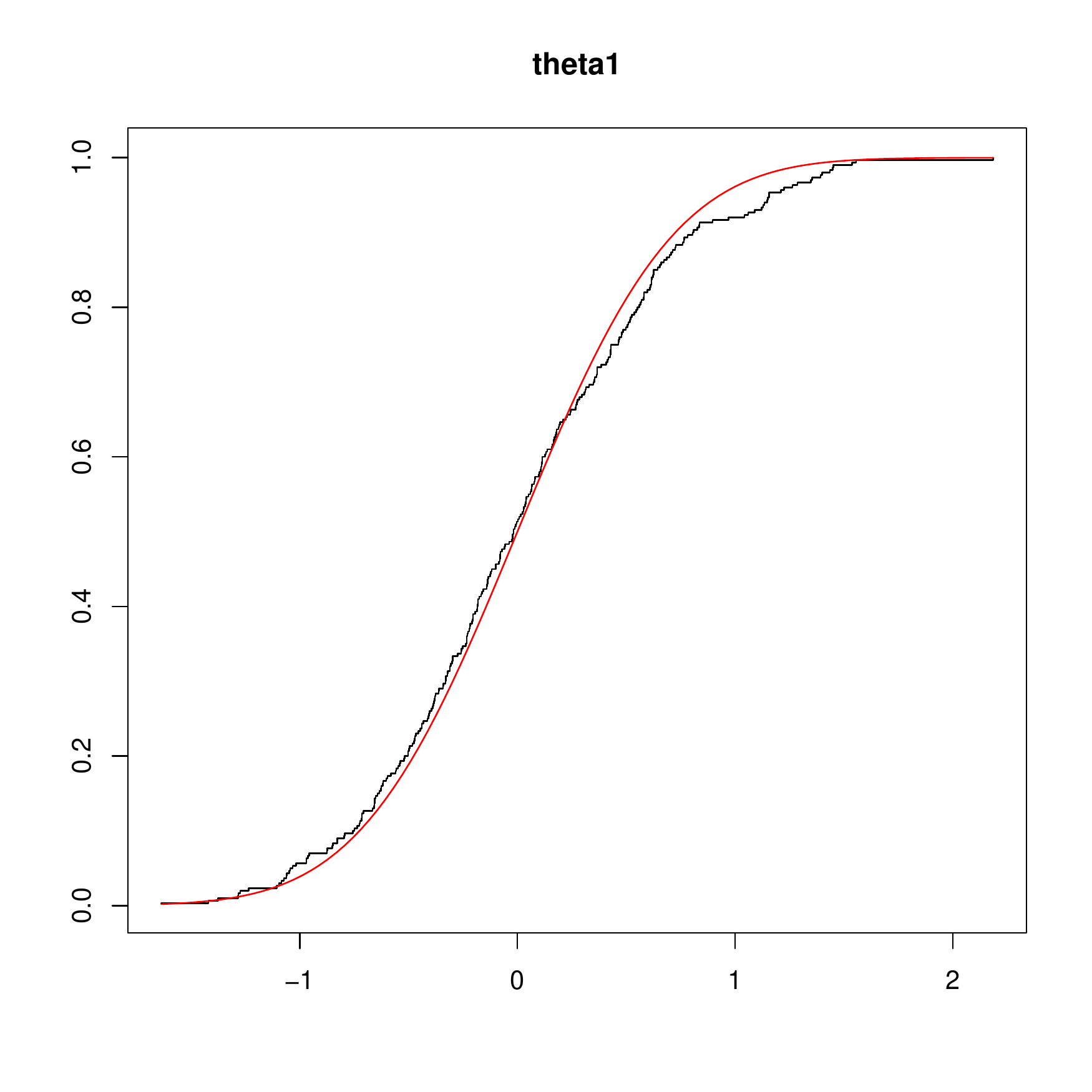}
\includegraphics[width=5cm,pagebox=cropbox,clip]{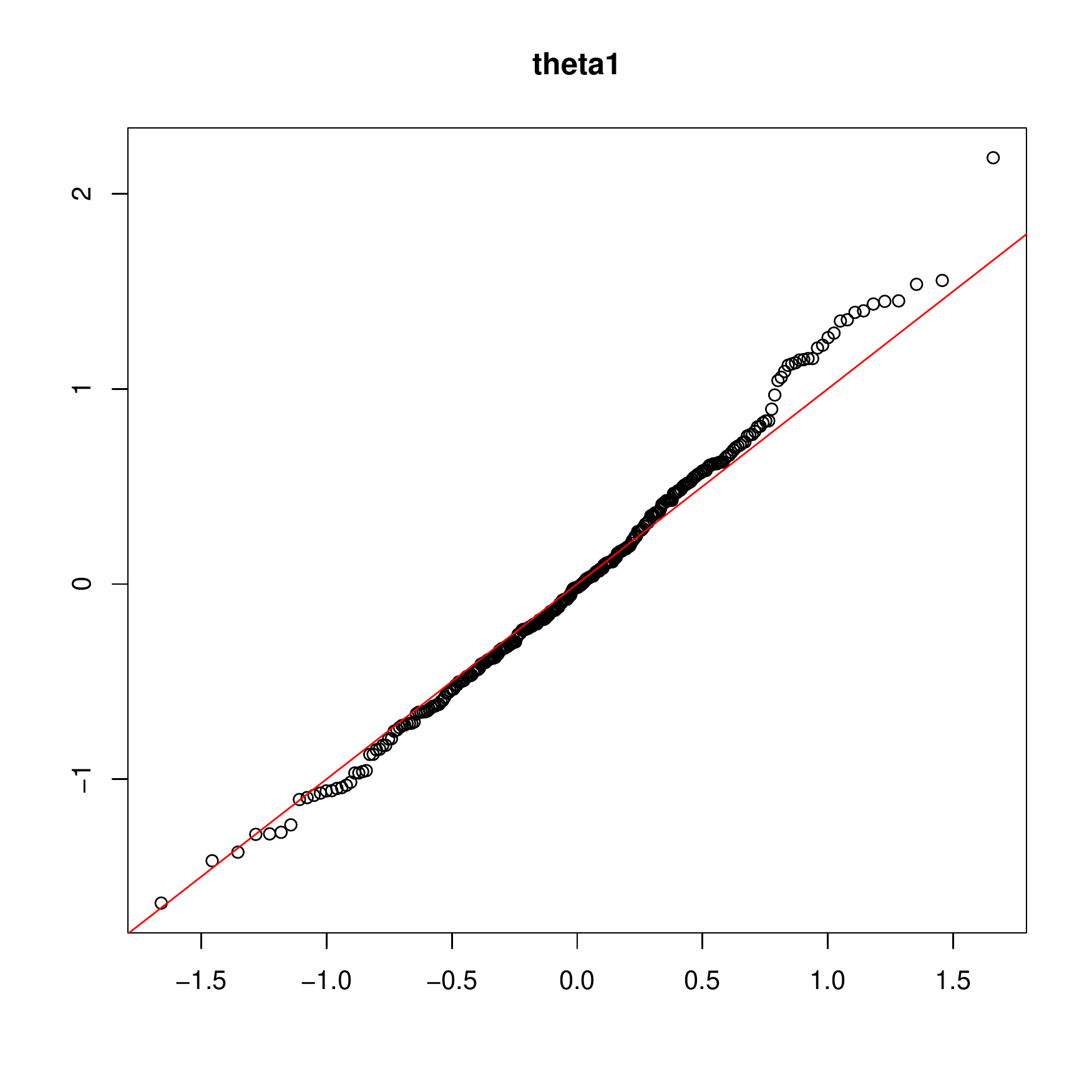}
\includegraphics[width=5cm,pagebox=cropbox,clip]{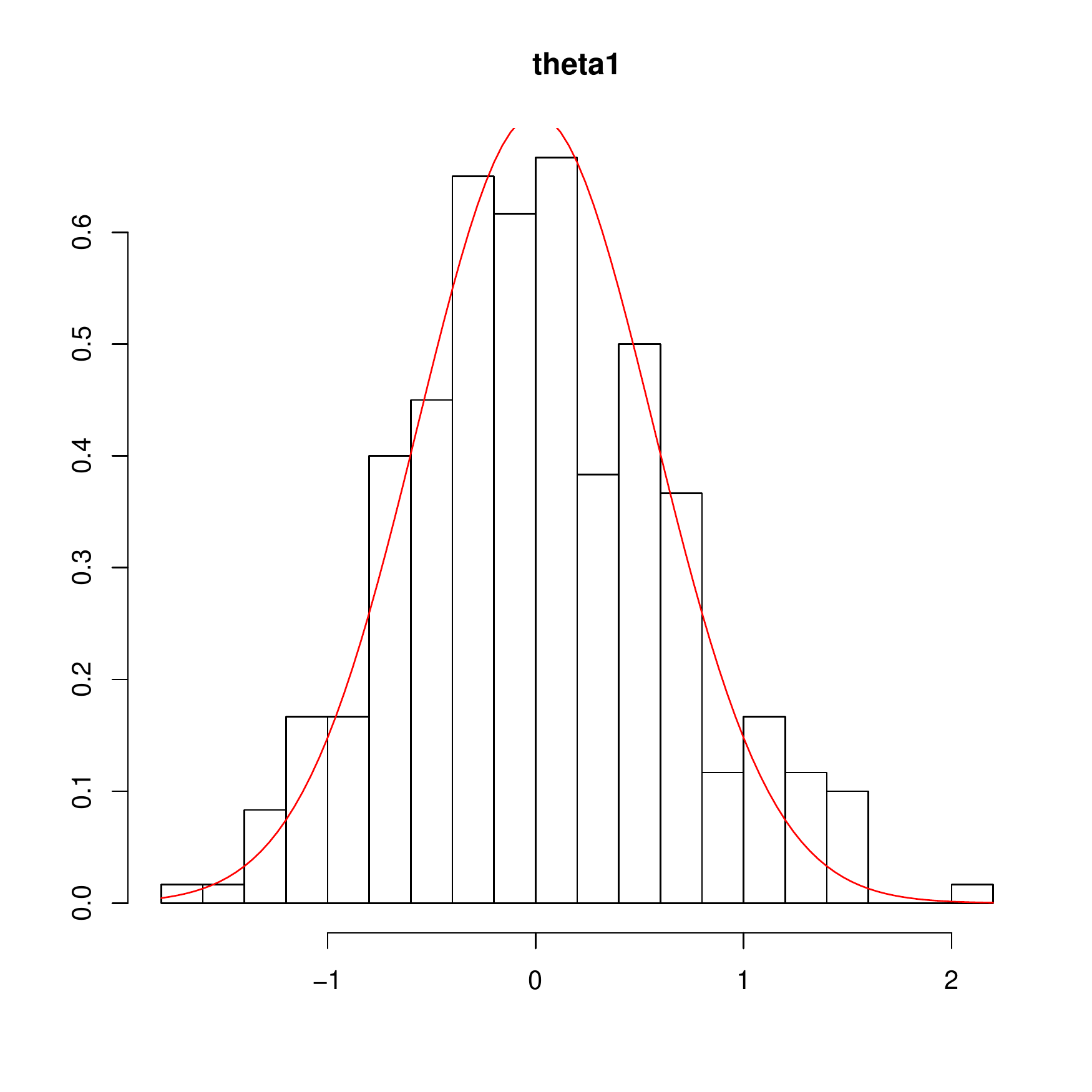}\\
\caption{Simulation results of $\hat{\theta}_1$ \label{fig9}}
\end{center}
\end{figure}

\begin{figure}[h]
\begin{center}
\includegraphics[width=5cm,pagebox=cropbox,clip]
{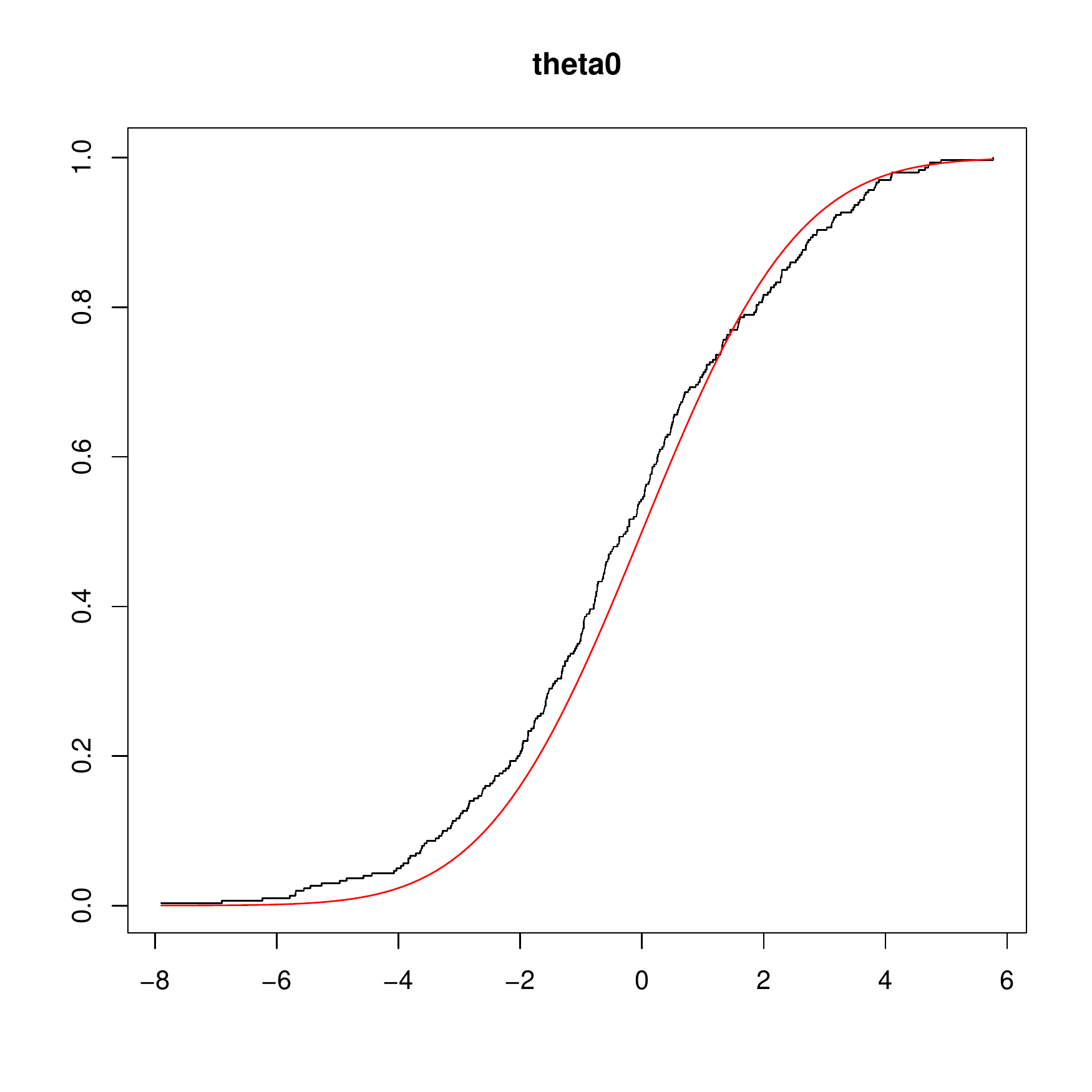}
\includegraphics[width=5cm,pagebox=cropbox,clip]{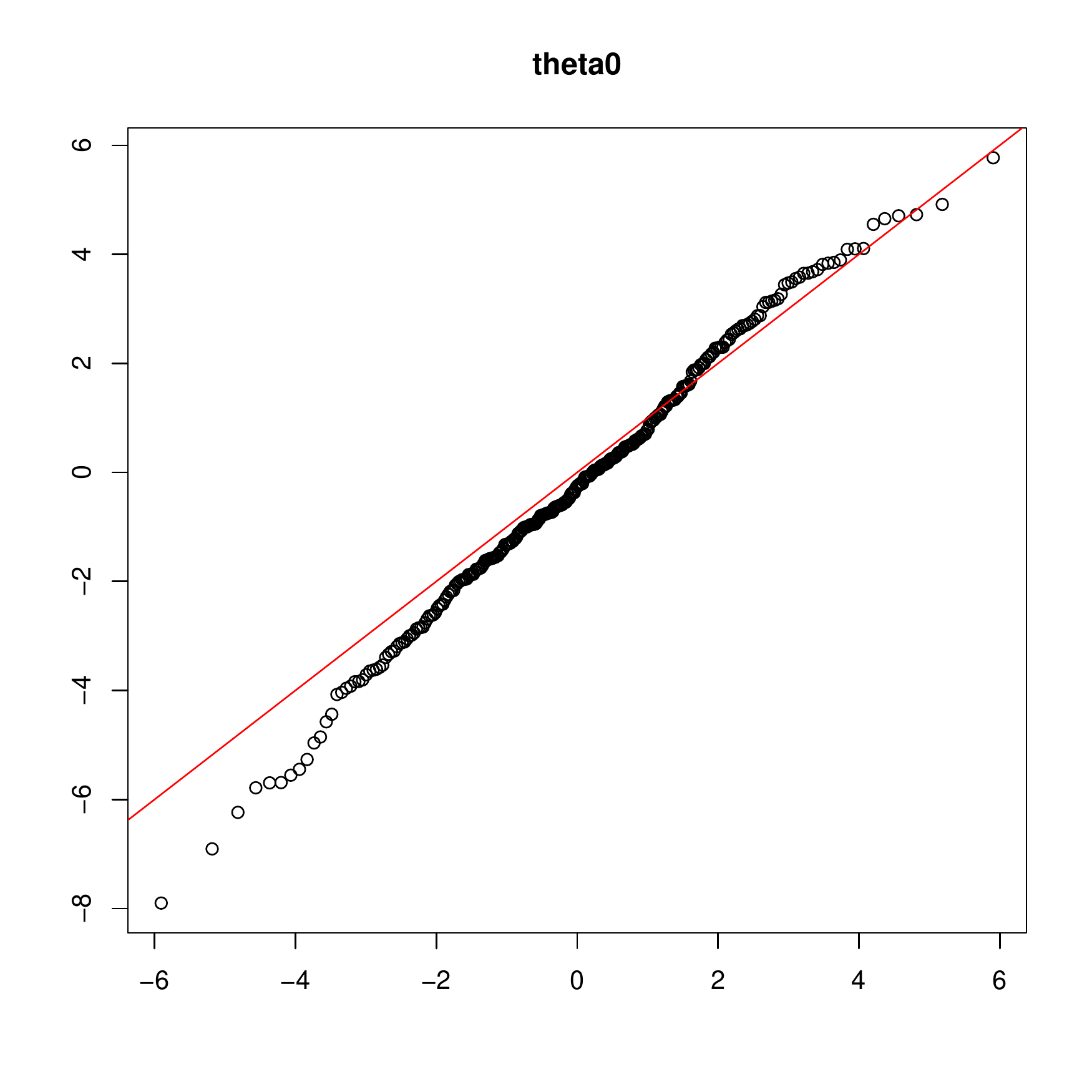}
\includegraphics[width=5cm,pagebox=cropbox,clip]{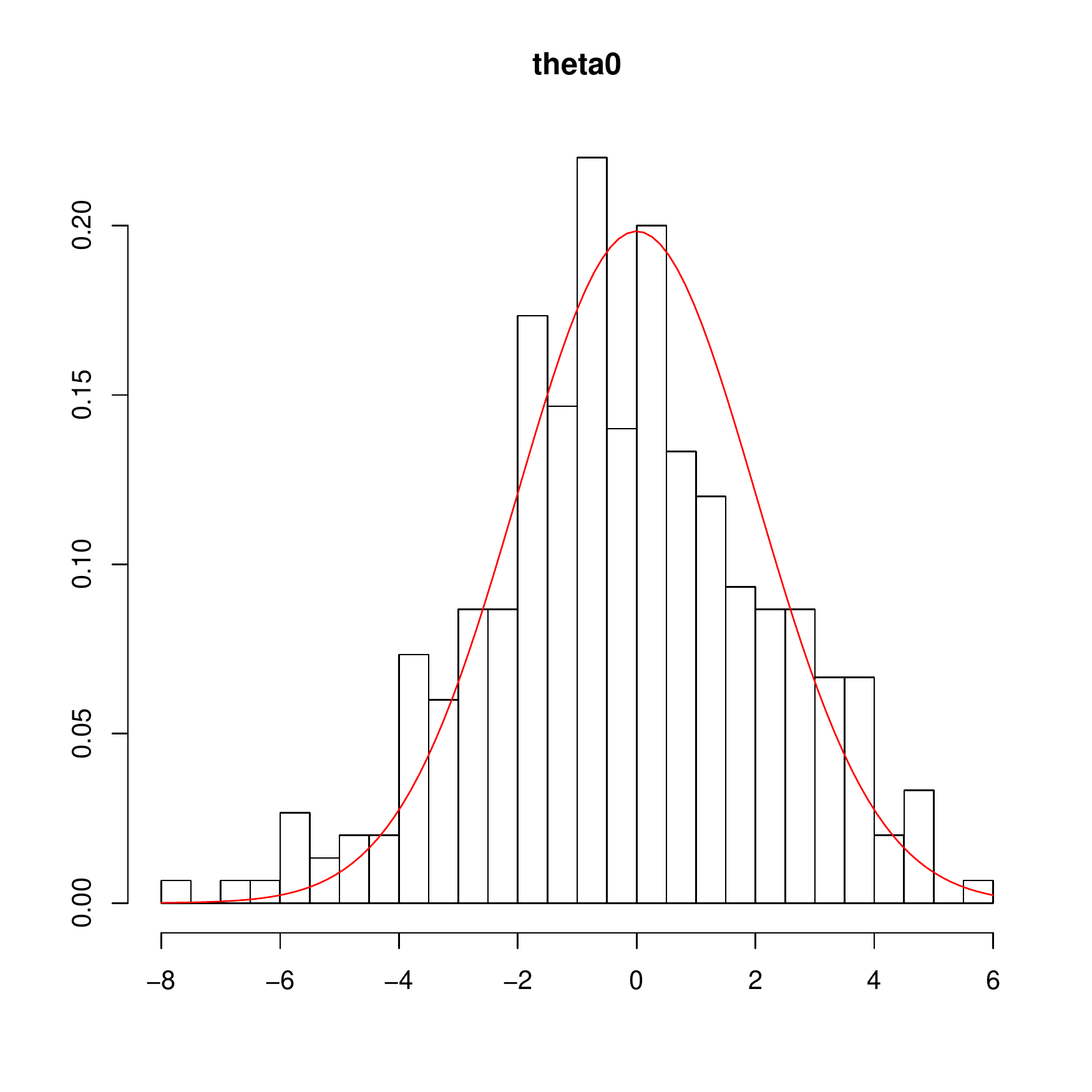}\\
\caption{Simulation results of $\hat{\theta}_0$ \label{fig10}}
\end{center}
\end{figure}

Figures \ref{fig7}-\ref{fig10} are the simulation results of $\hat{\sigma}^2$, $\hat{\theta}_2$, $\hat{\theta}_1$  
and  $\hat{\theta}_0$
in Theorem \ref{thm4}
with $(N, {\bar{m}}, T, {\bar{N}_2}) = (10^5, 24, 100, 800)$.
Note that 
{
$\frac{\bar{N}_2^{\frac{5}{2}}}{ T^{\frac{3}{2}} \bar{m} N} =0.0075$
and 
$\frac{\bar{N}_2^3}{T^2 M^{1-0.01}} =0.57$. }
%
From Figures \ref{fig7}-\ref{fig10}, we can see that 
these estimates have good performance.
When $T$ is large, we can estimate all parameters, $\theta_0$, $\theta_1$, $\theta_2$ and $\sigma$ 
of the SPDE (\ref{spde2}).

\clearpage

\section{Proofs}

{Proof of Theorem 2}. 
Let $s_i = s_{i:N_2}$. 
\begin{align}
	\check{\sigma}^2 =& \sum_{i=1}^{N_2} \left( \check{x}_k(s_i)-\check{x}_k(s_{i-1}) \right)^2 
			\nonumber\\		
		=& \sum_{i=1}^{N_2} \left\{ \check{x}_k(s_i)-x_k(s_i)-\{ \check{x}_k(s_{i-1})-x_k(s_{i-1})\}
			+x_k(s_i)-x_k(s_{i-1})\right\}^2  \nonumber \\
		=& \sum_{i=1}^{N_2} \left\{ (\check{x}_k(s_i)-x_k(s_i))^2 \right. \label{thm2-I} \\
		&+ \left(\check{x}_k(s_{i-1})-x_k(s_{i-1}) \right)^2 \label{thm2-II} \\
		&-2(\check{x}_k(s_i)-x_k(s_i))(\check{x}_k(s_{i-1})-x_k(s_{i-1})) \label{thm2-III} \\
		&+2\left\{ \check{x}_k(s_i)-x_k(s_i)-\{\check{x}_k(s_{i-1})-x_k(s_{i-1})\}\right\}
			(x_k(s_i)-x_k(s_{i-1})) \label{thm2-IV} \\
		&+\left. (x_k(s_i)-x_k(s_{i-1}))^2 \right\}.  \nonumber
\end{align}
First of all, we will show that  
\begin{equation}
	\sqrt{N_2}(\check{\sigma}^2-(\sigma^*)^2)-\sqrt{N_2}
		\left\{\sum_{i=1}^{N_2}(x_k(s_i)-x_k(s_{i-1}))^2-(\sigma^*)^2\right\} =o_p(1). \label{thm2-final}
\end{equation}

Let $g_k(t,y,\eta)=X_t(y)\sqrt{2}\sin(\pi ky)\exp\left\{ \frac{\eta}{2}y \right\}$. 
Note that
\begin{align*}
	x_k(t) =& \int_0^1 X_t(y)\sqrt{2}\sin(\pi ky) \exp\left\{\frac{\eta}{2}y \right\} dy
		= \sum_{j=1}^M \int_{\frac{j-1}{M}}^{\frac{j}{M}} g_k(t,y,\eta)dy, \\
	\check{x}_k(t) =& \frac{1}{M}\sum_{j=1}^M X_t(y_j)\sqrt{2}\sin(\pi ky_j) 
		\exp\left\{\frac{ \check{\eta} }{2} y_j \right\} 
		= \frac{1}{M}\sum_{j=1}^M g_k(t,y_j, \check{\eta}). 
\end{align*}
For the evaluation of (\ref{thm2-I}) and (\ref{thm2-II}), noting that  
\begin{align*}
	\left( \check{x}_k(s_i)-x_k(s_i) \right)^2 =& \left( M \frac{1}{M} \sum_{j=1}^M
		\int_{\frac{j-1}{M}}^{\frac{j}{M}} (g_k(s_i,y,\eta)-g_k(s_i,y_j,\check{\eta}) dy \right)^2 \\
		\leq& M^2 \frac{1}{M} \sum_{j=1}^M 
			\left(\int_{\frac{j-1}{M}}^{\frac{j}{M}} 1^2 dy \right)
			\int_{\frac{j-1}{M}}^{\frac{j}{M}} 
			\left(g_k(s_i,y,\eta)-g_k(s_i,y_j,\check{\eta}) \right)^2dy \\
		=& \sum_{j=1}^M \int_{\frac{j-1}{M}}^{\frac{j}{M}} 
			\left(g_k(s_i,y,\eta)-g_k(s_i,y_j,\check{\eta})\right)^2dy,
\end{align*}
one has that 
\begin{equation*}
	{Z_0} := \sqrt{N_2}\sum_{i=1}^{N_2} \left( \check{x}_k(s_i)-x_k(s_i) \right)^2
		\leq \sqrt{N_2}\sum_{i=1}^{N_2} \sum_{j=1}^M \int_{\frac{j-1}{M}}^{\frac{j}{M}} 
			\left(g_k(s_i,y,\eta)-g_k(s_i,y_j,\check{\eta})\right)^2dy.
\end{equation*}
It follows that 
\begin{align}
	g_k(t,y,\eta)-g_k(t,y_j,\check{\eta}) 
		=& X_t(y)\sqrt{2}\sin(\pi ky) \exp\left\{\frac{\eta}{2} y \right\}
		-X_t(y_j)\sqrt{2}\sin(\pi ky_j) \exp\left\{\frac{\check{\eta}}{2} y_j \right\} \nonumber\\
		=& (X_t(y)-X_t(y_j))\sqrt{2}\sin(\pi ky) \exp\left\{\frac{\eta}{2} y \right\} \label{thm2-1} \\
		&+ X_t(y_j)\left(\sqrt{2}\sin(\pi ky) \exp\left\{\frac{\eta}{2} y \right\}
		-\sqrt{2}\sin(\pi ky_j) \exp\left\{\frac{\eta}{2} y_j \right\} \right) \label{thm2-2} \\
		&+ X_t(y_j)\sqrt{2}\sin(\pi ky_j) \left\{ \exp\left\{\frac{\eta}{2} y_j \right\}
		-\exp\left\{\frac{\check{\eta}}{2} y_j \right\} \right\} \label{thm2-3} \\
		=:& g^{(1)}(t,y) + g^{(2)}(t,y) + g^{(3)}(t,y_j,\check{\eta}) \nonumber
\end{align}
and that 
\begin{equation*}
	Z_0 \leq C\sqrt{N_2}\sum_{i=1}^{N_2} \sum_{j=1}^M \int_{\frac{j-1}{M}}^{\frac{j}{M}} 
		\{
(g^{(1)}(s_i,y))^2 + (g^{(2)}(s_i,y))^2 + (g^{(3)}(s_i,y_j,\check{\eta}))^2 
		\} dy.
\end{equation*}
For the evaluation of (\ref{thm2-2}), noting that  
\begin{align*}
	E\left[ X_t^2(y) \right] =& \sum_{k=1}^\infty \sum_{l=1}^\infty 
		E[x_k(t)x_l(t)] e_k(y)e_l(y) \\
		=& \sum_{k=1}^\infty E[x_k^2(t)] e_k^2(y) 
		= \sum_{k=1}^\infty \frac{\sigma^2}{2\lambda_k} \left\{1-e^{-2\lambda_k t} \right\}e_k^2(y) \\
		\leq& C\sum_{k=1}^\infty \frac{1}{k^2} < \infty
\end{align*}
and
\begin{equation*}
	E[ (g^{(2)}(s_i,y))^2 ] \leq C(y_j-y_{j-1})^2 \leq \frac{C}{M^2},
\end{equation*}
one has that 
{under $\frac{N_2^{\frac{3}{2}}}{M^{1-\rho}} \rightarrow 0$, }
\begin{equation}
	\sqrt{N_2}N_2\frac{1}{N_2}\sum_{i=1}^{N_2} \sum_{j=1}^M 
		\int_{\frac{j-1}{M}}^{\frac{j}{M}} (g^{(2)}(s_i,y))^2 dy = o_p(1). \label{thm2-4}
\end{equation}
For the evaluation of (\ref{thm2-3}), 
\begin{align}
	& \sqrt{N_2}N_2 \frac{1}{N_2}\sum_{i=1}^{N_2} \frac{1}{M} \sum_{j=1}^M (g^{(3)}(s_i,y_j,\check{\eta}))^2 
		 \nonumber\\
		=& \sqrt{N_2}N_2\frac{1}{N_2}\sum_{i=1}^{N_2} \frac{1}{M} \sum_{j=1}^M X_{s_i}^2(y_j)
		{
		2\sin^2(\pi k y_j) 
		\left(\int_0^1\frac{y_j}{2}\exp \left\{\frac{y_j}{2}(\eta+u(\check{\eta}-\eta))\right\}du \right)^2
		}
		\frac{1}{N_m} \label{thm2-A}\\
		&\times (\sqrt{N_m}(\check{\eta}-\eta))^{2}. \nonumber
\end{align}
Let $\zeta>0$ and $\epsilon >0$. On $A=\{|\check{\eta}-\eta| < \zeta \}$,
\begin{equation*}
	( \ref{thm2-A}  ) \leq \frac{C}{N_2}\sum_{i=1}^{N_2}\frac{1}{M}\sum_{j=1}^M
	{ X_{s_i}^2(y_j) 2\sin^2(\pi k y_j) }\frac{N_2^{\frac{3}{2}}}{N_m},
\end{equation*}
and 
\begin{align*}
	P(|( \ref{thm2-A} )| > \varepsilon) =& P(\{| (\ref{thm2-A}) |>\varepsilon\}\cap A)
		+P(\{ | ( \ref{thm2-A} ) |>\varepsilon\} \cap A^c) \nonumber \\
	\leq&P\left(\frac{N_2^{\frac{3}{2}}}{N_m}\frac{1}{N_2}\sum_{i=1}^{N_2}
		\frac{1}{M}\sum_{j=1}^M X_{s_i}^2(y_j)2\sin^2(\pi y_j)>\varepsilon \right)
		+P(A^c) \nonumber \\
	\leq& \frac{C}{\varepsilon}\frac{N_2^{\frac{3}{2}}}{N_m}+o(1) \rightarrow 0, \nonumber \\
	(\ref{thm2-A})  =& o_p(1).
\end{align*}
Noting that $\sqrt{N_m}(\check{\eta}-\eta) = O_p(1)$,
one has that 
under $\frac{N_2^{\frac{3}{2}}}{N_m} \rightarrow 0$,
\begin{equation}
	\sqrt{N_2}\sum_{i=1}^{N_2}\frac{1}{M}\sum_{j=1}^M (g^{(3)}(s_i,y_j,\check{\eta}))^2  = o_p(1). 
		 \label{thm2-5}
\end{equation}
For the evaluation of (\ref{thm2-1}), noting that 
\begin{equation*}
	E[ (g^{(1)}(s_i,y))^2 ] \leq CE\left[\left(X_{s_i}(y)-X_{s_i}(y_j)\right)^2\right],
\end{equation*}
\begin{equation*}
	\left(X_t(y)-X_t(y_j)\right)^2 = \sum_{k,l} x_k(t)(e_k(y)-e_l(y_j))x_l(t)(e_l(y)-e_l(y_j)),
\end{equation*}
\begin{align*}
	&|e_k(y)-e_k(y_j)| \\
	=& \left| \int_0^1 \left(\sqrt{2}\pi k \cos(\pi k(y_j+u(y-y_j)))
		\exp\left\{-\frac{\eta}{2}(y_j+u(y-y_j) ) \right\} \right. \right. \\
		&+ \left.\left.\sqrt{2}\sin(\pi k(y_j+u(y-y_j)))\left(-\frac{\eta}{2}\right)
		\exp\left\{-\frac{\eta}{2}(y_j+u(y-y_j))\right\}\right)du \times (y-y_j)\right| \\
		\leq& C\frac{k}{M},
\end{align*}
\begin{align*}
	E\left[\left(X_t(y)-X_t(y_j)\right)^2\right] =& 
		\sum_{k=1}^\infty E[x_k^2(t)] (e_k(y)-e_k(y_j))^2 \\
		\leq& \sum_{k=1}^\infty \frac{C}{k^2}\left(\frac{k}{M}\right)^{1-\rho}C^{1+\rho}
		\leq \frac{C_1}{M^{1-\rho}}, 
\end{align*}
we obtain that 
under $\frac{N_2^{\frac{3}{2}}}{M^{1-\rho}} \rightarrow 0$,
\begin{equation}
	\sqrt{N_2}\sum_{i=1}^{N_2} 
	\sum_{j=1}^M \int_{\frac{j-1}{M}}^{\frac{j}{M}} 
(g^{(1)}(s_i,y))^2
		dy = O_p \left(\frac{N_2^{\frac{3}{2}}}{M^{1-\rho}}\right) = o_p(1).
                  \label{thm2-6}
\end{equation}
Therefore, 
it follows from (\ref{thm2-4}),  (\ref{thm2-5}) and (\ref{thm2-6}) that
{under $\frac{N_2^{\frac{3}{2}}}{N m} \rightarrow 0$} and 
$\frac{N_2^{\frac{3}{2}}}{M^{1-\rho}} \rightarrow 0$,
\begin{equation*}
	\sqrt{N_2}\sum_{i=1}^{N_2}\left( \check{x}_k(s_i)-x_k(s_i) \right)^2 = o_p(1)
\end{equation*}
and
\begin{equation*}
	\sqrt{N_2}\sum_{i=1}^{N_2}\left( \check{x}_k(s_{i-1})-x_k(s_{i-1}) \right)^2 = o_p(1).
\end{equation*}

For the evaluation of (\ref{thm2-III}), it follows that 
\begin{align*}
	& \left| \sqrt{N_2}\sum_{i=1}^{N_2}\left(\check{x}_k(s_i)-x_k(s_i)\right)
		\left( \check{x}_k(s_{i-1})-x_k(s_{i-1})\right)\right|^2 \\
	\leq& N_2\sum_{i=1}^{N_2}\left(\check{x}_k(s_i)-x_k(s_i)\right)^2
		\sum_{i=1}^{N_2}\left(\check{x}_k(s_{i-1})-x_k(s_{i-1})\right)^2\\
	=& \sqrt{N_2}\sum_{i=1}^{N_2}\left(\check{x}_k(s_i)-x_k(s_i)\right)^2
		\sqrt{N_2}\sum_{i=1}^{N_2}\left(\check{x}_k(s_{i-1})-x_k(s_{i-1})\right)^2\\
	=& o_p(1).
\end{align*}

For the evaluation of (\ref{thm2-IV}), 
setting that $\Delta X_{s_i}(y) = X_{s_i}(y)-X_{s_{i-1}}(y)$ and 
	$\Delta_{s_i}x  = x_1(s_i)-x_1(s_{i-1})$,
one has that 
\begin{align}
	U_1 :=& 2\sqrt{N_2} \sum_{i=1}^{N_2} 
		\{\check{x}_1(s_i)-\check{x}_1(s_{i-1})-(x_1(s_i)-x_1(s_{i-1}))\} (x_1(s_i)-x_1(s_{i-1})) \nonumber  \\
		=& 2\sqrt{N_2} \sum_{i=1}^{N_2}\frac{1}{M}\sum_{j=1}^M
			\Delta X_{s_i}(y_j)
			\sqrt{2}\sin(\pi y_j)
			\left(\exp\left\{\frac{\check{\eta}}{2}y_j\right\}
			-\exp\left\{\frac{\eta}{2}y_j\right\}\right) 
			\times  \Delta_{s_i}x  \label{thm2-i} \\
		&+ 2\sqrt{N_2} \sum_{i=1}^{N_2}\sum_{j=1}^M \int_{\frac{j-1}{M}}^{\frac{j}{M}} 
			\Delta X_{s_i}(y_j)
			(\sqrt{2}\sin(\pi y_j)-\sqrt{2}\sin(\pi y))
			\exp\left\{\frac{\eta}{2}y_j\right\}dy
			\times \Delta_{s_i}x \label{thm2-ii} \\
		&+ 2\sqrt{N_2} \sum_{i=1}^{N_2}\sum_{j=1}^M \int_{\frac{j-1}{M}}^{\frac{j}{M}} 
			\Delta X_{s_i}(y_j)
			\sqrt{2}\sin(\pi y)
			\left(\exp\left\{\frac{\eta}{2}y_j\right\}
			-\exp\left\{\frac{\eta}{2}y\right\}\right) dy		
			\times \Delta_{s_i}x  \label{thm2-iii} \\
		&+ 2\sqrt{N_2} \sum_{i=1}^{N_2}\sum_{j=1}^M \int_{\frac{j-1}{M}}^{\frac{j}{M}} 
			\left\{
			\Delta X_{s_i}(y_j) -\Delta X_{s_i}(y)
			\right\}\sqrt{2}\sin(\pi y)
			\exp\left\{\frac{\eta}{2}y\right\} dy
			\times \Delta_{s_i}x.  \label{thm2-iv} 
\end{align}

For the evaluation of (\ref{thm2-i}), noting that 
\begin{align*}
	(\ref{thm2-i}) =& 2\sqrt{N_2} \sum_{i=1}^{N_2}\frac{1}{M}\sum_{j=1}^M
			\Delta X_{s_i}(y_j) 
			\sqrt{2}\sin(\pi y_j)
			\int_0^1 \frac{y_j}{2}
			\exp\left\{\frac{y_j}{2}(\eta+u(\check{\eta}-\eta))\right\}du 
			(\Delta_{s_i}x)
			(\check{\eta}-\eta), 
\end{align*}
one has that 
\begin{align*}
	(\ref{thm2-i})^2 \leq& 4 N_2 \sum_{i=1}^{N_2}\frac{1}{Nm}\frac{1}{M}\sum_{j=1}^M
			\left(
			\Delta X_{s_i}(y_j) 
			\right)^2
			\left(
						\sqrt{2}\sin(\pi y_j)
			\int_0^1 \frac{y_j}{2}
			\exp\left\{\frac{y_j}{2}(\eta+u(\check{\eta}-\eta))\right\}du 
			\right)^2\\
		&\times \sum_{i=1}^{N_2} (\Delta_{s_i}x)^2
			\times (\sqrt{Nm}(\check{\eta}-\eta))^{2}.
\end{align*}
Set that
\begin{align*}			
	B :=& 4 N_2 \sum_{i=1}^{N_2}\frac{1}{Nm}\frac{1}{M}\sum_{j=1}^M
			\left(
			\Delta X_{s_i}(y_j) 
			\right)^2
			\left(	\sqrt{2}\sin(\pi y_j)
			\int_0^1 \frac{y_j}{2}
			\exp\left\{\frac{y_j}{2}(\eta+u(\check{\eta}-\eta))\right\}du 
			\right)^2.
\end{align*}
Since 
\begin{align*}
	&\sum_{i=1}^{N_2}  (\Delta_{s_i}x)^2 = O_p(1), \quad
	(\sqrt{Nm}(\check{\eta}-\eta))^{2} = O_p(1), 
\end{align*}
we obtain that 
\begin{align*}
	P(|B| > \varepsilon) &= P(\{|B|>\varepsilon\}\cap A)
		+P(\{|B|>\varepsilon\} \cap A^c) \\
		&\leq  C \frac{N_2 N_2}{Nm} 
		E\left[ 
			\left(	\Delta X_{s_i}(y_j) \right)^2
		\right]
			\frac{1}{\varepsilon} + P(A^c) \\
		&= C_1 \frac{N_2^{\frac{3}{2}}}{Nm}\frac{1}{\varepsilon}+o(1) \\
		&\rightarrow 0 
\end{align*}
under $\frac{N_2^{\frac{3}{2}}}{Nm} \rightarrow 0$.
Therefore $(\ref{thm2-i})^2 = o_p(1)$. 

For the evaluation of ({\ref{thm2-ii}),  noting that
\begin{align*}
		& N_2 \sum_{i=1}^{N_2} \frac{1}{M} \sum_{j=1}^M
		\left(
			\Delta X_{s_i}(y_j) 
			\right)^2
		|y-y_j|^2
			= O_p\left(\frac{N_2^2}{\sqrt{N_2}M^2}\right), \\
		&  \sum_{i=1}^{N_2} (\Delta_{s_i}x)^2  =O_p(1),
\end{align*}
one has that 
under $\frac{N_2^{\frac{3}{2}}}{M^{1-\rho}} \rightarrow 0$, 
\begin{align*}
	(\ref{thm2-ii})^2 \leq& N_2 \sum_{i=1}^{N_2} \frac{1}{M} \sum_{j=1}^M 
			\left(
			\Delta X_{s_i}(y_j) 
			\right)^2
			 |y-y_j|^2
			\sum_{i=1}^{N_2} (\Delta_{s_i}x)^2 = o_p(1).
\end{align*}

For the evaluation of (\ref{thm2-iii}), noting that
\begin{align*}
		& N_2 \sum_{i=1}^{N_2} M^2 \frac{1}{M} \sum_{j=1}^M \frac{1}{M}
		\int_{\frac{j-1}{M}}^{\frac{j}{M}} 
		\left(
			\Delta X_{s_i}(y_j) 
			\right)^2
			2\sin(\pi y)|y-y_j|^2 dy
		= O_p\left(\frac{N_2^2}{\sqrt{N_2}M^2} \right), \\
		&  \sum_{i=1}^{N_2} (\Delta_{s_i}x)^2  =O_p(1),
\end{align*}
we obtain that 
under $\frac{N_2^{\frac{3}{2}}}{M^{1-\rho}} \rightarrow 0$,
\begin{align*}
	({\ref{thm2-iii}})^2 \leq& N_2 \sum_{i=1}^{N_2} M^2 \frac{1}{M} \sum_{j=1}^M \frac{1}{M}
		\int_{\frac{j-1}{M}}^{\frac{j}{M}} 
		\left(
			\Delta X_{s_i}(y_j) 
			\right)^2 
			2\sin(\pi y)|y-y_j|^2 dy		
			\sum_{i=1}^{N_2} (\Delta_{s_i}x)^2  = o_p(1).
\end{align*}

For the evaluation of (\ref{thm2-iv}), setting that
\begin{align*}
	h(y, \Delta_{s_i}x) :=& \left(\sum_{k=1}^\infty (x_k(s_i)-x_k(s_{i-1}))(e_k(y_j)-e_k(y))\right)^2
		(x_1(s_i)-x_1(s_{i-1}))^2,
\end{align*}
one has that
\begin{align*}
	E[h(y, \Delta_{s_i}x)] =& \sum_{k=2}^\infty E[(x_k(s_i)-x_k(s_{i-1}))^2] E[(x_1(s_i)-x_1(s_{i-1}))^2]
		(e_k(y_j)-e_k(y))^2 \\
		&+ E[(x_1(s_i)-x_1(s_{i-1}))^4](e_1(y_j)-e_1(y))^2 \\
		\leq& \frac{C}{M^{1-\rho}}\frac{1}{N_2}+\frac{1}{N_2^2 M^2}
\end{align*}
and
\begin{equation*}
	({\ref{thm2-iv}})^2 \leq N_2 \sum_{i=1}^{N_2}M^2 \frac{1}{M}\sum_{j=1}^M \frac{1}{M}
		\int_{\frac{j-1}{M}}^{\frac{j}{M}} 
		h(y, \Delta_{s_i}x) 
		dy.
\end{equation*}
Since it follows that 
under $\frac{N_2^{\frac{3}{2}}}{M^{1-\rho}} \rightarrow 0$, 
\begin{equation*}
	E[(\ref{thm2-iv})^2] \leq \frac{N_2}{M^{1-\rho}}+\frac{1}{M^2} \rightarrow 0, 
	\end{equation*}
we obtain that 
$({\ref{thm2-iv}}) = o_p(1)$.
Hence, 
$U_1 = o_p(1)$.
Consequently, 
under $\frac{N_2^{\frac{3}{2}}}{N_m} \rightarrow 0$
and $\frac{N_2^{\frac{3}{2}}}{M^{1-\rho}} \rightarrow 0$, 
(\ref{thm2-final}) {holds} true, which yields that
\begin{equation}
\sqrt{N_2}(\check{\sigma}^2-(\sigma^*)^2) 
\stackrel{d}{\rightarrow} 
N \left( 
0,
2 (\sigma^*)^4 
\right). \label{thm2-AN_1}
\end{equation}
For the estimator of $\theta_2$, we obtain that
\begin{eqnarray*}
\sqrt{N_2} (\check{\theta}_2 -\theta_2^*) 
&=& \sqrt{N_2} \left( \left( \frac{\check{\sigma}^2}{\check{\sigma}_0^2} \right)^2  - \left( \frac{({\sigma^*})^2}{({\sigma_0^*})^2} \right)^2 \right) \\
&=&
\sqrt{N_2} 
\left(  
\check{\sigma}^4 
\left\{
\left( \frac{1}{\check{\sigma}_0^2} \right)^2  - \left( \frac{1}{({\sigma_0^*})^2} \right)^2 
\right\} 
+ \frac{1}{(\sigma_0^*)^4} \left\{ \check{\sigma}^4 - ({\sigma^*})^4 \right\}  
\right)  \\
&=&
\frac{\sqrt{N_2}}{\sqrt{m N}}   
\check{\sigma}^4
\sqrt{m N} 
\left\{
\left( \frac{1}{\check{\sigma}_0^2} \right)^2  - \left( \frac{1}{({\sigma_0^*})^2} \right)^2 
\right\} 
+ \sqrt{N_2}  \frac{1}{(\sigma_0^*)^4} \left\{ \check{\sigma}^4 - ({\sigma^*})^4 \right\}    \\
&=&
\sqrt{N_2}  \frac{1}{(\sigma_0^*)^4} \left\{ \check{\sigma}^4 - ({\sigma^*})^4 \right\} +o_p(1).
\end{eqnarray*}
For the estimator of $\theta_1$, one has that
\begin{eqnarray*}
\sqrt{N_2} (\check{\theta}_1 -\theta_1^*) 
&=& \sqrt{N_2} \left(  \check{\eta} \hat{\theta}_2 - \eta^* \theta_2^* \right) 
= \sqrt{N_2} \left(  \check{\theta}_2 \left( \check{\eta}  - \eta^* \right) + \eta^* \left( \check{\theta}_2 -\theta_2^* \right) \right) \\
&=& \sqrt{N_2} \eta^* \left( \check{\theta}_2 -\theta_2^* \right) + o_p(1) \\
&=& \sqrt{N_2}  \eta^* \frac{1}{(\sigma_0^*)^4} \left\{ \check{\sigma}^4 - ({\sigma^*})^4 \right\} +o_p(1).
\end{eqnarray*}

\noindent
By noting that
\begin{equation*}
\begin{pmatrix}
\sqrt{N_2}(\check{\sigma}^2-(\sigma^*)^2) \\
\sqrt{N_2}(\check{\theta}_2- \theta_2^*) \\
\sqrt{N_2}(\check{\theta}_1- \theta_1^*) 
\end{pmatrix}
=
\begin{pmatrix}
\sqrt{N_2}(\check{\sigma}^2-(\sigma^*)^2) \\
\sqrt{N_2}\frac{1}{(\sigma_0^*)^4}(\check{\sigma}^2-(\sigma^*)^2) \\
\sqrt{N_2}\frac{\eta^*}{(\sigma_0^*)^4}(\check{\sigma}^2-(\sigma^*)^2) 
\end{pmatrix}
+o_p(1),
\end{equation*}

\noindent
it follows from (\ref{thm2-AN_1}) and the delta method that
under $\frac{N_2^{\frac{3}{2}}}{N_m} \rightarrow 0$
and $\frac{N_2^{\frac{3}{2}}}{M^{1-\rho}} \rightarrow 0$, 
\begin{equation*}
\begin{pmatrix}
\sqrt{N_2}(\check{\sigma}^2-(\sigma^*)^2) \\
\sqrt{N_2}\frac{1}{(\sigma_0^*)^4}(\check{\sigma}^2-(\sigma^*)^2) \\
\sqrt{N_2}\frac{\eta^*}{(\sigma_0^*)^4}(\check{\sigma}^2-(\sigma^*)^2) 
\end{pmatrix}
\stackrel{d}{\rightarrow} 
N \left( 
\begin{pmatrix}
0 \\
0 \\
0 
\end{pmatrix}, 
\begin{pmatrix}
2 (\sigma^*)^4 & 4 \theta_2^* (\sigma^*)^2 & 4 \theta_1^* (\sigma^*)^2  \\
 4 \theta_2^* (\sigma^*)^2 &  8 (\theta_2^*)^2 & 8 \theta_1^* \theta_2^*  \\
 4 \theta_1^* (\sigma^*)^2 & 8 \theta_1^* \theta_2^* & 8 (\theta_1^*)^2  
\end{pmatrix}
\right), \label{AN_2}
\end{equation*}
which completes the proof.

\vspace{1cm}

\noindent
{Proof of Theorem 3}. 
By a similar way to the proof of Theorem 5.1 in Bibinger and Trabs (2007),
we can show the result
under $N h_{N:T}^2 \rightarrow 0$,  
$\bar{m} \rightarrow \infty$ and 
$\bar{m} = O(h_{N:T}^{-\rho})$ 
for
$\rho \in (0, 1/2)$.

\vspace{1cm}

\noindent
{Proof of Theorem 4}. 
Let $\delta =\delta_{\bar{N}_2:T}$
and $s_i = s_{i:\bar{N}_2:T}$.
The quasi log-likelihood function based on ${\bf \bar{x}}_k =\{\bar{x}_k(s_{i:\bar{N}_2:T})\}_{i=1,\ldots, \bar{N}_2}
 =\{\bar{x}_k(s_{i})\}_{i=1,\ldots, \bar{N}_2}$  
is as follows. 
\[
	l_{\bar{N}_2}(\lambda_k,\sigma^2 \ | \ {\bf \bar{x}}_k  ) = -\frac{1}{2} \sum_{i=1}^{\bar{N}_2} 
		\left\{\log \frac{\sigma^2(1-e^{-2\lambda_k \dlnt})}{2\lambda_k}
		+\frac{\left( \bar{x}_k(s_i)-e^{-\lambda_k \dlnt} \bar{x}_k(s_{i-1})\right)^2}{
		\frac{\sigma^2(1-e^{-2\lambda_k \dlnt})}{2\lambda_k}}
		\right\}.
\]
Set that $k=1$, $\lambda=\lambda_1$,
${\bf \bar{x}} = {\bf \bar{x}}_1 =\{\bar{x}_1(s_{i:\bar{N}_2:T})\}_{i=1,\ldots, \bar{N}_2}
=\{\bar{x}_1(s_{i})\}_{i=1,\ldots, \bar{N}_2}$    and
\[
         \Xi(\lambda) = \frac{(1-e^{-2\lambda \dlnt})}{2\lambda \dlnt}.  
\]

For the consistency of $\hat{\sigma}_2$ and $\hat{\lambda}$,
it is enough to show that 
under
{
$\frac{\bar{N}_2^{\frac{5}{2}}}{T^{\frac{5}{2}} N \bar{m}} \to 0$ and 
$\frac{\bar{N}_2^3}{T^3 M^{1-\rho_1}} \to 0$, }
\begin{eqnarray}
Z:= 
\frac{1}{T} 
\left\{
l_{\bar{N}_2}(\lambda, \sigma^2 \ | \ {\bf \bar{x}}  ) 
-
l_{\bar{N}_2}(\lambda, \sigma^2 \ | \ {\bf {x}}  ) 
\right\} = o_p(1) \label{consistency-1}
\end{eqnarray}
uniformly in $(\lambda, \sigma^2)$, where
$Z$ is the difference between the quasi log-likelihood functions based on 
${\bf \bar{x}}$ and  ${\bf {x}}$. 
Note that (\ref{consistency-1}) yields that 
\begin{eqnarray*} 
\frac{1}{\bN2} 
\left\{
 l_{\bar{N}_2}(\lambda, \sigma^2  \ | \ {\bf \bar{x}}  ) 
-
 l_{\bar{N}_2}(\lambda, \sigma^2 \ | \ {\bf {x}} ) 
\right\} = o_p(1)
\end{eqnarray*}
uniformly in $(\lambda, \sigma^2)$.

Since
\begin{align*}
Z=& \frac{1}{T} \frac{1}{2\sigma^2 \dlnt \Xi(\lambda)}
		\sum_{i=1}^{\bN2} \left\{\left(
			\bar{x}_1(s_i)-e^{-\lambda \dlnt} \bar{x}_1(s_{i-1})
		\right)^2 \right. 
	\left. -\left(x_1(s_i)-e^{-\lambda \dlnt} x_1(s_{i-1})\right)^2 \right\}
\end{align*}
and 
\begin{align*}
	& \left(	
		\bar{x}_1(s_i)-e^{-\lambda \dlnt} \bar{x}_1(s_{i-1})
	\right)^2 \\
	=& \left\{
		\bar{x}_1(s_i)-x_1(s_i)
		-e^{-\lambda \dlnt}(\bar{x}_1(s_i)-x_1(s_{i-1}))
		+(x_1(s_i)-e^{-\lambda \dlnt} x_1(s_{i-1})) 
	\right\}^2 \\
	=& \left\{
		(\bar{x}_1(s_i)-x_1(s_i))
		-e^{-\lambda \dlnt}(\bar{x}_1(s_{i-1})-x_1(s_{i-1})) \right\}^2 \\
	&+2\{ \bar{x}_1(s_i)-x_1(s_i)
		-e^{-\lambda \dlnt}(\bar{x}_1(s_{i-1})-x_1(s_{i-1})) \}
		(x_1(s_i)-e^{-\lambda \dlnt} x_1(s_{i-1})) \\
	&+(x_1(s_i)-e^{-\lambda \dlnt} x_1(s_{i-1}))^2,
\end{align*}
it follows that 
\begin{align}
	Z =& \frac{({\dlnt})^{-1}}{T} \frac{1}{2\sigma^2\Xi(\lambda)}
		\sum_{i=1}^{\bN2} \left[\left\{
		\left( \bar{x}_1(s_i)-x_1(s_i) \right)
		-e^{-\lambda \dlnt} \left( \bar{x}_1(s_{i-1})-x_1(s_{i-1}) \right)
		\right\}^2
	\right. \nonumber \\
	&\left. +2 \left( \bar{x}_1(s_i)-x_1(s_i) \right) 
		\left( x_1(s_i)-e^{-\lambda \dlnt} x_1(s_{i-1}) \right) \right. \nonumber\\
	&\left. -2e^{-\lambda \dlnt} 
		\left( \bar{x}_1(s_{i-1})-x_1(s_{i-1})\right)
		\left( x_1(s_i)-e^{-\lambda \dlnt} x_1(s_{i-1}) \right)
	\right] \nonumber\\
	=& \frac{({\dlnt})^{-1}}{T} \frac{1}{2\sigma^2 \Xi(\lambda)}
		\sum_{i=1}^{\bN2} \left[ \left\{ \left(
		\bar{x}_1(s_i)-x_1(s_i) \right)
		-e^{-\lambda \dlnt} \left(\bar{x}_1(s_{i-1})-x_1(s_{i-1}) \right)
		\right\}^2 \right.  \label{Z-Thm4-1}\\
	&\left. +2 \left\{\bar{x}_1(s_i)-\bar{x}_1(s_{i-1})
		-\left( x_1(s_i)-x_1(s_{i-1})\right)\right\}
		\left( x_1(s_i)-e^{-\lambda \dlnt} x_1(s_{i-1})\right) 
		\right. \label{Z-Thm4-2} \\
	&\left. +2(1-e^{-\lambda \dlnt})
		\left( \bar{x}_1(s_{i-1})-x_1(s_{i-1})\right)
		\left( x_1(s_i)-e^{-\lambda \dlnt} x_1(s_{i-1}) \right)
	\right]  \label{Z-Thm4-3} \\
	=:& W_1+W_2+W_3. \nonumber 
\end{align}

For the evaluation of (\ref{Z-Thm4-1}), 
we set that
\begin{eqnarray*}
g_1(t, y, \eta) &=& X_t(y)\sqrt{2} \sin (\pi y) 
	\exp\left\{\frac{\eta}{2}y \right\}.
\end{eqnarray*}
Noting that
\begin{eqnarray*}	
x_1(t) &=& \int_0^1 X_t(y)\sqrt{2} \sin (\pi y) 
	\exp\left\{\frac{\eta}{2}y \right\} dy =  \int_0^1 g_1(t, y, \eta) dy,
	\\
\bar{x}_1(t) &=& \frac{1}{M} \sum_{j=1}^M X_t(y_j) \sqrt{2} \sin (\pi y_j) \exp \left\{ \frac{\hat{\eta}}{2} y_j \right\} 
=  \frac{1}{M} \sum_{j=1}^M g_1(t, y_j, \hat{\eta}),
\end{eqnarray*}
we have that 
\begin{eqnarray*}
	Z_1 &:=& \frac{1}{T} \frac{1}{\dlnt} \sum_{i=1}^{\bN2} 
		\left( x_1(s_i)-\bar{x}_1(s_i) \right)^2  \\
	&=& \frac{1}{T} \frac{1}{\dlnt} \sum_{i=1}^{\bN2}
		\left\{
			M \frac{1}{M} \sum_{j=1}^M \int_{\frac{j-1}{M}}^{\frac{j}{M}}
			\{g_1(s_i,y,\eta)-g_1(s_i,y_j,\hat{\eta})\}dy
		\right\}^2  \\
	&\leq& \frac{1}{T} \frac{1}{\dlnt} \sum_{i=1}^{\bN2} M^2 \frac{1}{M}
		\sum_{j=1}^{M} \frac{1}{M} \int_{\frac{j-1}{M}}^{\frac{j}{M}}
		\{g_1(s_i,y,\eta)-g_1(s_i,y_j,\hat{\eta})\}^2 dy  \\
	&=& \frac{1}{T} \frac{1}{\dlnt} \sum_{i=1}^{\bN2} \sum_{j=1}^M
		\int_{\frac{j-1}{M}}^{\frac{j}{M}} \{g_1(s_i,y,\eta)-g_1(s_i,y_j,\hat{\eta})\}^2 dy.
\end{eqnarray*}
Moreover, 
\begin{eqnarray}
	& & g_1(s_i,y,\eta)-g_1(s_i,y_j,\hat{\eta}) \nonumber\\
	&=& X_{s_i}(y)\sqrt{2}\sin(\pi y) \exp \left\{ \frac{\eta}{2} y\right\}
		-X_{s_i}(y_j)\sqrt{2}\sin(\pi y_j) \exp \left\{ \frac{\hat{\eta}}{2} y_j
		\right\} \nonumber\\
	&=& \left( X_{s_i}(y)-X_{s_i}(y_j) \right) \sqrt{2} \sin(\pi y) 
		\exp \left\{ \frac{\eta}{2} y\right\}  \label{Z-thm4-i} \\
	& &+ X_{s_i}(y_j) \left( \sqrt{2}\sin(\pi y) 
		\exp\left\{ \frac{\eta}{2} y \right\} -\sqrt{2}\sin(\pi y_j) 
		\exp\left\{ \frac{\eta}{2} y_j\right\} \right) \label{Z-thm4-ii}\\
	& &+ X_{s_i}(y_j) \sqrt{2}\sin(\pi y_j)
		\left( \exp\left\{\frac{\eta}{2} y_j \right\} 
		-\exp\left\{\frac{\hat{\eta}}{2} y_j\right\}\right). \label{Z-thm4-iii} \\
	&=& \bar{g}^{(1)}(s_i,y) + \bar{g}^{(2)}(s_i,y) + \bar{g}^{(3)}(s_i,y_j,\hat{\eta}) \nonumber
\end{eqnarray}

Set 
$
R(y_j, \hat{\eta})
:=\int_0^1 \frac{y_j}{2} \exp \left\{ \frac{y_j}{2} (\eta + u ( \hat{\eta} - \eta)) \right\} du
$.
Let  $\delta_1>0$.
Since
on $J=\{|\hat{\eta}-\eta| < \delta_1 \}$
\begin{align*}
&	Z_2 :=	\frac{({\dlnt})^{-1}}{T} \sum_{i=1}^{\bN2} \frac{1}{M}
		\sum_{j=1}^M X_{s_i}^2(y_j)2\sin^2(\pi y_j) 
(R(y_j, \hat{\eta}))^2 \frac{1}{N \bar{m}}
		\leq C_1 \frac{{\bN2}^{2}}{T^2} \frac{1}{N \bar{m}},
\end{align*}
we obtain that 
\begin{align}
	P(|Z_2| > \varepsilon) &= P(|Z_2| > \varepsilon \cap J)+P(|Z_2| > 
		\varepsilon \cap J^c) 
		\leq C_1 \frac{{\bN2}^{2}}{T^2} \frac{1}{N \bar{m}} \frac{1}{\varepsilon}+o(1). \label{Z-Res1}
\end{align}

\noindent
It follows that
\begin{eqnarray*} 
  E\left[ (\bar{g}^{(1)}(s_i,y) )^2 \times  \frac{\bN2^2}{T^2} \right] &\leq& 
 C_1 E[ \left( X_{s_i}(y)-X_{s_i}(y_j)\right)^2]  \frac{\bN2^2}{T^2}
		\leq \frac{C_2}{M^{1-\rho_1}} \frac{\bN2^2}{T^2},
\\
   E \left[ (\bar{g}^{(2)}(s_i,y) )^2 \times \frac{\bN2^2}{T^2}  \right] 
   &\leq& C_1 (y-y_j)^2 \frac{\bN2^2}{T^2}  
		\leq \frac{C_1}{M^2} \frac{\bN2^2}{T^2},
\\
\frac{({\dlnt})^{-1}}{T} \sum_{i=1}^{\bN2} ( \bar{g}^{(3)}(s_i,y_j,\hat{\eta}))^2 
&=&
		\frac{({\dlnt})^{-1}}{T} \sum_{i=1}^{\bN2} \frac{1}{M}
		\sum_{j=1}^M X_{s_i}^2(y_j)2\sin^2(\pi y_j) 
(R(y_j, \hat{\eta}))^2 \frac{1}{N \bar{m}} 
		\left(\sqrt{N \bar{m}}(\hat{\eta}-\eta) \right)^{2} 
		\\
	&=& Z_2 \left(\sqrt{N \bar{m}}(\hat{\eta}-\eta) \right)^{2}  = O_p\left( \frac{{\bN2}^{2}}{T^2} \frac{1}{N \bar{m}} \right),
\end{eqnarray*}
where we use  (\ref{Z-Res1}) for the last estimate.  

Hence
\begin{eqnarray*}
	Z_1 &=& O_p \left(\frac{{\bN2}^2}{T^2 M^{1-\rho_1}} \right)
		+O_p\left(\frac{{\bN2}^2}{T^2 N \bar{m}}\right), \\
	W_1 &=& O_p \left(\frac{{\bN2}^2}{T^2 M^{1-\rho_1}} \right)
		+O_p\left(\frac{{\bN2}^2}{T^2 N \bar{m}}\right). 
\end{eqnarray*}

For the evaluation of \eqref{Z-Thm4-3},  one has that
\begin{align*}
	W_3^2 \leq& C \frac{({\dlnt})^{-1}}{T} ({\dlnt})^2
		\frac{({\dlnt})^{-1}}{T}\sum_{i=1}^{\bN2}
		\left( \bar{x}_1(s_i)-x_1(s_{i-1}) \right)^2 \sum_{i=1}^{\bN2} 
		\left( x_1(s_i)-e^{-\lambda \dlnt}x_1(s_{i-1}) \right)^2 \\
	=& \frac{T}{\bN2} \times \left(  
O_p \left(\frac{{\bN2}^2}{T^2 M^{1-\rho_1}} \right)
		+O_p\left(\frac{{\bN2}^2}{T^2 N \bar{m}}\right)
		 \right) 
		\\
	=& O_p \left(\frac{{\bN2}}{T M^{1-\rho_1}} \right)
		+O_p\left(\frac{{\bN2}}{T N \bar{m}}\right). 
\end{align*}

For the evaluation of \eqref{Z-Thm4-2}, 
setting that $\Delta X_{s_i}(y) = X_{s_i}(y)-X_{s_{i-1}}(y)$, 
we obtain that 
\begin{align*}
	& \frac{({\dlnt})^{-1}}{T} \sum_{i=1}^{\bN2}
		\left\{\bar{x}_1(s_i)-\bar{x}_1(s_{i-1})
		-\left(x_1(s_i)-x_1(s_{i-1})\right) \right\} 
	\left(x_1(s_i)-e^{-\lambda \dlnt}x_1(s_{i-1})\right) \nonumber \\
	&= \frac{({\dlnt})^{-1}}{T} \sum_{i=1}^{\bN2} \frac{1}{M}
		\sum_{j=1}^M 
                  \Delta X_{s_i}(y_j) 
		\sqrt{2}\sin(\pi y_j)
		\left(\exp\left\{\frac{\hat{\eta}}{2}y_j \right\}
		-\exp\left\{\frac{\eta}{2} y_j\right\}\right) 
		\left(x_1(s_i)-e^{-\lambda \dlnt}x_1(s_{i-1})\right) \nonumber \\
		&+ \frac{({\dlnt})^{-1}}{T} 
		\sum_{i=1}^{\bN2} \sum_{j=1}^{M}
		\int_{\frac{j-1}{M}}^{\frac{j}{M}}
                  \Delta X_{s_i}(y_j) 
		\left(\sqrt{2}\sin(\pi y_j)-\sqrt{2}\sin(\pi y)\right)
		\exp\left\{\frac{\eta}{2}y_j\right\} dy 
                 \left(x_1(s_i)-e^{-\lambda \dlnt}x_1(s_{i-1})\right) \nonumber \\
	&+ \frac{({\dlnt})^{-1}}{T}
		\sum_{i=1}^{\bN2} \sum_{j=1}^{M}
		\int_{\frac{j-1}{M}}^{\frac{j}{M}}
                  \Delta X_{s_i}(y_j) 
		\sqrt{2}\sin(\pi y) \left(\exp\left\{\frac{\eta}{2}y_j\right\}
		-\exp\left\{\frac{\eta}{2}y\right\} \right)dy 
		\left(x_1(s_i)-e^{-\lambda \dlnt}x_1(s_{i-1})\right) \nonumber \\
	&+ \frac{({\dlnt})^{-1}}{T}
		\sum_{i=1}^{\bN2} \sum_{j=1}^{M}
		\int_{\frac{j-1}{M}}^{\frac{j}{M}}
		\left\{ \Delta X_{s_i}(y_j) -\Delta X_{s_i}(y) \right\}		
		\sqrt{2}\sin(\pi y)\exp\left\{\frac{\eta}{2}y\right\}dy 
		\left(x_1(s_i)-e^{-\lambda \dlnt}x_1(s_{i-1})\right) \nonumber \\
	&=: \ ({\rm I})+({\rm II})+({\rm III})+({\rm IV}). \nonumber
\end{align*}

For the evaluation of (I), one has that
\begin{align}
	({\rm I})^2 \leq& \left(\frac{({\dlnt})^{-1}}{T}\right)^2 \sum_{i=1}^{\bN2}\frac{1}{M}
		\sum_{j=1}^{M} 
                \left( \Delta X_{s_i}(y_j) \right)^2 
                \left(\sqrt{2}\sin(\pi y_j) R(y_j,\hat{\eta} )  \right)^2 
		\frac{T}{N \bar{m}}  \label{Z-C1} \\
	& \times \frac{1}{T} \sum_{i=1}^{\bN2} \left(x_1(s_i)-e^{-\lambda \dlnt}x_1(s_{i-1})\right)^2 
		\times \left( \sqrt{N \bar{m}}(\hat{\eta}-\eta)\right)^{2}.  \nonumber
\end{align}
Let  $\eta_1 >0$ and $\epsilon>0$.
On $J=\{|\hat{\eta}-\eta| < \eta_1\}$,
\begin{align*}
	(\ref{Z-C1}) \leq&  C_1 \frac{({\dlnt})^{-2}}{N \bar{m}} \frac{1}{T }
	\sum_{i=1}^{\bN2}\frac{1}{M}
		\sum_{j=1}^{M} \left( X_{s_i}(y_j)-X_{s_{i-1}}(y_j) \right)^2 
	= O_p \left( {\frac{1}{ \dlnt^{\frac{5}{2}} N \bar{m}} } \right)
	= O_p \left( {\frac{{\bN2}^{\frac{5}{2}}}{T^{\frac{5}{2}} N \bar{m}}  } \right)
\end{align*}
because 
$E \left[ 
\left( 
\Delta X_{s_i}(y_j)
\right)^2 
\right] = \sqrt{{\dlnt}}=\sqrt{T/\bN2}$.
It follows that 
under
{
$\frac{\bar{N}_2^{\frac{5}{2}}}{T^{\frac{5}{2}} N \bar{m}} \to 0$, 
}

\begin{align*}
	P(|(\ref{Z-C1})| > \varepsilon) &= P(|(\ref{Z-C1})| > \varepsilon \cap J)+P(|(\ref{Z-C1})| > 
		\varepsilon \cap J^c) 
		\leq  {\frac{{\bN2}^{\frac{5}{2}}}{T^{\frac{5}{2}} N \bar{m}}  }  \frac{1}{\varepsilon}+o(1) 
		\rightarrow 0.
\end{align*}
Therefore, $({\rm I})=o_p(1)$.

For the evaluation of (II), we obtain that
\begin{align*}
	({\rm II})^2 \leq& C_1 \left( \frac{{\dlnt}^{-1}}{T}\right)^2
		\sum_{i=1}^{\bN2} \frac{1}{M}
		\sum_{j=1}^{M}
		 \left( \Delta X_{s_i}(y_j) \right)^2 
		  |y-y_j|^2 
		\sum_{i=1}^{\bN2} \left(x_1(s_i)-e^{-\lambda \dlnt}x_1(s_{i-1})\right)^2 \\
	=& O_p \left( \left(\frac{\bN2}{T^2}\right)^2
		\bN2 \sqrt{\frac{T}{\bN2}} \frac{1}{M^2} T \right)
		=O_p \left(\frac{{\bN2}^{\frac{5}{2}}}{T^{\frac{5}{2}}M^2} \right). 
\end{align*}

For the evaluation of (III), one has that 
\begin{align*}
	({\rm III})^2 \leq& C_1 \left( \frac{{\dlnt}^{-1}}{ T }\right)^2
		\sum_{i=1}^{\bN2} \frac{1}{M} \sum_{j=1}^{M}
		 \left( \Delta X_{s_i}(y_j) \right)^2 
		|y-y_j|^2
		\sum_{i=1}^{\bN2} \left(x_1(s_i)-e^{-\lambda \dlnt}x_1(s_{i-1})\right)^2 \\
		=& O_p \left(\frac{{\bN2}^{\frac{5}{2}}}{T^{\frac{5}{2}}M^2} \right). 
\end{align*}

For the evaluation of (IV),  we obtain that 
\begin{eqnarray*}
	({\rm IV})^2 &\leq&  \left( \frac{{\dlnt}^{-1}}{T}\right)^2 {\bar{N}_2}
		\sum_{i=1}^{\bN2} M^2 \frac{1}{M} \\
		& & \times \sum_{j=1}^{M} \frac{1}{M}
		\int_{\frac{j-1}{M}}^{\frac{j}{M}}
		\left( \left( \sum_{k=1}^\infty (x_k(s_i)-x_k(s_{i-1}))(e_k(y_j) -e_k(y)) \right)^2
		(x_1(s_i)-{e^{-\lambda \delta} }x_1(s_{i-1}))^2 \right) dy.
\end{eqnarray*}
By setting that
\begin{equation*}
	Z_3 := \left( \sum_{k=1}^\infty (x_k(s_i)-x_k(s_{i-1}))(e_k(y_j) -e_k(y)) \right)^2
		(x_1(s_i)-{e^{-\lambda \delta} }x_1(s_{i-1}))^2, 
\end{equation*}
it follows that 
\begin{align*}
	E[Z_3] =& \sum_{k=2}^\infty E[(x_k(s_i)-x_k(s_{i-1}))^2]
		E[(x_1(s_i)-{e^{-\lambda \delta} }x_1(s_{i-1}))^2](e_k(y_j) -e_k(y))^2 \\
	&+ E[{(x_1(s_i)-x_1(s_{i-1}))^2 (x_1(s_i)- e^{-\lambda \delta} x_1(s_{i-1}))^2} ](e_1(y_j) -e_1(y))^2 \\
	\leq& C_1 \left(  \frac{1}{M^{1-\rho_1}} \dlnt+\frac{{\dlnt}^2}{M^2} \right).
\end{align*}
Therefore, 
under
{
$\frac{\bar{N}_2^3}{T^3 M^{1-\rho_1}} \to 0$, }

\begin{align*}
	E[({\rm IV})^2] \leq& \frac{{\dlnt}^{-2}}{T^2} {\bN2^2} \cdot \frac{\dlnt}{M^{1-\rho_1}}
		= \frac{{\bN2^3}}{T^3 M^{1-\rho_1}} \rightarrow 0. 
\end{align*}
Consequently, 
under
{
$\frac{\bar{N}_2^{\frac{5}{2}}}{T^{\frac{5}{2}} N \bar{m}} \to 0$ and 
$\frac{\bar{N}_2^3}{T^3 M^{1-\rho_1}} \to 0$, }
one has that $W_2 =o_p(1)$ and 
\begin{equation*}
 \frac{1}{T} \left\{ 
		 l_{\bN2}(\lambda,\sigma^2 \ | \ {\bf \bar{x}}  )
		- l_{\bN2}(\lambda,\sigma^2 \ | \ {\bf {x}}  ) \right\} 	
= o_p(1)
\end{equation*}
uniformly in $(\lambda, \sigma^2)$, which completes the proof of consistency of 
$(\hat{\lambda}, \check{\sigma}^2)$.

Next, we will show the asymptotic normality of $(\hat{\lambda}, \check{\sigma}^2)$.
The derivatives of the quasi  log-likelihood function with respect to the parameters are as follows.
\begin{align*}
	\partial_\lambda l_{\bar{N}_2}(\lambda, \sigma^2 \ | \ {\bf \bar{x}}  ) 
		=& -\frac{1}{2} \sum_{i=1}^{\bar{N}_{2}}
		\left\{
			\frac{\partial_{\lambda}\Xi(\lambda)}{\Xi(\lambda)}
			-\frac{(\partial_\lambda \Xi(\lambda)) \left(\bar{x}_1(s_i)-e^{-\lambda \dlnt} \bar{x}_1(s_{i-1})\right)^2}{\sigma^2 \Xi(\lambda)^2 \dlnt} \right.\\
			& \left. +\frac{2\dlnt e^{-\lambda \dlnt} \bar{x}_1(s_{i-1})\left(\bar{x}_1(s_i)-e^{-\lambda \dlnt} \bar{x}_1(s_{i-1})\right)}{\sigma^2 \Xi(\lambda) \dlnt}
    \right\}.
\end{align*}
\[
	\partial_{\sigma^2}l_{\bar{N}_2}(\lambda, \sigma^2 \ | \ {\bf \bar{x}}  ) 
		=-\frac{1}{2} \sum_{i=1}^{\bar{N}_2} \left\{
		\frac{1}{\sigma^2}
		-\frac{1}{\sigma^4}
		\frac{\left(\bar{x}_1(s_i)-e^{-\lambda \dlnt}\bar{x}_1(s_{i-1})\right)^2}{\Xi(\lambda)\dlnt}
	\right\}.
\]
\begin{align*}
	\partial_{\lambda}^2 l_{\bar{N}_2}(\lambda, \sigma^2 \ | \ {\bf \bar{x}}  ) 
		=&-\frac{1}{2}\sum_{i=1}^{\bN2} \left\{ \partial_\lambda\left(
		\frac{\partial_\lambda\Xi(\lambda)}{\Xi(\lambda)}\right)
		-\partial_\lambda \left(
			\frac{\partial_\lambda \Xi(\lambda)}{\Xi(\lambda)^2}
		\right)
		\frac{\left(\bar{x}_1(s_i)-e^{-\lambda \dlnt} \bar{x}_1(s_{i-1})\right)^2}{\sigma^2 \dlnt} 
	\right. \\
		&\left. 
		+\left( \frac{\partial_\lambda \Xi(\lambda)}{\Xi(\lambda)^2} \right)
		\frac{2 \dlnt e^{-\lambda \dlnt} \bar{x}_1(s_{i-1})
		\left(\bar{x}_1(s_i)-e^{-\lambda \dlnt}\bar{x}_1(s_{i-1})\right)} {\sigma^2 \dlnt} 
	\right.\\
		&\left. -\frac{\partial_\lambda\Xi(\lambda)}{\Xi(\lambda)^2}
		\frac{2 \dlnt e^{-\lambda \dlnt} \bar{x}_1(s_{i-1})
		\left(\bar{x}_1(s_i)-e^{-\lambda \dlnt}\bar{x}_1(s_{i-1})\right)}{\sigma^2 \dlnt}
	\right.\\
		&\left. -\frac{1}{\Xi(\lambda)}
		\frac{2 {\dlnt}^2 e^{-\lambda \dlnt} \bar{x}_1(s_{i-1}) 
		\left(\bar{x}_1(s_i)-e^{-\lambda \dlnt} \bar{x}_1(s_{i-1}) \right)}{\sigma^2 \dlnt}
	\right.\\
		&\left. -\frac{1}{\Xi(\lambda)} 
		\frac{2 {\dlnt}^2 e^{-2 \lambda \dlnt}\bar{x}_1^2 (s_{i-1})}{\sigma^2 \dlnt}	
\right\}.
\end{align*}
\[
	\partial_{\sigma^2}^2 l_{\bar{N}_2}(\lambda, \sigma^2 \ | \ {\bf \bar{x}}  ) 
		=-\frac{1}{2} \sum_{i=1}^{\bN2} \left\{\frac{-1}{\sigma^4}
		+\frac{2}{\sigma^6}
		\frac{\left( \bar{x}_1(s_i)-e^{-\lambda \dlnt} \bar{x}_1(s_{i-1}) \right)^2}{\Xi(\lambda) \dlnt}
\right\}.
\]
\begin{align*}
	\partial_\lambda \partial_{\sigma^2} l_{\bar{N}_2}(\lambda, \sigma^2 \ | \ {\bf \bar{x}}  ) 
		=&-\frac{1}{2} \sum_{i=1}^{\bN2} \left\{\frac{-1}{\sigma^4}
		\frac{-\partial_\lambda \Xi(\lambda)}{\Xi(\lambda)^2}
		\frac{\left(\bar{x}_1(s_i)-e^{-\lambda \dlnt} \bar{x}_1(s_{i-1}) \right)^2}{\dlnt}
	\right.\\
		&\left. +\frac{-1}{\sigma^4}
		\frac{1}{\Xi(\lambda)}
		\frac{2 \dlnt e^{-\lambda \dlnt}\bar{x}_1(s_{i-1}) \left(\bar{x}_1(s_i)-e^{-\lambda \dlnt} \bar{x}_1(s_{i-1}) \right)}{\dlnt}
\right\}.
\end{align*}

\begin{en-text}
Let the discrete data of the coordinate process ${\bf x} =\{ x_1(s_{i:\bar{N}_2:T})\}_{i=1,\ldots, \bar{N}_2}$. 
The difference between the score functions of the volatility parameter $\sigma^2$ 
based on ${\bf \bar{x}}$ and  ${\bf {x}}$ 
is as follows.
\begin{align*}
	(A) :=& \frac{1}{\sqrt{\bN2}} \left\{ 
		\partial_{\sigma^2} l_{\bN2}(\lambda,\sigma^2 \ | \ {\bf \bar{x}}  )
		-\partial_{\sigma^2} l_{\bN2}(\lambda,\sigma^2 \ | \ {\bf {x}}  ) \right\} \\
	=& \frac{1}{\sqrt{\bN2}} \frac{1}{2\sigma^4 \dlnt \Xi(\lambda)}
		\sum_{i=1}^{\bN2} \left\{\left(
			\bar{x}_1(s_i)-e^{-\lambda \dlnt} \bar{x}_1(s_{i-1})
		\right)^2 \right. 
	\left. -\left(x_1(s_i)-e^{-\lambda \dlnt} x_1(s_{i-1})\right)^2 \right\}. 
\end{align*}

One has that 
\begin{align*}
	& \left(	
		\bar{x}_1(s_i)-e^{-\lambda \dlnt} \bar{x}_1(s_{i-1})
	\right)^2 \\
	=& \left\{
		\bar{x}_1(s_i)-x_1(s_i)
		-e^{-\lambda \dlnt}(\bar{x}_1(s_i)-x_1(s_{i-1}))
		+(x_1(s_i)-e^{-\lambda \dlnt} x_1(s_{i-1})) 
	\right\}^2 \\
	=& \left\{
		(\bar{x}_1(s_i)-x_1(s_i))
		-e^{-\lambda \dlnt}(\bar{x}_1(s_{i-1})-x_1(s_{i-1})) \right\}^2 \\
	&+2\{ \bar{x}_1(s_i)-x_1(s_i)
		-e^{-\lambda \dlnt}(\bar{x}_1(s_{i-1})-x_1(s_{i-1})) \}
		(x_1(s_i)-e^{-\lambda \dlnt} x_1(s_{i-1})) \\
	&+(x_1(s_i)-e^{-\lambda \dlnt} x_1(s_{i-1}))^2.
\end{align*}
It follows that 
\begin{align}
	(A) =& \frac{({\dlnt})^{-1}}{\sqrt{\bN2}} \frac{1}{2\sigma^2\Xi(\lambda)}
		\sum_{i=1}^{\bN2} \left[\left\{
		\left( \bar{x}_1(s_i)-x_1(s_i) \right)
		-e^{-\lambda \dlnt} \left( \bar{x}_1(s_{i-1})-x_1(s_{i-1}) \right)
		\right\}^2
	\right. \nonumber \\
	&\left. +2 \left( \bar{x}_1(s_i)-x_1(s_i) \right) 
		\left( x_1(s_i)-e^{-\lambda \dlnt} x_1(s_{i-1}) \right) \right. \nonumber\\
	&\left. -2e^{-\lambda \dlnt} 
		\left( \bar{x}_1(s_{i-1})-x_1(s_{i-1})\right)
		\left( x_1(s_i)-e^{-\lambda \dlnt} x_1(s_{i-1}) \right)
	\right] \nonumber\\
	=& \frac{({\dlnt})^{-1}}{\sqrt{\bN2}} \frac{1}{2\sigma^2 \Xi(\lambda)}
		\sum_{i=1}^{\bN2} \left[ \left\{ \left(
		\bar{x}_1(s_i)-x_1(s_i) \right)
		-e^{-\lambda \dlnt} \left(\bar{x}_1(s_{i-1})-x_1(s_{i-1}) \right)
		\right\}^2 \right.  \label{Thm4-1}\\
	&\left. +2 \left\{\bar{x}_1(s_i)-\bar{x}_1(s_{i-1})
		-\left( x_1(s_i)-x_1(s_{i-1})\right)\right\}
		\left( x_1(s_i)-e^{-\lambda \dlnt} x_1(s_{i-1})\right) 
		\right. \label{Thm4-2} \\
	&\left. +2(1-e^{-\lambda \dlnt})
		\left( \bar{x}_1(s_{i-1})-x_1(s_{i-1})\right)
		\left( x_1(s_i)-e^{-\lambda \dlnt} x_1(s_{i-1}) \right)
	\right].  \label{Thm4-3}
\end{align}

For the evaluation of \eqref{Thm4-1}, 
we set that
\begin{eqnarray*}
g_1(t, y, \eta) &=& X_t(y)\sqrt{2} \sin (\pi y) 
	\exp\left\{\frac{\eta}{2}y \right\}.
\end{eqnarray*}
Noting that
\begin{eqnarray*}	
x_1(t) &=& \int_0^1 X_t(y)\sqrt{2} \sin (\pi y) 
	\exp\left\{\frac{\eta}{2}y \right\} dy =  \int_0^1 g_1(t, y, \eta) dy,
	\\
\bar{x}_1(t) &=& \frac{1}{M} \sum_{j=1}^M X_t(y_j) \sqrt{2} \sin (\pi y_j) \exp \left\{ \frac{\hat{\eta}}{2} y_j \right\} 
=  \frac{1}{M} \sum_{j=1}^M g_1(t, y_j, \hat{\eta}),
\end{eqnarray*}
we have that 
\begin{eqnarray*}
	(B) &:=& \frac{1}{\sqrt{\bN2}} \frac{1}{\dlnt} \sum_{i=1}^{\bN2} 
		\left( x_1(s_i)-\bar{x}_1(s_i) \right)^2  \\
	&=& \frac{1}{\sqrt{\bN2}} \frac{1}{\dlnt} \sum_{i=1}^{\bN2}
		\left\{
			M \frac{1}{M} \sum_{j=1}^M \int_{\frac{j-1}{M}}^{\frac{j}{M}}
			\{g_1(s_i,y,\eta)-g_1(s_i,y_j,\hat{\eta})\}dy
		\right\}^2  \\
	&\leq& \frac{1}{\sqrt{\bN2}} \frac{1}{\dlnt} \sum_{i=1}^{\bN2} M^2 \frac{1}{M}
		\sum_{j=1}^{M} \frac{1}{M} \int_{\frac{j-1}{M}}^{\frac{j}{M}}
		\{g_1(s_i,y,\eta)-g_1(s_i,y_j,\hat{\eta})\}^2 dy  \\
	&=& \frac{1}{\sqrt{\bN2}} \frac{1}{\dlnt} \sum_{i=1}^{\bN2} \sum_{j=1}^M
		\int_{\frac{j-1}{M}}^{\frac{j}{M}} \{g_1(s_i,y,\eta)-g_1(s_i,y_j,\hat{\eta})\}^2 dy.
\end{eqnarray*}
Moreover, 
\begin{eqnarray}
	& & g_1(s_i,y,\eta)-g_1(s_i,y_j,\hat{\eta}) \nonumber\\
	&=& X_{s_i}(y)\sqrt{2}\sin(\pi y) \exp \left\{ \frac{\eta}{2} y\right\}
		-X_{s_i}(y_j)\sqrt{2}\sin(\pi y_j) \exp \left\{ \frac{\hat{\eta}}{2} y_j
		\right\} \nonumber\\
	&=& \left( X_{s_i}(y)-X_{s_i}(y_j) \right) \sqrt{2} \sin(\pi y) 
		\exp \left\{ \frac{\eta}{2} y\right\}  \label{thm4-i} \\
	& &+ X_{s_i}(y_j) \left( \sqrt{2}\sin(\pi y) 
		\exp\left\{ \frac{\eta}{2} y \right\} -\sqrt{2}\sin(\pi y_j) 
		\exp\left\{ \frac{\eta}{2} y_j\right\} \right) \label{thm4-ii}\\
	& &+ X_{s_i}(y_j) \sqrt{2}\sin(\pi y_j)
		\left( \exp\left\{\frac{\eta}{2} y_j \right\} 
		-\exp\left\{\frac{\hat{\eta}}{2} y_j\right\}\right). \label{thm4-iii}
\end{eqnarray}

Set 
$
R(y_j, \hat{\eta})
:=\int_0^1 \frac{y_j}{2} \exp \left\{ \frac{y_j}{2} (\eta + u ( \hat{\eta} - \eta)) \right\} du
$.
Let  $\delta_1>0$.
Since
on $J=\{|\hat{\eta}-\eta| < \delta_1 \}$
\begin{align*}
&	(D) :=	\frac{({\dlnt})^{-1}}{\sqrt{\bN2}} \sum_{i=1}^{\bN2} \frac{1}{M}
		\sum_{j=1}^M X_{s_i}^2(y_j)2\sin^2(\pi y_j) 
(R(y_j, \hat{\eta}))^2
		\leq C_1 \frac{{\bN2}^{\frac{3}{2}}}{T},
\end{align*}
we obtain that 
\begin{align}
	P(|({\rm D})| > \varepsilon) &= P(|({\rm D})| > \varepsilon \cap J)+P(|({\rm D})| > 
		\varepsilon \cap J^c) 
		\leq C_1 \frac{{\bN2}^{\frac{3}{2}}}{T} \frac{1}{\varepsilon}+o(1). \label{Res1}
\end{align}

\noindent
It follows that
\begin{eqnarray*} 
  E[(\ref{thm4-i})^2 \times \sqrt{\bN2} ({\dlnt})^{-1}] &\leq& 
 C_1 E[ \left( X_{s_i}(y)-X_{s_i}(y_j)\right)^2] \sqrt{\bN2}  ({\dlnt})^{-1} 
		\leq \frac{C_2}{M^{1-\rho_1}} \frac{\bN2^{\frac{3}{2}}}{T},
\\
   E[(\ref{thm4-ii})^2 \times \sqrt{\bN2} ({\dlnt})^{-1}] &\leq& C_1 (y-y_j)^2 \sqrt{\bN2} ({\dlnt})^{-1} 
		\leq \frac{C_1}{M^2} \frac{\bN2^{\frac{3}{2}}}{T},
\\
\frac{({\dlnt})^{-1}}{\sqrt{\bN2}} \sum_{i=1}^{\bN2} (\ref{thm4-iii})^2 
&=&
		\frac{({\dlnt})^{-1}}{\sqrt{\bN2}} \sum_{i=1}^{\bN2} \frac{1}{M}
		\sum_{j=1}^M X_{s_i}^2(y_j)2\sin^2(\pi y_j) 
(R(y_j, \hat{\eta}))^2 \frac{1}{N \bar{m}} 
		\left(\sqrt{N \bar{m}}(\hat{\eta}-\eta) \right)^{2} 
		\\
	&=& O_p\left( \frac{{\bN2}^{\frac{3}{2}}}{T} \frac{1}{N \bar{m}} \right),
\end{eqnarray*}
where we use  (\ref{Res1}) for the last estimate.  

Hence
\begin{equation*}
	(B) = O_p \left(\frac{{\bN2}^{\frac{3}{2}}}{T M^{1-\rho_1}} \right)
		+O_p\left(\frac{{\bN2}^{\frac{3}{2}}}{T N \bar{m}}\right), \quad 
	\eqref{Thm4-1} = O_p \left(\frac{{\bN2}^{\frac{3}{2}}}{T M^{1-\rho_1}} \right)
		+O_p\left(\frac{{\bN2}^{\frac{3}{2}}}{T N \bar{m}}\right). 
\end{equation*}

For the evaluation of \eqref{Thm4-3},  one has that
\begin{align*}
	\eqref{Thm4-3}^2 \leq& \frac{({\dlnt})^{-1}}{\sqrt{\bN2}} ({\dlnt})^2
		\frac{({\dlnt})^{-1}}{\sqrt{\bN2}}\sum_{i=1}^{\bN2}
		\left( \bar{x}_1(s_i)-x_1(s_{i-1}) \right)^2 \sum_{i=1}^{\bN2} 
		\left( x_1(s_i)-e^{-\lambda \dlnt}x_1(s_{i-1}) \right)^2 \\
	=& \frac{T^2}{\sqrt{\bN2} \bN2} \times \left(  O_p \left(\frac{{\bN2}^{\frac{3}{2}}}{T M^{1-\rho_1}} \right)
		+O_p\left(\frac{{\bN2}^{\frac{3}{2}}}{T N \bar{m}}\right) \right) 
		\\
	=& O_p\left(\frac{T}{M^{1-\rho_1}}\right)+O_p\left(\frac{T}{N \bar{m}}\right). 
\end{align*}

For the evaluation of \eqref{Thm4-2}, 
setting that $\Delta X_{s_i}(y) = X_{s_i}(y_j)-X_{s_{i-1}}(y_j)$, 
we obtain that 
\begin{align*}
	\eqref{Thm4-2} =& \frac{({\dlnt})^{-1}}{\sqrt{\bN2}} \sum_{i=1}^{\bN2}
		\left\{\bar{x}_1(s_i)-\bar{x}_1(s_{i-1})
		-\left(x_1(s_i)-x_1(s_{i-1})\right) \right\} 
	\left(x_1(s_i)-e^{-\lambda \dlnt}x_1(s_{i-1})\right) \nonumber \\
	=& \frac{({\dlnt})^{-1}}{\sqrt{\bN2}} \sum_{i=1}^{\bN2} \frac{1}{M}
		\sum_{j=1}^M 
                  \Delta X_{s_i}(y_j) 
		\sqrt{2}\sin(\pi y_j)
		\left(\exp\left\{\frac{\hat{\eta}}{2}y_j \right\}
		-\exp\left\{\frac{\eta}{2} y_j\right\}\right) 
		\left(x_1(s_i)-e^{-\lambda \dlnt}x_1(s_{i-1})\right) \nonumber \\
		&+ \frac{({\dlnt})^{-1}}{\sqrt{\bN2}} 
		\sum_{i=1}^{\bN2} \sum_{j=1}^{M}
		\int_{\frac{j-1}{M}}^{\frac{j}{M}}
                  \Delta X_{s_i}(y_j) 
		\left(\sqrt{2}\sin(\pi y_j)-\sqrt{2}\sin(\pi y)\right)
		\exp\left\{\frac{\eta}{2}y_j\right\} dy 
                 \left(x_1(s_i)-e^{-\lambda \dlnt}x_1(s_{i-1})\right) \nonumber \\
	&+ \frac{({\dlnt})^{-1}}{\sqrt{\bN2}}
		\sum_{i=1}^{\bN2} \sum_{j=1}^{M}
		\int_{\frac{j-1}{M}}^{\frac{j}{M}}
                  \Delta X_{s_i}(y_j) 
		\sqrt{2}\sin(\pi y) \left(\exp\left\{\frac{\eta}{2}y_j\right\}
		-\exp\left\{\frac{\eta}{2}y\right\} \right)dy 
		\left(x_1(s_i)-e^{-\lambda \dlnt}x_1(s_{i-1})\right) \nonumber \\
	&+ \frac{({\dlnt})^{-1}}{\sqrt{\bN2}}
		\sum_{i=1}^{\bN2} \sum_{j=1}^{M}
		\int_{\frac{j-1}{M}}^{\frac{j}{M}}
		\left\{ \Delta X_{s_i}(y_j) -\Delta X_{s_i}(y) \right\}		
		\sqrt{2}\sin(\pi y)\exp\left\{\frac{\eta}{2}y\right\}dy 
		\left(x_1(s_i)-e^{-\lambda \dlnt}x_1(s_{i-1})\right) \nonumber \\
	&=: \ ({\rm I})+({\rm II})+({\rm III})+({\rm IV}). \nonumber
\end{align*}

For the evaluation of (I), one has that
\begin{align}
	({\rm I})^2 \leq& \left(\frac{({\dlnt})^{-1}}{\sqrt{\bN2}}\right)^2 \sum_{i=1}^{\bN2}\frac{1}{M}
		\sum_{j=1}^{M} 
                \left( \Delta_{s_i}(y_j) \right)^2 
                \left(\sqrt{2}\sin(\pi y_j) R(y_j,\hat{\eta} )  \right)^2 
		\frac{T}{N \bar{m}}  \label{C1} \\
	& \times \frac{1}{T} \sum_{i=1}^{\bN2} \left(x_1(s_i)-e^{-\lambda \dlnt}x_1(s_{i-1})\right)^2 
		\times \left( \sqrt{N \bar{m}}(\hat{\eta}-\eta)\right)^{2}.  \nonumber
\end{align}
Let  $\delta_1>0$.
On $J=\{|\hat{\eta}-\eta| < \delta_1\}$,
\begin{align*}
	(\ref{C1}) \leq&  C_1 \frac{({\dlnt})^{-2}}{N \bar{m}} \frac{T}{\bN2 }
	\sum_{i=1}^{\bN2}\frac{1}{M}
		\sum_{j=1}^{M} \left( X_{s_i}(y_j)-X_{s_{i-1}}(y_j) \right)^2 
	=& O_p \left( \frac{T}{({\dlnt})^{\frac{3}{2}} N \bar{m}} \right)
	= O_p \left( \frac{{\bN2}^{\frac{3}{2}}}{T^{\frac{1}{2}} N \bar{m}} \right)
\end{align*}
because 
$E \left[ 
\left( 
\Delta X_{s_i}(y_j)
\right)^2 
\right] = \sqrt{{\dlnt}}=\sqrt{T/\bN2}$.
It follows that under 
$\frac{{\bN2}^{\frac{3}{2}}}{T^{\frac{1}{2}} N \bar{m}} \rightarrow 0$, 
\begin{align*}
	P(|(\ref{C1})| > \varepsilon) &= P(|(\ref{C1})| > \varepsilon \cap J)+P(|(\ref{C1})| > 
		\varepsilon \cap J^c) 
		\leq \frac{{\bN2}^{\frac{3}{2}}}{T^{\frac{1}{2}}N \bar{m}} \frac{1}{\varepsilon}+o(1) 
		\rightarrow 0.
\end{align*}
Therefore, $({\rm I})=o_p(1)$.

For the evaluation of (II), we obtain that
\begin{align*}
	({\rm II})^2 \leq& C_1 \left( \frac{{\dlnt}^{-1}}{\sqrt{\bN2}}\right)^2
		\sum_{i=1}^{\bN2} \frac{1}{M}
		\sum_{j=1}^{M}
		 \left( \Delta X_{s_i}(y_j) \right)^2 
		  |y-y_j|^2 
		\sum_{i=1}^{\bN2} \left(x_1(s_i)-e^{-\lambda \dlnt}x_1(s_{i-1})\right)^2 \\
	=& O_p \left( \left(\frac{\sqrt{\bN2}}{T}\right)^2
		\bN2 \sqrt{\frac{T}{\bN2}} \frac{1}{M^2} T \right)
		=O_p \left(\frac{{\bN2}^{\frac{3}{2}}}{T^{\frac{1}{2}}M^2} \right). 
\end{align*}

For the evaluation of (III), one has that 
\begin{align*}
	({\rm III})^2 \leq& C_1 \left( \frac{{\dlnt}^{-1}}{\sqrt{\bN2}}\right)^2
		\sum_{i=1}^{\bN2} \frac{1}{M} \sum_{j=1}^{M}
		 \left( \Delta X_{s_i}(y_j) \right)^2 
		|y-y_j|^2
		\sum_{i=1}^{\bN2} \left(x_1(s_i)-e^{-\lambda \dlnt}x_1(s_{i-1})\right)^2 \\
		=& O_p \left(\frac{{\bN2}^{\frac{3}{2}}}{T^{\frac{1}{2}}M^2} \right). 
\end{align*}

For the evaluation of (IV),  we obtain that 
\begin{eqnarray*}
	({\rm IV})^2 &\leq&  \left( \frac{{\dlnt}^{-1}}{\sqrt{\bN2}}\right)^2
		\sum_{i=1}^{\bN2} M^2 \frac{1}{M} \\
		& & \times \sum_{j=1}^{M} \frac{1}{M}
		\int_{\frac{j-1}{M}}^{\frac{j}{M}}
		\left( \left( \sum_{k=1}^\infty (x_k(s_i)-x_k(s_{i-1}))(e_k(y_j) -e_k(y)) \right)^2
		(x_1(s_i)-x_1(s_{i-1}))^2 \right) dy.
\end{eqnarray*}
By setting that
\begin{equation*}
	({\rm E}) := \left( \sum_{k=1}^\infty (x_k(s_i)-x_k(s_{i-1}))(e_k(y_j) -e_k(y)) \right)^2
		(x_1(s_i)-x_1(s_{i-1}))^2, 
\end{equation*}
it follows that 
\begin{align*}
	E[({\rm E})] =& \sum_{k=2}^\infty E[(x_k(s_i)-x_k(s_{i-1}))^2]
		E[(x_1(s_i)-x_1(s_{i-1}))^2](e_k(y_j) -e_k(y))^2 \\
	&+ E[(x_1(s_i)-x_1(s_{i-1}))^4](e_1(y_j) -e_1(y))^2 \\
	\leq& C_1 \left(  \frac{1}{M^{1-\rho_1}} \dlnt+\frac{{\dlnt}^2}{M^2} \right).
\end{align*}
Therefore, 
under 
$\frac{{\bN2}^{\frac{3}{2}}} {T M^{1-\rho_1}} \rightarrow 0$, 
\begin{align*}
	E[({\rm IV})^2] \leq& \frac{{\dlnt}^{-2}}{\bN2} \bN2 \cdot \frac{\dlnt}{M^{1-\rho_1}}
		=\frac{1}{\dlnt M^{1-\rho_1}} = \frac{\bN2}{T M^{1-\rho_1}} \rightarrow 0. 
\end{align*}
\end{en-text}

The difference between 
{the score function of the volatility parameter $\sigma^2$ 
based on ${\bf \bar{x}}$ 
and  that based on ${\bf {x}}$} 
is as follows.
\begin{align*}
	& \frac{1}{\sqrt{\bN2}} \left\{ 
		\partial_{\sigma^2} l_{\bN2}(\lambda,\sigma^2 \ | \ {\bf \bar{x}}  )
		-\partial_{\sigma^2} l_{\bN2}(\lambda,\sigma^2 \ | \ {\bf {x}}  ) \right\} \\
	=& \frac{1}{\sqrt{\bN2}} \frac{1}{2\sigma^4 \dlnt \Xi(\lambda)}
		\sum_{i=1}^{\bN2} \left\{\left(
			\bar{x}_1(s_i)-e^{-\lambda \dlnt} \bar{x}_1(s_{i-1})
		\right)^2 \right. 
	\left. -\left(x_1(s_i)-e^{-\lambda \dlnt} x_1(s_{i-1})\right)^2 \right\}. 
\end{align*}

By an analogous manner to (\ref{consistency-1}),
it is shown that
under
{
$\frac{\bar{N}_2^{\frac{3}{2}}}{T^{\frac{1}{2}} N \bar{m}} \to 0$ and 
$\frac{\bar{N}_2^2}{T M^{1-\rho_1}} \to 0$, }
\begin{equation*}
 \frac{1}{\sqrt{\bN2}} \left\{ 
		\partial_{\sigma^2} l_{\bN2}(\lambda,\sigma^2 \ | \ {\bf \bar{x}}  )
		-\partial_{\sigma^2} l_{\bN2}(\lambda,\sigma^2 \ | \ {\bf {x}}  ) \right\} 	
= o_p(1)
\end{equation*}
uniformly in $(\lambda, \sigma^2)$.

The difference between 
{the score function of the drift parameter $\lambda$ 
based on ${\bf \bar{x}}$ and  
that based on ${\bf {x}}$} 
is as follows.
\begin{align*}
	F :=& \frac{1}{\sqrt{T}} \left\{
		\partial_\lambda l_{\bN2}(\lambda,\sigma^2 \ | \ {\bf \bar{x}}  )
		-\partial_\lambda l_{\bN2}(\lambda,\sigma^2 \ | \ {\bf {x}}  ) \right\} \\
	=&\left(-\frac{1}{2}\right)\frac{1}{\sqrt{T}}\sum_{i=1}^{\bN2}
		\left[\left(-\frac{\partial_\lambda \Xi(\lambda)}{\sigma^2\Xi(\lambda)^2}\right)
		\left\{\left(\bar{x}_1(s_i)-e^{-\lambda \dlnt}\bar{x}_1(s_{i-1})\right)^2
		-\left(x_1(s_i)-e^{-\lambda \dlnt}x_1(s_{i-1}) \right)^2 \right\} {\dlnt}^{-1} \right.\\
	&+2 \left. \frac{e^{-\lambda \dlnt}}{\sigma^2 \Xi(\lambda)}
		\left\{\bar{x}_1(s_i) \left(\bar{x}_1(s_i)
		-e^{-\lambda \dlnt}  \bar{x}_1(s_{i-1})\right)
		-x_1(s_i)\left(x_1(s_i)-e^{-\lambda \dlnt}x_1(s_{i-1}) \right)\right\}\right] \\
	=:& F_1+F_2.
\end{align*}
By a similar way to (\ref{consistency-1}),
one has that
under
{
$\frac{\bar{N}_2^{\frac{5}{2}}}{T^{\frac{3}{2}} N \bar{m}} \to 0$ and 
$\frac{\bar{N}_2^3}{T^2 M^{1-\rho_1}} \to 0$, }
\begin{align*}
		F_1 &= o_p(1).
\end{align*}
For the evaluation of $F_2$, one has that 
\begin{align}
	F_2 =& \left( \frac{-1}{\sqrt{T}} \right) 
		\frac{e^{-\lambda \dlnt}}{\sigma^2\Xi(\lambda)}
		\sum_{i=1}^{\bN2} \left[ \left( \bar{x}_1(s_i)-x_1(s_i)  \right)
		\left( \bar{x}_1(s_i)-e^{-\lambda \dlnt}  \bar{x}_1(s_{i-1}) 
		\right) \right. \nonumber\\
	&+ \left. x_1(s_i) \left\{ \bar{x}_1(s_i)-x_1(s_i)
		-e^{-\lambda \dlnt}( \bar{x}_1(s_{i-1})-x_1(s_{i-1})) \right\}\right] \nonumber\\
	=& \left( \frac{-1}{\sqrt{T}} \right) 
		\frac{e^{-\lambda \dlnt}}{\sigma^2\Xi(\lambda)}
		\sum_{i=1}^{\bN2} \left[ \left( \bar{x}_1(s_i)-x_1(s_i)  \right) 
		\left\{ \bar{x}_1(s_i)-x_1(s_i)-e^{-\lambda \dlnt}( {\bar{x}_1(s_{i-1})}-x_1(s_{i-1})
		\right\} \right. \nonumber\\
	&+  \left( \bar{x}_1(s_i)-x_1(s_i)  \right)
		\left(x_1(s_i)-e^{-\lambda \dlnt}x_1(s_{i-1})  \right)  \nonumber\\
	&+ \left. x_1(s_i) \left\{ \bar{x}_1(s_i)-x_1(s_i)
		-e^{-\lambda \dlnt}( \bar{x}_1(s_{i-1})-x_1(s_{i-1})) \right\}\right] \nonumber
\end{align}
\begin{align}
=&  \left( \frac{-1}{\sqrt{T}} \right) 
		\frac{e^{-\lambda \dlnt}}{\sigma^2\Xi(\lambda)}
		\sum_{i=1}^{\bN2} \left[ \left( \bar{x}_1(s_i)-x_1(s_i)  \right)^2 \right.  
		\label{est-F-i}  \hspace{6.5cm} \\
	&- \left. e^{-\lambda \dlnt}\left( \bar{x}_1(s_i)-x_1(s_i)\right)
		\left( \bar{x}_1(s_{i-1})-x_1(s_{i-1})\right) \right.
		\label{est-F-ii}  \\
	&+ \left. \left( \bar{x}_1(s_i)-x_1(s_i)  \right)
		\left( x_1(s_i)-e^{-\lambda \dlnt}x_1(s_{i-1})\right)
		\right\}. \label{est-F-iii}  \\
	&+ \left. x_1(s_i) \left\{ \bar{x}_1(s_i)-x_1(s_i)
		-e^{-\lambda \dlnt}( \bar{x}_1(s_{i-1})-x_1(s_{i-1})) \right\}\right] 	\label{est-F-iv} \\
	=:& H_1 +H_2+H_3+H_4.  \nonumber 
\end{align}

For the evaluation of (\ref{est-F-i}), 
it follows from the evaluation of 
(\ref{Z-Thm4-1})
that
\begin{equation*}
	\frac{1}{\sqrt{T}\sqrt{\dlnt}} \sum_{i=1}^{\bN2}
		\left( \bar{x}_1(s_i)-x_1(s_i) \right)^2
	= O_p\left(\frac{{\bN2}^{\frac{3}{2}}}{T M^{1-\rho_1}} \right)
		+O_p\left( \frac{{\bN2}^{\frac{3}{2}}}{ T N \bar{m}} \right)
\end{equation*}
and that 
under
{
$\frac{\bar{N}_2^{\frac{5}{2}}}{T^{\frac{3}{2}} N \bar{m}} \to 0$ and 
$\frac{\bar{N}_2^3}{T^2 M^{1-\rho_1}} \to 0$, }
\begin{align*}
	|H_1| &\leq C_1 \frac{1}{\sqrt{T}}\sum_{i=1}^{\bN2}
		\left( \bar{x}_1(s_i)-x_1(s_i) \right)^2
		= O_p\left(\frac{{\bN2}}{\sqrt{T}M^{1-\rho_1}} \right)
		+O_p\left( \frac{{\bN2}}{\sqrt{T}N \bar{m}} \right) \stackrel{p}{\rightarrow} 0. 
\end{align*}

For the evaluation of (\ref{est-F-ii}), we obtain that 
under
{
$\frac{\bar{N}_2^{\frac{5}{2}}}{T^{\frac{3}{2}} N \bar{m}} \to 0$ and 
$\frac{\bar{N}_2^3}{T^2 M^{1-\rho_1}} \to 0$, }
\begin{align*}
	H_2^2 &\leq \frac{1}{T} \sum_{i=1}^{\bN2}
		\left( \bar{x}_1(s_i)-x_1(s_i) \right)^2
		\sum_{i=1}^{\bN2} \left( \bar{x}_1(s_{i-1})-x_1(s_{i-1}) \right)^2 \\
	&= O_p\left(\left(\frac{\bN2}{\sqrt{T}M^{1-\rho_1}} \right)^2\right)
		+O_p\left(\left( \frac{\bN2}{\sqrt{T}N \bar{m}} \right)^2\right) 
		\stackrel{p}{\rightarrow} 0.
\end{align*}

For the evaluation of (\ref{est-F-iii}), 
setting that 
\begin{equation*}
	\Delta x ( s_i, s_{i-1})  = x_1(s_i)-e^{-\lambda \dlnt}x_1(s_{i-1}),
\end{equation*}
one has that 
\begin{align}
	& -H_3 
	=
	 \frac{1}{\sqrt{T}}\frac{e^{-\lambda \dlnt}}{\sigma^2 \Xi(\lambda)}
		\sum_{i=1}^{\bN2} \left\{\left( \bar{x}_1(s_i)-x_1(s_i)\right)
		\left(x_1(s_i)-e^{-\lambda \dlnt}x_1(s_{i-1})\right)\right\} \nonumber
		\\
	 &= \frac{1}{\sqrt{T}}\frac{e^{-\lambda \dlnt}}{\sigma^2 \Xi(\lambda)}
	 \sum_{i=1}^{\bN2}\frac{1}{M}
		\sum_{j=1}^M X_{s_i}(y_j)\sqrt{2}\sin(\pi y_j)
		\left(\exp\left\{\frac{\hat{\eta}}{2}y_j\right\}
		-\exp\left\{\frac{\eta}{2}y_j\right\}\right)
		\Delta x ( s_i, s_{i-1})  \label{est-1} \\
	&+\frac{1}{\sqrt{T}}\frac{e^{-\lambda \dlnt}}{\sigma^2 \Xi(\lambda)}
	\sum_{i=1}^{\bN2}\sum_{j=1}^M\int_{\frac{j-1}{M}}^{\frac{j}{M}}
		X_{s_i}(y_j)\left(\sqrt{2}\sin(\pi y_j)-\sqrt{2}\sin(\pi y)\right)
		\exp\left\{\frac{\eta}{2}y_j\right\}dy  \Delta x ( s_i, s_{i-1})   \label{est-2} \\
	&+\frac{1}{\sqrt{T}}\frac{e^{-\lambda \dlnt}}{\sigma^2 \Xi(\lambda)}
	\sum_{i=1}^{\bN2}\sum_{j=1}^M\int_{\frac{j-1}{M}}^{\frac{j}{M}}
		X_{s_i}(y_j)\sqrt{2}\sin(\pi y)\left(\exp\left\{\frac{\eta}{2}y_j\right\}
		-\exp\left\{\frac{\eta}{2}y\right\}\right)dy \Delta x ( s_i, s_{i-1})   \label{est-3} \\
	&+\frac{1}{\sqrt{T}}\frac{e^{-\lambda \dlnt}}{\sigma^2 \Xi(\lambda)}
	\sum_{i=1}^{\bN2}\sum_{j=1}^M \int_{\frac{j-1}{M}}^{\frac{j}{M}}
		\left(X_{s_i}(y_j)-X_{s_i}(y)\right)\sqrt{2}\sin(\pi y)
		\exp\left\{\frac{\eta}{2}y\right\}dy \Delta x ( s_i, s_{i-1}) .  \label{est-4} 
\end{align}

For the evaluation of (\ref{est-1}), it follows that 
under
{
$\frac{\bar{N}_2^{\frac{5}{2}}}{T^{\frac{3}{2}} N \bar{m}} \to 0$, }
\begin{align*}
	(\ref{est-1})^2 \leq& \frac{1}{T}\sum_{i=1}^{\bN2}\frac{1}{M}\sum_{j=1}^M
		\left(X_{s_i}(y_j)\sqrt{2}\sin(\pi y_j) 
                 R(y_j, \hat{\eta}) 
		 (\hat{\eta}-\eta)\right)^2
		\sum_{i=1}^{\bN2} (\Delta x ( s_i, s_{i-1}) )^2 \\
	=& O_p\left(\frac{\bN2}{N \bar{m}}\right)\stackrel{p}{\rightarrow} 0.
\end{align*}

For the evaluation of (\ref{est-2}), one has that 
under
{
$\frac{\bar{N}_2^3}{T^2 M^{1-\rho_1}} \to 0$, }

\begin{align*}	
	(\ref{est-2})^2 \leq& \frac{1}{T}\sum_{i=1}^{\bN2} \sum_{j=1}^M  \int_{\frac{j-1}{M}}^{\frac{j}{M}} (X_{s_i}(y_j))^2(y-y_j)^2 dy
		\sum_{i=1}^{\bN2} (\Delta x ( s_i, s_{i-1}) )^2 \\
	=& O_p\left(\frac{\bN2}{M^2} \right) \stackrel{p}{\rightarrow} 0.
\end{align*}
It is shown that  $(\ref{est-3}) \stackrel{p}{\rightarrow} 0$ in the same way as (\ref{est-2}).

For the evaluation of (\ref{est-4}), setting that
\begin{equation*}
	G_i := \left( \sum_{k=1}^{\infty}x_k(s_i)(e_k(y)-e_k(y_j))\right)^2
		(\Delta x ( s_i, s_{i-1}) )^2, 
\end{equation*}
one has that
\begin{align*}
	E[G_i] =& \sum_{k=2}^{\infty} E[x_k^2(s_i)](e_k(y)-e_k(y_j))^2
		E\left[(x_1(s_i)-e^{-\lambda \dlnt}x_1(s_{i-1})^2\right] \\
	&+ E\left[x_1^2(s_i)(x_1(s_i)-e^{-\lambda \dlnt}x_1(s_{i-1}))^2\right](e_1(y)-e_1(y_j))^2\\
	\leq& C_1 \left( \frac{1}{M^{1-\rho_1}} \dlnt+\frac{\dlnt}{M^2} \right).
\end{align*}
Noting that
\begin{eqnarray*}
	(\ref{est-4})^2 &\leq& \frac{C_1}{T}{\bN2}^2\frac{1}{\bN2}
		\sum_{i=1}^{\bN2}\sum_{j=1}^M \int_{\frac{j-1}{M}}^{\frac{j}{M}}
		\left( \sum_{k=1}^{\infty}x_k(s_i)(e_k(y)-e_k(y_j))\right)^2 dy
		(\Delta x ( s_i, s_{i-1}) )^2 \\
		&=&  \frac{C_1}{T}{\bN2}^2\frac{1}{\bN2}
		\sum_{i=1}^{\bN2}\sum_{j=1}^M \int_{\frac{j-1}{M}}^{\frac{j}{M}}
		G_i  dy, 
\end{eqnarray*}
we obtain that 
under
{
$\frac{\bar{N}_2^3}{T^2 M^{1-\rho_1}} \to 0$, }
\begin{equation*}
	E[(\ref{est-4})^2]  \leq C_1 \frac{\bN2}{\dlnt}\frac{\dlnt}{M^{1-\rho_1}}
		= {C_1  \frac{\bN2}{M^{1-\rho_1}}} 
		\rightarrow 0.
\end{equation*}

\begin{en-text}
\begin{equation*}
	({\rm iii}) \stackrel{p}{\rightarrow}  0 \quad\mbox{under} \quad
		\frac{{\bN2}^2}{T^{\frac{3}{2}}M^{1-\rho_1}} \rightarrow 0, 
		\quad \frac{{\bN2}^2}{T^{\frac{3}{2}} N \bar{m}} \rightarrow 0
\end{equation*}
\begin{equation*}
	({\rm **})_2 \stackrel{p}{\rightarrow} 0 \quad\mbox{under} \quad
		\frac{{\bN2}^2}{T^{\frac{3}{2}}M^{1-\rho_1}} \rightarrow 0, 
		\quad\frac{{\bN2}^2}{T^{\frac{3}{2}} N \bar{m}} \rightarrow 0
\end{equation*}
\end{en-text}

\noindent
Hence,
under
{
$\frac{\bar{N}_2^{\frac{5}{2}}}{T^{\frac{3}{2}} N \bar{m}} \to 0$ and 
$\frac{\bar{N}_2^3}{T^2 M^{1-\rho_1}} \to 0$, }
$$H_3 \stackrel{p}{\rightarrow}  0.
$$

For the evaluation of (\ref{est-F-iv}), 
noting that
\begin{align}
         -H_4 =&  \frac{1}{\sqrt{T}}\frac{e^{-\lambda \dlnt}}{\sigma^2 \Xi(\lambda)}
	\sum_{i=1}^{\bN2}
		\left\{\bar{x}_1(s_i)-\bar{x}_1(s_{i-1})
		-\left(x_1(s_i)-x_1(s_{i-1})\right) \right\} 
	x_1(s_i) \label{est-F-iv-a} \\
	+&  \frac{1}{\sqrt{T}}\frac{e^{-\lambda \dlnt}}{\sigma^2 \Xi(\lambda)}
	\sum_{i=1}^{\bN2}
		\left(1-e^{\lambda \delta} \right) \left(\bar{x}_1(s_{i-1})-x_1(s_{i-1})\right)  
	x_1(s_i), \label{est-F-iv-b} 
\end{align}
one has that
under
{
$\frac{\bar{N}_2^{\frac{5}{2}}}{T^{\frac{3}{2}} N \bar{m}} \to 0$ and 
$\frac{\bar{N}_2^3}{T^2 M^{1-\rho_1}} \to 0$, }
\begin{align*}
(\ref{est-F-iv-b})^2 \leq \frac{\delta^2}{T}  \sum_{i=1}^{\bar{N}_2} (x_1(s_i))^2 
\sum_{i=1}^{\bar{N}_2} \left(\bar{x}_1(s_{i-1})-x_1(s_{i-1})\right)^2
=O_p \left( \frac{T}{M^{1-\rho_1}} \right)+ O_p \left( \frac{T}{N \bar{m}} \right)
 \stackrel{p}{\rightarrow}  0.
\end{align*}
Moreover, we set that 
\begin{align*}
	-(\ref{est-F-iv-a}) =&  \frac{1}{\sqrt{T}}\frac{e^{-\lambda \dlnt}}{\sigma^2 \Xi(\lambda)}
	 \sum_{i=1}^{\bN2} \frac{1}{M}
		\sum_{j=1}^M 
                  \Delta X_{s_i}(y_j) 
		\sqrt{2}\sin(\pi y_j)
		\left(\exp\left\{\frac{\hat{\eta}}{2}y_j \right\}
		-\exp\left\{\frac{\eta}{2} y_j\right\}\right) 
		x_1(s_i) \nonumber \\
		&+ \frac{1}{\sqrt{T}}\frac{e^{-\lambda \dlnt}}{\sigma^2 \Xi(\lambda)}
		\sum_{i=1}^{\bN2} \sum_{j=1}^{M}
		\int_{\frac{j-1}{M}}^{\frac{j}{M}}
                  \Delta X_{s_i}(y_j) 
		\left(\sqrt{2}\sin(\pi y_j)-\sqrt{2}\sin(\pi y)\right)
		\exp\left\{\frac{\eta}{2}y_j\right\} dy 
                 x_1(s_i) \nonumber \\
	&+  \frac{1}{\sqrt{T}}\frac{e^{-\lambda \dlnt}}{\sigma^2 \Xi(\lambda)}
		\sum_{i=1}^{\bN2} \sum_{j=1}^{M}
		\int_{\frac{j-1}{M}}^{\frac{j}{M}}
                  \Delta X_{s_i}(y_j) 
		\sqrt{2}\sin(\pi y) \left(\exp\left\{\frac{\eta}{2}y_j\right\}
		-\exp\left\{\frac{\eta}{2}y\right\} \right)dy 
		x_1(s_i) \nonumber \\
	&+  \frac{1}{\sqrt{T}}\frac{e^{-\lambda \dlnt}}{\sigma^2 \Xi(\lambda)}
		\sum_{i=1}^{\bN2} \sum_{j=1}^{M}
		\int_{\frac{j-1}{M}}^{\frac{j}{M}}
		\left\{ \Delta X_{s_i}(y_j) -\Delta X_{s_i}(y) \right\}		
		\sqrt{2}\sin(\pi y)\exp\left\{\frac{\eta}{2}y\right\}dy 
		x_1(s_i) \nonumber \\
	&=: \ ({\rm V})+({\rm VI})+({\rm VII})+({\rm VIII}). \nonumber
\end{align*}

For the evaluation of (V), one has that
\begin{align}
	({\rm V})^2 \leq& \frac{1}{T} \sum_{i=1}^{\bN2}\frac{1}{M}
		\sum_{j=1}^{M} 
                \left( \Delta X_{s_i}(y_j) \right)^2 
                \left(\sqrt{2}\sin(\pi y_j) R(y_j,\hat{\eta} )  \right)^2 
		\frac{\bar{N}_2}{N \bar{m}}  \label{thm4-Z-C1} \\
	& \times  \frac{1}{\bar{N}_2} \sum_{i=1}^{\bN2} \left(x_1(s_i) \right)^2 
		\times \left( \sqrt{N \bar{m}}(\hat{\eta}-\eta)\right)^{2}.  \nonumber
\end{align}
Let  $\delta_1>0$.
On $J=\{|\hat{\eta}-\eta| < \delta_1\}$,
\begin{align*}
	(\ref{thm4-Z-C1}) \leq&  C_1 \frac{\bN2}{N \bar{m}} \frac{1}{T }
	\sum_{i=1}^{\bN2}\frac{1}{M}
		\sum_{j=1}^{M} \left( X_{s_i}(y_j)-X_{s_{i-1}}(y_j) \right)^2 
	= O_p \left( \frac{{\bN2}^{\frac{3}{2}}}{T^{\frac{1}{2}} N \bar{m}} \right)
\end{align*}
because 
$E \left[ 
\left( 
\Delta X_{s_i}(y_j)
\right)^2 
\right] = \sqrt{{\dlnt}}=\sqrt{T/\bN2}$.
It follows that 
under
{
$\frac{\bar{N}_2^{\frac{5}{2}}}{T^{\frac{3}{2}} N \bar{m}} \to 0$, } 
\begin{align*}
	P(|(\ref{thm4-Z-C1})| > \varepsilon) &= P(|(\ref{thm4-Z-C1})| > \varepsilon \cap J)+P(|(\ref{thm4-Z-C1})| > 
		\varepsilon \cap J^c) 
		\leq \frac{{\bN2}^{\frac{3}{2}}}{T^{\frac{1}{2}}N \bar{m}} \frac{1}{\varepsilon}+o(1) 
		\rightarrow 0.
\end{align*}
Therefore, $({\rm V})=o_p(1)$.

For the evaluation of (VI), we obtain that
under
{
$\frac{\bar{N}_2^3}{T^2 M^{1-\rho_1}} \to 0$, }
\begin{align*}
	|({\rm VI})| \leq& C_1 \frac{1}{\sqrt{T}}
		\sum_{i=1}^{\bN2} \frac{1}{M}
		\sum_{j=1}^{M}
		 \left| \Delta X_{s_i}(y_j) x_1(s_i)  \right| 
		  |y_j-y_{j-1}| \\
	=& O_p \left( 
	        \frac{\bN2}{\sqrt{T}} 
		\left( \frac{T}{\bN2} \right)^{1/4} 
		\frac{1}{M}
		\right)
		=o_p(1). 
\end{align*}

For the evaluation of (VII), one has that 
under
{
$\frac{\bar{N}_2^3}{T^2 M^{1-\rho_1}} \to 0$, }
\begin{align*}
	|({\rm VII})| \leq& C_1 \frac{1}{\sqrt{T}}
		\sum_{i=1}^{\bN2} \frac{1}{M}
		\sum_{j=1}^{M}
		 \left| \Delta X_{s_i}(y_j) x_1(s_i)  \right| 
		  |y_j-y_{j-1}| \\
	=& O_p \left( 
	        \frac{\bN2}{\sqrt{T}} 
		\left( \frac{T}{\bN2} \right)^{1/4} 
		\frac{1}{M}
		\right)
		=o_p(1). 
\end{align*}

For the evaluation of (VIII),  we obtain that 
\begin{eqnarray*}
	({\rm VIII})^2 &\leq&  \frac{1}{T} {\bar{N}_2}
		\sum_{i=1}^{\bN2} M^2 \frac{1}{M} \\
		& & \times \sum_{j=1}^{M} \frac{1}{M}
		\int_{\frac{j-1}{M}}^{\frac{j}{M}}
		\left( \left( \sum_{k=1}^\infty (x_k(s_i)-x_k(s_{i-1}))(e_k(y_j) -e_k(y)) \right)^2
		(x_1(s_i))^2 \right) dy.
\end{eqnarray*}
By setting that
\begin{equation*}
	Z_4 := \left( \sum_{k=1}^\infty (x_k(s_i)-x_k(s_{i-1}))(e_k(y_j) -e_k(y)) \right)^2
		(x_1(s_i))^2, 
\end{equation*}
it follows that 
\begin{align*}
	E[Z_4] =& \sum_{k=2}^\infty E[(x_k(s_i)-x_k(s_{i-1}))^2]
		E[(x_1(s_i))^2](e_k(y_j) -e_k(y))^2 \\
	&+ E[(x_1(s_i)-x_1(s_{i-1}))^2 (x_1(s_i))^2 ](e_1(y_j) -e_1(y))^2 \\
	\leq& C_1 \left(  \frac{1}{M^{1-\rho_1}} +\frac{{\dlnt}}{M^2} \right).
\end{align*}
Therefore, 
under
{
$\frac{\bar{N}_2^3}{T^2 M^{1-\rho_1}} \to 0$, }
\begin{align*}
	E[({\rm VIII})^2] \leq& \frac{{\bN2^2}}{T} \cdot \frac{1}{M^{1-\rho_1}}
		= \frac{{\bN2^2}}{T M^{1-\rho_1}} \rightarrow 0. 
\end{align*}

\begin{en-text}
Consequently, 
under
$\frac{\bN2^{2}}{T^2 M^{1-\rho_1}} \rightarrow 0$ and
$\frac{\bN2^{2}}{T^{\frac{3}{2}} N \bar{m}} \rightarrow 0$,
one has that 
\begin{equation*}
 \frac{1}{T} \left\{ 
		 l_{\bN2}(\lambda,\sigma^2 \ | \ {\bf \bar{x}}  )
		- l_{\bN2}(\lambda,\sigma^2 \ | \ {\bf {x}}  ) \right\} 	
= o_p(1)
\end{equation*}
uniformly in $(\lambda, \sigma^2)$, which completes the proof pf consistency of 
$(\hat{\lambda}, \check{\sigma}^2)$.
\end{en-text}

\noindent
We obtain that
under
{
$\frac{\bar{N}_2^{\frac{5}{2}}}{T^{\frac{3}{2}} N \bar{m}} \to 0$ and 
$\frac{\bar{N}_2^3}{T^2 M^{1-\rho_1}} \to 0$, }
$$
(\ref{est-F-iv-a})  \stackrel{p}{\rightarrow} 0,  \quad  H_4  \stackrel{p}{\rightarrow} 0,
\quad F_2 \stackrel{p}{\rightarrow} 0
$$
and  
\begin{equation*}
F= \frac{1}{\sqrt{\bN2}} \left\{ 
		\partial_{\sigma^2} l_{\bN2}(\lambda,\sigma^2 \ | \ {\bf \bar{x}}  )
		-\partial_{\sigma^2} l_{\bN2}(\lambda,\sigma^2 \ | \ {\bf {x}}  ) \right\} 	
= o_p(1).
\end{equation*}

\noindent 
Furthermore,
under
{
$\frac{\bar{N}_2^{\frac{5}{2}}}{T^{\frac{3}{2}} N \bar{m}} \to 0$ and 
$\frac{\bar{N}_2^3}{T^2 M^{1-\rho_1}} \to 0$, }
\begin{eqnarray*}
& & 
\frac{1}{\bN2} 
\left\{
\partial_{\sigma^2}^2 l_{\bar{N}_2}(\lambda, \sigma^2 \ | \ {\bf \bar{x}}  ) 
-
\partial_{\sigma^2}^2 l_{\bar{N}_2}(\lambda, \sigma^2 \ | \ {\bf {x}}  ) 
\right\} = o_p(1), 
\\
& & 
\frac{1}{T} 
\left\{
\partial_{\lambda}^2 l_{\bar{N}_2}(\lambda, \sigma^2 \ | \ {\bf \bar{x}}  ) 
-
\partial_{\lambda}^2 l_{\bar{N}_2}(\lambda, \sigma^2 \ | \ {\bf {x}}  ) 
\right\} = o_p(1), 
\\
& & 
\frac{1}{\sqrt{\bN2 T}} 
\left\{
\partial_\lambda \partial_{\sigma^2} l_{\bar{N}_2}(\lambda, \sigma^2 \ | \ {\bf \bar{x}}  ) 
-
\partial_\lambda \partial_{\sigma^2} l_{\bar{N}_2}(\lambda, \sigma^2 \ | \ {\bf {x}}  ) 
\right\} = o_p(1)
\end{eqnarray*}
uniformly in $(\lambda, \sigma^2)$.
These results imply that
\begin{equation}
\begin{pmatrix}
\sqrt{\bN2}(\hat{\sigma}^2-(\sigma^*)^2) \\
\sqrt{T}(\hat{\lambda}-\lambda^*)
\end{pmatrix}
\stackrel{d}{\rightarrow} 
N \left( 
\begin{pmatrix}
0 \\
0
\end{pmatrix}, 
\begin{pmatrix}
2 (\sigma^*)^4 & 0 \\
0 & 2 \lambda^*
\end{pmatrix}
\right). \label{AN_1}
\end{equation}
For the estimator of $\theta_2$, we obtain that
\begin{eqnarray*}
\sqrt{\bN2} (\hat{\theta}_2 -\theta_2^*) 
&=& \sqrt{\bN2} \left( \left( \frac{\hat{\sigma}^2}{\hat{\sigma}_0^2} \right)^2  
- \left( \frac{({\sigma^*})^2}{({\sigma_0^*})^2} \right)^2 \right) \\
&=&
\sqrt{\bN2} 
\left(  
\hat{\sigma}^4 
\left\{
\left( \frac{1}{\hat{\sigma}_0^2} \right)^2  - \left( \frac{1}{({\sigma_0^*})^2} \right)^2 
\right\} 
+ \frac{1}{(\sigma_0^*)^4} \left\{ \hat{\sigma}^4 - ({\sigma^*})^4 \right\}  
\right)  \\
&=&
\frac{\sqrt{\bN2}}{\sqrt{\bar{m} N}}   
(\hat{\sigma}^2)^2 
\sqrt{\bar{m} N} 
\left\{
\left( \frac{1}{\hat{\sigma}_0^2} \right)^2  - \left( \frac{1}{({\sigma_0^*})^2} \right)^2 
\right\} 
+ \sqrt{\bN2}  \frac{1}{(\sigma_0^*)^4} \left\{ \hat{\sigma}^4 - ({\sigma^*})^4 \right\}    \\
&=&
\sqrt{\bN2}  \frac{1}{(\sigma_0^*)^4} \left\{ \hat{\sigma}^4 - ({\sigma^*})^4 \right\} +o_p(1).
\end{eqnarray*}
For the estimator of $\theta_1$, one has that
\begin{eqnarray*}
\sqrt{\bN2} (\hat{\theta}_1 -\theta_1^*) 
&=& \sqrt{\bN2} \left(  \hat{\eta} \hat{\theta}_2 - \eta^* \theta_2^* \right) 
= \sqrt{\bN2} \left(  \hat{\theta}_2 \left( \hat{\eta}  - \eta^* \right) + \eta^* \left( \hat{\theta}_2 -\theta_2^* \right) \right) \\
&=& \sqrt{\bN2} \eta^* \left( \hat{\theta}_2 -\theta_2^* \right) + o_p(1) \\
&=& \sqrt{\bN2}  \eta^* \frac{1}{(\sigma_0^*)^4} \left\{ (\hat{\sigma}^2)^2 - (({\sigma^*})^2)^2 \right\} +o_p(1).
\end{eqnarray*}
For the estimator of $\theta_0$, one has that
\begin{eqnarray*}
\sqrt{T} (\hat{\theta}_0 -\theta_0^*) 
&=& \sqrt{T} \left( \hat{\lambda}_1 - \lambda_1^* 
+  \frac{\left( \hat{\theta}_1 \right)^2}{  4 \hat{\theta}_2}  - \frac{\left( \theta_1 \right)^2}{  4 \theta_2^*}
+ \pi^2 (\hat{\theta}_2 -\theta_2^*) \right)
\\
&=&
 \sqrt{T} \left( \hat{\lambda}_1 - \lambda_1^* \right)
+ \frac{ \sqrt{T} \left( (\hat{\theta}_1)^2 - ({\theta_1^*})^2 \right)}{4 \hat{\theta}_2}  
+ 
\frac{(\theta_1^*)^2}{4} 
\sqrt{T}
\left( \frac{1}{\hat{\theta}_2}   - \frac{1}{\theta_2^*} \right) 
+ \pi^2 \sqrt{T} \left( \hat{\theta}_2 -\theta_2^* \right)
\\
&=&
\sqrt{T} \left( \hat{\lambda}_1 - \lambda_1^* \right) +o_p(1).
\end{eqnarray*}

\noindent
By noting that
\begin{equation*}
\begin{pmatrix}
\sqrt{\bN2}(\hat{\sigma}^2-(\sigma^*)^2) \\
\sqrt{\bN2}(\hat{\theta}_2- \theta_2^*) \\
\sqrt{\bN2}(\hat{\theta}_1- \theta_1^*) \\
\sqrt{T}(\hat{\theta}_0-\theta_0^*)
\end{pmatrix}
=
\begin{pmatrix}
\sqrt{\bN2}(\hat{\sigma}^2-(\sigma^*)^2) \\
\sqrt{\bN2}\frac{1}{(\sigma_0^*)^4}(\hat{\sigma}^4 - ({\sigma^*})^4) \\
\sqrt{\bN2}\frac{\eta^*}{(\sigma_0^*)^4}(\hat{\sigma}^4 - ({\sigma^*})^4) \\
\sqrt{T}(\hat{\lambda}-\lambda^*)
\end{pmatrix}
+o_p(1),
\end{equation*}

\noindent
it follows from (\ref{AN_1}) and the delta method that
\begin{equation*}
\begin{pmatrix}
\sqrt{\bN2}(\hat{\sigma}^2-(\sigma^*)^2) \\
\sqrt{\bN2}\frac{1}{(\sigma_0^*)^4}(\hat{\sigma}^4 - ({\sigma^*})^4) \\
\sqrt{\bN2}\frac{\eta^*}{(\sigma_0^*)^4}(\hat{\sigma}^4 - ({\sigma^*})^4) \\
\sqrt{T}(\hat{\lambda}-\lambda^*)
\end{pmatrix}
\stackrel{d}{\rightarrow} 
N \left( 
\begin{pmatrix}
0 \\
0 \\
0 \\
0
\end{pmatrix}, 
\begin{pmatrix}
2 (\sigma^*)^4 & 4 \theta_2^* (\sigma^*)^2 & 4 \theta_1^* (\sigma^*)^2 & 0 \\
 4 \theta_2^* (\sigma^*)^2 &  8 (\theta_2^*)^2 & 8 \theta_1^* \theta_2^* & 0 \\
 4 \theta_1^* (\sigma^*)^2 & 8 \theta_1^* \theta_2^* & 8 (\theta_1^*)^2 & 0 \\
0 & 0 & 0 & 2 \lambda_1^*
\end{pmatrix}
\right), \label{AN_2}
\end{equation*}
which completes the proof.

\begin{en-text}
Therefore,
\begin{align*}
	F &= \frac{1}{\sqrt{T}} \{\partial_\lambda l_{\bN2}
		(\hat{x}^{(\mu)}|\lambda,\sigma^2)
		-\partial_\lambda l_{\bN2}(x|\lambda,\sigma)\} \\
		&\stackrel{p}{\rightarrow} 0 \quad \mbox{under} 
		\frac{{\bN2}^2}{T^{\frac{3}{2}}M^{1-\rho_1}},
		\mbox{under} \frac{{\bN2}^2}{T^{\frac{3}{2}}{1-\rho_1}}\rightarrow 0,
		\frac{{\bN2}^2}{T^{\frac{3}{2}}N \bar{m}}
\end{align*}
\end{en-text}

\section{Appendix} \label{appendix}

{In order to help us understand the characteristics of the parameters
$\theta_0$, $\theta_1$, $\theta_2$ and $\sigma$ of the SPDE (\ref{spde0}), 
we can refer {some} sample paths
with different values of the parameters as follows.
}

Figure \ref{sigma11} are {the} sample paths,
{ where} $\theta_0$, $\theta_1$, and $\theta_2$ are fixed and only $\sigma$ is changed.
The shape of the sample paths does not change
{ and} only {their} height changes.

\begin{figure}[h]
\begin{minipage}{0.32\hsize}
\begin{center}
\includegraphics[width=4.5cm]{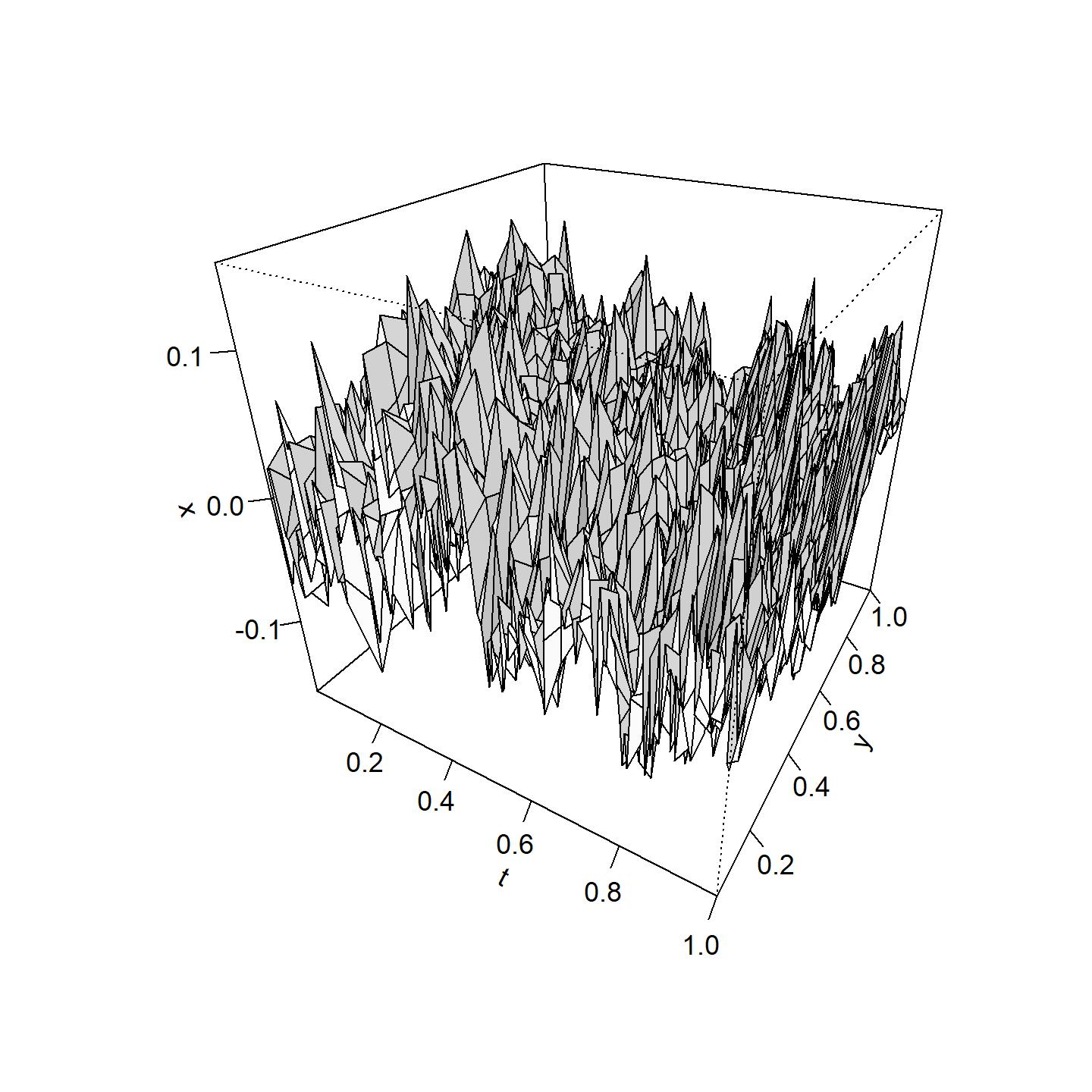}
\captionsetup{labelformat=empty,labelsep=none}
\subcaption{$\theta=$(0,0.1,0.1,0.1)}
\end{center}
\end{minipage}
\begin{minipage}{0.32\hsize}
\begin{center}
\includegraphics[width=4.5cm]{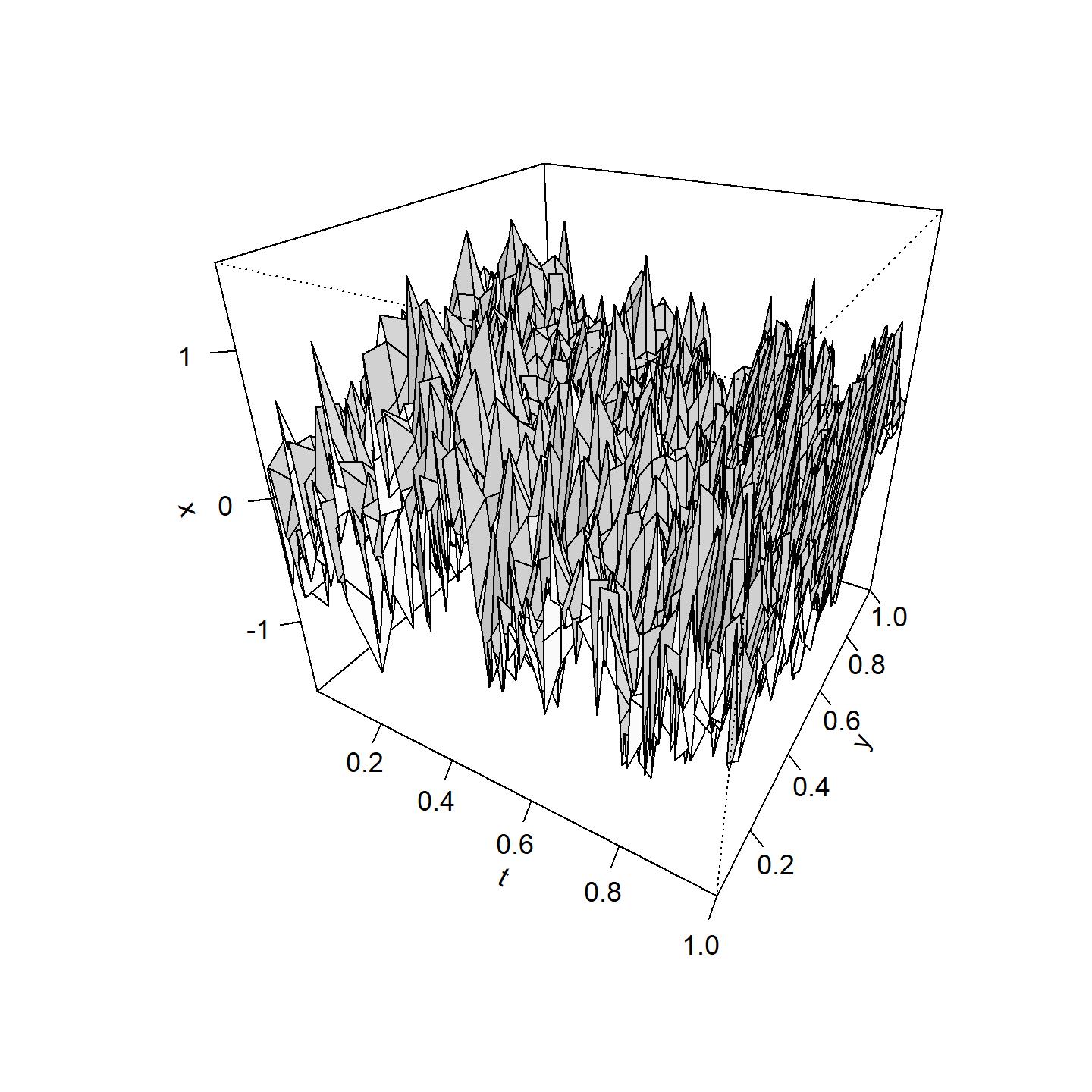}
\captionsetup{labelformat=empty,labelsep=none}
\subcaption{$\theta=$(0,0.1,0.1,1)}
\end{center}
\end{minipage}
\begin{minipage}{0.32\hsize}
\begin{center}
\includegraphics[width=4.5cm]{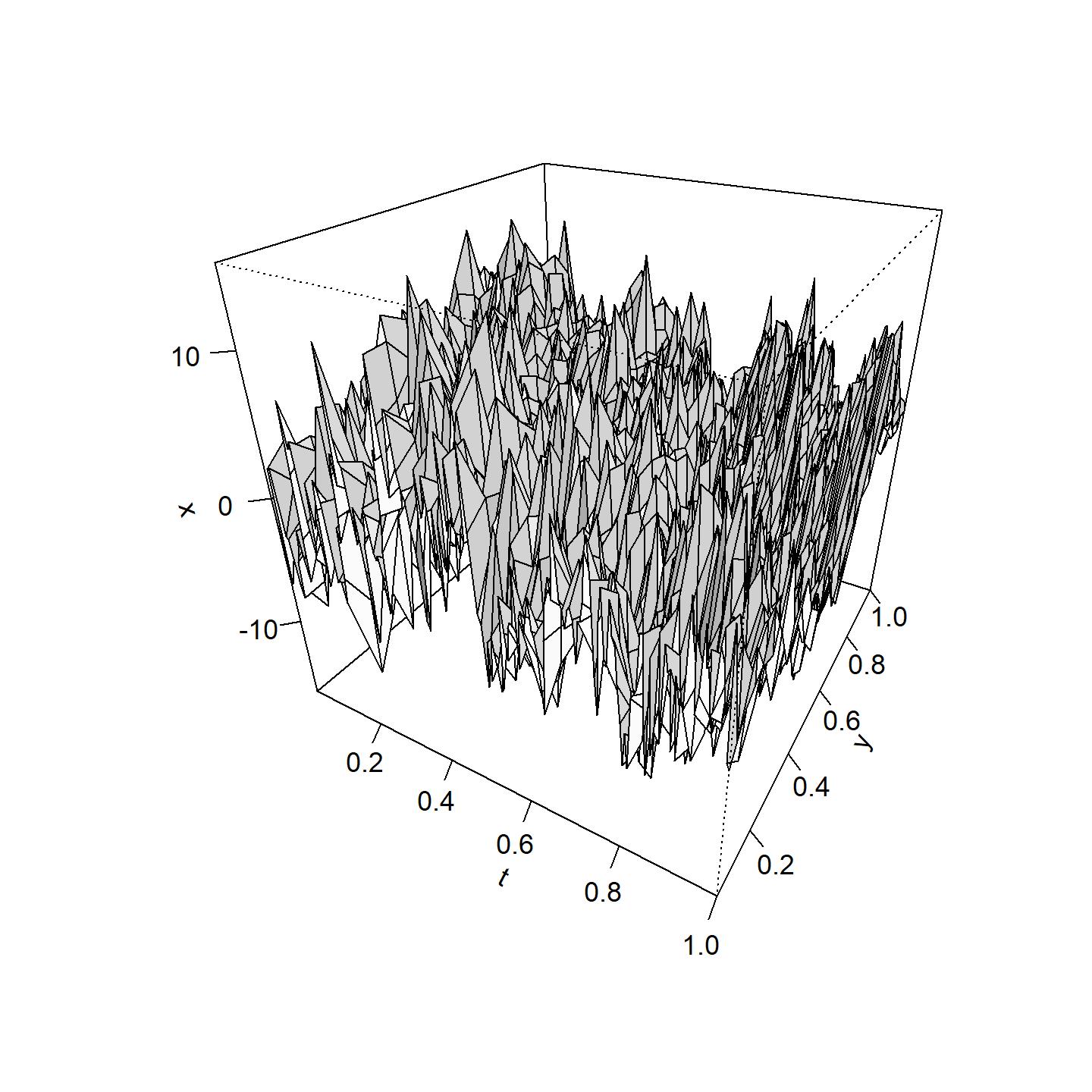}
\captionsetup{labelformat=empty,labelsep=none}
\subcaption{$\theta=$(0,0.1,0.1,10)}
\end{center}
\end{minipage}
\caption{Sample paths with $\sigma=0.1$, 1, 10}
\label{sigma11}
\end{figure}

\clearpage

Figures \ref{t1-11}-\ref{t1-16} {the} sample paths,
{ where}  $\theta_0$, $\theta_2$ and $\sigma$ are fixed and only $\theta_1$ is changed.
Figures \ref{t1-11}-\ref{t1-13} show that {the} variation of the sample path  
is large near {$y = 0$ }
and small near {$y = 1$} when $\theta_1 > 0$.
This trend increases as the value of $ \theta_1 $ increases.


\begin{figure}[H]
\begin{minipage}{0.32\hsize}
\begin{center}
\includegraphics[width=4.5cm]{p2all.jpeg}
\captionsetup{labelformat=empty,labelsep=none}
\subcaption{$\theta=$(0,0.1,0.1,1)}
\end{center}
\end{minipage}
\begin{minipage}{0.32\hsize}
\begin{center}
\includegraphics[width=4.5cm]{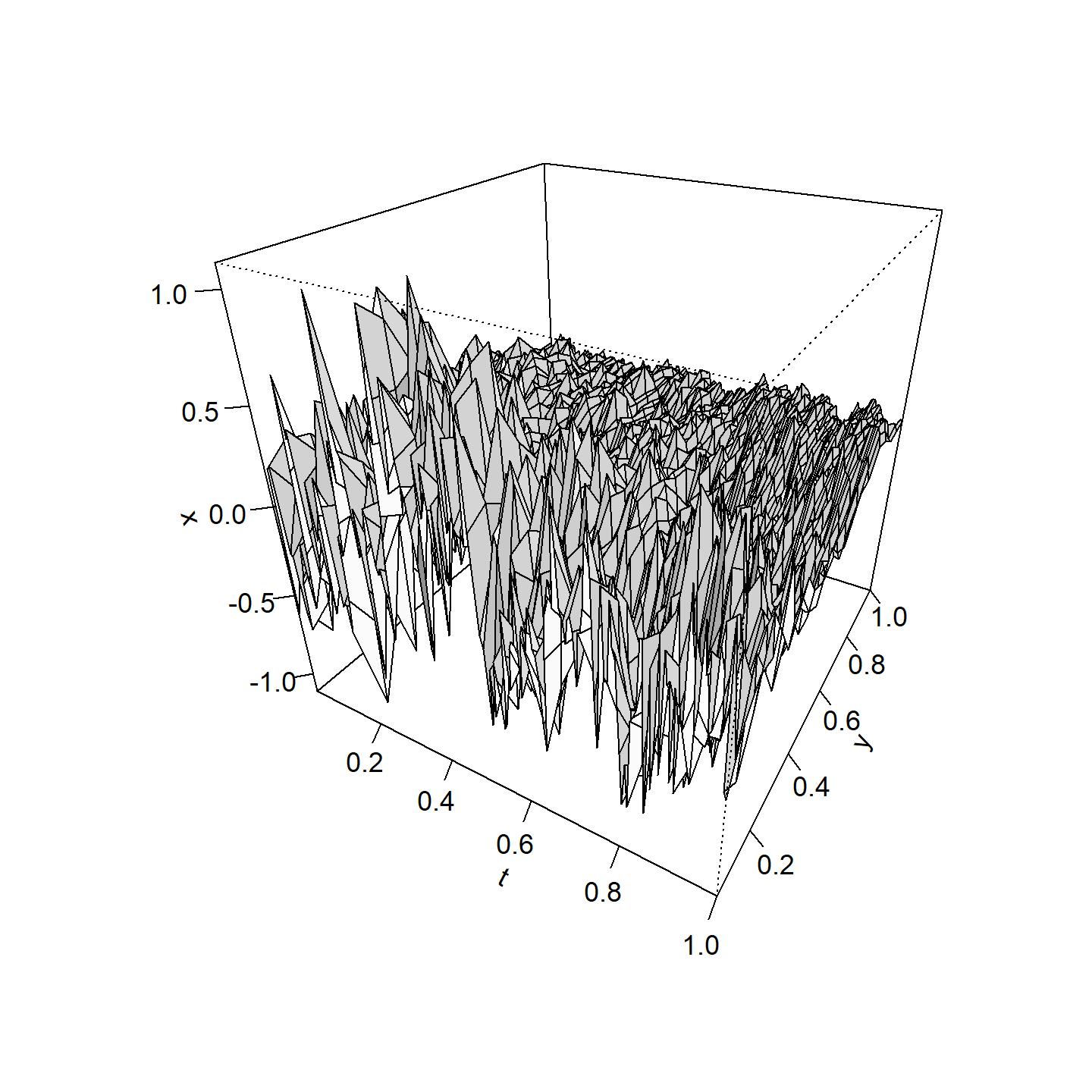}
\captionsetup{labelformat=empty,labelsep=none}
\subcaption{$\theta=$(0,0.5,0.1,1)}
\end{center}
\end{minipage}
\begin{minipage}{0.32\hsize}
\begin{center}
\includegraphics[width=4.5cm]{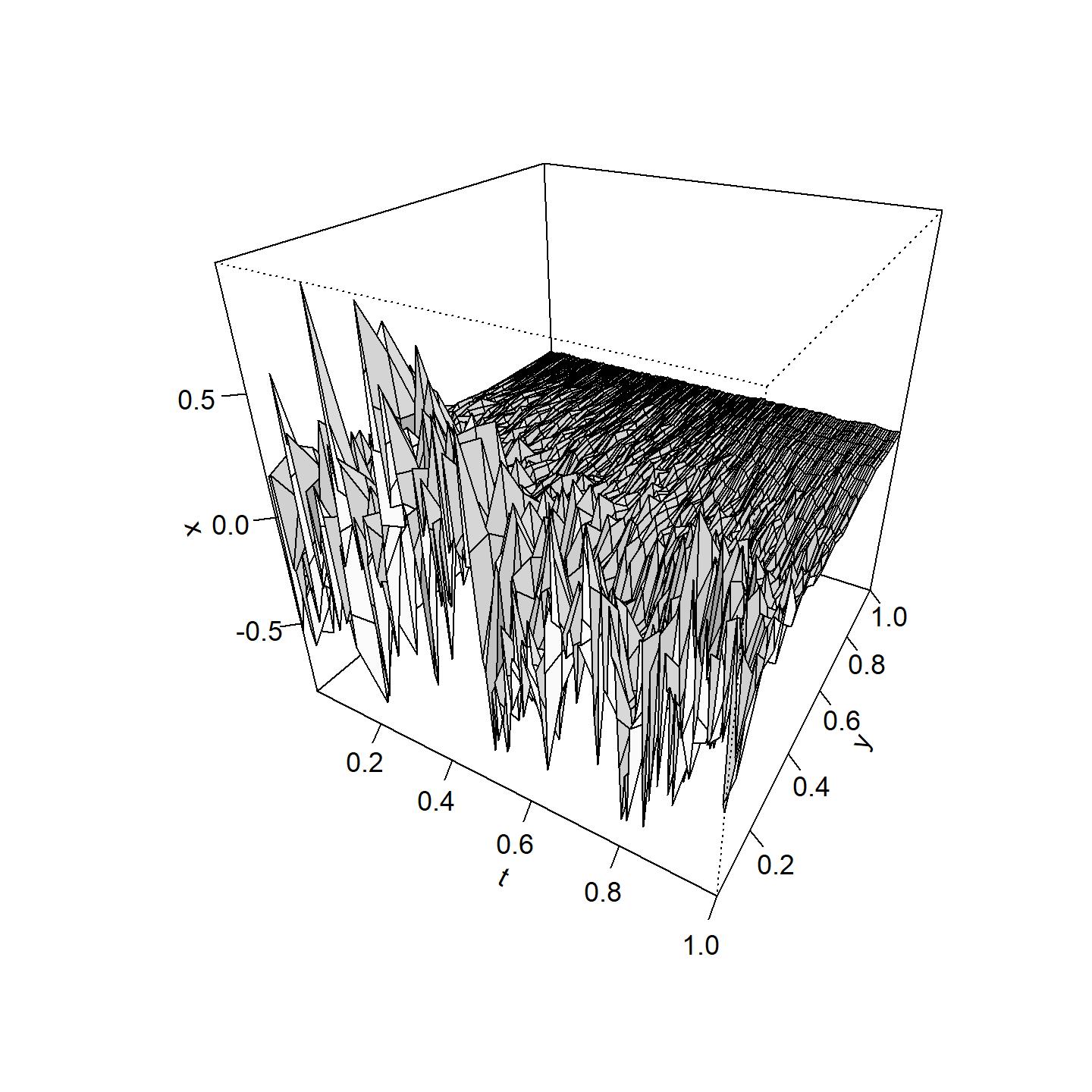}
\captionsetup{labelformat=empty,labelsep=none}
\subcaption{$\theta=$(0,1,0.1,1)}
\end{center}
\end{minipage}
\caption{Sample paths with $\theta_1=0.1, 0.5, 1$}
\label{t1-11}

\begin{minipage}{0.32\hsize}
\begin{center}
\includegraphics[width=4.3cm]{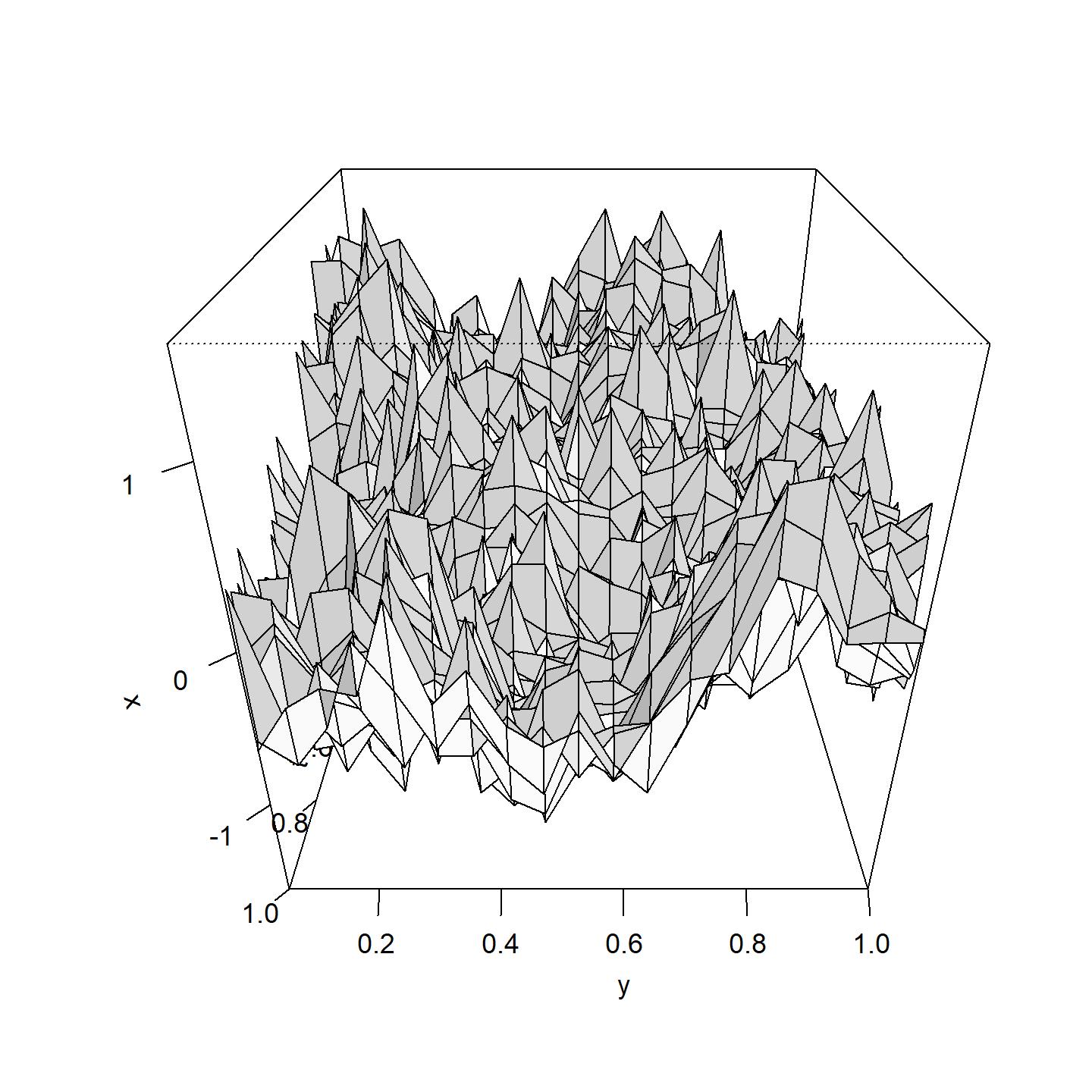}
\captionsetup{labelformat=empty,labelsep=none}
\subcaption{$\theta=$(0,0.1,0.1,1)}
\end{center}
\end{minipage}
\begin{minipage}{0.32\hsize}
\begin{center}
\includegraphics[width=4.3cm]{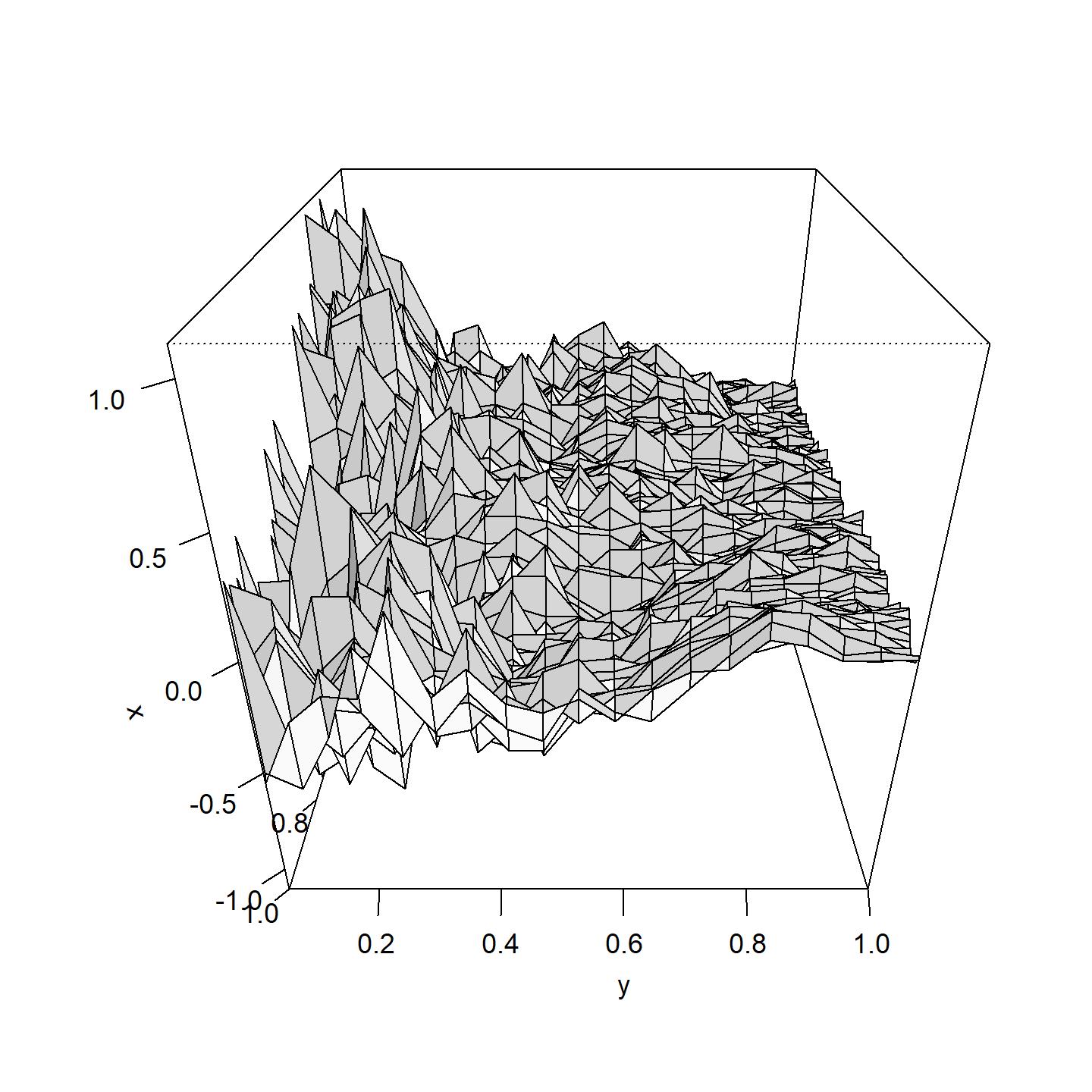}
\captionsetup{labelformat=empty,labelsep=none}
\subcaption{$\theta=$(0,0.5,0.1,1)}
\end{center}
\end{minipage}
\begin{minipage}{0.32\hsize}
\begin{center}
\includegraphics[width=4.3cm]{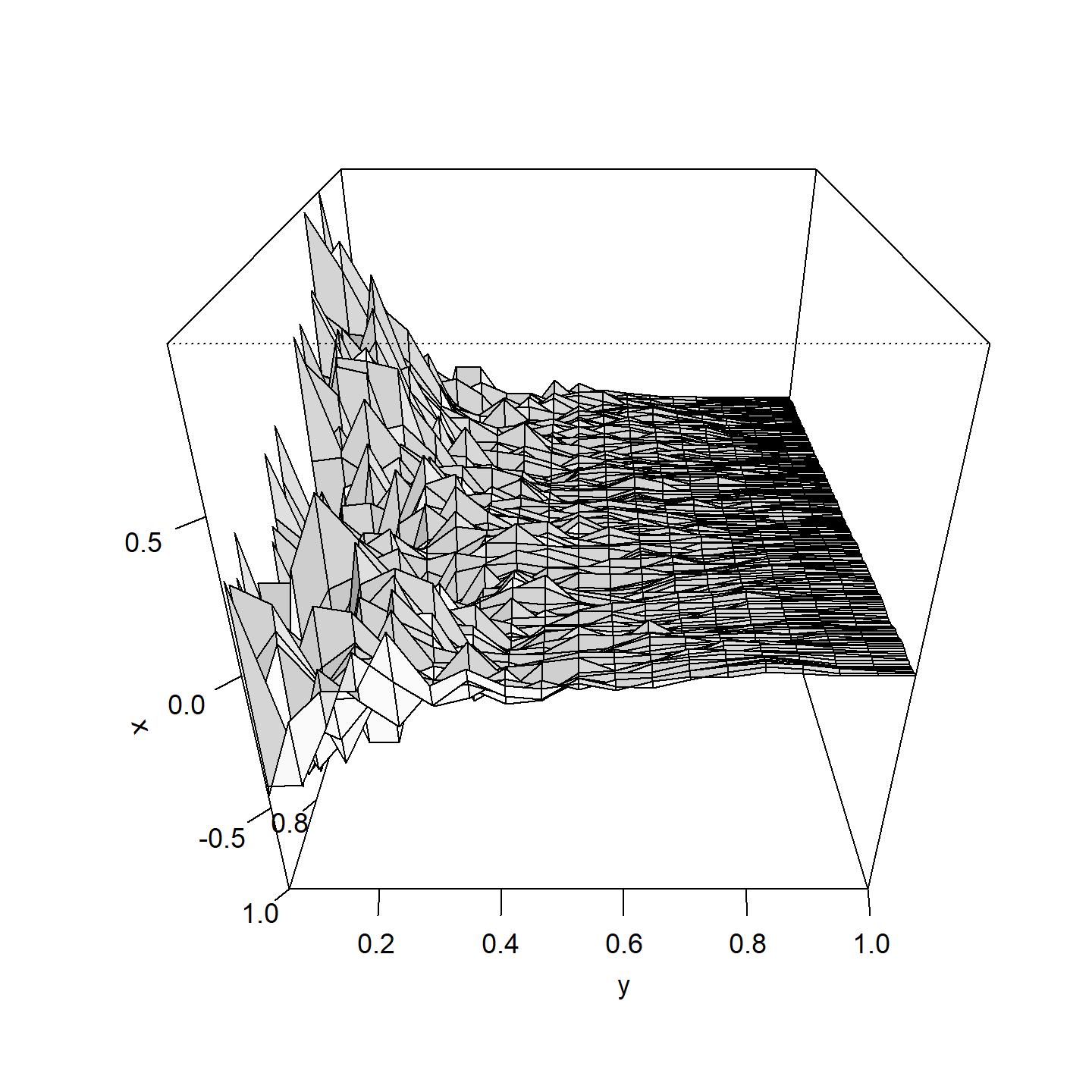}
\captionsetup{labelformat=empty,labelsep=none}
\subcaption{$\theta=$(0,1,0.1,1)}
\end{center}
\end{minipage}
\caption{Sample paths with $\theta_1=0.1, 0.5, 1$ (y-axis side)}
\label{t1-12}

\begin{minipage}{0.32\hsize}
\begin{center}
\includegraphics[width=4.3cm]{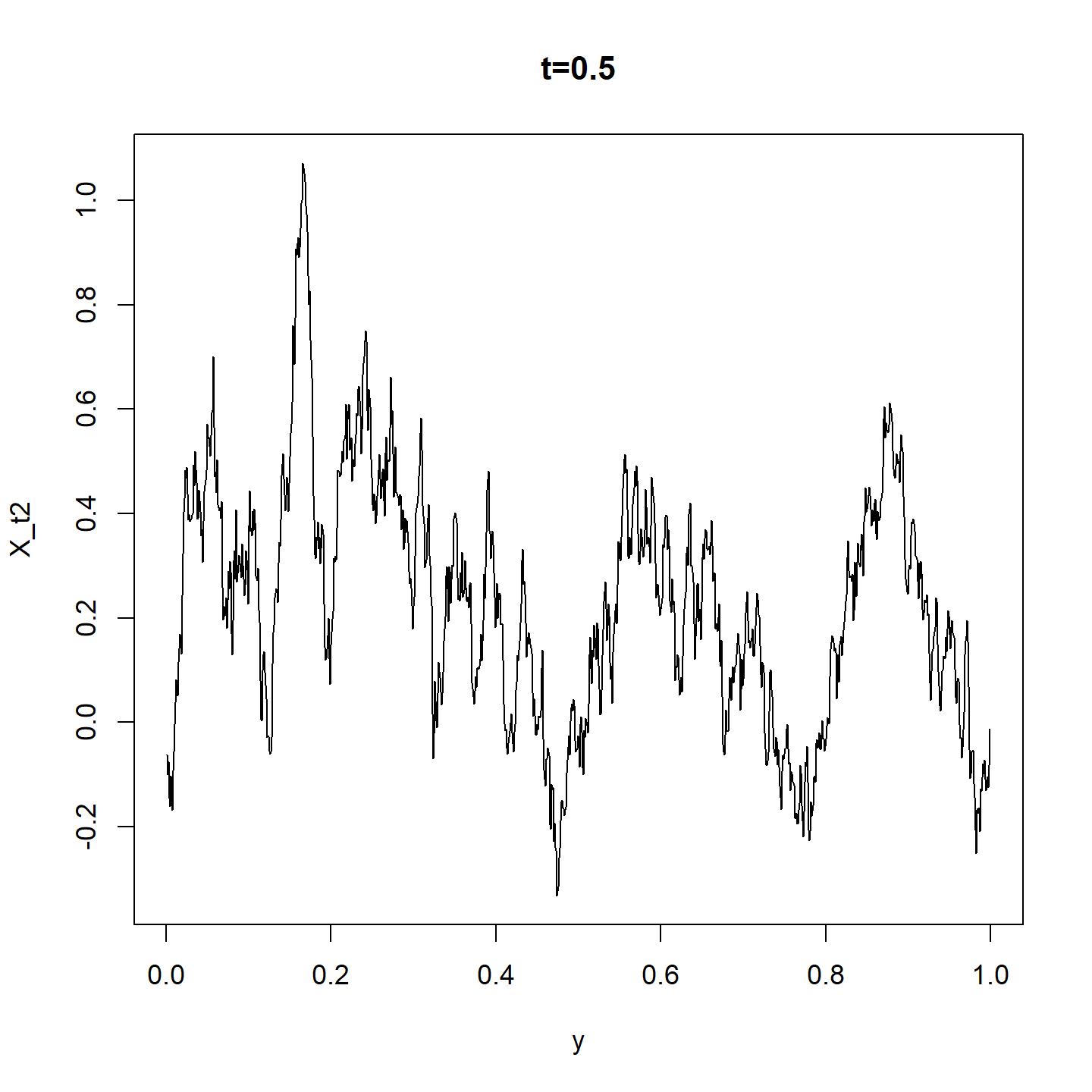}
\captionsetup{labelformat=empty,labelsep=none}
\subcaption{$\theta=$(0,0.1,0.1,1)}
\end{center}
\end{minipage}
\begin{minipage}{0.32\hsize}
\begin{center}
\includegraphics[width=4.3cm]{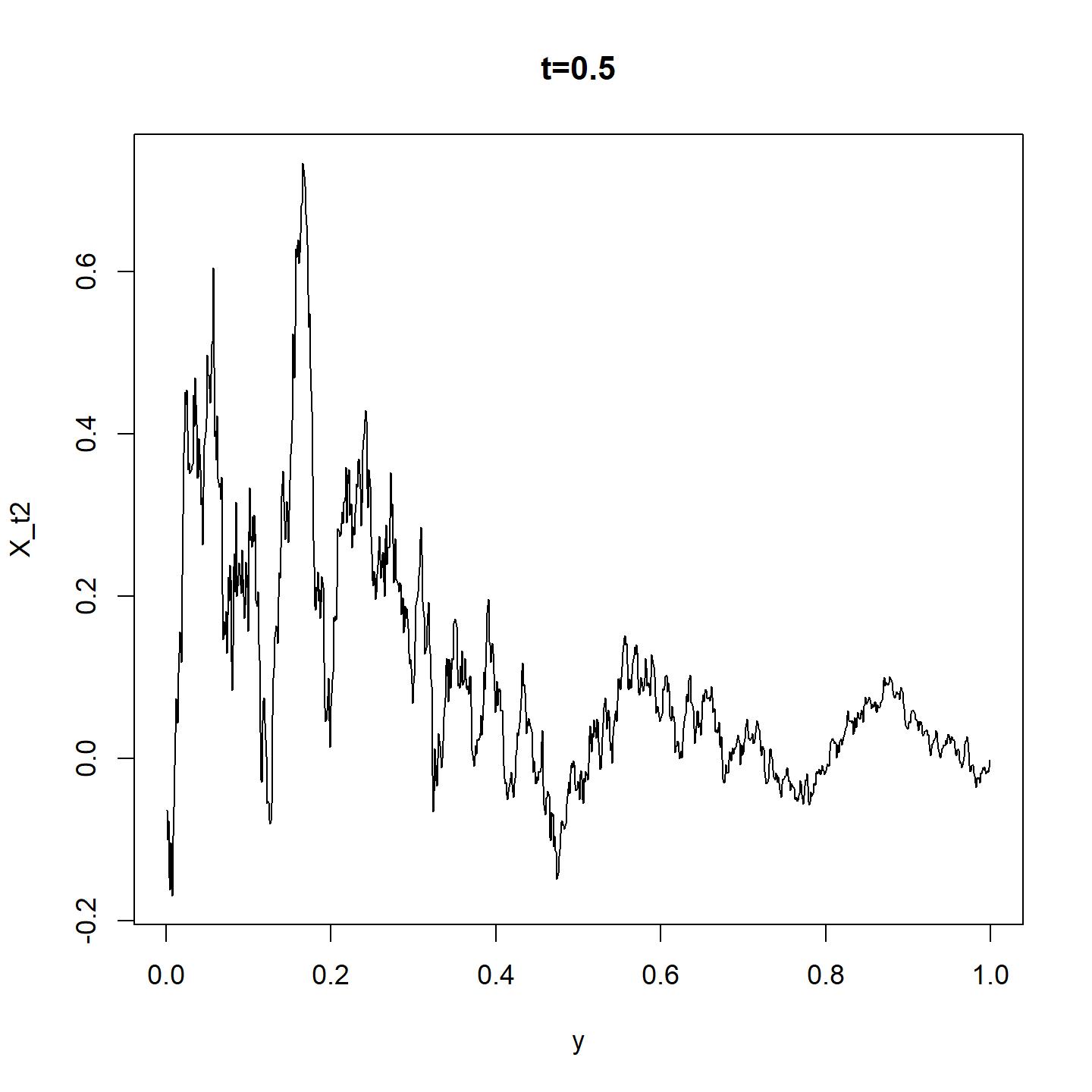}
\captionsetup{labelformat=empty,labelsep=none}
\subcaption{$\theta=$(0,0.5,0.1,1)}
\end{center}
\end{minipage}
\begin{minipage}{0.32\hsize}
\begin{center}
\includegraphics[width=4.3cm]{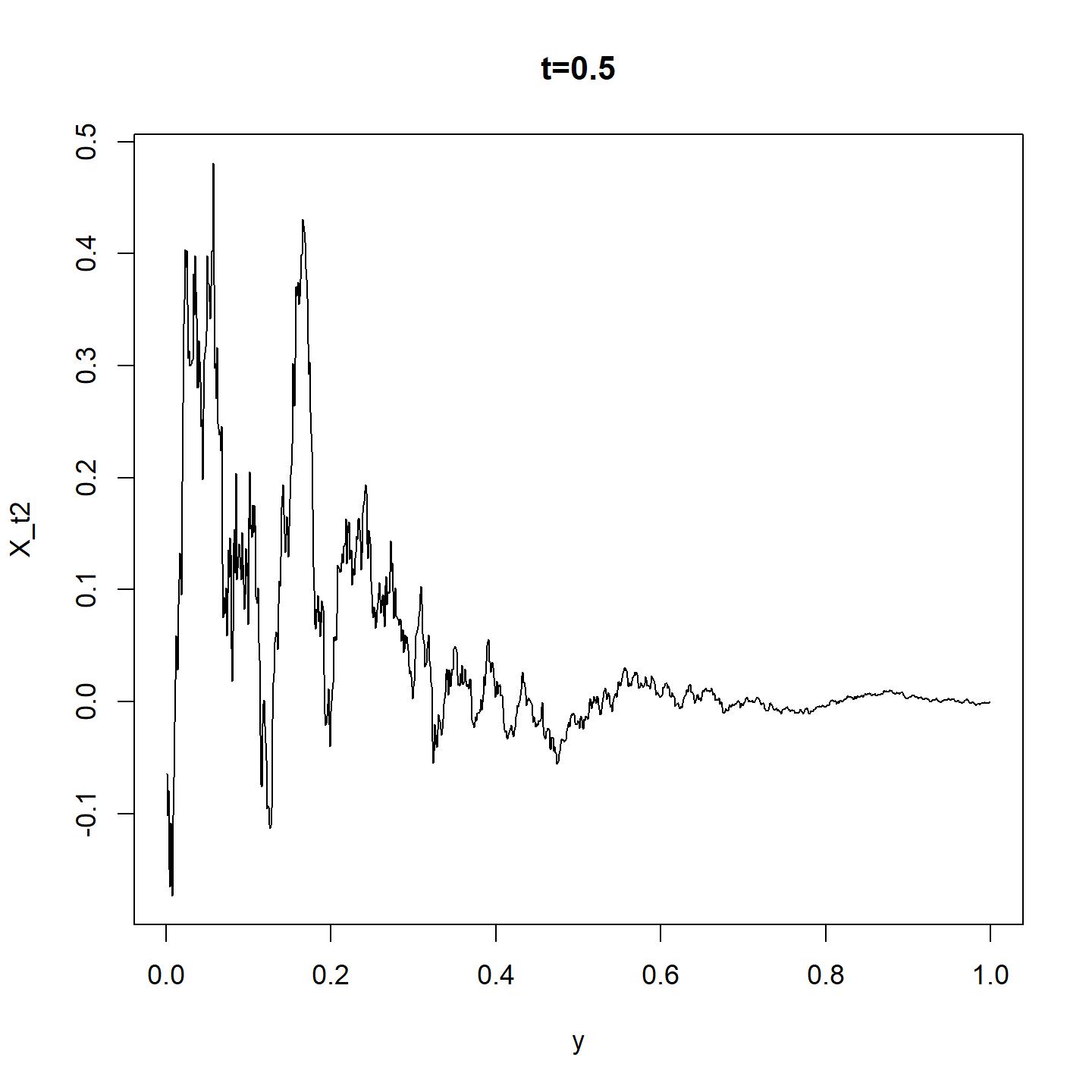}
\captionsetup{labelformat=empty,labelsep=none}
\subcaption{$\theta=$(0,1,0.1,1)}
\end{center}
\end{minipage}
\caption{Sample paths with $\theta_1=0.1, 0.5, 1$ (cross section at $t = 0.5$)}
\label{t1-13}
\end{figure}

\clearpage

Figures \ref{t1-14}-\ref{t1-16} show that {the} variation of the sample path  is small near 
 {$y = 0$ } 
and large near  {$y = 1$ } when $\theta_1 < 0$.
This trend increases as the value of $ \theta_1 $ decreases.

\begin{figure}[H]
\begin{minipage}{0.32\hsize}
\begin{center}
\includegraphics[width=4.5cm]{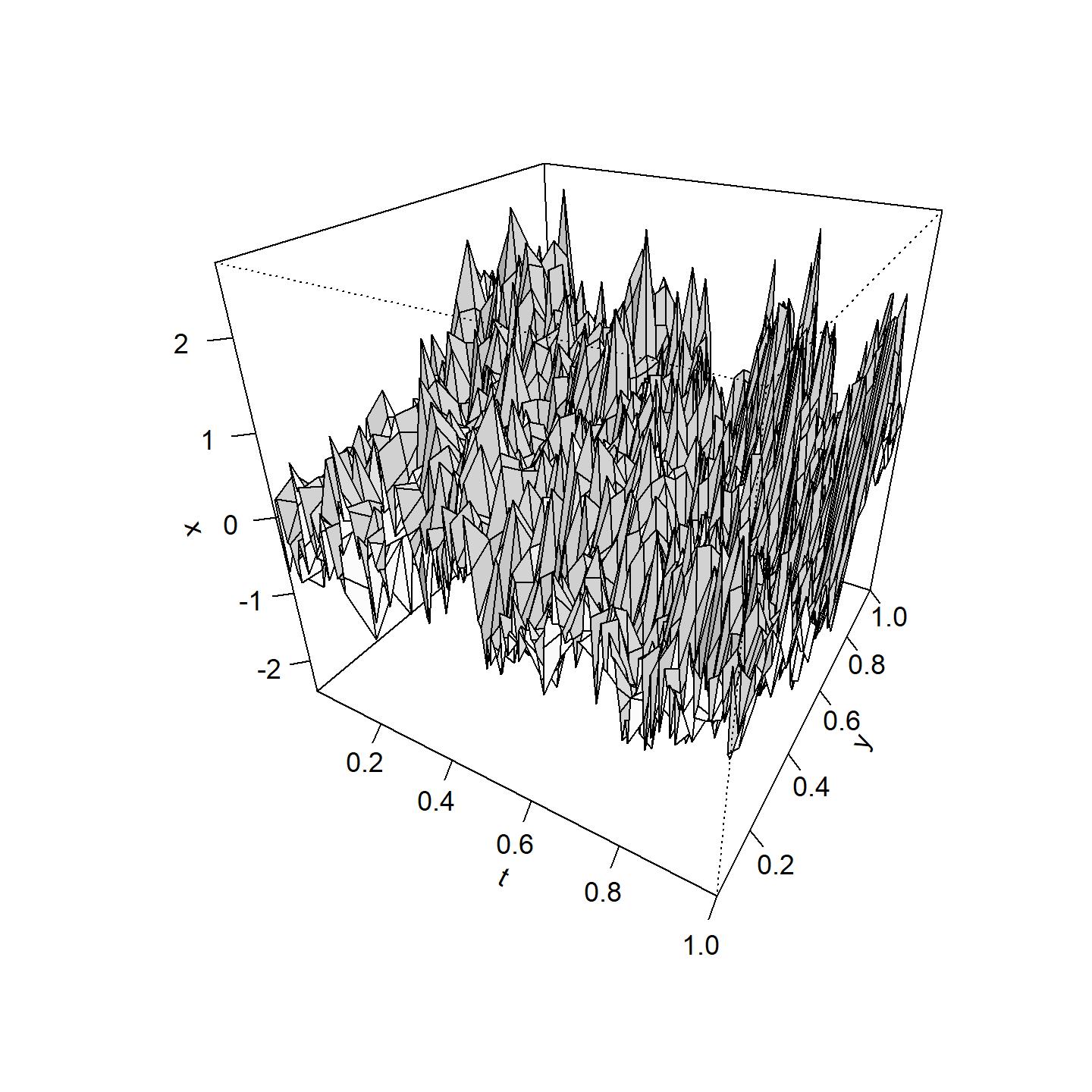}
\captionsetup{labelformat=empty,labelsep=none}
\subcaption{$\theta=$(0,-0.1,0.1,1)}
\end{center}
\end{minipage}
\begin{minipage}{0.32\hsize}
\begin{center}
\includegraphics[width=4.5cm]{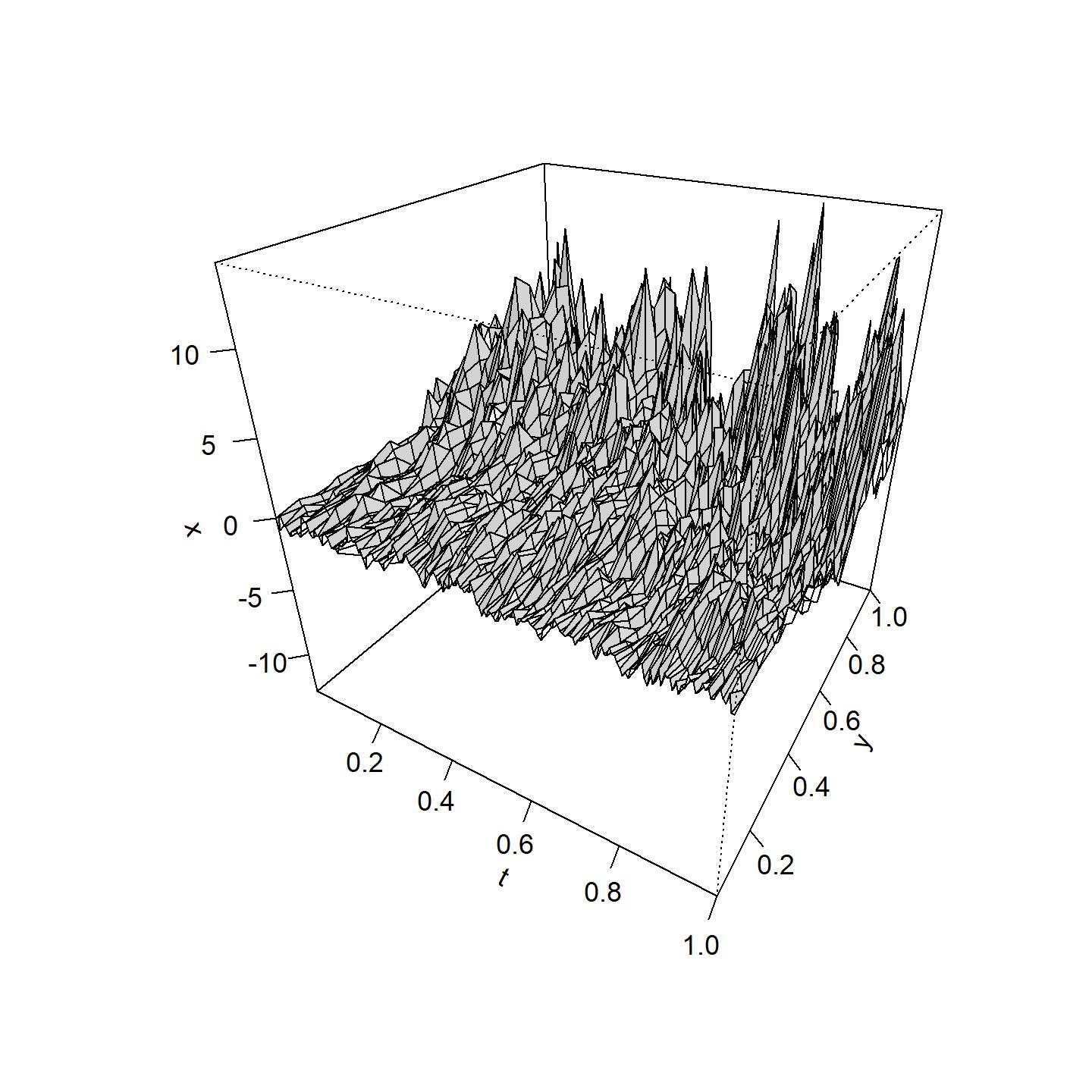}
\captionsetup{labelformat=empty,labelsep=none}
\subcaption{$\theta=$(0,-0.5,0.1,1)}
\end{center}
\end{minipage}
\begin{minipage}{0.32\hsize}
\begin{center}
\includegraphics[width=4.5cm]{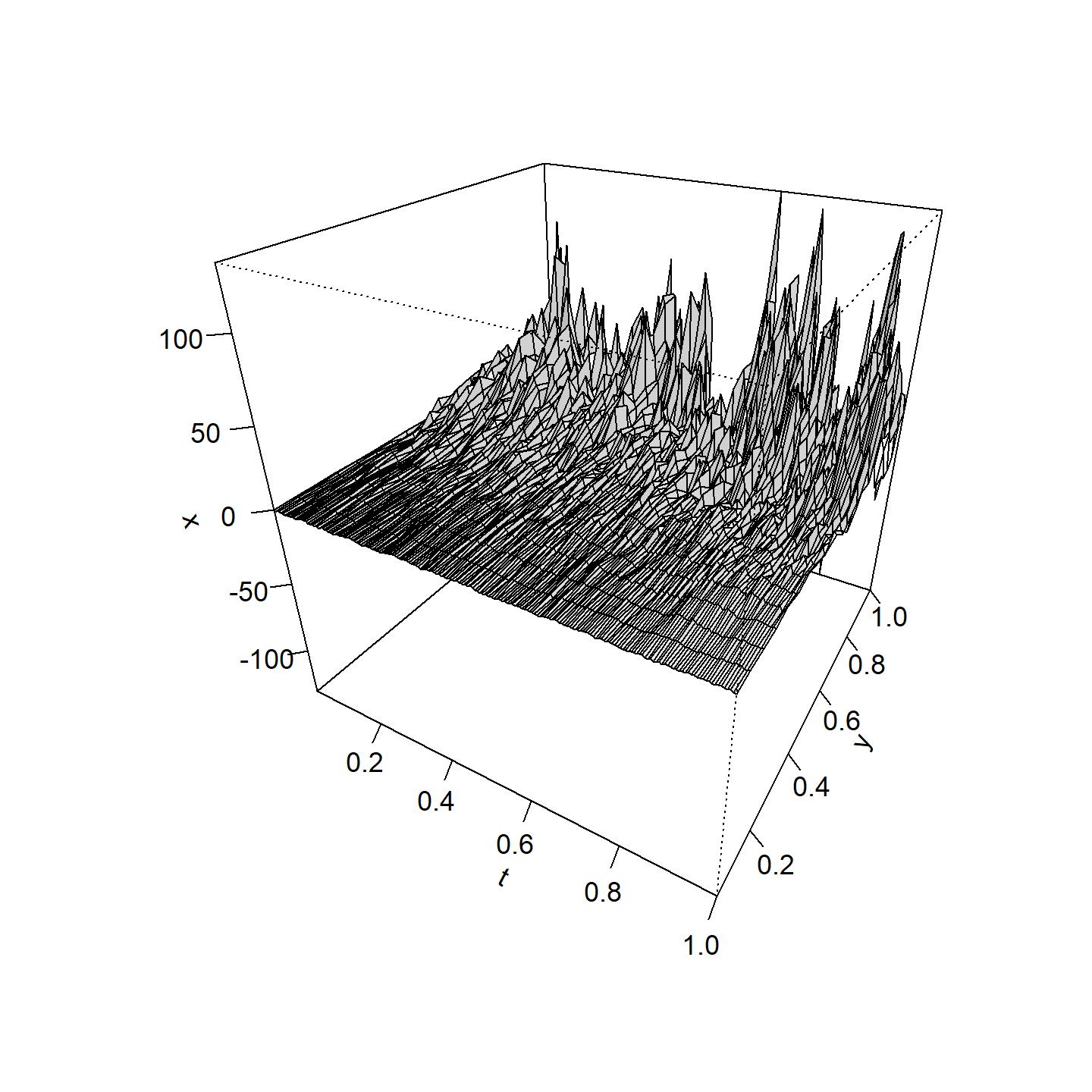}
\captionsetup{labelformat=empty,labelsep=none}
\subcaption{$\theta=$(0,-1,0.1,1)}
\end{center}
\end{minipage}
\caption{Sample paths with $\theta_1=-0.1, -0.5, -1$}
\label{t1-14}

\begin{minipage}{0.32\hsize}
\begin{center}
\includegraphics[width=4.3cm]{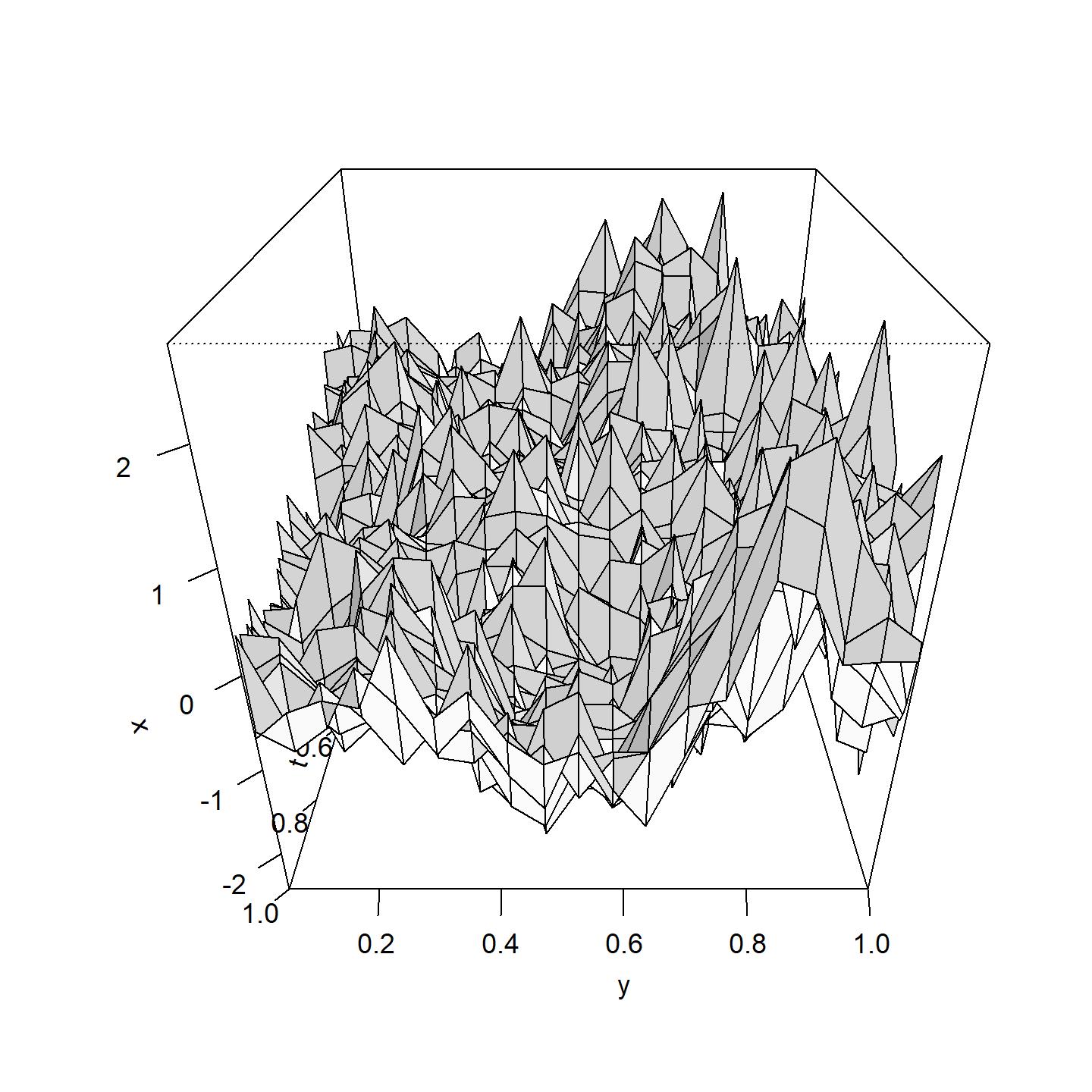}
\captionsetup{labelformat=empty,labelsep=none}
\subcaption{$\theta=$(0,-0.1,0.1,1)}
\end{center}
\end{minipage}
\begin{minipage}{0.32\hsize}
\begin{center}
\includegraphics[width=4.3cm]{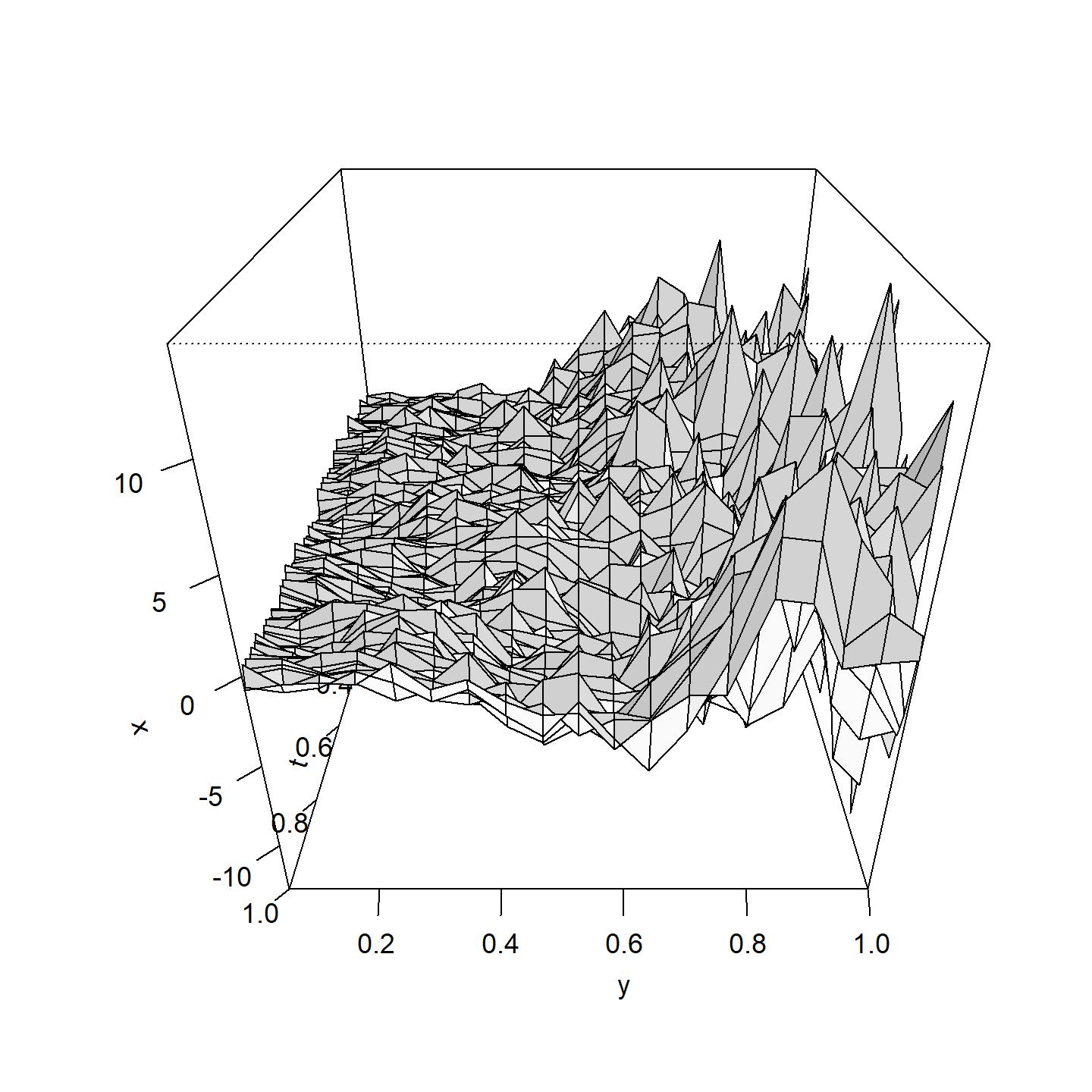}
\captionsetup{labelformat=empty,labelsep=none}
\subcaption{$\theta=$(0,-0.5,0.1,1)}
\end{center}
\end{minipage}
\begin{minipage}{0.32\hsize}
\begin{center}
\includegraphics[width=4.3cm]{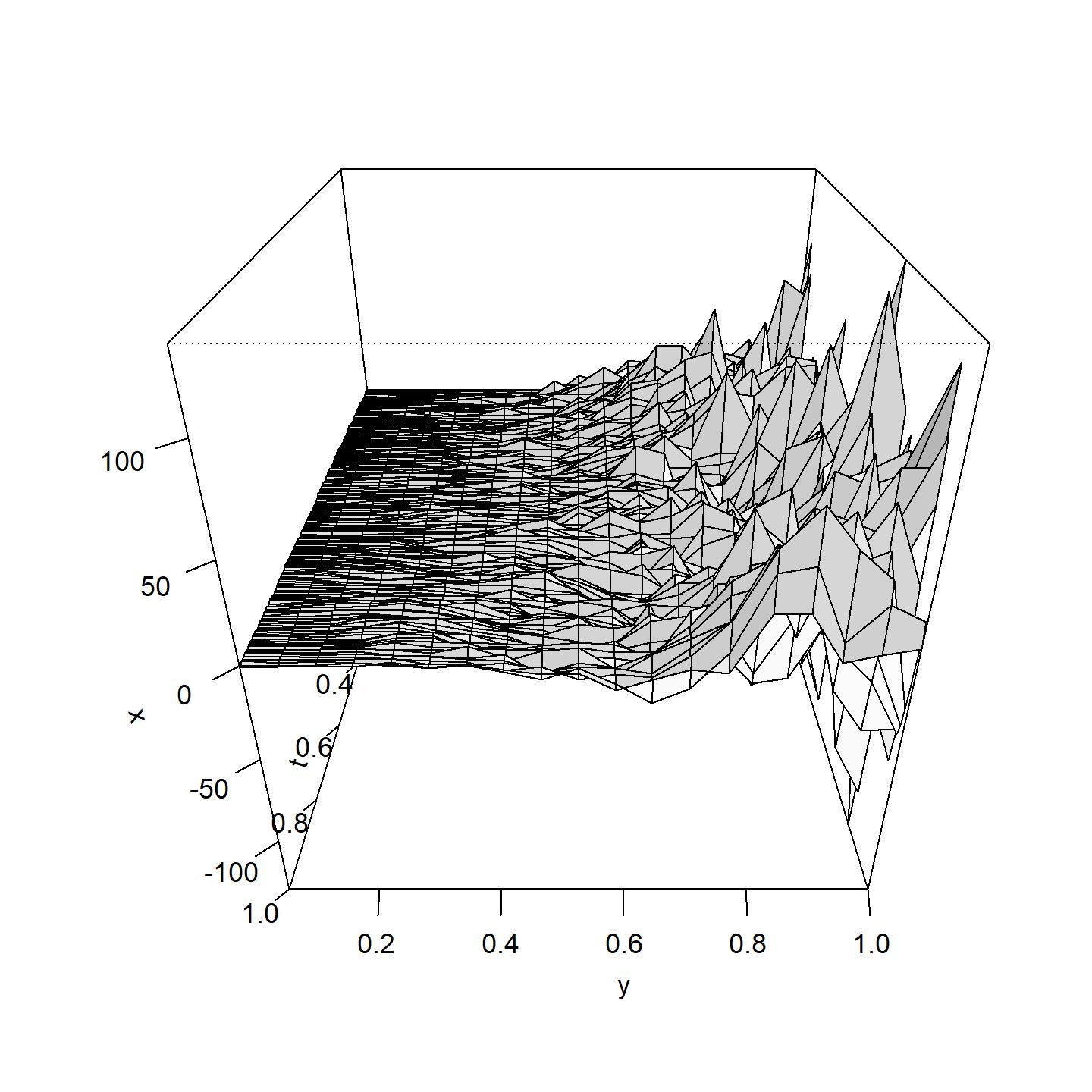}
\captionsetup{labelformat=empty,labelsep=none}
\subcaption{$\theta=$(0,-1,0.1,1)}
\end{center}
\end{minipage}
\caption{Sample paths with $\theta_1=-0.1, -0.5, -1$ (y-axis side)}
\label{t1-15}

\begin{minipage}{0.32\hsize}
\begin{center}
\includegraphics[width=4.3cm]{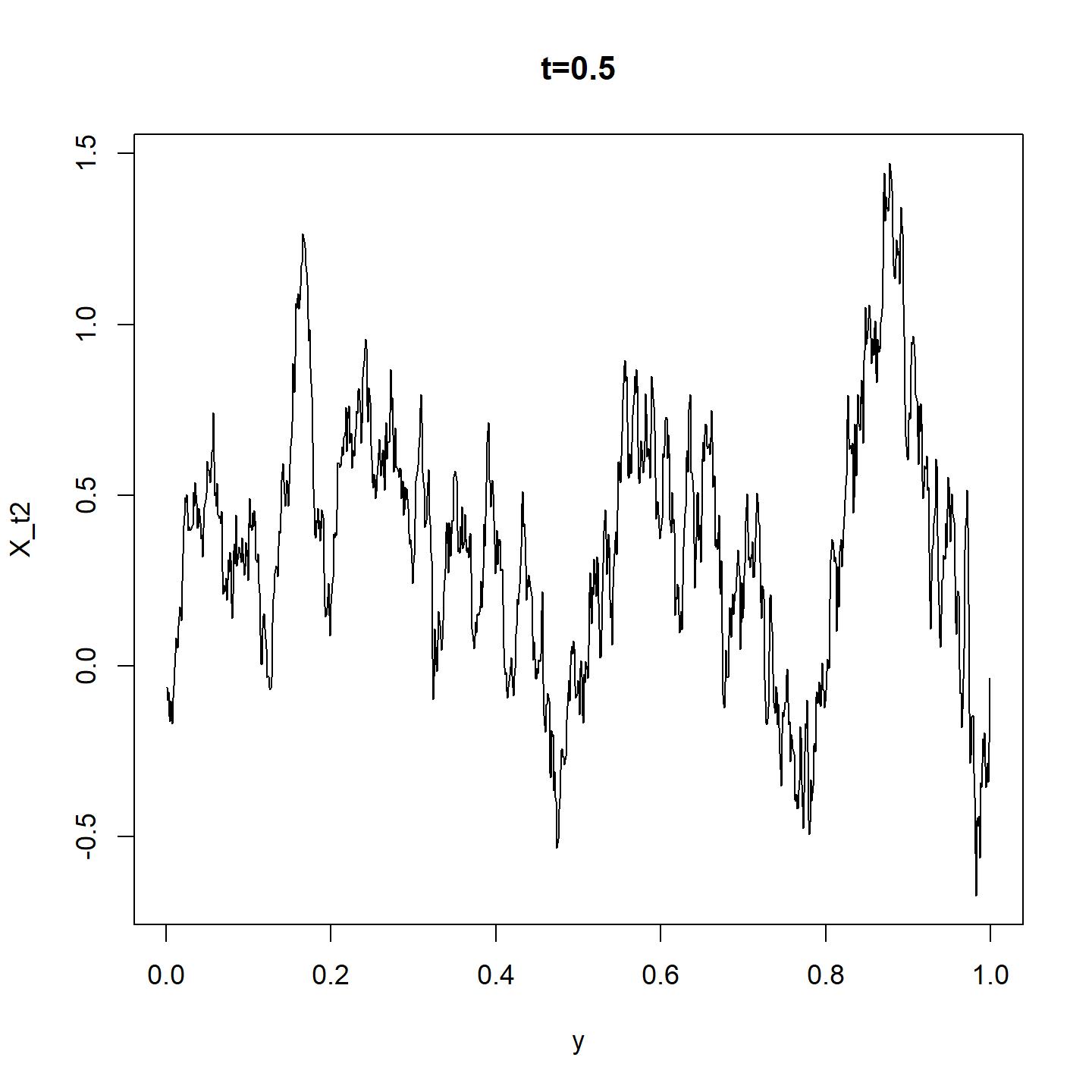}
\captionsetup{labelformat=empty,labelsep=none}
\subcaption{$\theta=$(0,-0.1,0.1,1)}
\end{center}
\end{minipage}
\begin{minipage}{0.32\hsize}
\begin{center}
\includegraphics[width=4.3cm]{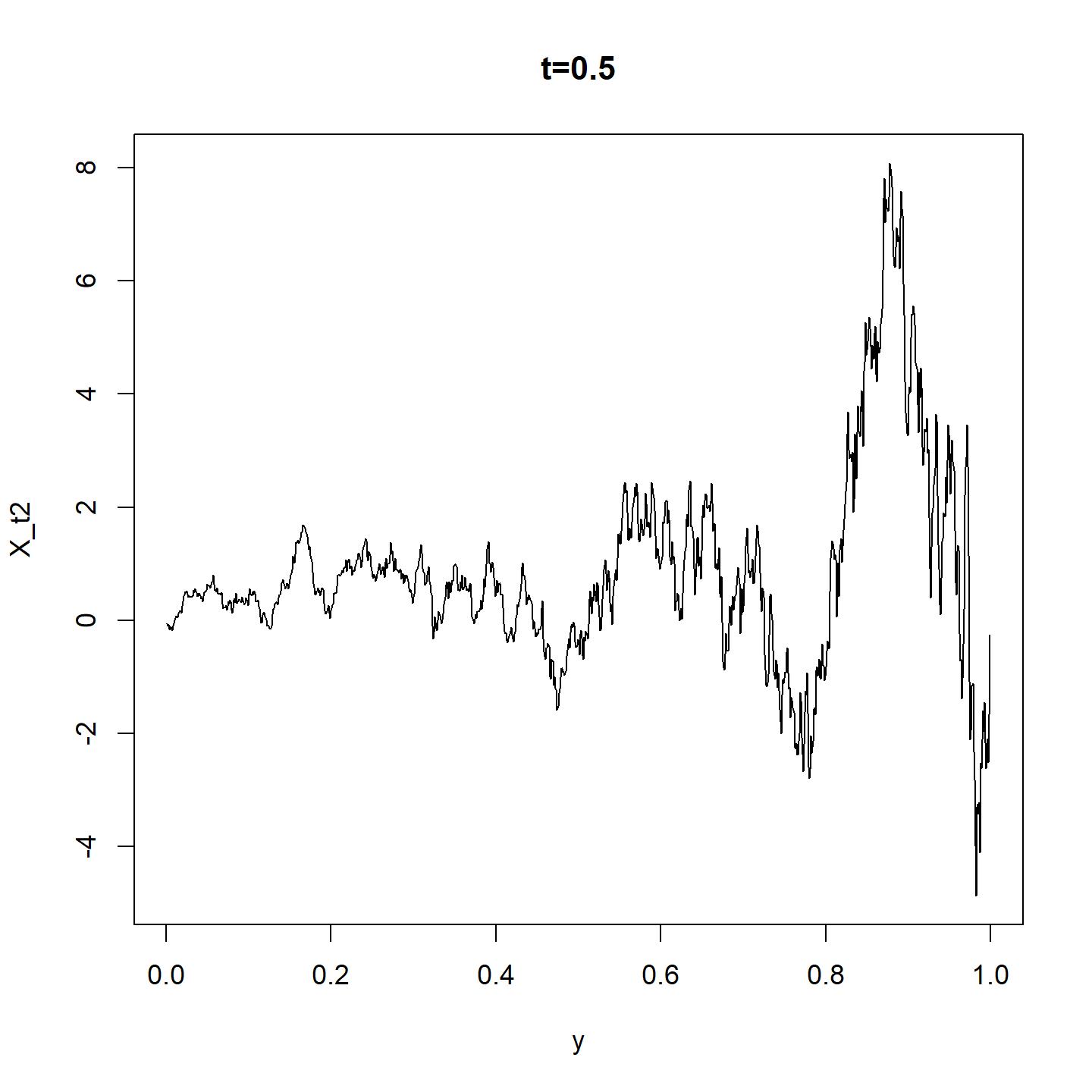}
\captionsetup{labelformat=empty,labelsep=none}
\subcaption{$\theta=$(0,-0.5,0.1,1)}
\end{center}
\end{minipage}
\begin{minipage}{0.32\hsize}
\begin{center}
\includegraphics[width=4.3cm]{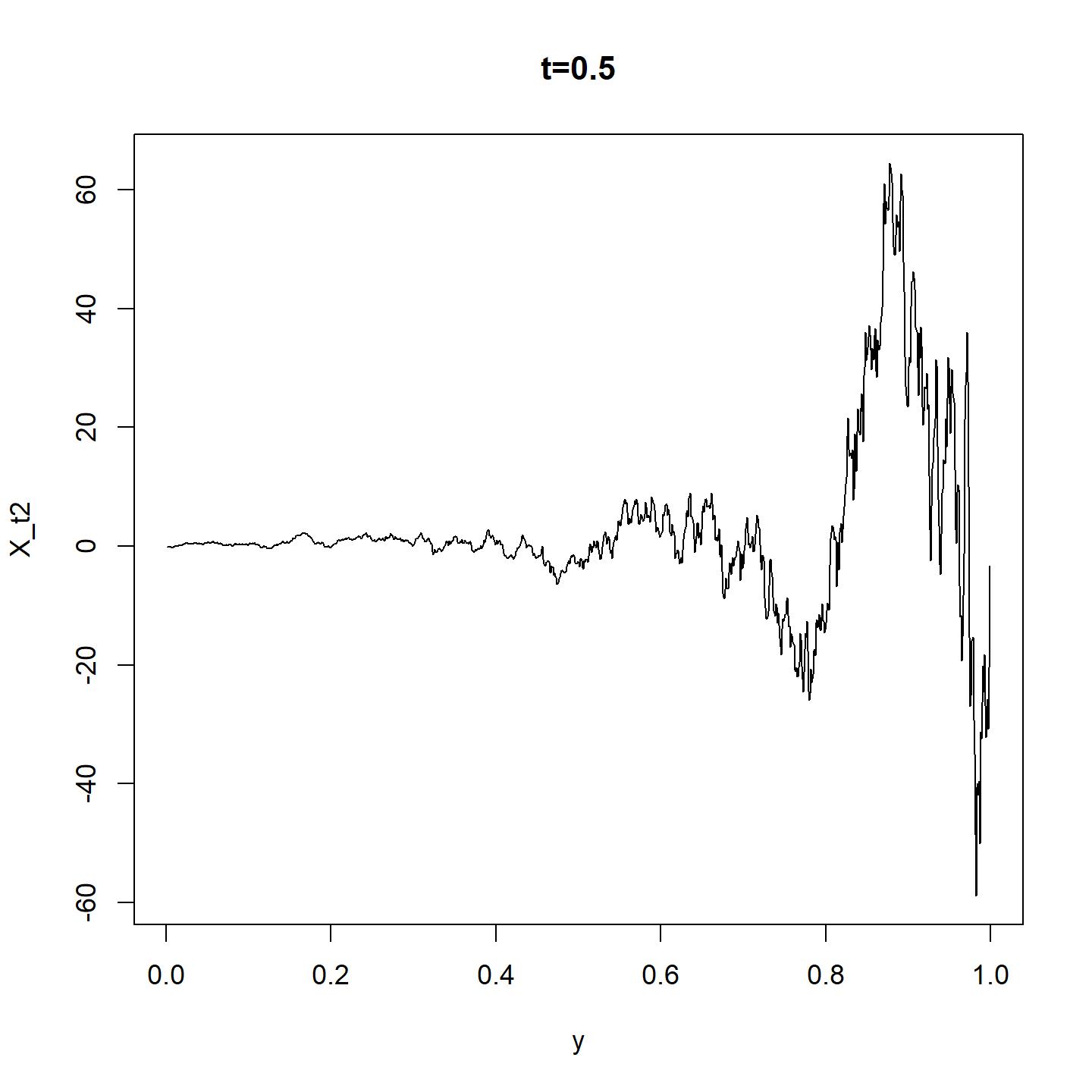}
\captionsetup{labelformat=empty,labelsep=none}
\subcaption{$\theta=$(0,-1,0.1,1)}
\end{center}
\end{minipage}
\caption{Sample paths with $\theta_1=-0.1, -0.5, -1$ (cross section at  {$t = 0.5$})}
\label{t1-16}
\end{figure}

\clearpage
Figures \ref{t2-1}-\ref{t2-3} are {the} sample paths,
{ where} $\theta_0$, $\theta_1$ and $\sigma$ are fixed and only $\theta_2$ is changed.
These show that {the} variation of the sample path  is large near 
 {$y = 0$ } and small near 
 {$y = 1$ } and
this trend increases as the value of $ \theta_2 $ decreases.


\begin{figure}[H]
\begin{minipage}{0.32\hsize}
\begin{center}
\includegraphics[width=4.5cm]{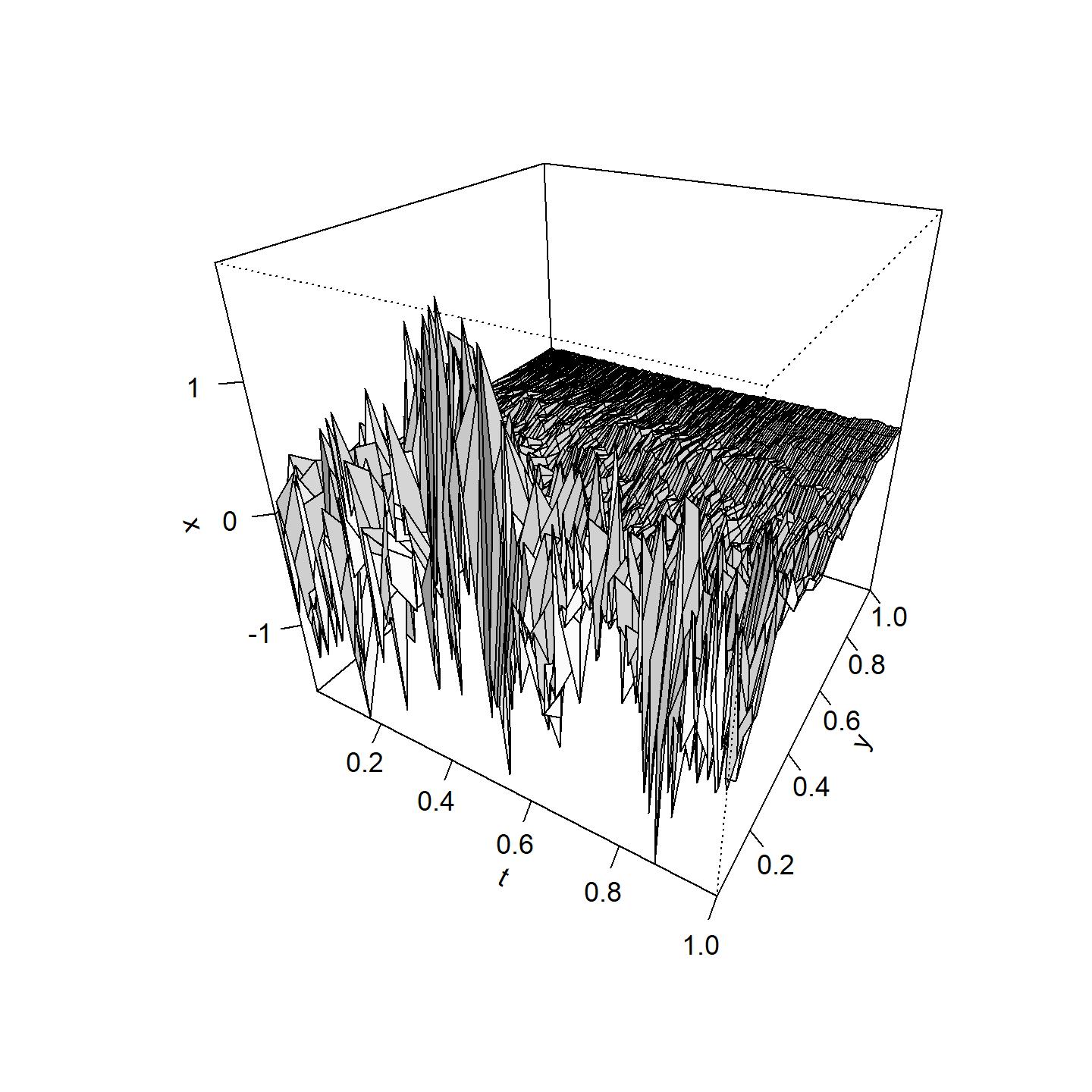}
\captionsetup{labelformat=empty,labelsep=none}
\subcaption{$\theta=$(0,0.1,0.01,1)}
\end{center}
\end{minipage}
\begin{minipage}{0.32\hsize}
\begin{center}
\includegraphics[width=4.5cm]{p2all.jpeg}
\captionsetup{labelformat=empty,labelsep=none}
\subcaption{$\theta=$(0,0.1,0.1,1)}
\end{center}
\end{minipage}
\begin{minipage}{0.32\hsize}
\begin{center}
\includegraphics[width=4.5cm]{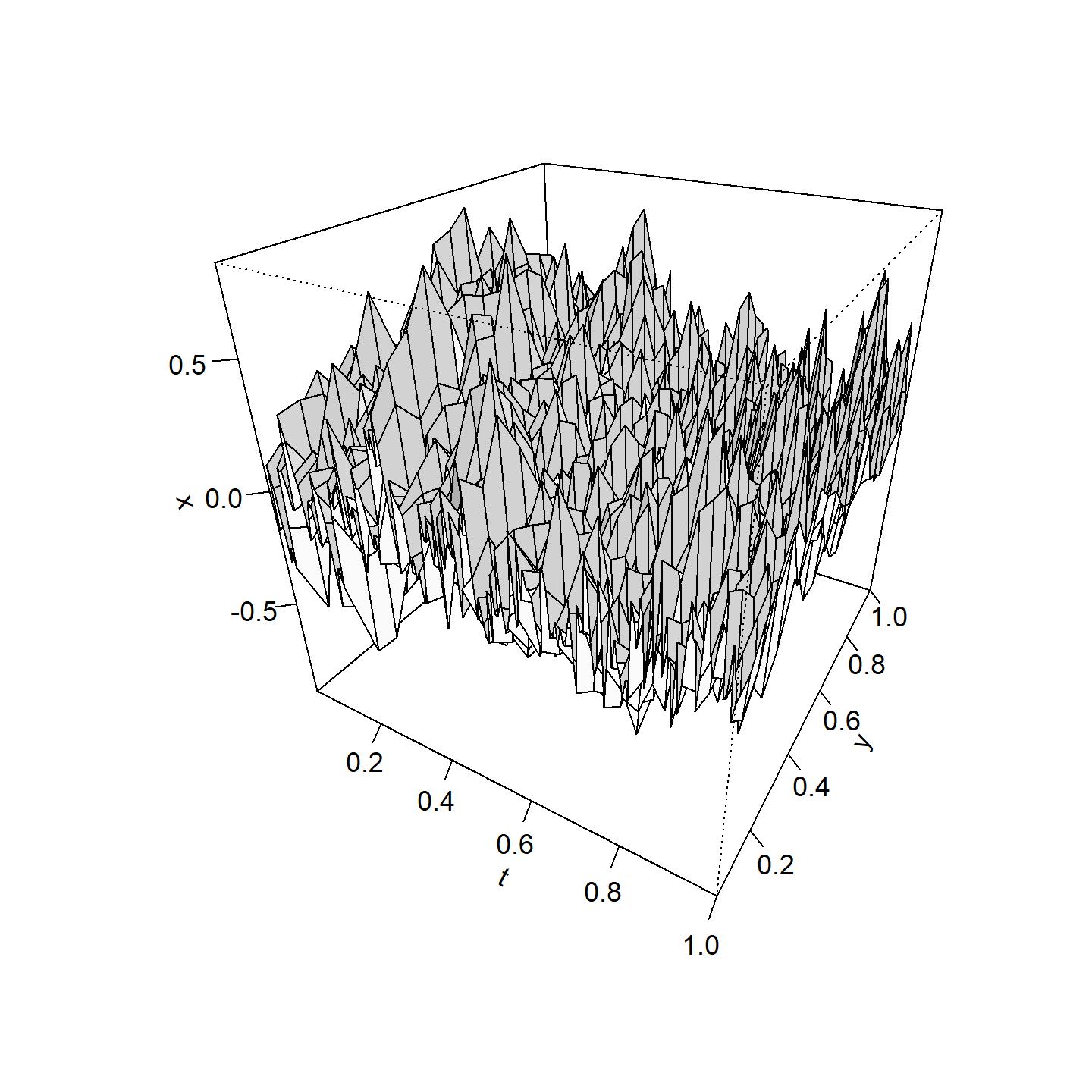}
\captionsetup{labelformat=empty,labelsep=none}
\subcaption{$\theta=$(0,0.1,1,1)}
\end{center}
\end{minipage}
\caption{Sample paths with $\theta_2=0.01, 0.1, 1$}\label{t2-1}

\begin{minipage}{0.32\hsize}
\begin{center}
\includegraphics[width=4.3cm]{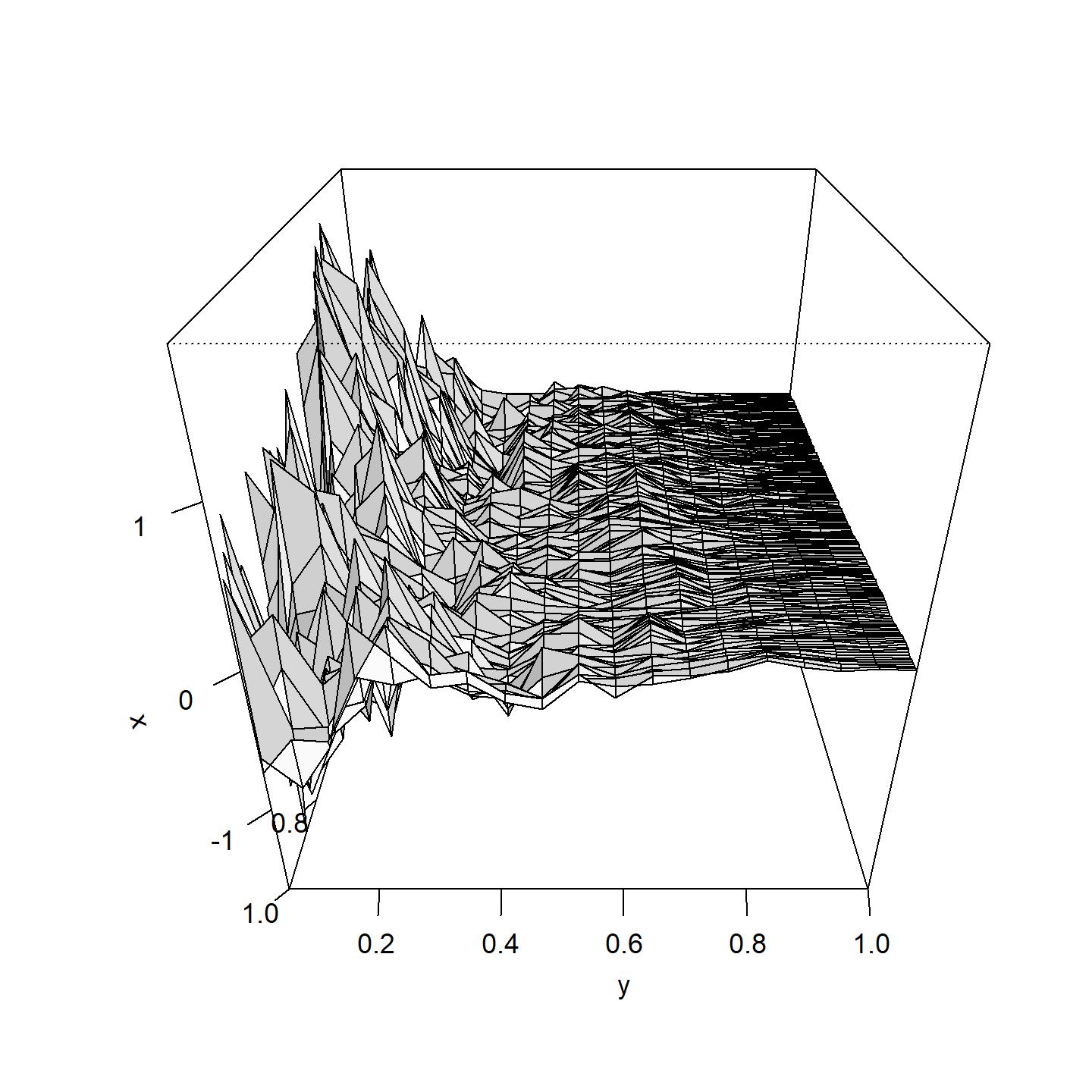}
\captionsetup{labelformat=empty,labelsep=none}
\subcaption{$\theta=$(0,0.1,0.01,1)}
\end{center}
\end{minipage}
\begin{minipage}{0.32\hsize}
\begin{center}
\includegraphics[width=4.3cm]{p2yside.jpeg}
\captionsetup{labelformat=empty,labelsep=none}
\subcaption{$\theta=$(0,0.1,0.1,1)}
\end{center}
\end{minipage}
\begin{minipage}{0.32\hsize}
\begin{center}
\includegraphics[width=4.3cm]{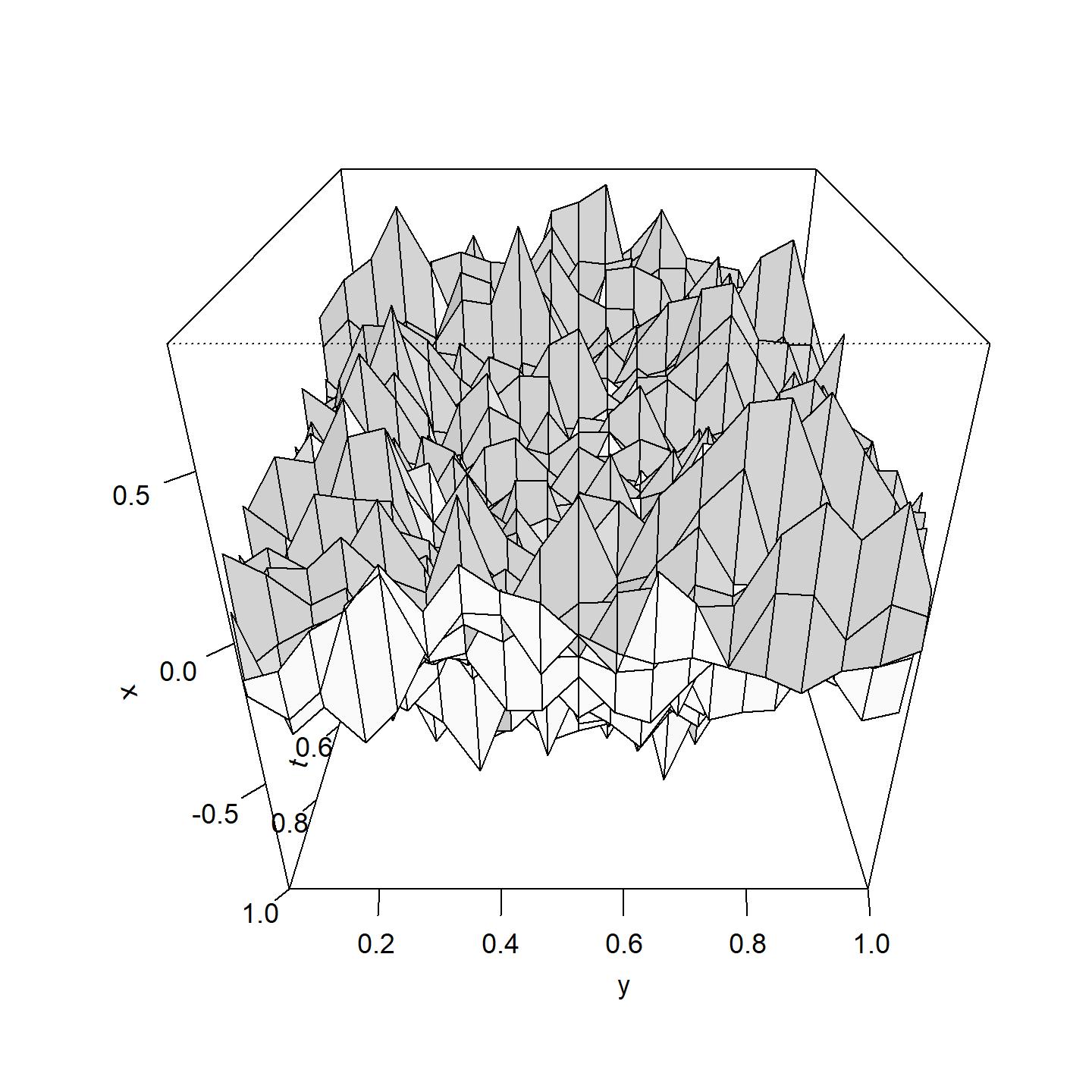}
\captionsetup{labelformat=empty,labelsep=none}
\subcaption{$\theta=$(0,0.1,1,1)}
\end{center}
\end{minipage}
\caption{Sample paths with $\theta_2=0.01, 0.1, 1$ (y-axis side)}\label{t2-2}

\begin{minipage}{0.32\hsize}
\begin{center}
\includegraphics[width=4.3cm]{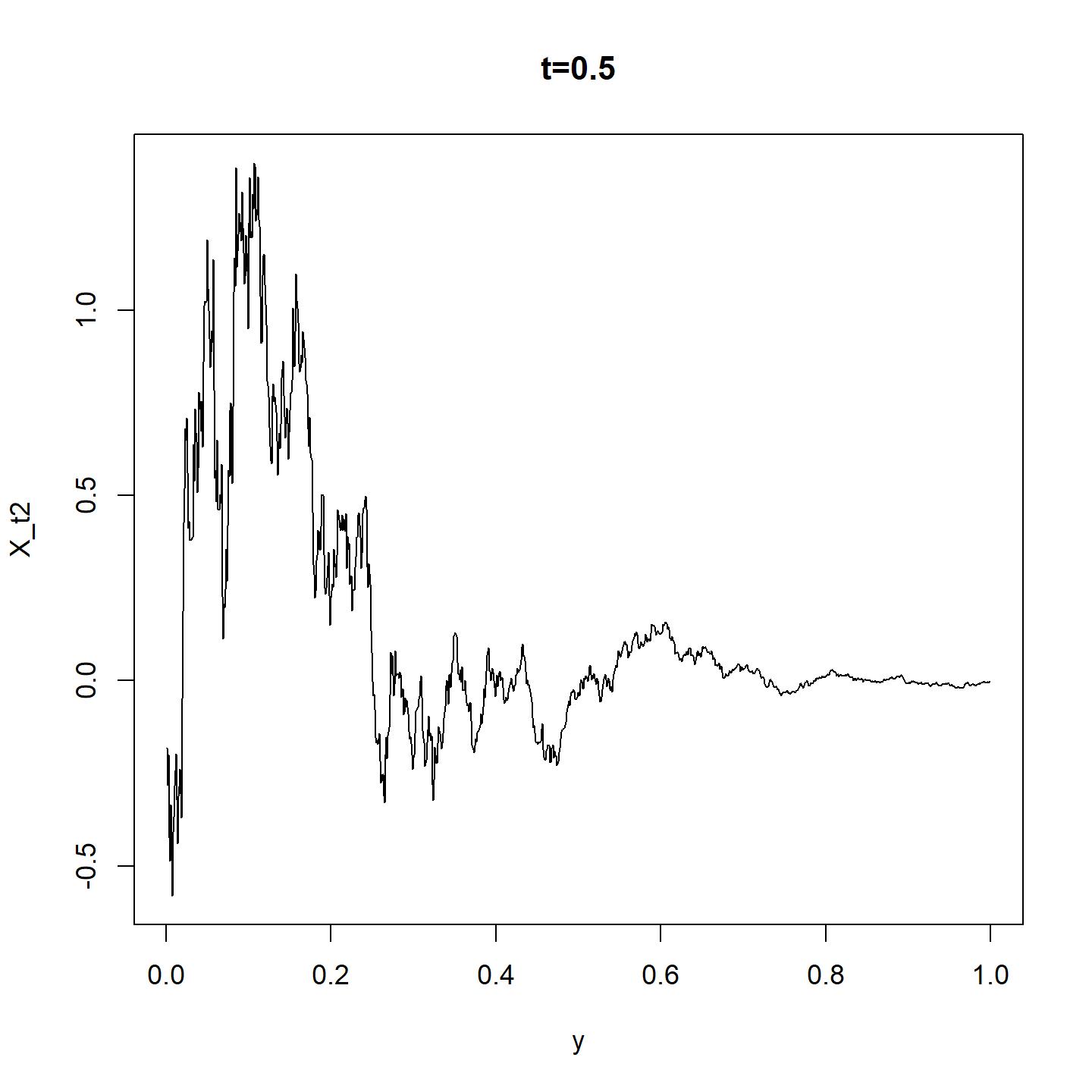}
\captionsetup{labelformat=empty,labelsep=none}
\subcaption{$\theta=$(0,0.1,0.01,1)}
\end{center}
\end{minipage}
\begin{minipage}{0.32\hsize}
\begin{center}
\includegraphics[width=4.3cm]{p2t=05.jpeg}
\captionsetup{labelformat=empty,labelsep=none}
\subcaption{$\theta=$(0,0.1,0.1,1)}
\end{center}
\end{minipage}
\begin{minipage}{0.32\hsize}
\begin{center}
\includegraphics[width=4.3cm]{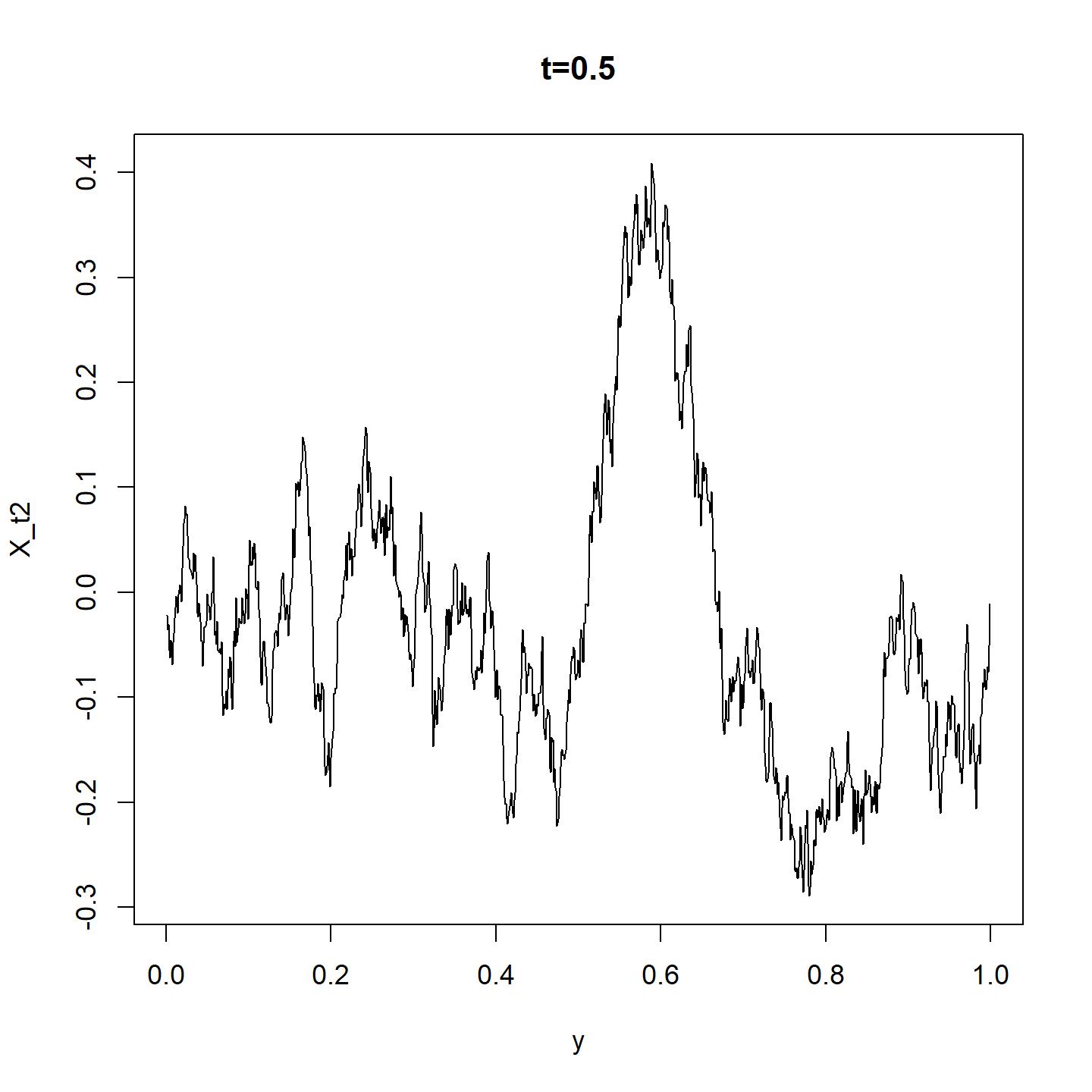}
\captionsetup{labelformat=empty,labelsep=none}
\subcaption{$\theta=$(0,0.1,1,1)}
\end{center}
\end{minipage}
\caption{Sample paths with $\theta_2=0.01, 0.1, 1$ (cross section at $t=0.5$)}\label{t2-3}
\end{figure}

\clearpage
Figures \ref{t0-1}-\ref{t0-6} are {the} sample paths,
{ where}  $\theta_1$, $\theta_2$ and $\sigma$ are fixed and only $\theta_0$ is changed.
Figures \ref{t0-1}-\ref{t0-3} show  {the} sample paths with {$T=1$} and 
Figures \ref{t0-4}-\ref{t0-6} show  {the} sample paths with {$T=100$}. 
$ \theta_0 $ affects the shape of  {the} sample path when $y$ is fixed and $t$ is changed.
When  {$T$} is large, the effect of $\theta_0$ is large.

\begin{figure}[H]
\begin{minipage}{0.32\hsize}
\begin{center}
\includegraphics[width=4.5cm]{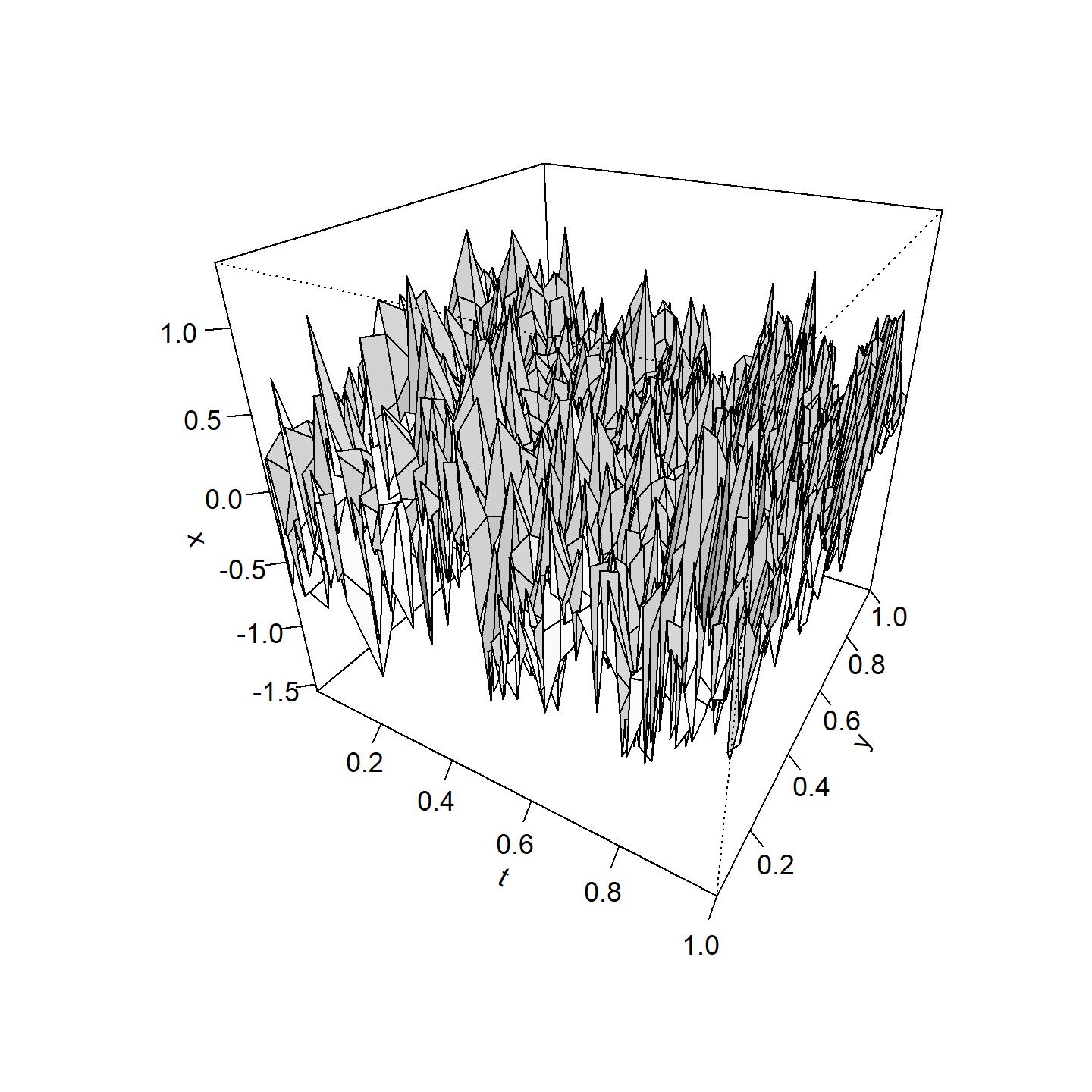}
\captionsetup{labelformat=empty,labelsep=none}
\subcaption{$\theta=$(-5,0.1,0.1,1)}
\end{center}
\end{minipage}
\begin{minipage}{0.32\hsize}
\begin{center}
\includegraphics[width=4.5cm]{p2all.jpeg}
\captionsetup{labelformat=empty,labelsep=none}
\subcaption{$\theta=$(0,0.1,0.1,1)}
\end{center}
\end{minipage}
\begin{minipage}{0.32\hsize}
\begin{center}
\includegraphics[width=4.5cm]{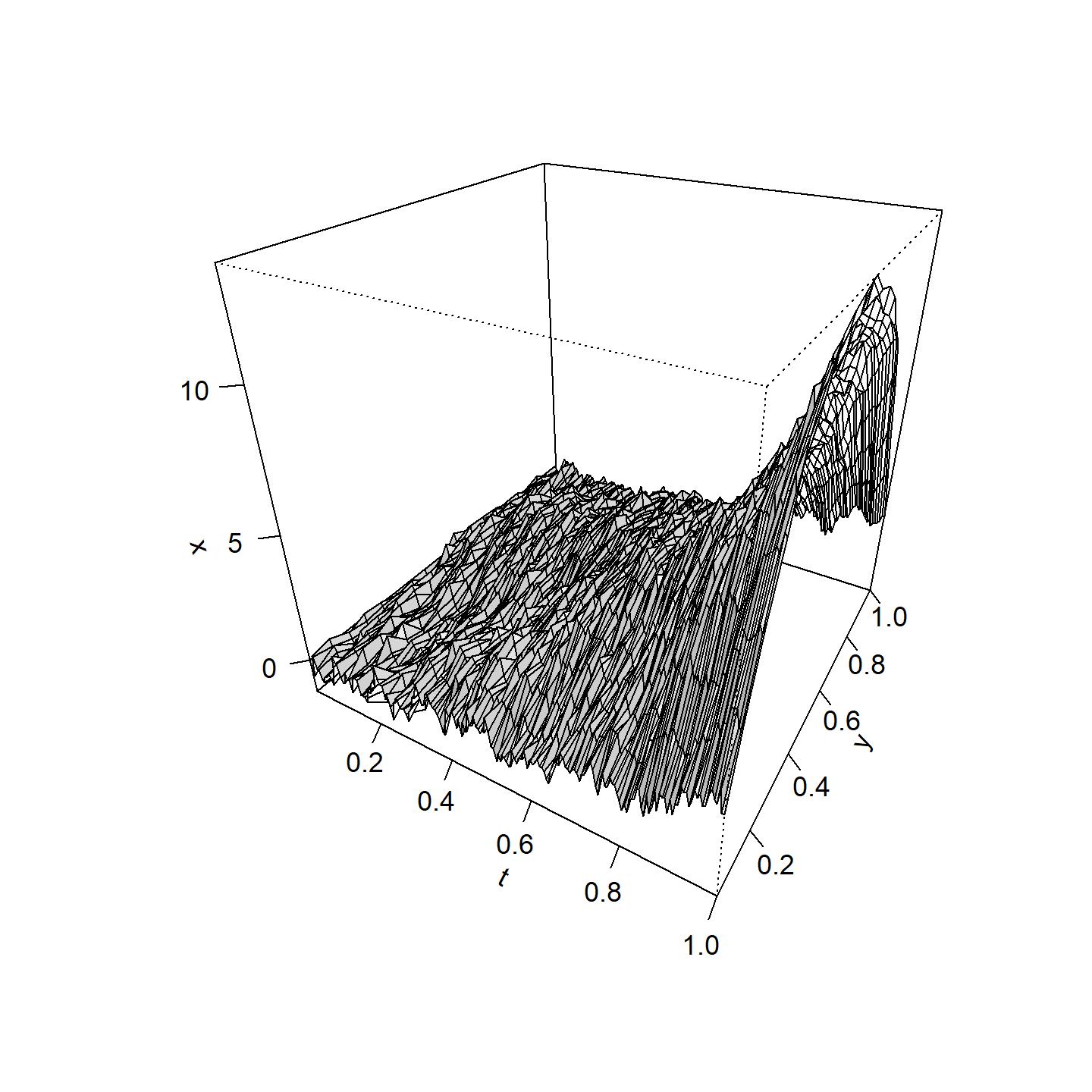}
\captionsetup{labelformat=empty,labelsep=none}
\subcaption{$\theta=$(5,0.1,0.1,1)}
\end{center}
\end{minipage}
\caption{Sample paths with $\theta_0=-5, 0, 5$ and $T=1$}\label{t0-1}

\begin{minipage}{0.32\hsize}
\begin{center}
\includegraphics[width=4.3cm]{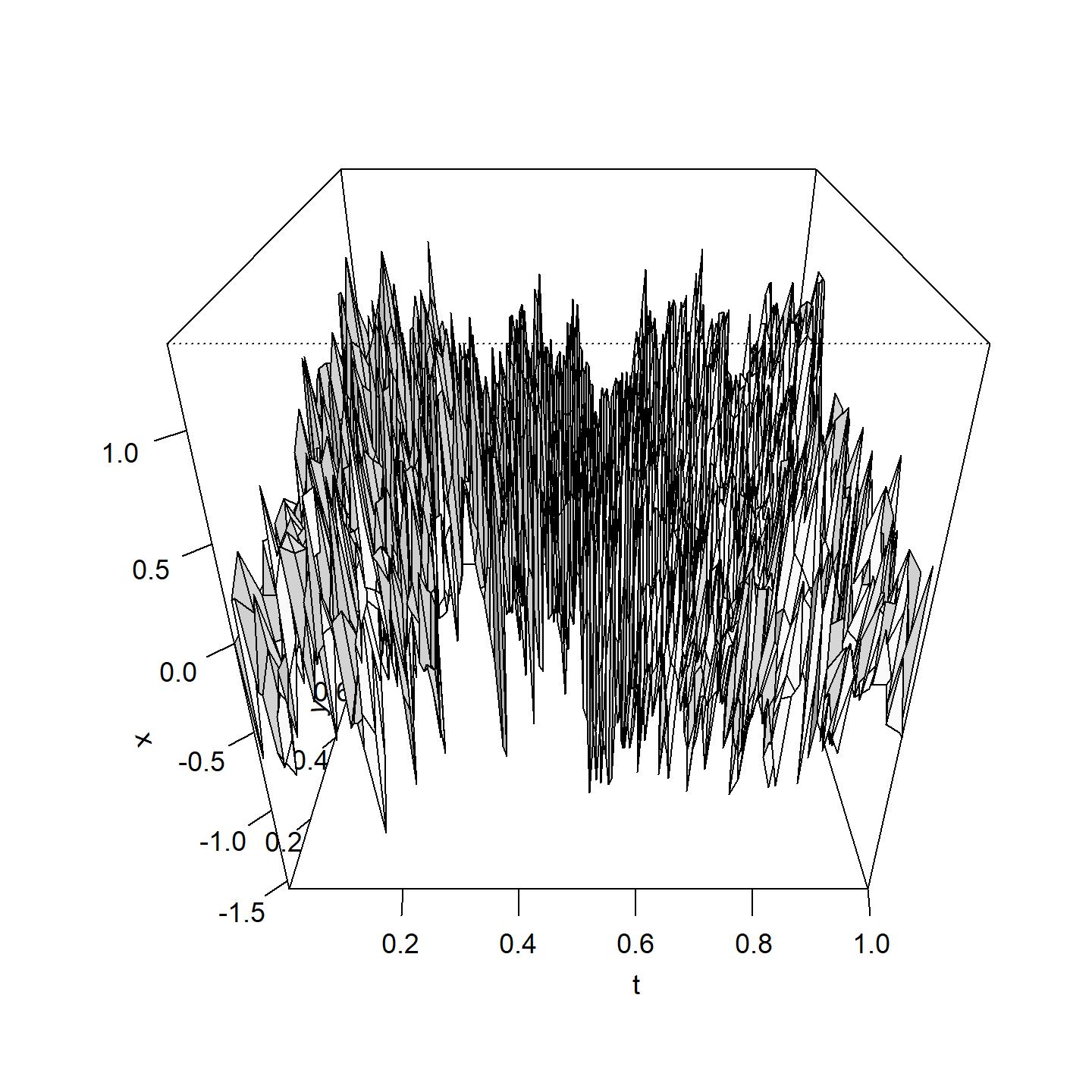}
\captionsetup{labelformat=empty,labelsep=none}
\subcaption{$\theta=$(-5,0.1,0.1,1)}
\end{center}
\end{minipage}
\begin{minipage}{0.32\hsize}
\begin{center}
\includegraphics[width=4.3cm]{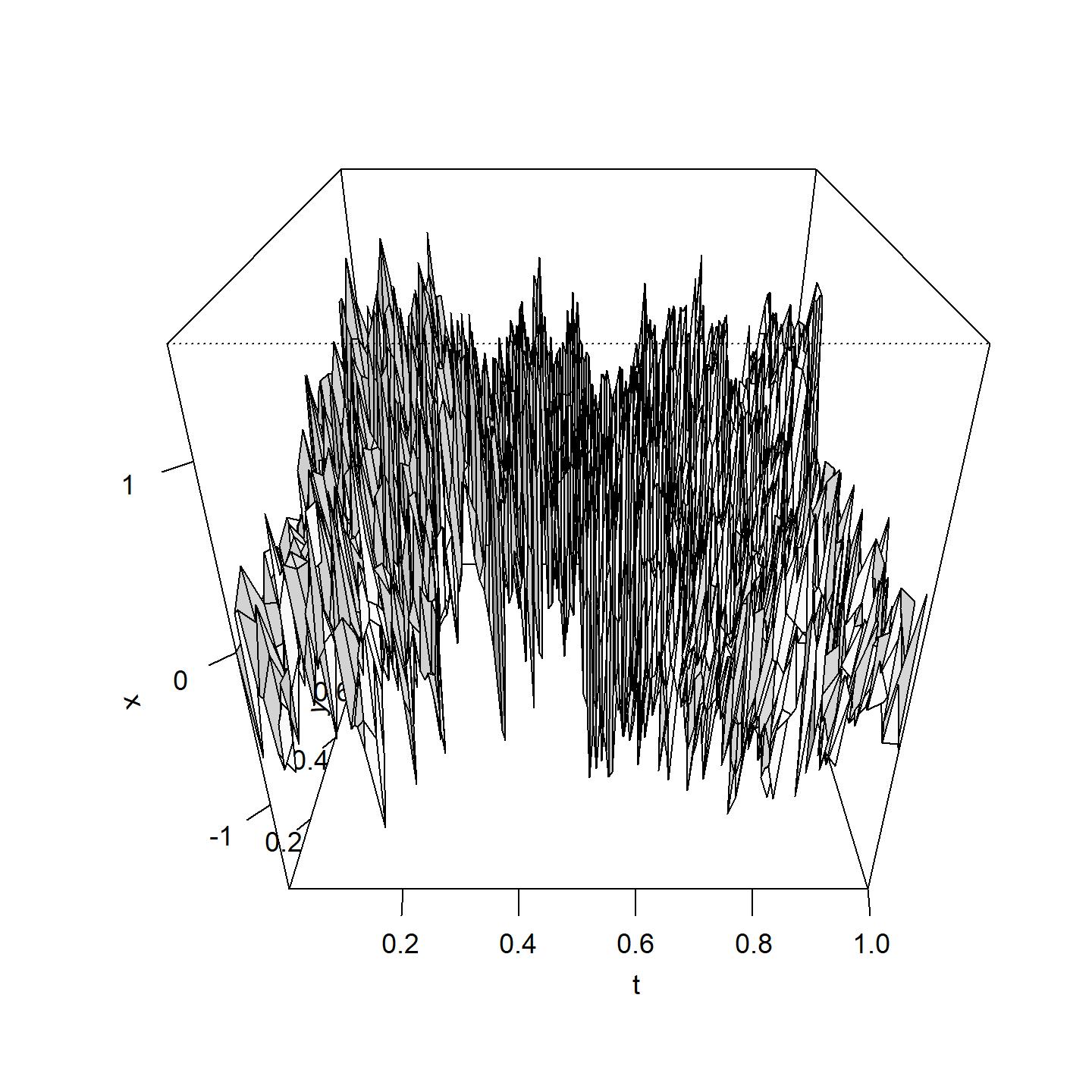}
\captionsetup{labelformat=empty,labelsep=none}
\subcaption{$\theta=$(0,0.1,0.1,1)}
\end{center}
\end{minipage}
\begin{minipage}{0.32\hsize}
\begin{center}
\includegraphics[width=4.3cm]{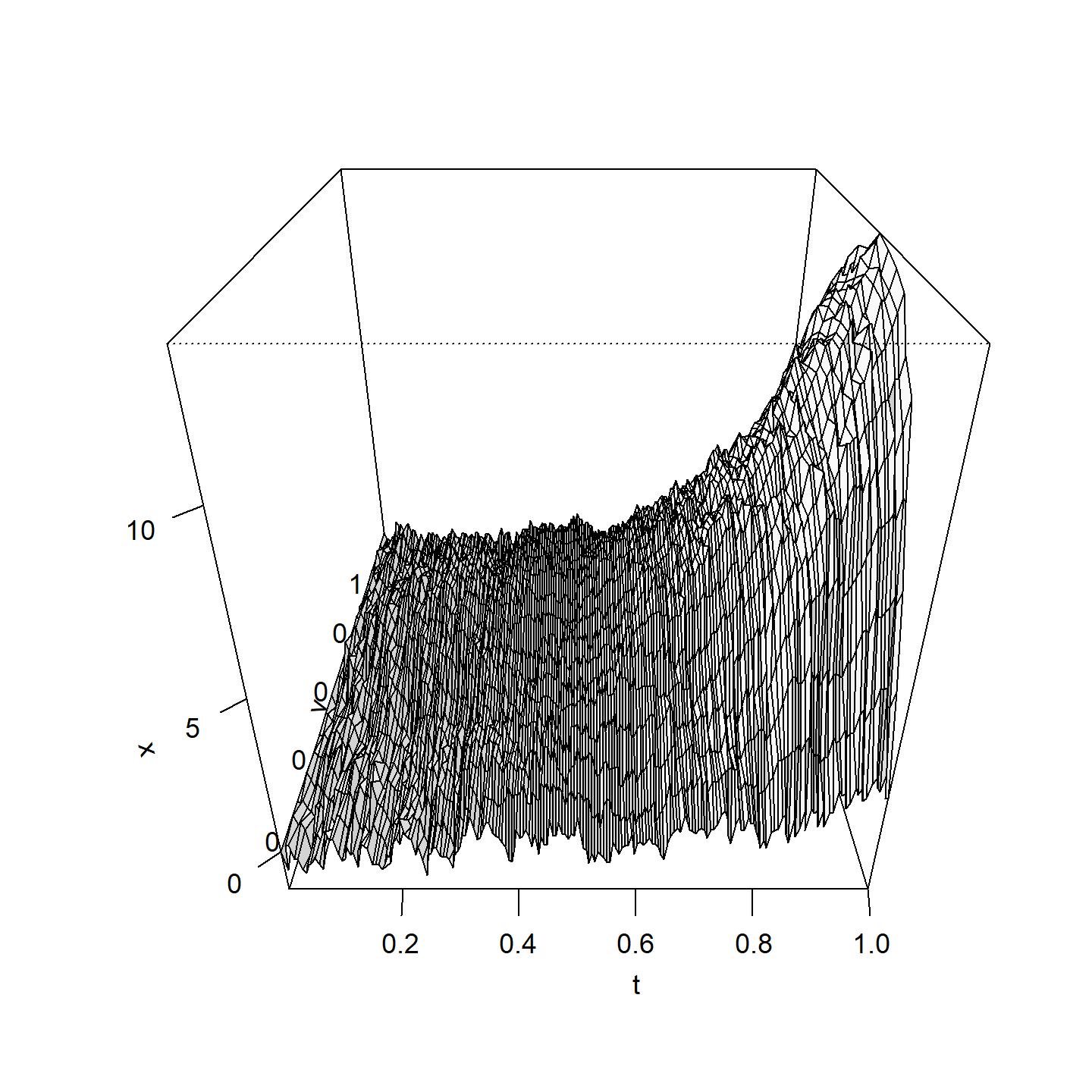}
\captionsetup{labelformat=empty,labelsep=none}
\subcaption{$\theta=$(5,0.1,0.1,1)}
\end{center}
\end{minipage}
\caption{Sample paths with $\theta_0=-5, 0, 5$ and $T=1$ (t-axis side)}\label{t0-2}

\begin{minipage}{0.32\hsize}
\begin{center}
\includegraphics[width=4.3cm]{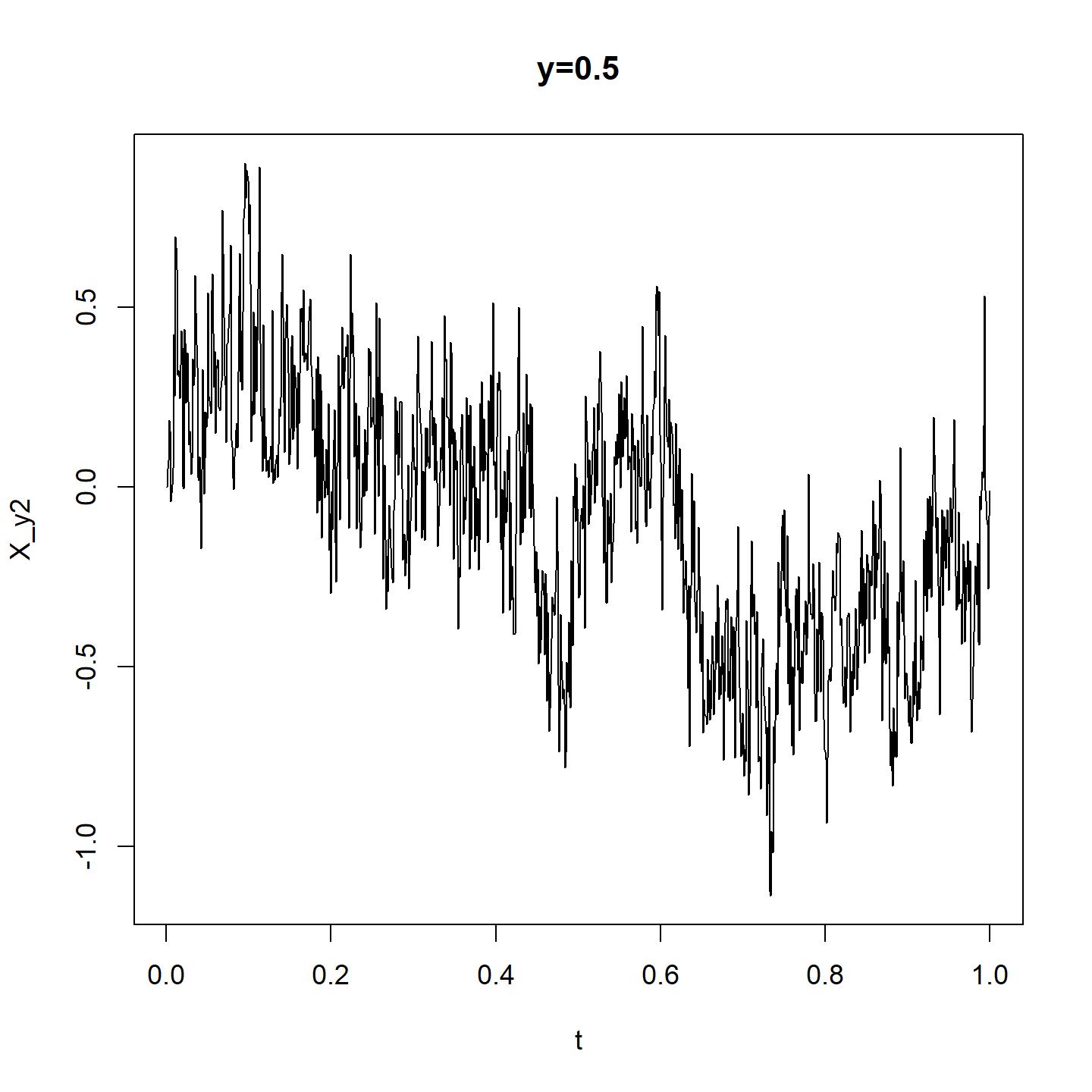}
\captionsetup{labelformat=empty,labelsep=none}
\subcaption{$\theta=$(-5,0.1,0.1,1)}
\end{center}
\end{minipage}
\begin{minipage}{0.32\hsize}
\begin{center}
\includegraphics[width=4.3cm]{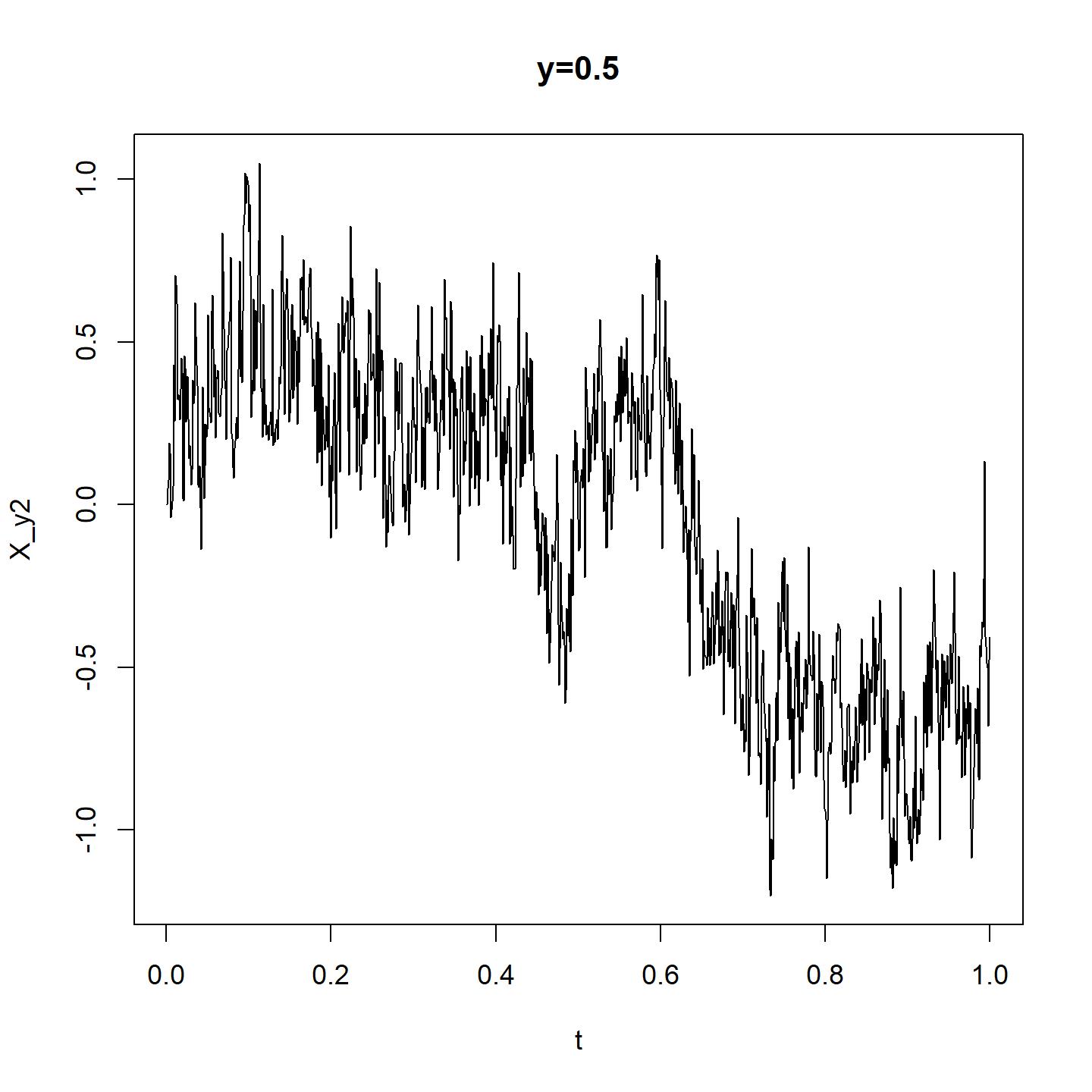}
\captionsetup{labelformat=empty,labelsep=none}
\subcaption{$\theta=$(0,0.1,0.1,1)}
\end{center}
\end{minipage}
\begin{minipage}{0.32\hsize}
\begin{center}
\includegraphics[width=4.3cm]{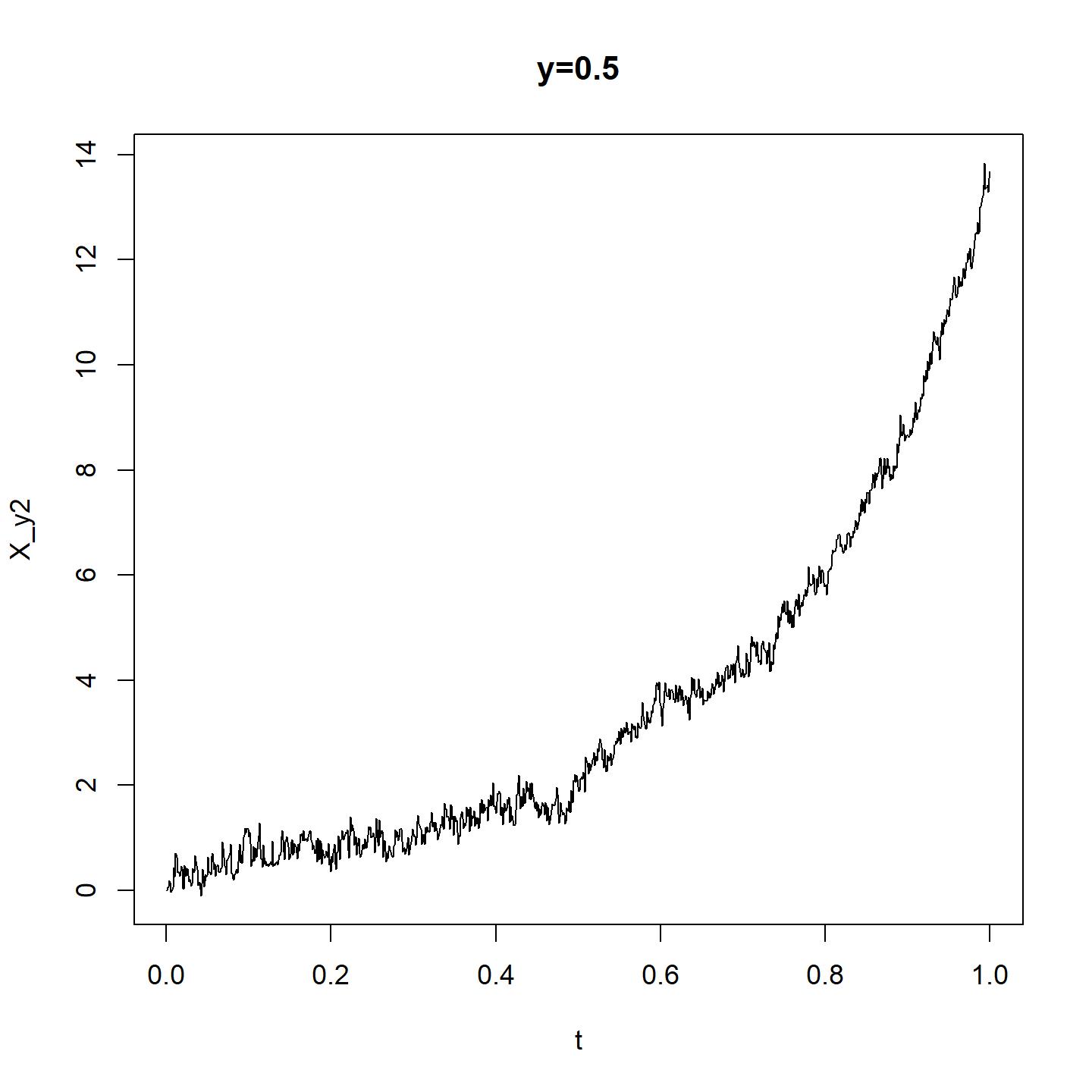}
\captionsetup{labelformat=empty,labelsep=none}
\subcaption{$\theta=$(5,0.1,0.1,1)}
\end{center}
\end{minipage}
\caption{Sample paths with $\theta_0=-5, 0, 5$ and $T=100$ (cross section at $y=0.5$)}\label{t0-3}
\end{figure}

\begin{figure}[H]
\begin{minipage}{0.32\hsize}
\begin{center}
\includegraphics[width=4.5cm]{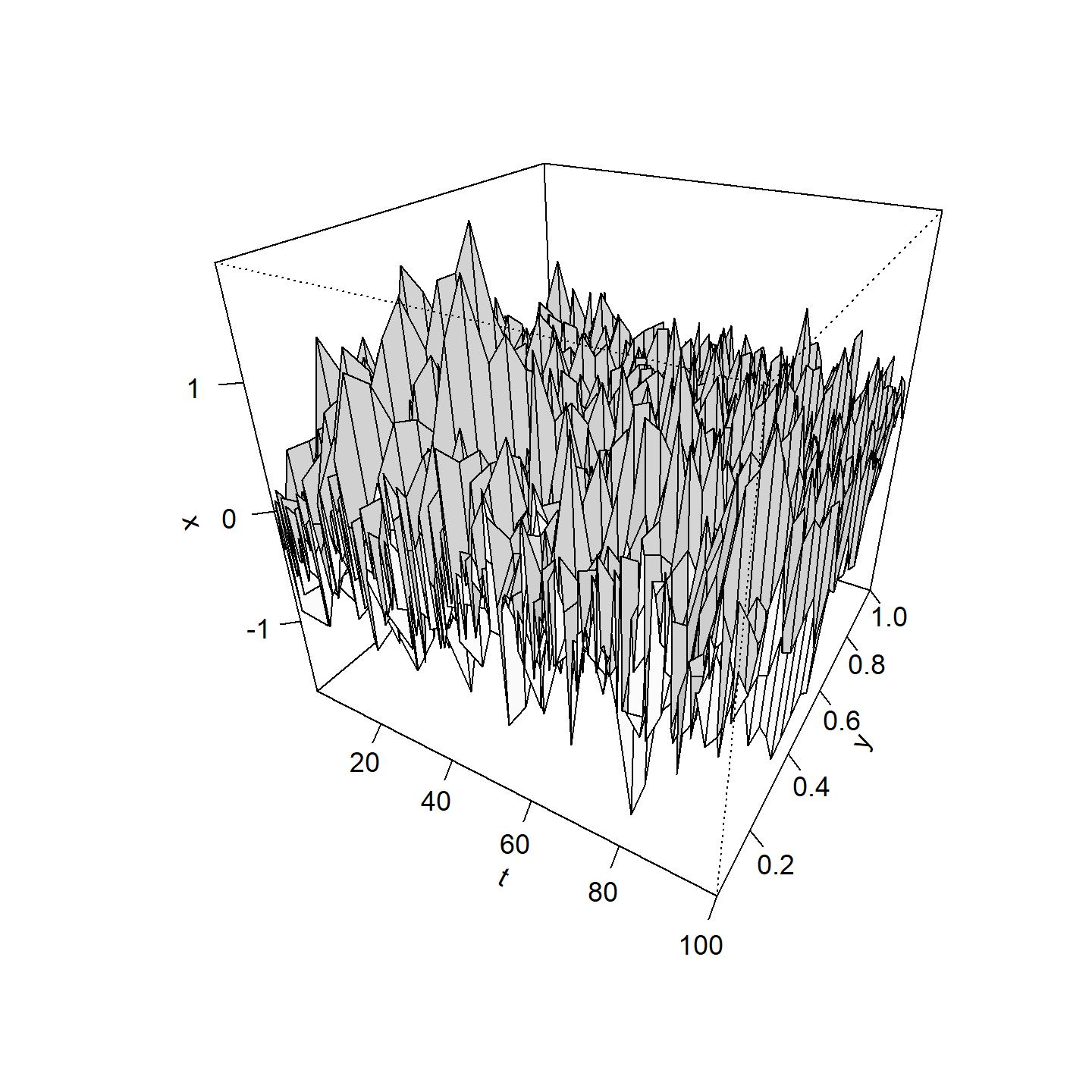}
\captionsetup{labelformat=empty,labelsep=none}
\subcaption{$\theta=$(-5,0.1,0.1,1)}
\end{center}
\end{minipage}
\begin{minipage}{0.32\hsize}
\begin{center}
\includegraphics[width=4.5cm]{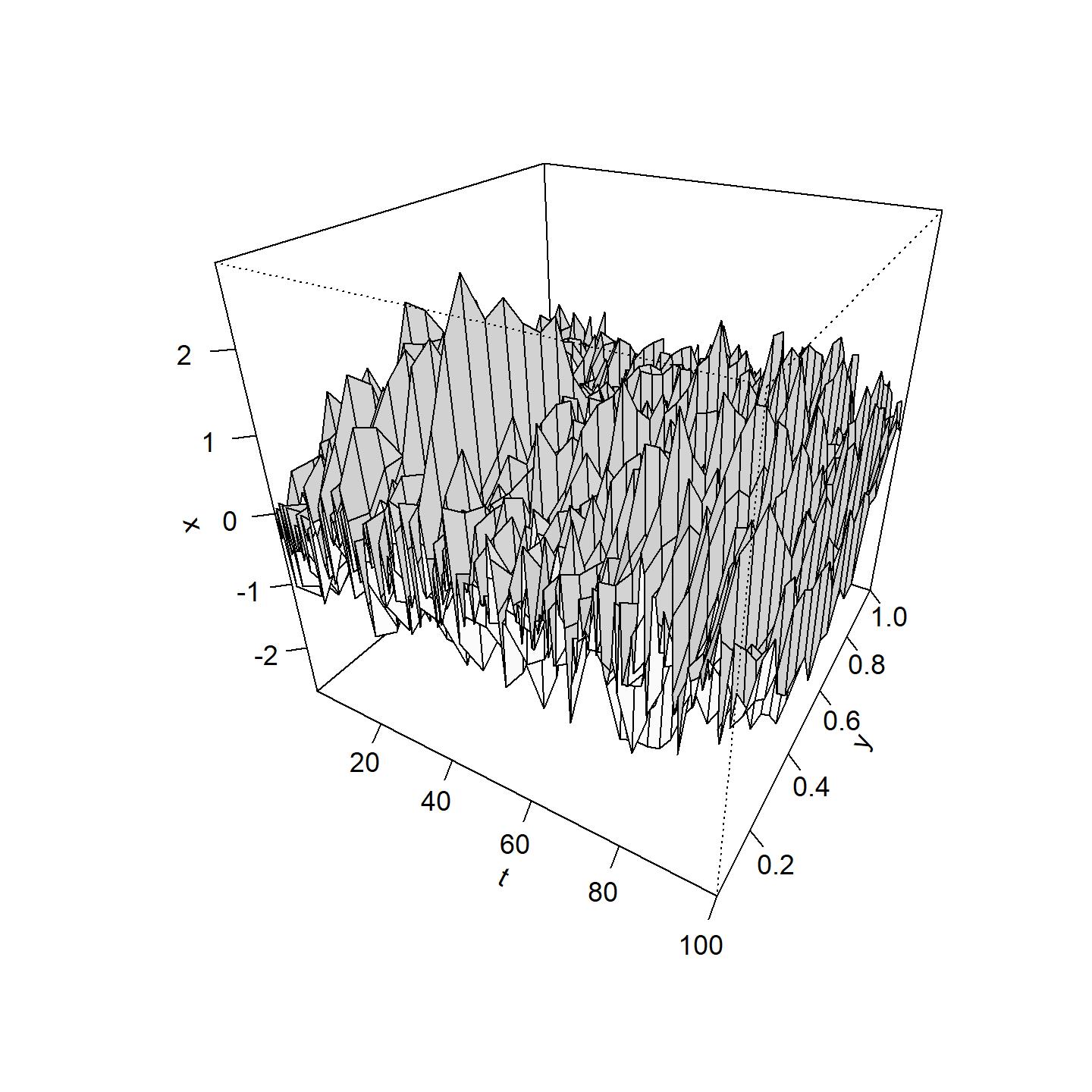}
\captionsetup{labelformat=empty,labelsep=none}
\subcaption{$\theta=$(0,0.1,0.1,1)}
\end{center}
\end{minipage}
\begin{minipage}{0.32\hsize}
\begin{center}
\includegraphics[width=4.5cm]{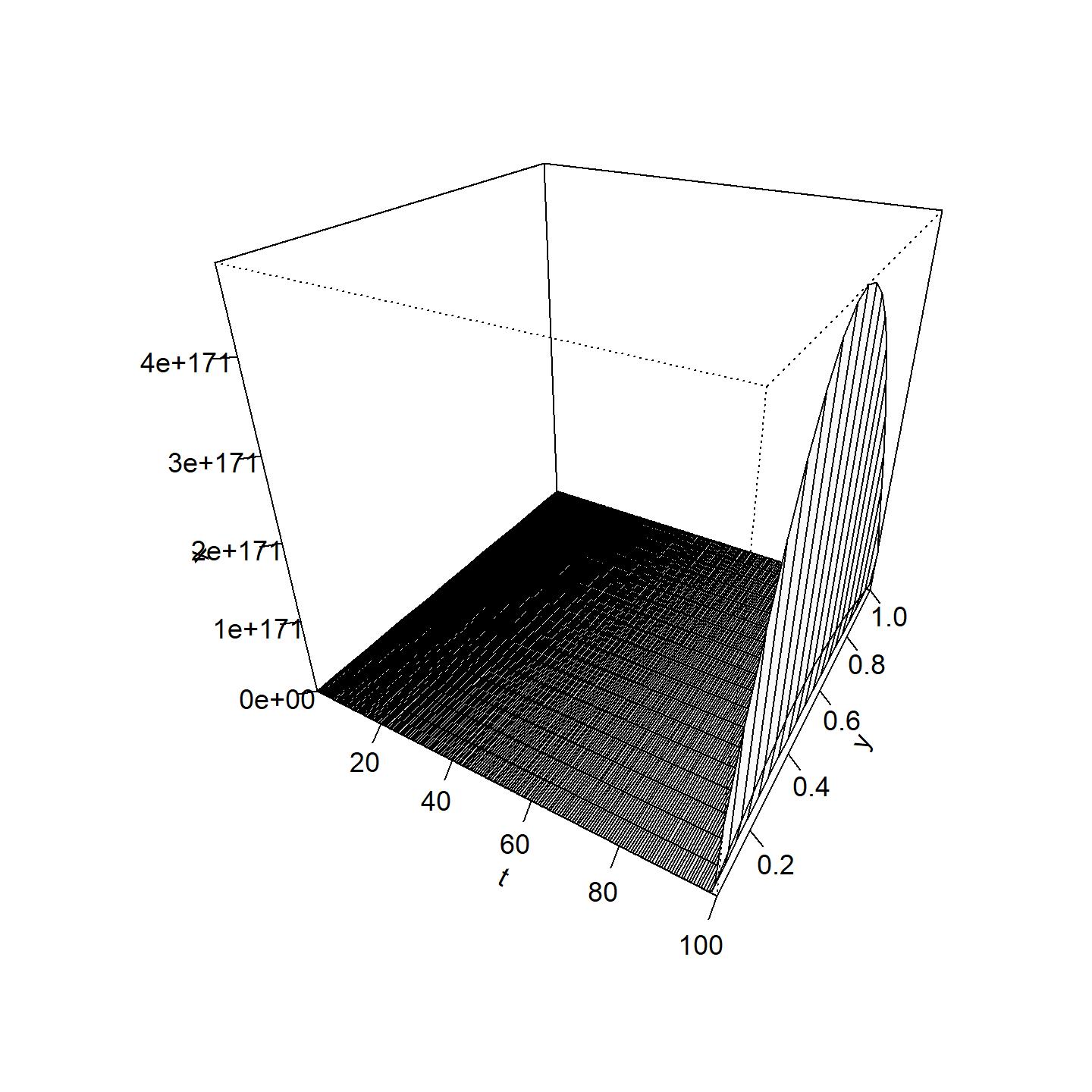}
\captionsetup{labelformat=empty,labelsep=none}
\subcaption{$\theta=$(5,0.1,0.1,1)}
\end{center}
\end{minipage}
\caption{Sample paths with $\theta_0=-5, 0, 5$ and $T=100$}\label{t0-4}
\end{figure}

\begin{figure}[H]
\begin{minipage}{0.32\hsize}
\begin{center}
\includegraphics[width=4.3cm]{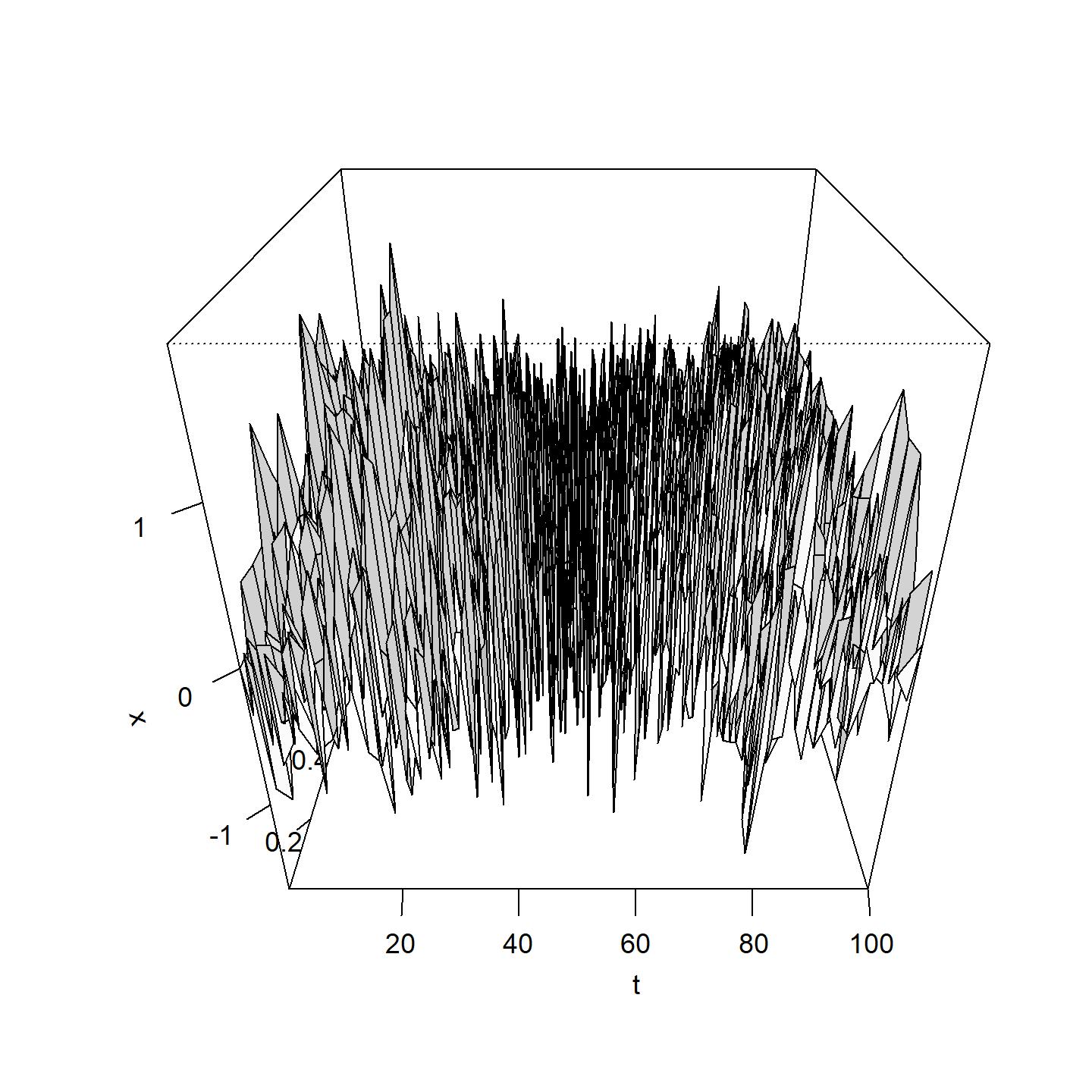}
\captionsetup{labelformat=empty,labelsep=none}
\subcaption{$\theta=$(-5,0.1,0.1,1)}
\end{center}
\end{minipage}
\begin{minipage}{0.32\hsize}
\begin{center}
\includegraphics[width=4.3cm]{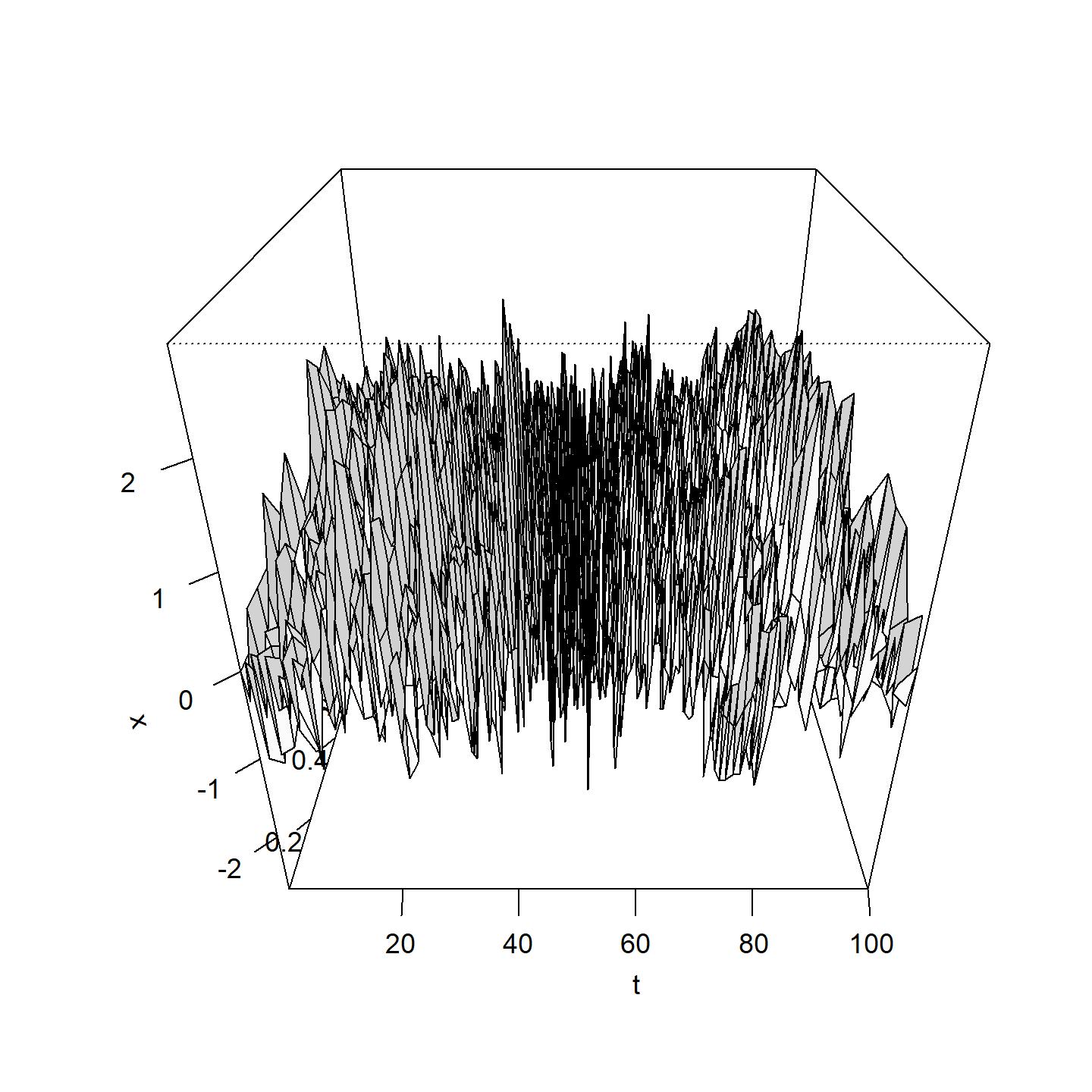}
\captionsetup{labelformat=empty,labelsep=none}
\subcaption{$\theta=$(0,0.1,0.1,1)}
\end{center}
\end{minipage}
\begin{minipage}{0.32\hsize}
\begin{center}
\includegraphics[width=4.3cm]{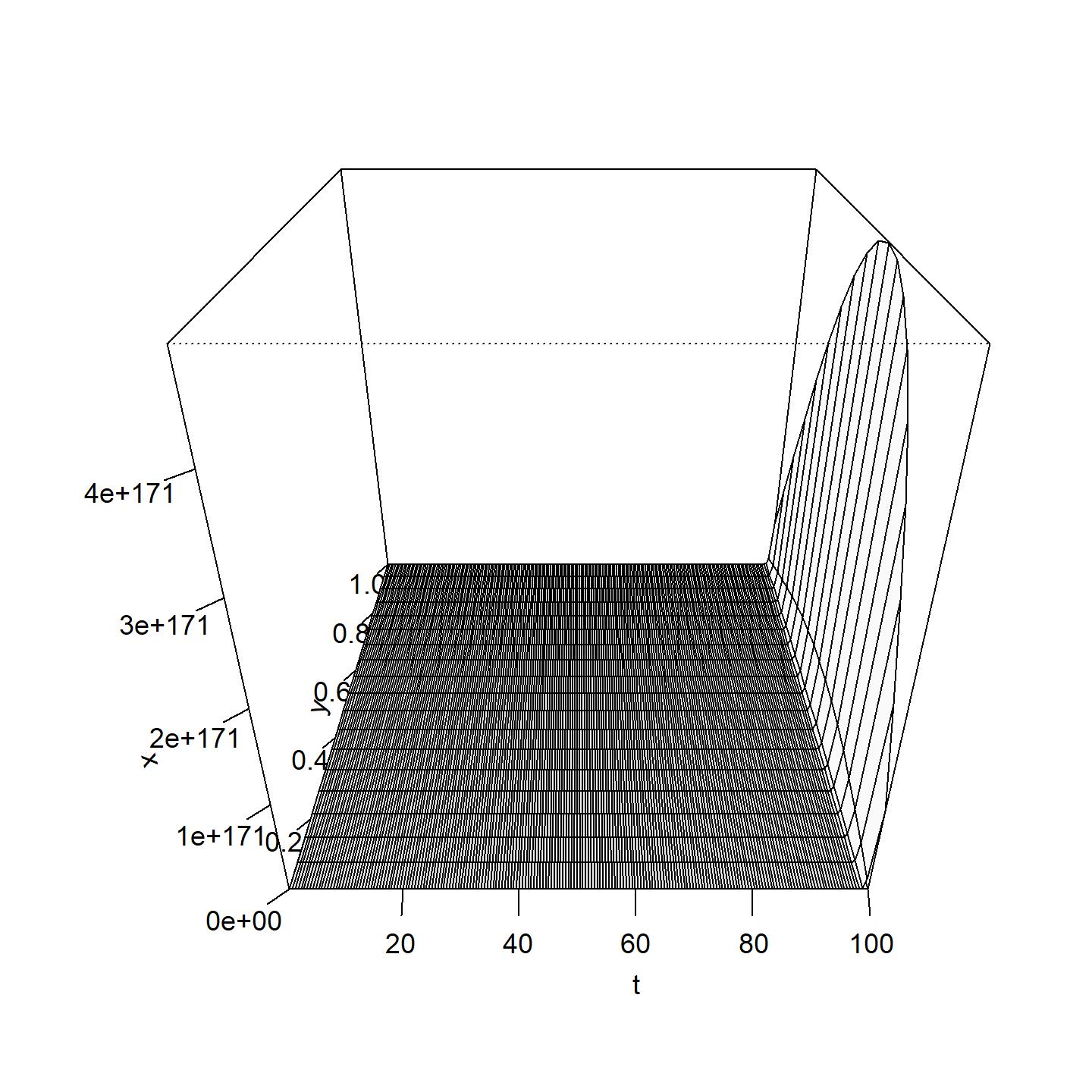}
\captionsetup{labelformat=empty,labelsep=none}
\subcaption{$\theta=$(5,0.1,0.1,1)}
\end{center}
\end{minipage}
\caption{Sample paths with $\theta_0=-5, 0, 5$ and $T=100$ (t-axis side)}\label{t0-5}
\end{figure}

\begin{figure}[H]
\begin{minipage}{0.32\hsize}
\begin{center}
\includegraphics[width=4.3cm]{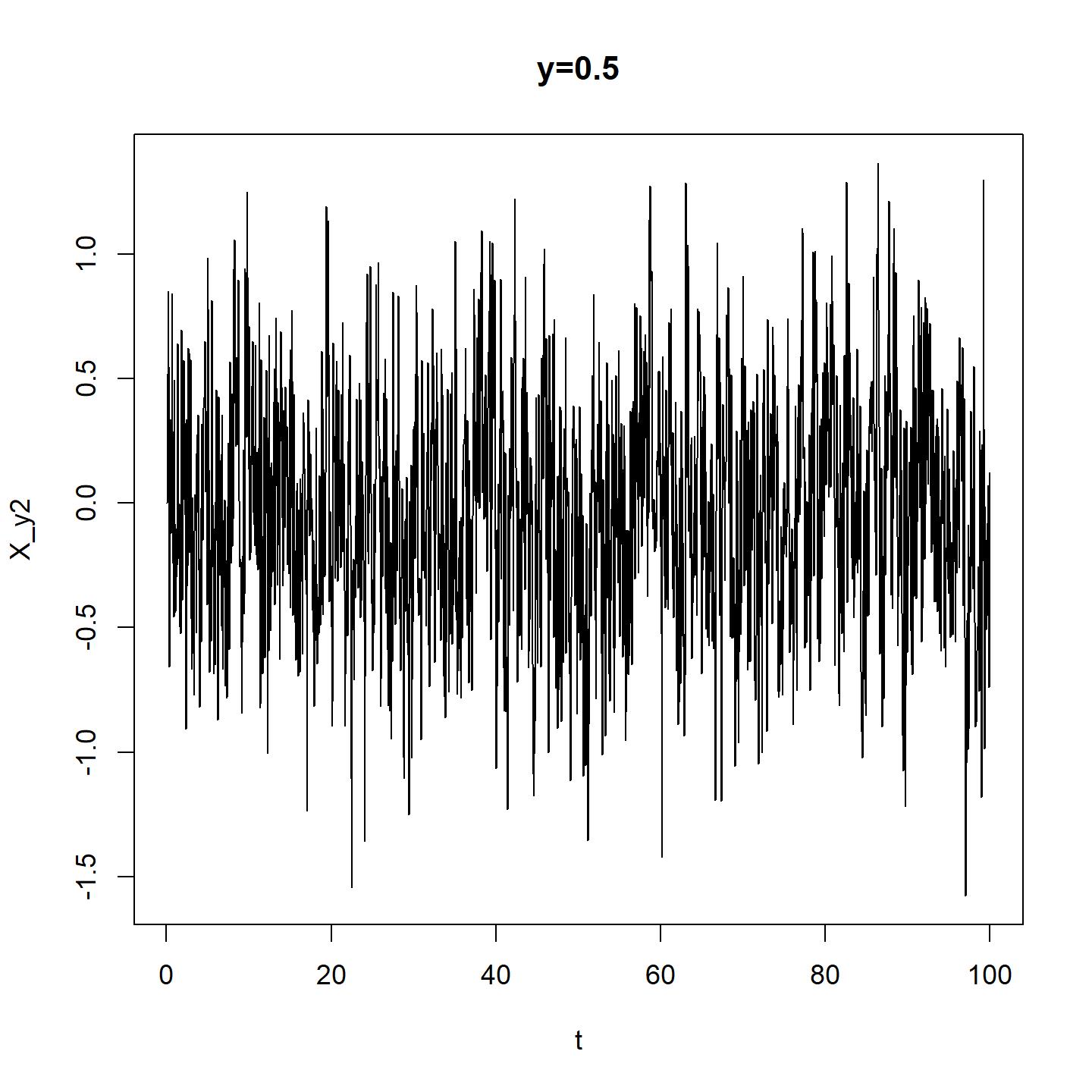}
\captionsetup{labelformat=empty,labelsep=none}
\subcaption{$\theta=$(-5,0.1,0.1,1)}
\end{center}
\end{minipage}
\begin{minipage}{0.32\hsize}
\begin{center}
\includegraphics[width=4.3cm]{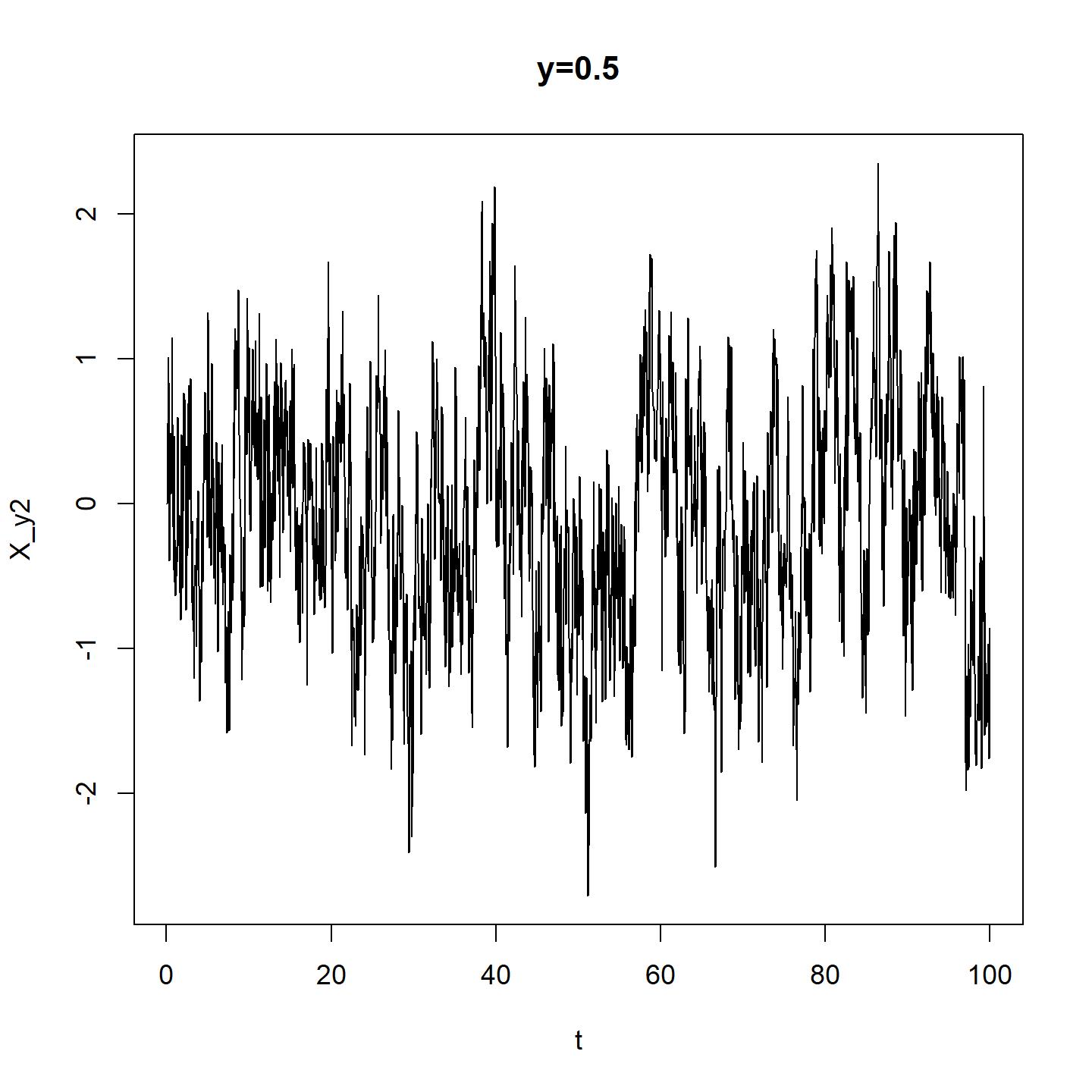}
\captionsetup{labelformat=empty,labelsep=none}
\subcaption{$\theta=$(0,0.1,0.1,1)}
\end{center}
\end{minipage}
\begin{minipage}{0.32\hsize}
\begin{center}
\includegraphics[width=4.3cm]{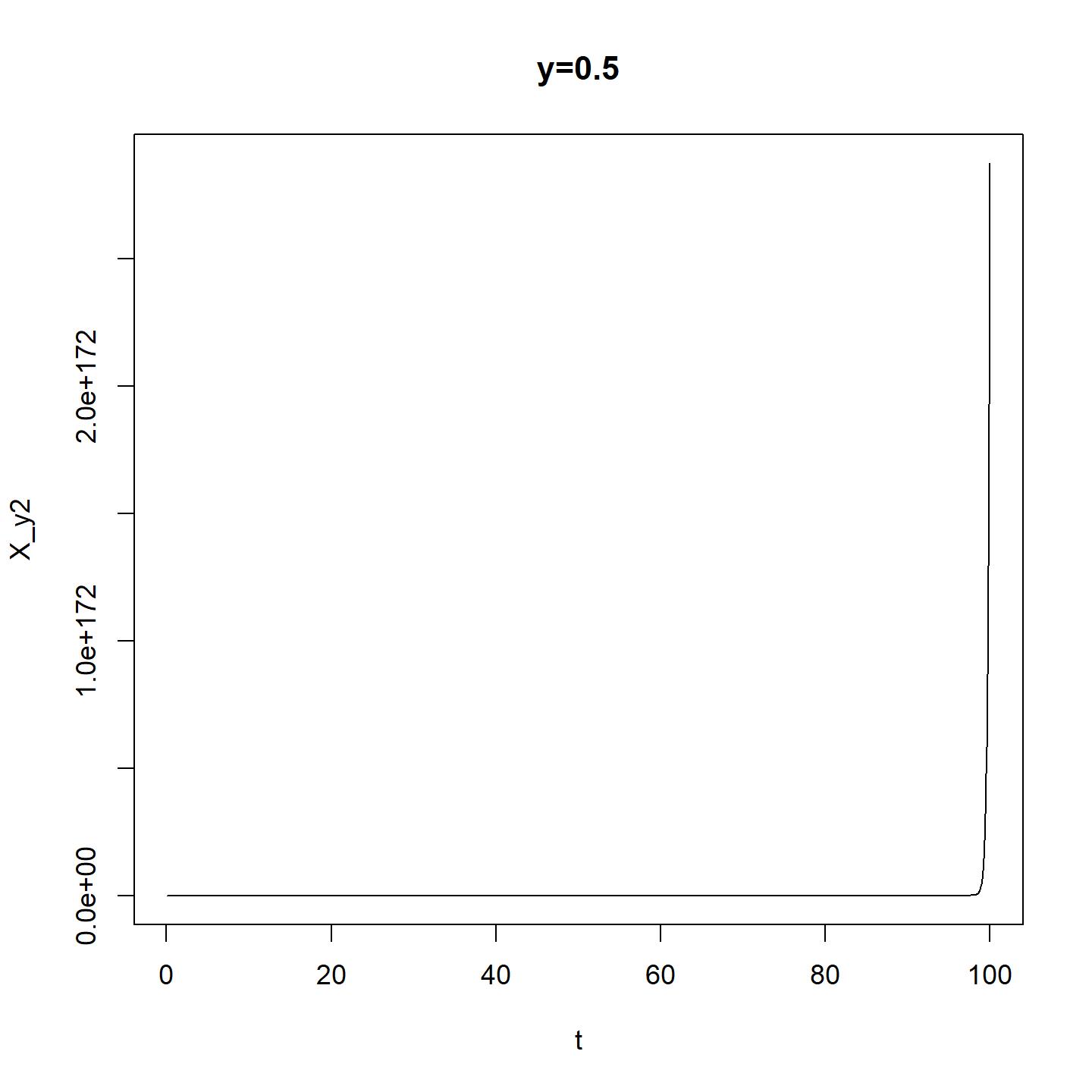}
\captionsetup{labelformat=empty,labelsep=none}
\subcaption{$\theta=$(5,0.1,0.1,1)}
\end{center}
\end{minipage}
\caption{Sample paths with $\theta_0=-5, 0, 5$ and $T=100$ (cross section at $y=0.5$)}\label{t0-6}
\end{figure}

\newpage

\begin{en-text}
\section*{Acknowledgement}
The authors would like to thank the editor, the associate editor, and the two reviewers
for their valuable comments.
This work 
was partially supported by 
JST CREST,
JSPS KAKENHI Grant Number 
JP17H01100 
and Cooperative Research Program
of the Institute of Statistical Mathematics.
\end{en-text}

\section*{References}
\begin{description}{}{
}

\item
Bibby, B.\ M.\ and S{\o}rensen, M.\ 
(1995).\ 
Martingale estimating functions for discretely 
observed diffusion processes.\
{\it Bernoulli}, 
{\bf 1}, 17--39.

\item
Bibinger, M.\ and Trabs, M.\ 
(2017).\ 
Volatility estimation for stochastic pdes using high-frequency observations.\
arXiv:1710.03519.


\item
Chong, C.\ 
(2019).\ 
High-frequency analysis of parabolic stochastic PDEs.\  
Forthcoming in The Annals of Statistics.

\item
Cialenco, I.\
 (2018).\ 
Statistical inference for SPDEs:\ an overview.
{\it Statistical Inference for Stochastic Processes},
{\bf 21}, 309--329.

\item
Cialenco, I., Delgado-Vences, F.\ and Kim, H.-J.\
(2019).\
Drift Estimation for Discretely Sampled SPDEs.\
arXiv:1904.10884.

\item
Cialenco, I.\ and Glatt-Holtz, N.\ 
(2011).\ 
Parameter estimation for the stochastically perturbed Navier-Stokes equations.\ 
{\it Stochastic Processes and their Applications.\ }
{\bf 121}, 701--724.

\item
Cialenco, I.,  Gong, R.\ and Y. Huang, Y.\ 
(2018).\ 
Trajectory fitting estimators for SPDEs driven by additive noise.\ 
{\it Statistical Inference for Stochastic Processes.\ } 
{\bf 21}, 1--19.

\item
Cialenco, I.\ and Huang, Y.\ 
(2017).\ 
A note on parameter estimation for discretely sampled SPDEs.\
arXiv:1710.01649.



\item
Cont, R. (2005). 
Modeling term structure dynamics: an infinite dimensional approach.\ 
{\it International Journal of Theoretical and Applied Finance}, 
{\bf 8}, 357--380.

\item
De Gregorio, A.\ and Iacus, S.\ M.\
(2013).\ 
On a family of test statistics for discretely observed diffusion processes.\ 
{\it Journal of Multivariate Analysis}, {\bf 122}, 292--316.

\item	
Dohnal, G.\
(1987).\ 
On Estimating the Diffusion Coefficient.\
{\it Journal of Applied Probability},
{\bf  24}, 105--114.


\item
Florens-Zmirou, D.\ 
(1989).\ 
Approximate discrete time schemes for statistics 
of diffusion processes.\ 
{\it Statistics}, {\bf 20}, 547--557.

\item
Genon-Catalot, V.\ and Jacod, J.\ 
(1993).\ 
On the estimation of the diffusion coefficient 
for multidimensional diffusion processes.\
{\it Annales de l’Institut Henri Poincar{\'e} Probabilit{\'e}s et Statistiques, }{\bf 29}, 119--151.

\item
Genon-Catalot, V.\  and Jacod, J.\ 
(1994).\ 
Estimation of the Diffusion Coefficient for Diffusion Processes: Random Sampling.\
{\it Scandinavian Journal of Statistics}, 
{\bf 21}, 193--221.

\item
Kaino, Y.\ and Uchida, M.\ (2018).\ 
Hybrid estimators for stochastic differential equations from reduced data.\ 
{\it Statistical Inference for Stochastic Processes}, {\bf 21}, 435--454.

\item
Kamatani, K.\ and Uchida, M.\ 
(2015).\ 
Hybrid multi-step estimators for stochastic differential equations based on sampled data.\
{\it Statistical Inference for Stochastic Processes}, 
{\bf 18}, 177--204.

\item 
Kessler, M.\ 
(1995).\ 
Estimation des param{\`e}tres d'une diffusion par des contrastes corrig{\'e}s.\
{\it Comptes Rendus de l'Acad{\'e}mie des Sciences - Series I - Mathematics}, 
{\bf 320}, 359--362. 

\item
Kessler, M.\ 
(1997).\ 
Estimation of an ergodic diffusion from discrete observations.\
{\it Scandinavian Journal of Statistics}, {\bf 24}, 211--229.

\item
Kutoyants, Yu.\ A.\ 
(1984).\ 
Parameter estimation for stochastic processes.\ 
Prakasa Rao, B.L.S.\ (ed.\ )
Heldermann, Berlin.

\item
Kutoyants, Yu.\ A.\ 
(2004).\ 
{\it Statistical inference for ergodic diffusion processes.\ }
Springer-Verlag, London.\

\item 
Markussen, B.\
(2003).\
Likelihood inference for a discretely observed stochastic partial differential equation.\
{\it Bernoulli},
{\bf 9}, 745--762.

\item
Masuda, H.\
(2013a).\
Asymptotics for functionals of self-normalized residuals of discretely observed stochastic processes.\ 
{\it Stochastic Processes and their Applications}, {\bf 123}, 2752--2778.

\item
Masuda, H.\
(2013b).\
Convergence of Gaussian quasi-likelihood random fields for ergodic Levy driven SDE observed at high frequency.\ 
{\it The Annals of Statistics}, {\bf 41}, 1593--1641.

\item
Nakakita, S.\ H. and Uchida, M.\ (2019).\ 
Inference for ergodic diffusions plus noise.\ 
{\it Scandinavian Journal of Statistics}, {\bf 46}, 470--516.


\item
Ogihara, T.\  (2018).\
Parametric inference for nonsynchronously observed diffusion
processes in the presence of market microstructure noise.\
{\it Bernoulli},
{\bf 24}, 3318--3383.

\item 
Ogihara, T.\ and Yoshida, N.\
(2011).\ 
Quasi-likelihood analysis for the stochastic differential equation with jumps.\
{\it Statistical Inference for Stochastic Processes}, {\bf 14}, 189--229.

\item 
Ogihara, T.\ and Yoshida, N.\
(2014).\ 
Quasi-likelihood analysis for nonsynchronously observed diffusion processes.\
{\it Stochastic Processes and their Applications}, {\bf 124}, 2954--3008.

\item
Prakasa Rao, B.\ L.\ S.\ 
(1983).\
Asymptotic theory for nonlinear least squares estimator for diffusion processes.\ 
{\it Math.\ Operationsforsch.\ Statist.\ Ser.\ Statist.,} {\bf 14}, 195--209.

\item
Prakasa Rao, B.\ L.\ S.\ 
(1988).\
Statistical inference from sampled data for stochastic processes.\
{\it Contemporary Mathematics}, 
{\bf 80}, 249--284.\ 
Amer.\ Math.\ Soc., Providence, RI.

\item
Shimizu, Y. 
(2006).\ 
M-estimation for discretely observed ergodic diffusion processes with infinitely many jumps.\
{\it Statistical Inference for Stochastic Processes},  {\bf 9}, 179--225.

\item 
Shimizu, Y.\ and Yoshida, N.\ 
(2006).\ 
Estimation of parameters for diffusion processes with jumps from discrete observations.\ 
{\it Statistical Inference for Stochastic Processes}, {\bf 9}, 227--277.

\item
Uchida, M.\ 
(2010).\ 
Contrast-based information criterion for ergodic diffusion processes from discrete observations.\ 
{\it Annals of the Institute of Statistical Mathematics}, {\bf 62}, 161--187. 

\item
Uchida, M.\ and Yoshida, N.\ 
(2012).\ 
Adaptive estimation of an ergodic diffusion process
based on sampled data.\ 
{\it Stochastic Processes and their Applications}, {\bf 122}, 2885--2924.

\item
Uchida, M. and Yoshida, N. (2013). 
 Quasi likelihood analysis of volatility and nondegeneracy of statistical random field.\ 
{\it Stochastic Processes and their Applications}, {\bf 123},  2851--2876.

\item
Uchida, M.\ and Yoshida, N.\ 
(2014).\ 
Adaptive Bayes type estimators of ergodic diffusion processes 
from discrete observations.\
{\it Statistical Inference for Stochastic Processes}, 
{\bf 17}, 181--219.

\item
Yoshida, N.\ 
(1992).\ 
Estimation for diffusion processes from discrete observation.\
{\it Journal of Multivariate Analysis}, 
{\bf 41}, 220--242.

\item
Yoshida, N.\
(2011).\
Polynomial type large deviation inequalities and quasi-likelihood analysis for stochastic differential equations.\
{\it Annals of the Institute of Statistical Mathematics}, {\bf 63}, 431--479.

\end{description}


\end{document}